\documentclass[titlepage,11pt]{article}
\usepackage[USenglish]{babel}
\usepackage{multicol}
\usepackage{a4wide}
\usepackage[all]{xy}
\usepackage{amssymb}
\usepackage{amsmath}
\usepackage{latexsym}
\usepackage{pifont}
\usepackage{pdfsync}
\usepackage{color}
\usepackage{mathrsfs}
\usepackage{enumerate}
\usepackage{float,fancyvrb,verbatim,multicol}
\usepackage{ifthen}
\usepackage[parfill]{parskip}    
\usepackage{xr,url}
\usepackage{makeidx}
\usepackage[utf8]{inputenc} 
\title{Topologies and all that --- A Tutorial}
\author{Ernst-Erich Doberkat\\Math ++ Software, Bochum\\\texttt{eed@doberkat.de}}
\date{\today}
\newcommand{\labelImpl}[2]{\ensuremath{\ref{#1}~\Rightarrow~\ref{#2}}}

\newcommand{\Klasse}[2]{\left[#1\right]_{#2}}
\newcommand{\Faktor}[2]{{#1}/{#2}}
\newcommand{\fMap}[1]{\eta_{#1}}
\newcommand{\Bild}[2]{{#1}\left[#2\right]}
\newcommand{\InvBild}[2]{\Bild{#1^{-1}}{#2}}

\newcommand{\Folge}[1]{(#1_n)_{n \in \Nat}}

\newcommand{\spaceFont}[1]{\mathfrak{#1}}

\newcommand{\Category}[1]{\ensuremath{\mathbf{#1}}}
\newcommand{\SubProbSenza}{\spaceFont{S}}
\newcommand{\ProbSenza}{\spaceFont{P}}
\newcommand{\PowerSet}[1]{\ensuremath{\mathcal{P}\left(#1\right)}}

\newcommand{\FunctorSenza}[1]{\ensuremath{\spaceFont{#1}}}

\newcommand{\Borel}[1]{\ensuremath{{\mathcal B}(#1)}}

\edef\LinkeKlammer{\lbrack\!\lbrack}
\edef\RechteKlammer{\rbrack\!\rbrack}
\newcommand{\Gilt}[1][\phi]{\ensuremath{\LinkeKlammer#1\RechteKlammer}}

\newcommand{\Closure}[1]{\ensuremath{\mathsf{cl}(#1)}}

\newtheorem{definition}{Definition}[section]
\newcommand{\BeginDefinition}[1]{%
  \begin{definition}\label{#1}
}
\newcommand{\EndDefinition}{\end{definition}}

\newtheorem{example}[definition]{Example}
\newcommand{\BeginExample}[1]{%
  \begin{example}\label{#1}\rm
}
\newcommand{\EndExample}{--- \end{example}}

\newtheorem{observation}[definition]{Observation}
\newcommand{\BeginObservation}[1]{
  \begin{observation}\label{#1}\rm
}
\newcommand{\EndObservation}{--- \end{observation}}

\newtheorem{theorem}[definition]{Theorem}
\newcommand{\BeginTheorem}[1]{%
  \begin{theorem}\label{#1}
}
\newcommand{\EndTheorem}{\end{theorem}}

\newtheorem{corollary}[definition]{Corollary}
\newcommand{\BeginCorollary}[1]{
  \begin{corollary}\label{#1}
}

\newtheorem{proposition}[definition]{Proposition}
\newcommand{\BeginProposition}[1]{%
  \begin{proposition}\label{#1}
}
\newcommand{\EndProposition}{\end{proposition}}

\newcommand{\EndCorollary}{\end{corollary}}
\newtheorem{lemma}[definition]{Lemma}
\newcommand{\BeginLemma}[1]{%
  \begin{lemma}\label{#1}
}
\newcommand{\EndLemma}{\end{lemma}}

\newtheorem{claim}{Claim}
\newcommand{\BeginClaim}[1]{%
  \begin{claim}\label{#1}
}
\newcommand{\EndClaim}{\end{claim}}

\newenvironment{proof}{\textbf{Proof\ }}{\ensuremath{\QED}}
\newcommand{\BeginProof}{\begin{proof}}
\newcommand{\EndProof}{\end{proof}}

\newenvironment{remark}{\textbf{Remark:\ }}{}
\newcommand{\BeginRemark}{\begin{remark}}
\newcommand{\EndRemark}{\QED\end{remark}}
\newcommand{\QED}{%
\ensuremath{\dashv}
}

\newcommand{\Real}{\mathbb{R}}
\newcommand{\pReal}{\mathbb{R}_{+}}
\newcommand{\Nat}{\mathbb{N}}
\newcommand{\Rational}{\mathbb{Q}}

\def\theta{\vartheta}
\def\dots{\ldots}

\def\Ganz{\mathbb{Z}}
\makeatletter
\def\@axx#1{\ensuremath{\mathbb{(#1)}}}
\def\AC{\@axx{AC}}
\def\WO{\@axx{WO}}
\def\ZL{\@axx{ZL}}
\def\MP{\@axx{MP}}
\def\MI{\@axx{MI}}
\def\AD{\@axx{AD}}
\makeatother

\newcommand{\isEquiv}[3]{\ensuremath{{#1}\ {#3}\ {#2}}}
\newenvironment{theExercises}
{\begin{exercise}\rm}
{\end{exercise}}
\newtheorem{exercise}{Exercise}
\newcommand{\BeginExercise}[1]{%
  \begin{theExercises}\label{#1}
}
\newcommand{\EndExercise}{\end{theExercises}}
\newcommand{\Exercise}[2]{\BeginExercise{#1}{#2}\EndExercise}

\def\endEx{{\Large\ding{44}}}
\renewcommand{\EndExample}{\endEx \end{example}}

\def\CatFont{\mathbf}
\def\Category#1{\ensuremath{\CatFont{#1}}}
\renewcommand{\spaceFont}[1]{\CatFont{#1}}

\def\catTS{\Category{TS}}
\def\catTop{\Category{Top}}

\def\catCha{\Category{cHA}}

\def\hom#1{\ensuremath{\mathrm{hom_{#1}}}}

\def\funV{\FunctorSenza{V}}
\def\funSP{\FunctorSenza{SP}}
\def\funLoc{\FunctorSenza{Loc}}

\makeatletter
\def\@theProbSenza#1{\ensuremath{\mathbb{#1}}}
\def\SubProbSenza{\@theProbSenza{S}}
\def\ProbSenza{\@theProbSenza{P}}
\def\FinSenza{\@theProbSenza{M}}
\def\SigmaSenza{\@theProbSenza{M}_{\ensuremath{\sigma}}}
\makeatother

\makeatletter
\def\@beta{\pmb{\beta}}
\def\betaSenza#1{\@beta_{#1}}

\makeatother

\renewcommand{\Borel}[1]{\ensuremath{{\mathcal B}(#1)}}

\def\o{\varnothing}
\def\half{\ensuremath{\frac{1}{2}}}

\def\DefSect{
\ifthenelse{\boolean{isBook}}{
\def\Section{\chapter}
\def\Subsection{\section}
\def\Subsubsection{\subsection}
\def\Subsubsubsection{\subsubsection}
}{
\def\Section{\section}
\def\Subsection{\subsection}
\def\Subsubsection{\subsubsection}
\def\Subsubsubsection{\paragraph}
}
}

\def\Closure#1{\ensuremath{{#1}^{a}}}
\def\Interior#1{\ensuremath{{#1}^{o}}}
\def\closOp#1{\ensuremath{{#1}}^{\textbf{c}}}
\makeatletter
\def\@Tut#1#2{\cite[#1]{#2}}
\def\CategCite#1{\@Tut{#1}{EED-Categs}}
\def\SetCite#1{\@Tut{#1}{EED-Tut_sets}}
\def\TopCite#1{\@Tut{#1}{EED-Topologies}}
\def\MeasCite{\cite{EED-Meas}}
\makeatother

\makeatletter
\def\@norm#1#2#3{\ensuremath{||{#1}||_{#2}^{#3}}}
\def\infNorm#1#2{\@norm{#1}{\infty}{#2}}
\def\aNorm#1{\@norm{#1}{}{}}
\makeatother

\makeatletter
\def\@conv#1{\stackrel{#1}{\longrightarrow}} 
\def\aeC{\@conv{a.e.}} 
\def\nmC{\@conv{i.m.}}
\makeatother

\makeatletter
\newcommand{\@theL}[3]{\ensuremath{{#3}_{#1}{#2}}}
\newcommand{\cLp}[2][p]{\@theL{#1}{(#2)}{{\cal L}}}
\newcommand{\rLp}[2][p]{\@theL{#1}{(#2)}{L}}
\newcommand\cLpS[1][p]{\@theL{#1}{}{{\cal L}}}
\newcommand\rLpS[1][p]{\@theL{#1}{}{L}}
\makeatother

\renewcommand{\Bild}[2]{{#1}\bigl[#2\bigr]}
\renewcommand{\Folge}[2][n]{\ensuremath{({#2}_{#1})_{#1\in\Nat}}}
\def\partMap#1#2{\ensuremath{{#1}\rightharpoonup{#2}}}

\def\allCompact{\mathfrak{C}}

\makeatletter
\newcommand{\@imText}[2][b]{\marginpar{\parbox[#1]{\marginparwidth}{\raggedleft{#2}}}}
\definecolor{@leichtgrau}{gray}{.60}
\def\Rot#1{\begin{color}{red}#1\end{color}}
\newcommand{\@randMarkierung}[1]{
\mbox{}\marginpar{\raggedleft\hspace{0pt}\textsf{\Rot{#1}}}}
\newcommand{\mmP}[2][b]{\@imText[#1]{\Rot{#2}}}
\makeatother
\newcommand{\MMP}[2][b]{\mmP[#1]{#2}}
\def\zZ#1{{#1}\times{#1}}

\makeatletter
\def\@pt{\ensuremath{\flat}}
\def\@om{\ensuremath{\sharp}}
\def\@comb#1#2{\ensuremath{{#1}^{#2}}}
\newcommand{\pTf}[1][f]{\@comb{#1}{\@pt}}
\newcommand{\pTX}[1][X]{\@comb{#1}{\@pt}}
\newcommand{\oMf}[1][f]{\@comb{#1}{\@om}}
\newcommand{\oMX}[1][X]{\@comb{#1}{\@om}}
\makeatother

\newcommand{\Cont}[1][X]{\ensuremath{{\cal C}(#1)}}
\edef\LinksOffen{(\!|}
\edef\RechtsOffen{|\!)}
\newcommand{\eXT}[1][a]{\ensuremath{\LinksOffen#1\RechtsOffen}}
\def\Zwei{2\!\!2}
\renewcommand{\upsilon}{\ensuremath{{\filterFont U}}}
\def\filterFont{\mathfrak}
\def\fiF{\ensuremath{{\filterFont F}}}
\def\theLang{\ensuremath{{\filterFont L}}} 
\newboolean{isBook}
\setboolean{isBook}{false}
\DefSect

\def\phi{\varphi}

\usepackage{fancyhdr}
\makeindex
\begin{document}
\maketitle
\begin{abstract} This is a brief introduction to the basic concepts of topology. It includes the basic constructions, discusses separation properties, metric and pseudometric spaces, and gives some applications arising from the use of topology in computing. 
\end{abstract}

\tableofcontents\newpage
{
\Section{Topological Spaces}
A topology formalizes the notion of an open set; call a set open iff each of its members leaves a little room like a breathing space around it. This gives immediately a hint at the structure of the collection of  open sets --- they should be closed under finite intersections, but under arbitrary unions, yielding the base for a calculus of observable properties, as outlined in~\cite[Chapter 1]{Smyth} or in~\cite{Vickers}. It makes use of properties of topological spaces, but puts its emphasis subtly away from the classic approach, e.g., in mathematical analysis or probability theory, by stressing different properties of a space. The traditional approach, for example, stresses separation properties like being able to separate two distinct points through an open set. Such a strong emphasis is not necessarily observed in the computationally oriented use of topologies, where for example pseudometrics for measuring the conceptual distance between objects are important, when it comes to find an approximation between Markov transition systems. 

We give in this short treatise a brief introduction to some of the main properties of topological spaces, given that we have touched upon topologies already in the context of the Axiom of Choice~\SetCite{Sect. 1.5.8}. The objective is to provide the tools and methods offered by set-theoretic topology to an application oriented reader, thus we introduce the very basic notions of topology, and hint at applications of these tools. Some connections to logic and set theory are indicated, but as Moschovakis writes ``General (pointset) topology is to set theory like parsley to Greek food: some of it gets in almost every dish, but there are no 'parsley recipes' that the good Greek cook needs to know.''~\cite[6.27, p. 79]{Moschovakis-Notes}. In this metaphor, we study the parsley here, so that it can get into the dishes which require it. 

The goal of making topology useful suggests the following core areas: one should first discuss the  \emph{basic notion} of a topology and its construction, including bases and subbases. Since compactness has been made available very early, compact spaces serve occasionally as an exercise ground. Continuity is an important topic in this context, and the basic constructions like product or quotients which are enabled by it. Since some interesting and important topological constructions are tied to filters, we study \emph{filters and convergence}, comparing in examples the sometimes more easily handled nets to the occasionally more cumbersome filters, which, however, offer some conceptual advantages. Talking about convergence, \emph{separation properties} suggest themselves; they are studied in detail, providing some classic results like Urysohn's Theorem. It happens so often that one works with a powerful concept, but that this concept requires assumptions which are too strong, hence one has to weaken it in a sensible way. This is demonstrated in the transition from compactness to local compactness; we discuss local compact spaces, and we give an example of a compactification. Quantitative aspects enter when one measures openness through a pseudometric; here many concepts are seen in a new, sharper light, in particular the problem of completeness comes up --- you have a sequence the elements of which are eventually very close to each other, and you want to be sure that a limit exists. This is possible on complete spaces, and, even better, if a space is not complete, then you can complete it. Complete spaces have some very special properties, for example the intersection of countably many open dense sets is dense again. This is Baire's Theorem, we show through a Banach-Mazur game played on a topological space that being of first category can be determined through Demon having a winning strategy. 

This completes the round trip of basic properties of topological spaces. We then present a small gallery in which topology is in action. The reason for singling out some topics is that we want to demonstrate the techniques developed with topological spaces for some interesting applications. For example, Gödel's Completeness Theorem for (countable) first order logic has been proved by Rasiowa and Sikorski through a combination of Baire's Theorem and Stone's topological representation of Boolean algebras. This topic is discussed. The calculus of observations, which is mentioned above, leads to the notion of topological systems, as demonstrated by Vickers. This hints at an interplay of topology and order, since a topology is after all a complete Heyting algebra. Another important topic is the approximation of continuous functions by a given class of functions, like the polynomials on an interval, leading quickly to the Stone-Weierstraß Theorem on a compact topological space, a topic with a rich history. Finally, the relationship of pseudometric spaces to general topological spaces is reflected again, we introduce uniform spaces as a rich class of spaces which is more general than pseudometric spaces, but less general than their topological cousins. Here we find concepts like completeness or uniform continuity, which are formulated for metric spaces, but which cannot be realized in general topological ones. This gallery could be extended, for example, Polish spaces could be discussed here with considerable relish, but it seemed to be more adequate to discuss these spaces in the context of their measure theoretic use. 

\begin{center}
  \fbox{We assume throughout that the Axiom of Choice is valid.}
\end{center}


\Subsection{Defining Topologies} 
\label{sec:top-deftops}
\def\cal{\mathcal}
\def\mathit{}
Recall that a topology $\tau$ on a carrier set $X$ is a collection of subsets
which contains both $\emptyset$ and $X$, and which is closed under
finite intersections and arbitrary unions. The elements of $\tau$ are
called the \emph{open sets}. Usually a topology is not written down as
one set, but it is specified what an open set looks like. This is done
through a base or a subbase. Recall\MMP{Base, subbase} that a
\emph{\index{topology!base}base} $\beta$ for $\tau$ is a set of subsets of $\tau$
such that for any $x\in G$ there exists $B\in\beta$ with $x\in
B\subseteq G$. A subbase is a family of sets for which the finite
intersections form a base.

Not every family of subsets qualifies as a subbase or a base. We have the following characterization of a base.

\BeginProposition{when-is-a-base}
A family $\beta$ of sets is the base for a topology on $X = \bigcup\beta$ iff given $U, V\in\beta$ and $x\in U\cap V$, there exists $W\in \beta$ with $x\in W\subseteq U\cap V$, and if $X $.
\EndProposition

Kelley~\cite[p. 47]{Kelley} gives the following example: Put $X := \{0, 1, 2\}$, $A := \{0, 1\}$ and $B := \{1, 2\}$, then $\beta := \{X, A, B, \emptyset\}$ cannot be the base for a topology. Assume it is, then the topology must be $\beta$ itself, but $A\cap B\not\in \beta$. So we have to be a bit careful. Let us have a look at the proof.
 
\BeginProof
Checking the properties for a base shows that the condition is certainly necessary. Suppose that the condition holds, and define 
\begin{equation*}\textstyle
  \tau :=\{\bigcup \beta_{0}\mid \beta_{0}\subseteq\beta\}.
\end{equation*}
Then $\emptyset, X\in \tau$, and $\tau$ is closed under arbitrary unions, so that we have to check whether $\tau$ is closed under finite intersections. In fact, let $x\in U\cap V$ with $U, V\in \tau$, then we can find $U_{0}, V_{0}\in \beta$ with $x\in U_{0}\cap V_{0}$. By assumption there exists $W\in \beta$ with $x\in W\subseteq U_{0}\cap V_{0}\subseteq U\cap V$, so that $U\cap V$ can be written as union of elements in $\beta$. 
\EndProof


We perceive a base and a subbase, resp., relative to a topology,
but it is usually clear what the topology looks like, once a basis is
given. Let us have a look at some examples to clarify things.

\BeginExample{ex-topol-bases-real}
Consider the real numbers $\Real$ with the Euclidean topology
$\tau$. We say that a set $G$ is open iff given $x\in G$, there exists
an open interval $]a, b[$ with $x\in ]a, b[\ \subseteq G$. Hence the set
$\bigl\{]a, b[\mid a, b\in \Real, a < b\bigr\}$ forms a base for $\tau$;
actually, we could have chosen $a$ and $b$ as rational numbers, so
that we have even a countable base for $\tau$. Note that although we
can find a closed interval $[v, w]$ such that $x\in [v, w]\ \subseteq\ ]a,
b[\ \subseteq G$, we could not have used the closed intervals for a
description of $\tau$, since otherwise the singleton sets $\{x\} = [x,
x]$ would be open as well. This is both undesirable and counter
intuitive: in an open set we expect that each element has some
breathing space around it.
\EndExample

The next example looks at Euclidean spaces; here we do not have
intervals directly at our disposal, but we can measure distances as
well, which is a suitable generalization, given that the interval
$]x-r, x+r[$ equals $\{y\in \Real \mid |x-y|< r\}$.

\BeginExample{ex-topol-bases-eucl}
Consider the three dimensional space $\Real^{3}$, and define for $x,
y\in \Real^{3}$ their distance 
\begin{equation*}
   d(x, y) := \sum_{i=1}^{3}|x_{i}- y_{i}|.
\end{equation*}
Call $G\subseteq\Real^{3}$ open iff given $x\in G$, there exists $r >
0$ such that $\{y\in \Real^{3}\mid d(x, y) < r\}\subseteq G$. Then it
is clear that the set of all open sets form a topology:
\begin{itemize}
\item Both the empty set and $\Real^{3}$ are open.
\item The union of an arbitrary collection of open sets is open again.
\item Let $G_{1}, \dots, G_{k}$ be open, and $x\in G_{1}\cap\dots\cap
  G_{k}$. Take an index $i$; since $x\in G_{i}$, there exists
  $r_{i}>0$ such that $K(d, x, r) := \{y\in \Real^{3}\mid d(x, y) < r_{i}\}\subseteq
  G_{i}$. Let $r := \min\{r_{1}, \dots, r_{k}\}$, then 
  \begin{equation*}
    \{y\in \Real^{3}\mid d(x, y) < r\}=\bigcap_{i=1}^{k}\{y\in
    \Real^{3}\mid d(x, y) < r_{i}\}
\subseteq\bigcap_{i=1}^{k}G_{i}.
  \end{equation*}
Hence the intersection of a finite number of open sets is open again.
\end{itemize}
This argument would not work with a countable number of open sets, by
the way. 

We could have used other measures for the distance, e.g.,
\begin{align*}
  d'(x, y) & := \sqrt{\sum_{i}|x_{i}- y_{i}|^{2}},\\
d''(x, y) & := \max_{1\leq i \leq 3}|x_{i}- y_{i}|.
\end{align*}
Then it is not difficult to see that all three describe the same
collection of open sets. This is so because we can find for $x$ and
$r>0$ some $r'>0$ and $r'' > 0$ with $K(d', x, r') \subseteq K(d, x,
r)$ and $K(d'', x, r'') \subseteq K(d, x,
r)$, similarly for the other combinations. 

It is noted that $3$ is not a magical number here, we can safely
replace it with any positive $n$, indicating an arbitrary finite dimension. Hence we have
shown that $\Real^{n}$ is for each $n\in \Nat$ a topological space in
the Euclidean topology. 
\EndExample

The next example uses also some notion of distance between two
elements, which are given through evaluating real valued
functions. Think of $f(x)$ as the numerical value of attribute $f$ for object
$x$, then $|f(x) - f(y)|$ indicates how far apart $x$ and $y$ are with
respect to their attribute values.

\BeginExample{ex-topol-bases-weak}
Let $X$ be an arbitrary non-empty set, and ${\cal E}$ be a non-empty
collections of functions $f: X\to \Real$. Define for the finite
collection ${\cal F}\subseteq {\cal E}$, for
$r>0$, and for $x\in X$  the base set 
\begin{equation*}
  W_{{\cal F}; r} (x) := \{y\in X \mid 
  |f(x)-f(y)|<r\text{ for all }f\in {\cal F}\}.
\end{equation*}
We define as a base $\beta := \{W_{{\cal F}; r} (x) \mid x\in X, r >
0, {\cal F}\subseteq{\cal G}\text{ finite}\}$, and hence
call $G\subseteq X$ open iff given $x\in G$, there exists ${\cal
  F}\subseteq{\cal E}$ finite and $r>0$ such that $W_{{\cal F};
  r}(x)\subseteq G$.

It is immediate that the finite intersection of open sets is open
again. Since the other properties are checked easily as well, we have
defined a topology\MMP{Weak topology}, which is sometimes called the \emph{weak
  \index{topology!weak}topology} on $X$ induced by ${\cal E}$.

It is clear that in the last example the argument would not work if
we restrict ourselves to elements of ${\cal G}$ for defining the base,
i.e., to sets of the form $W_{\{g\}; r}$. These sets, however, have
the property that they form a subbase, since finite intersections
of these sets form a base. 
\EndExample

The next example shows that a topology may be defined on the set of
all partial functions from some set to another one. In contrast to the
previous example, we do without any numerical evaluations.

\BeginExample{partial-maps-top}
{
\def\partMap#1#2{\ensuremath{{#1}\rightharpoonup{#2}}}
Let $A$ and $B$ be non-empty sets, define $$\partMap{A}{B} :=
\{f\subseteq A\times B \mid f\text{ is a partial map}\}.$$ A set
$G\subseteq\partMap{A}{B}$ is called open iff given $f\in G$ there
exists a finite $f_{0}\in\partMap{A}{B}$ such that 
\begin{equation*}
  f \in N(f_{0}) := \{g\in \partMap{A}{B} \mid f_{0}\subseteq g\}
  \subseteq G.
\end{equation*}
Thus we can find for $f$ a finite partial map $f_{0}$ which is
extended by $f$ such that all extensions to $f_{0}$ are contained in
$G$. 

Then this is in fact a topology. The collection of open sets is
certainly closed under arbitrary unions, and both the empty set and
the whole set $\partMap{A}{B}$ are open. Let $G_{1}, \dots, G_{n}$ be
open, and $f\in G := G_{1}\cap\dots\cap G_{n}$, then we can find
finite partial maps $f_{1}, \dots, f_{n}$ which are extended by $f$
such that $N(f_{i})\subseteq G_{i}$ for $1 \leq i \leq n$. Since $f$
extends all these maps, $f_{0} := f_{1}\cup\dots\cup f_{n}$ is a well
defined finite partial map which is extended by $f$, and 
\begin{equation*}
  f\in N(f_{0}) = N(f_{1})\cap\dots\cup N(f_{n})\subseteq G.
\end{equation*}
Hence the finite intersection of open sets is open again. 

A base for this topology is the set $\{N(f)\mid f\text{ is finite}\}$,
a subbase is the set $\bigl\{N(\{\langle a, b\rangle\})\mid a\in A, b\in B\bigr\}$
}
\EndExample

The next example deals with a topology which is induced by an order
structure. Recall that a chain in a partially ordered set is a
non-empty totally ordered subset, and that in an inductively ordered
set each chain has an upper bound.

\BeginExample{scott-open}
Let $(P, \leq)$ be a inductively ordered set. Call $G\subseteq P$
\emph{Scott \index{open!Scott}open} iff 
\begin{enumerate}
\item $G$ is upward closed (hence $x\in G$ and $x\leq y$ imply
  $y\in G$).
\item If $S\subseteq P$ is a chain with $\sup S \in G$, then $S\cap
  G\not=\emptyset$. 
\end{enumerate}
Again, this defines a \index{topology!Scott}topology on $P$. In fact,
it is enough to show that $G_{1}\cap G_{2}$ is open, if $G_{1}$ and
$G_{2}$ are. Let $S$ be a chain with $\sup S\in G_{1}\cap G_{2}$, then
we find $s_{i}\in S$ with $s_{i}\in G_{i}$. Since $S$ is a chain, we
may and do assume that $s_{1}\leq s_{2}$, hence $s_{2}\in G_{1}$,
because $G_{1}$ is upward closed. Thus
$s_{2}\in S\cap (G_{1}\cap G_{2})$. Because $G_{1}$ and $G_{2}$ are
upward closed, so is $G_{1}\cap G_{2}$.

As an illustration, we show that the set
$F := \{x\in P\mid x \leq t\}$ is Scott closed for each $t\in P$. Put
$G := P\setminus F$. Let $x\in G$, and $x\leq y$, then obviously
$y \not\in F$, so $y\in G$. If $S$ is a chain with $\sup S\in G$, then
there exists $s\in S$ such that $s\not\in F$, hence
$S\cap G\not=\emptyset$.
\EndExample

\subsubsection{Continuous Functions}

A continuous map between topological spaces is compatible with the topological structure. This is familiar from real functions, but we cannot copy the definition, since we have no means of measuring the distance between points in a topological space. All we have is the notion of an open set. So the basic idea is to say that given an open neighborhood $U$ of the image, we want to be able to find an open neighborhood $V$ of the inverse image so that all element of $V$ are mapped to $U$. This is a direct translation of the familiar $\epsilon$-$\delta$-definition from calculus. Since we are concerned with continuity as a global concept (as opposed to one which focusses on a given point), we arrive at this definition, and show in the subsequent example the it is really a faithful translation.

\BeginDefinition{fnct-continuous}
Let $(X, \tau)$ and $(Y, \theta)$ be topological spaces. A map
$f: X\to Y$ is called 
$\tau$-$\theta$-\emph{\index{continuous}continuous} iff
$\InvBild{f}{H}\in\tau$ for all $H\in\theta$ holds, we write this also
as $f: (X, \tau)\to (Y, \theta)$. 
\EndDefinition

If the context is clear, we omit the reference to the topologies. Hence we say that the inverse image of an open set under a continuous map is an open set again.

Let us have a look at real functions.

\BeginExample{ex-cont-fncts-real}
Endow the reals with the Euclidean topology, and let
$f: \Real\to \Real$ be a map. Then the definition of continuity given
above coincides with the usual $\epsilon$-$\delta$-definition.

1.
Assuming the \MMP{$\epsilon$-$\delta$?}$\epsilon$-$\delta$-definition, we want to show that the
inverse image of an open set is open. In fact, let $G\subseteq\Real$
be open, and pick $x\in \InvBild{f}{G}$. Since $f(x)\in G$, we can
find $\epsilon>0$ such that $]f(x)-\epsilon,
f(x)+\epsilon[\ \subseteq G$. Pick $\delta>0$ for this $\epsilon$, hence
$x'\in\ ]x-\delta, x+\delta[$ implies $f(x')\in\ ]f(x)-\epsilon, f(x)+
\epsilon[\ \subseteq G$. Thus $x\in\ ]x-\delta,
x+\delta[\ \subseteq\InvBild{f}{G}$. 

2.  Assuming that the inverse image of an open set is open, we want to
establish the $\epsilon$-$\delta$-definition. Given $x\in \Real$, let
$\epsilon>0$ be arbitrary, we want to show that there exists
$\delta>0$ such that $|x-x'|<\delta$ implies
$|f(x)-f(x')|<\epsilon$. Now $]f(x)-\epsilon, f(x')+\epsilon[$ is an
open set hence $H := \InvBild{f}{]f(x)-\epsilon, f(x')+\epsilon[}$ is
open by assumption, and $x\in H$, Select $\delta>0$ with
$]x-\delta, x+\delta[\ \subseteq H$, then $|x-x'|<\delta$ implies
$x'\in H$, hence $f(x')\in\ ]f(x)-\epsilon, f(x')+\epsilon[$.
\EndExample

Thus we work on familiar ground, when it comes to the reals. 
Continuity may be tested on a subbase:

\BeginLemma{cont-for-subbase}
Let $(X, \tau)$ and $(Y, \theta)$ be topological spaces, $f: X\to Y$
be a map. Then $f$ is  $\tau$-$\theta$-continuous iff
$\InvBild{f}{S}\in \tau$ for each $S\in \sigma$ with
$\sigma\subseteq\theta$ a subbase. 
\EndLemma

\BeginProof
Clearly, the inverse image of a subbase element is open, whenever $f$
is continuous. Assume, conversely, that the $\InvBild{f}{S} \in \tau$
for each $S\in\sigma$. Then $\InvBild{f}{B}\in \tau$ for each element
$B$ of the base $\beta$ generated from $\sigma$, because $B$ is the
intersection of a finite number of subbase elements. Now, finally, if
$H\in\theta$, then $H = \bigcup\{B\mid B\in\beta, B\subseteq H\}$, so
that $\InvBild{f}{H} = \bigcup\{\InvBild{f}{B}\mid B\in\beta,
B\subseteq H\}\in\tau$. Thus the inverse image of an open set is open.
\EndProof

\BeginExample{partial-maps-top-cont} {
Take the
  topology from Example~\ref{partial-maps-top} on the space
  $\partMap{A}{B}$ of all partial maps. A map
  $q: (\partMap{A}{B})\to (\partMap{C}{D})$ is continuous in this topology
  iff the following condition holds: whenever $q(f)(c)=d$,
then there exists $f_{0}\subseteq f$ finite such that $q(f_{0})(c) =
d$. 

In fact, let $q$ be continuous, and $q(f)(c) = d$, then
$G := \InvBild{q}{N(\{\langle c, d\rangle\})}$ is open and contains
$f$, thus there exists $f_{0}\subseteq f$ with
$f\in N(f_{0})\subseteq G$, in particular $q(f_{0})(c) =
d$.
Conversely, assume that $H\subseteq \partMap{C}{D}$ is open, and we
want to show that $G := \InvBild{q}{H}\subseteq \partMap{A}{B}$ is
open. Let $f\in G$, thus $q(f)\in H$, hence there exists
$g_{0}\subseteq q(f)$ finite with $q(f)\in N(g_{0})\subseteq
H$.
$g_{0}$ is finite, say
$g_{0} = \{\langle c_{1}, d_{1}\rangle, \dots, \langle c_{n},
d_{n}\rangle\}$.
By assumption there exists $f_{0}\in\partMap{A}{B}$ with
$q(f_{0})(c_{i}) = d_{i}$ for $1\leq i\leq n$, then $f\in N(f_{0})
\subseteq G$, so that the latter set is open.
}
\EndExample

Let us have a look at the Scott topology.

\BeginExample{scott-continuous-map}
Let $(P, \leq)$ and $(Q, \leq)$ be inductively ordered sets, then $f:
P\to Q$ is Scott continuous (i.e., continuous, when both ordered sets
carry their respective Scott topology) iff $f$ is monotone, and if
$f(\sup S) = \sup\Bild{f}{S}$ holds for every chain $S$. 

Assume that $f$ is Scott continuous. If $x\leq x'$, then every open
set which contains $x$ also contains $x'$, so if $x\in\InvBild{f}{H}$
then $x'\in\InvBild{f}{H}$ for every Scott open $H\subseteq Q$; thus
$f$ is monotone. If $S\subseteq P$ is a chain, then $\sup S$ exists in
$P$, and $f(s)\leq f(\sup S)$ for all $s\in S$, so that
$\sup\Bild{f}{S}\leq f(\sup S)$. For the other inequality, assume that
$f(\sup S)\not\leq \sup\Bild{f}{S}$. We note that
$G := \InvBild{f}{\{q\in Q\mid q\not\leq \sup\Bild{f}{S}\}}$ is open
with $\sup S\in G$, hence there exists $s\in S$ with $s\in G$. But
this is impossible. On the other hand, assume that $H\subseteq Q$ is
Scott open, we want to show that $G := \InvBild{f}{H}\subseteq P$ is
Scott open. $G$ is upper closed, since $x\in G$ and $x\leq x'$ implies
$f(x)\in H$ and $f(x)\leq f(x')$, thus $f(x')\in H$, so that
$x'\in G$. Let $S\subseteq P$ be a chain with $\sup S\in G$, hence
$f(\sup S)\in H$. Since $\Bild{f}{S}$ is a chain, and
$f(\sup S) = \sup\Bild{f}{S}$, we infer that there exists $s\in S$
with $f(s)\in H$, hence there is $s\in S$ with $s\in G$. Thus $G$ is
Scott open in $P$, and $f$ is Scott continuous.
\EndExample

The interpretation of modal logics in a topological space is interesting, when we interpret the transition which is associated with the diamond operator through a continuous map; thus the next step of a transition is uniquely determined, and it depends continuously on its argument.

\BeginExample{interpr-modal-logics-top}
The syntax of our modal logics is given through
\begin{equation*}
  \phi ::= \top~\mid~p~\mid~\phi_{1}\vee \phi_{2}~\mid~\phi_{1}\wedge\phi_{2}~\mid~\neg \phi~\mid \Diamond\phi
\end{equation*} 
with $p\in\Phi$ an atomic proposition. The logic has the usual
operators, viz., disjunction and negation, and $\Diamond$ as the modal
operator. 

For interpreting the logic, we take a topological state space $(S, \tau)$
and a continuous map $f: X\to X$, and we  associate with each atomic
proposition $p$ an open set $V_{p}$ as the set of all states in which
$p$ is true. We want the validity set $\Gilt$ of all those states in which formula $\phi$ holds to be open, and define inductively the validity of a formula in a
state in the following way. 
\begin{align*}
\Gilt[\top] & := S\\
\Gilt[p] & := V_{p}, \text{ if $p$ is atomic}\\
\Gilt[\phi_{1}\vee\phi_{2}] & := \Gilt[\phi_{1}]\cup\Gilt[\phi_{2}]\\
\Gilt[\phi_{1}\wedge\phi_{2}] & := \Gilt[\phi_{1}]\cap\Gilt[\phi_{2}]\\
\Gilt[\neg \phi] & := \Interior{(S\setminus \Gilt)}\\
\Gilt[\Diamond \phi] & := \InvBild{f}{\Gilt}
\end{align*}
All definitions but the last two are self explanatory. The interpretation of $\Gilt[\Diamond\gamma]$ through $\InvBild{f}{\Gilt}$ suggests itself when considering the graph of $f$ in the usual interpretation of the diamond in modal logics, see~\CategCite{Sect. 2.7}.

Since we want $\Gilt[\neg\phi]$ be open, we cannot take the complement of $\Gilt$ and declare it as the validity set for $\phi$, because the complement of an open set is not necessarily open. Instead, we take the largest open set which is contained in $S\setminus\Gilt$ (this is the best we can do), and assign it to $\neg \phi$. One shows easily through induction on the structure of formula $\phi$ that $\Gilt$ is an open set.

But now look at this. Assume that $X := \Real$ in the usual topology,
$V_{p} = \Gilt[p] = ]0, +\infty[$, then
$\Gilt[\neg p] = \Interior{]-\infty, 0]} = ]-\infty, 0[$, thus
$\Gilt[p\vee\neg p] = ¸\Real\setminus\{0\} \not= \Gilt[\top]$. Thus
the law of the excluded middle does not hold in this model.
\EndExample

Returning to the general discussion, the following fundamental property is immediate.

\BeginProposition{cont-closed-under-compos}
The identity $(X, \tau)\to (X, \tau)$ is continuous, and continuous maps
are closed under composition. Consequently, topological spaces with
continuous maps form a category. \QED
\EndProposition

Continuous maps can be used to define topologies.

\BeginDefinition{initial-and-final-tops}
Given a family ${\cal F}$ of maps $f: A\to X_{f}$, where
$(X_{f}, \tau_{f})$ is a topological space for each $f\in {\cal F}$,
the \emph{initial \index{topology!initial}topology}
$\tau_{\mathit{in}
, {\cal F}}$ on $A$
with respect to ${\cal F}$ is the smallest topology on $A$ so that $f$
is $\tau_{\mathit{in}, {\cal F}}$-$\tau_{f}$-continuous for every $f\in{\cal F}$.  Dually,
given a family ${\cal G}$ of maps $g: X_{g}\to Z$, where
$(X_{g}, \tau_{g})$ is a topological space for each $g\in{\cal G}$,
the \emph{final \index{topology!final}topology} $\tau_{\mathit{fi}, {\cal G}}$ on $Z$ is
the largest topology on $Z$ so that $g$ is
$\tau$-$\tau_{\mathit{fi},{\cal G}}$-continuous for every $g\in{\cal G}$.
\EndDefinition

In the case of the initial topology for just one map $f: A\to X_{f}$, note that $\PowerSet{A}$ is a topology which renders $f$ continuous, so there exists in fact a smallest topology on $A$ with the desired property; because $\{\InvBild{f}{G}\mid G\in \tau_{f}\}$ is a topology that satisfies the requirement, and because each such topology must contain it, this is in fact the smallest one. If we have a family ${\cal F}$ of maps $A\to X_{f}$, then each topology making all $f\in{\cal F}$ continuous must contain $$\xi := \bigcup_{f\in {\cal F}}\{\InvBild{f}{G}\mid G\in \tau_{f}\},$$ so the initial topology with respect to ${\cal F}$ is just the smallest topology on $A$ containing $\xi$. Similarly, being the largest topology rendering each $g\in {\cal G}$ continuous, the final topology with respect to ${\cal G}$ must contain the set $\bigcup_{g\in{\cal G}}\{H\mid \InvBild{g}{H}\in\tau_{g}\}$. An easy characterization of the initial resp. the final
topology is proposed here:

\BeginProposition{initial-final}
Let $(Z, \tau)$ be a topological space, and ${\cal F}$ be a family of
maps $A\to X_{f}$ with $(X_{f}, \tau_{f})$ topological spaces; $A$ is
endowed with the initial topology $\tau_{\mathit{in}, {\cal F}}$ with
respect to ${\cal F}$. A map $h: Z\to A$ is
$\tau$-$\tau_{\mathit{in}, {\cal F}}$-continuous iff
$h\circ f: Z\to X_{f}$ is $\tau$-$\tau_{f}$-continuous for every
$f\in {\cal F}$.
\EndProposition

\BeginProof
1.
Certainly, if $h: Z\to A$ is $\tau$-$\tau_{\mathit{in}, {\cal F}}$ continuous,
then $h\circ f: Z\to X_{f}$ is $\tau$-$\tau_{f}$-continuous for every
$f\in {\cal F}$ by Proposition~\ref{cont-closed-under-compos}. 

2.
Assume, conversely, that $h\circ f$ is continuous for every $f\in{\cal
  F}$; we want to show that $h$ is continuous. Consider
\begin{equation*}
  \zeta := \{G\subseteq A\mid \InvBild{h}{G}\in\tau\}.
\end{equation*}
Because $\tau$ is a topology, $\zeta$ is; because $h\circ f$ is
continuous, $\zeta$ contains the sets $\{\InvBild{f}{H}\mid H\in
\tau_{f}\}$ for every $f\in{\cal F}$. But this implies that $\zeta$
contains $\tau_{\mathit{in}, {\cal F}}$, hence $\InvBild{h}{G}\in\tau$
for every $G\in\tau_{\mathit{in}, {\cal F}}$. This establishes the assertion. 
\EndProof

There is a dual characterization for the final topology, see Exercise~\ref{ex-char-final-top}. 

These are the most popular examples for initial and final topologies.
\begin{enumerate}
\item Given a family $(X_{i}, \tau_{i})_{i\in I}$ of topological spaces, let $X := \prod_{i\in I}X_{i}$ be the Cartesian product of the carrier sets\footnote{This works only if $X \not=\emptyset$, recall that we assume here that the Axiom of Choice is valid}. The\MMP{Product} \emph{product \index{topology!product}topology} $\prod_{i\in I}\tau_{i}$ is the initial topology on $X$ with respect to the projections $\pi_{i}: X\to X_{i}$. The product topology has as a base 
  \begin{equation*}\textstyle
    \{\prod_{i\in I}A_{i}\mid A_{i}\in \tau_{i}\text{ and }A_{i}\not= X_{i}\text{ only for finitely many indices}\}
  \end{equation*}
\item Let $(X, \tau)$ be a topological space, $A\subseteq X$. The
  \emph{\index{topology!trace}trace} $(A, \tau\cap A)$ of $\tau$ on
  $A$ is the initial topology on $A$ with respect to the embedding
  $i_{A}: A\to X$\MMP{Subspace}. It has the open sets $\{G\cap A\mid  G\in \tau\}$;
  this is sometimes called the \emph{subspace topology}~\SetCite{p. 40}. We do not assume that $A$ is open. 
\item Given the family of spaces as above, let $X := \sum_{i\in I}X_{i}$ be the direct sum. The \emph{sum \index{topology!sum}topology}\MMP{Sum} $\sum_{i\in I}\tau_{i}$ is the final topology on $X$ with respect to the injections $\iota_{i}: X_{i}\to X$. Its open sets are described through 
  \begin{equation*}
    \bigl\{\sum_{i\in I}\Bild{\iota_{i}}{G_{i}}\mid  G_{i}\in \tau_{i}\text{ for all } i\in I\bigr\}.
  \end{equation*}
 \item Let $\rho$ be an equivalence relation on $X$ with $\tau$ a
   topology on the base space.  The factor space $\Faktor{X}{\rho}$ is
   equipped with the final topology $\Faktor{\tau}{\rho}$ with respect
   to the factor map $\fMap{\rho}$ which sends each element to its
   $\rho$-class\MMP{Factor}. This topology is called the \emph{\index{topology!quotient}quotient topology} (with
   respect to $\tau$ and $\rho$). If a set $G\subseteq \Faktor{X}{\rho}$
   is open then its inverse image $\InvBild{\fMap{\rho}}{G}=\bigcup G
   \subseteq X$ is open in $X$. But the converse holds as well: assume
   that $\bigcup G$ is open in $X$ for some
   $G\subseteq\Faktor{X}{\rho}$, then $G = \Bild{\fMap{\rho}}{\bigcup
     G}$, and, because $\bigcup G$ is the union if equivalence
   classes, one shows that $\InvBild{\fMap{\rho}}{G} = \InvBild{\fMap{\rho}}{\Bild{\fMap{\rho}}{\bigcup
     G}} = \bigcup G$. But this means that $G$ is open in
 $\Faktor{X}{\rho}$. 
\end{enumerate}

Just to gain some familiarity with the concepts involved, we deal with
an induced map on a product space, and with the subspace coming from
the image of a map. The properties we find here will be useful later on as well.

The product space first. We will use that a map into a topological
product is continuous iff all its projections are; this follows from
the characterization of an initial topology. It goes like this.

\BeginLemma{into-unitcube}
Let $M$ and $N$ be non-empty sets, $f: M\to N$ be a map. Equip both
$[0, 1]^{M}$ and $[0, 1]^{N}$ with the product topology. Then 
\begin{equation*}
f^{*}:
  \begin{cases}
    [0, 1]^{N}& \to [0, 1]^{M}\\
g & \mapsto g\circ f
  \end{cases}
\end{equation*}
is continuous. 
\EndLemma

\BeginProof
Note the reversed order; we have $f^{*}(g)(m) = (g\circ f)(m) =
g(f(m))$ for $g\in [0, 1]^{N}$ and $m\in M$. 

Because $f^{*}$ maps $[0, 1]^{N}$ into $[0, 1]^{M}$, and the latter
space carries the initial topology with respect to the projections
$(\pi_{M, m})_{m\in N}$ with $\pi_{M, m}: q \mapsto q(m)$, it is by
Proposition~\ref{initial-final} sufficient to show that $\pi_{M,
  m}\circ f^{*}: [0, 1]^{N}\to [0, 1]$ is continuous for every $m\in
M$. But $\pi_{M, m}\circ f^{*} = \pi_{N, f(m)}$; this is a projection,
which is continuous by definition. Hence $f^{*}$ is continuous.
\EndProof

Hence an application of the projection defuses a seemingly complicated
map. Note in passing the neither $M$ nor $N$ are assumed to carry a
topology, they are simply plain sets. 

The next observation displays an example of a subspace topology. Each
continuous map $f: X\to Y$ of one topological space to another one
induces a subspace $\Bild{f}{X}$ of $Y$, which may or may not have
interesting properties. In the case considered, it inherits
compactness from its source. 

\BeginProposition{image-is-compact}
Let $(X, \tau)$ and $(Y, \theta)$ be topological spaces, $f: X\to Y$
be $\tau$-$\theta$-continuous. If $(X, \tau)$ is compact, so is
$(\Bild{f}{X}, \theta\cap\Bild{f}{X})$, the subspace of $(Y, \theta)$ induced by $f$. 
\EndProposition

\BeginProof
We take on open cover of $\Bild{f}{X}$ and show that it contains a
finite cover of this space. So let $(H_{i})_{i\in I}$ be an open cover
of $\Bild{f}{X}$. There exists open sets $H_{i}'\in\theta$ such
that $H_{i}' = H_{i}\cap\Bild{f}{X}$, since $(\Bild{f}{X},
\tau\cap\Bild{f}{X})$ carries the subspace topology. Then
$(\InvBild{f}{H_{i}'})_{i\in I}$ is an open cover of $X$, so there
exists a finite subset $J\subseteq I$ such that $X = \bigcup_{i\in
  J}\InvBild{f}{H_{i}'}$, since $X$ is compact. But then
$(H_{i}'\cap\Bild{f}{X})_{i\in J}$ is an open cover of
$\Bild{f}{X}$. Hence this space is compact. 
\EndProof

Before continuing, we introduce the notion of homeomorphism (as an
isomorphism in the category of topological spaces with continuous
maps).

\BeginDefinition{homeomorphism}
Let $X$ and $Y$ be topological spaces. A bijection $f: X\to Y$ is
called a \emph{\index{homeomorphism}homeomorphism} iff both $f$ and $f^{-1}$ are continuous.
\EndDefinition

It is clear that continuity and bijectivity alone do not make a
homeomorphism. Take as a trivial example the identity $(\Real,
\PowerSet{\Real})\to (\Real, \tau)$ with $\tau$ as the Euclidean
topology. It is continuos and bijective, but its inverse is not
continuous.

Let us have a look at some examples, first one for the quotient topology.

\BeginExample{example-quotient-top}
Let $U := [0, 2\cdot \pi]$, and identify the endpoints of the
interval, i.e., consider the equivalence relation
\begin{equation*}
  \rho := \{\langle x, x\rangle\mid x\in U\}\cup\{\langle 0, 2\cdot
  \pi\rangle, \langle 2\cdot\pi, 0\rangle\}.
\end{equation*}
Let $K := \Faktor{U}{\rho}$, and endow $K$ with the quotient
topology. 

A set $G\subseteq K$ is open iff $\InvBild{\fMap{\rho}}{G}\subseteq U$
is open, thus iff we can find an open set $H\subseteq\Real$ such that
$\InvBild{\fMap{\rho}}{G} = H\cap U$, since $U$ carries the trace of
$\Real$. Consequently, if $\Klasse{0}{\rho}\not\in G$, we find that 
$\InvBild{\fMap{\rho}}{G} = \{x\in U\mid \{x\}\in G\}$, which is open
by construction. If, however, $\Klasse{0}{\rho}\in G$, then
$\InvBild{\fMap{\rho}}{G} = \{x\in U\mid \{x\}\in G\}\cup\{0, 2\cdot
\pi\}$, which is open in $U$. 

We claim that $K$ and  the unit circle 
$ S := \{\langle s, t\rangle \mid 0\leq s, t\leq 1, s^{2}+t^{2} = 1\rangle\}$, are homeomorphic under the map
$\psi: \Klasse{x}{\rho}\mapsto \langle \sin x, \cos
x\rangle$. Because $\langle \sin 0, cos 0\rangle = \langle \sin 2\cdot
\pi, \cos 2\cdot \pi\rangle$, the map is well defined. Since we can write $S = \{\langle \sin x, \cos x\rangle\mid
0\leq x \leq 2\cdot \pi\}$, it is clear that $\psi$ is onto. The
topology on $S$ is inherited from the Cartesian plane, so open arcs
are a subbasis for it. Because the old Romans Sinus and Cosinus both are continuous, we find that
$\psi\circ \fMap{\rho}$ is continuous. We infer from
Exercise~\ref{ex-char-final-top} that $\psi$ is continous, since $K$
has the quotient topology, which is final. We want to show that
$\psi^{-1}$ is continuous. The argumentation is geometrical. Given an
open arc on $K$, we may describe it through $(P_{1}, P_{2})$ with a
clockwise movement. If the arc does not contain the critical point
$P := \langle0, 1\rangle$, we find an open interval $I := ]a, b[$ with $0 < a <
b < 2\cdot \pi$ such that $\Bild{\psi}{(P_{1}, P_{2})} =  \{\Klasse{x}{\rho}\mid x\in
  I\}$, which is open in $K$. If, however, $P$ is
on this arc, we decompose it into two parts $(P_{1}, P)\cup(P,
P_{2})$. Then $(P_{1}, P)$ is the image of some interval $]a, 2\cdot
\pi]$, and $(P, P_{2})$ is the image of an interval $[0, b[$, so that
$\Bild{\psi}{(P_{1}, P_{2})} = \Bild{\fMap{\rho}}{[0, b[\ \cup\ ]a,
    2\cdot \pi]}$, which is open in $K$ as well (note that $[0, b[$ as
  well as $]a, 2\cdot\pi]$ are open in $U$).

\EndExample

While we have described so far direct methods to describe a topology
by saying when a set is open, we turn now to an observation due to
Kuratowski which yields an indirect way. It describes axiomatically
what properties the closure of a set should have. Assume that we have
a \emph{\index{closure operator}closure operator}\MMP{Closure operator}, i.e., a map $A\mapsto \closOp{A}$ on the
powerset of a set $X$ with these properties:
\begin{enumerate}
\item $\closOp{\emptyset} = \emptyset$ and $\closOp{X} = X$.
\item $A\subseteq \closOp{A}$, and $\closOp{(A\cup B)} = \closOp{A}\cup\closOp{B}$.
\item $\closOp{(\closOp{A})} = \closOp{A}$.
\end{enumerate}
Thus the operator leaves the empty set and the whole set alone, the
closure of the union is the union of the closures, and the operator is
idempotent. One sees immediately that the operator which assigns to
each set its closure with respect to a given topology is such a
closure operator. It is also quite evident that the closure operator
is monotone. Assume that $A\subseteq B$, then $B = A \cup (B\setminus
A)$, so that $\closOp{B} = \closOp{A}\cup\closOp{(B\setminus
  A)}\supseteq\closOp{A}$.

\BeginExample{finite-ordered-for-closure}
Let $(D, \leq)$ be a finite partially ordered set. We put
$\closOp{\emptyset} := \emptyset$ and $\closOp{D} := D$, moreover,
\begin{equation*}
  \closOp{\{x\}} := \{y\in D \mid y \leq x\}
\end{equation*}
is defined for $x\in D$, and $\closOp{A} :=
\bigcup_{x\in X}\closOp{\{x\}}$ for subsets $A$ of $D$. Then this is a
closure operator. It is enough to check whether
$\closOp{(\closOp{\{x\}})} = \closOp{\{x\}}$ holds. In fact, we have
\begin{align*}
  z \in \closOp{(\closOp{\{x\}})} 
& \Leftrightarrow
z\in \closOp{\{y\}}\text{ for some }y\in \closOp{\{x\}}\\
& \Leftrightarrow
\text{there exists }y\leq x\text{ with }z\leq y\\
& \Leftrightarrow
z\leq x\\
& \Leftrightarrow
z\in\closOp{\{x\}}.
\end{align*}
Thus we associate with each finite partially ordered set a closure
operator, which assigns to each $A\subseteq D$ its down set. The map
$x\mapsto\closOp{\{x\}}$ embeds $D$ into a distributive lattice, see
the discussion in~\SetCite{Example 1.72}.
\EndExample

We will show now that we can obtain a topology by calling open all those sets the complements of which are fixed under the closure operator; in addition, it turns out that the topological closure and the one from the closure operator are the same.

\BeginTheorem{kuratowski-closure}
Let $\closOp{\cdot}$ be a closure operator. Then 
\begin{enumerate}
\item The set $\tau := \{X\setminus F \mid F\subseteq X, \closOp{F} = F\}$ is a topology.
\item For each set $\Closure{A} = \closOp{A}$ with $\Closure{\cdot}$ as the closure in $\tau$. 
\end{enumerate}
\EndTheorem

\BeginProof
1.  For establishing that $\tau$ is a topology, it is enough to show
that $\tau$ is closed under arbitrary unions, since the other
properties are evident. Let ${\cal G}\subseteq \tau$, and put $G :=
\bigcup {\cal G}$, so we want to know whether $\closOp{X\setminus G} =
X\setminus G$. If $H\in {\cal G}$, then $X\setminus G \subseteq
X\setminus H$, so $\closOp{(X\setminus G)}\subseteq
\closOp{(X\setminus H)} = X\setminus H$, thus $\closOp{(X\setminus G)}
\subseteq X\setminus G$. Since the operator is monotone, it follows
that $\closOp{(X\setminus G)} = X\setminus G$, hence $\tau$ is in fact
closed under arbitrary unions, hence it is a topology.

2.
Given $A\subseteq X$, 
\begin{equation*}\textstyle
 \Closure{A} = \bigcap\{F\subseteq X \mid F\text{ is closed, and }A\subseteq F\},
\end{equation*}
and $\closOp{A}$ takes part in the intersection, so that $\Closure{A}\subseteq \closOp{A}$. On the other hand, $A\subseteq \Closure{A}$, thus $\closOp{A}\subseteq \closOp{(\Closure{A})} = \Closure{A}$ by part 1. Consequently, $\Closure{A}$ and $\closOp{A}$ are the same.
\EndProof

It is on first sight a bit surprising that a topology can be described by finitary means, although arbitrary unions are involved. But we should not forget that we have also the subset relation at our disposal. Nevertheless, a rest of surprise remains. 

\Subsubsection{Neighborhood Filters}
\label{sec:nghb-filters}
The last method for describing a topology we are discussing here deals also with some order properties. Assume that we assign to each $x\in X$, where $X$ is a given carrier set, a filter $\upsilon(x)\subseteq \PowerSet{X}$ with the property that $x\in U$ holds for each $U\in \upsilon(x)$. Thus $\upsilon(x)$ has these properties:
\begin{enumerate}
\item $x\in U$ for all $U\in \upsilon(x)$.
\item If $U, V\in \upsilon(x)$, then $U\cap V\in\upsilon(x)$.
\item If $U\in\upsilon(x)$ and $U\subseteq V$, then $V\in\upsilon(x)$. 
\end{enumerate}
It is fairly clear that, given a topology $\tau$ on $X$, the \emph{neighborhood \index{filter!neighborhood}filter}\MMP{$\upsilon_{\tau}(x)$}
\begin{equation*}
  \upsilon_{\tau}(x) := \{V\subseteq X \mid \text{ there exists $U\in\tau$ with $x\in U$ and $U\subseteq V$}\}
\end{equation*}
for $x$ has these properties. It has also an additional property, which we will discuss shortly~---~for dramaturgical reasons.

Such a system of special filters defines a topology. We declare all
those sets as open which belong to the neighborhoods of their
elements. So if we take all balls in Euclidean $\Real^{3}$ as the
basis for a filter and assign each point the balls which it centers,
then the sphere of radius $1$ around the origin would not be open
(intuitively, it does not contain an open ball). So this appears to be
an appealing idea. In fact:

\BeginProposition{def-through-nbh-filters}
Let $\{\upsilon(x) \mid x\in X\}$ be a family of filters such that $x\in U$ for all $U\in \upsilon(x)$. Then 
\begin{equation*}
  \tau := \{U\subseteq X \mid U\in\upsilon(x)\text{ whenever }x\in U\}
\end{equation*}
defines a topology on $X$. 
\EndProposition

\BeginProof
We have to establish that $\tau$ is closed under finite intersections,
since the other properties are fairly straightforward. Now, let $U$
and $V$ be open, and take $x\in U\cap V$. We know that $U\in
\upsilon(x)$, since $U$ is open, and we have $V\in\upsilon(x)$ for the
same reason. Since $\upsilon(x)$ is a filter, it is closed under
finite intersections, hence $U\cap V\in\upsilon(x)$, thus $U\cap V$ is
open.
\EndProof

We cannot, however, be sure that the neighborhood filter
$\upsilon_{\tau}(x)$ for this new topology is the same as the given
one. Intuitively, the reason is that we do not know if we can find for
$U\in\upsilon(x)$ an open $V\in\upsilon(x)$ with $V\subseteq U$ such
that $V\in\upsilon(y)$ for all $y\in V$. To illustrate, look at
$\Real^{3}$, and take the neighborhood filter for, say, $0$ in the
Euclidean topology. Put for simplicity
\begin{equation*}
  \|x\| := \sqrt{x_{1}^{2}+x_{2}^{2}+x_{3}^{2}}.
\end{equation*}
Let $U\in\upsilon(0)$, then we can find an open
ball $V\in\upsilon(0)$ with $V\subseteq U$. In fact, assume $U = \{a\mid
\|a\|< q\}.$ Take $z\in U$, then we
can find $r>0$ such that the ball $V := \{y\mid
\|y-z\|<r\}$
is entirely contained in $U$ (select $\|z\|< r < q$), thus
$V\in\upsilon(0)$. Now let $y\in V$, let $0 < t < r-\|z-y\|$, then
$\{a\mid \|a-y\|<t\}\subseteq V$, since $\|a-z\|\leq
\|a-y\|+\|z-y\|<r$. Hence 
$U\in\upsilon(y)$ for all $y\in V$.

We obtain now as a simple corollary

\BeginCorollary{cor-def-through-nbh-filters}
Let $\{\upsilon(x)\mid x\in X\}$ be a family of filters such
that $x\in U$ for all $U\in \upsilon(x)$, and assume that for any
$U\in\upsilon(x)$ there exists $V\in\upsilon(x)$ with $V\subseteq U$
and $U\in\upsilon(y)$ for all $y\in V$. Then
$\{\upsilon(x)\mid x\in X\}$ coincides with the neighborhood
filter for the topology defined by this family.  \QED
\EndCorollary

In what follows, unless otherwise stated, \index{$\upsilon(x)$}$\upsilon(x)$ will denote
the neighborhood filter of a point $x$ in a topological space $X$.

\BeginExample{simple-ordered-set}
Let $L := \{1, 2, 3, 6\}$ be the set of all divisors of $6$, and
define $\isEquiv{x}{y}{\leq}$ iff $x$ divides $y$, so that we obtain
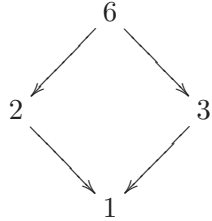
\begin{equation*}
\xymatrix{
&6\ar[dl]\ar[dr]\\
2\ar[dr]&&3\ar[dl]\\
&1
}
\end{equation*}
Let us compute ---~just for fun~--- the topology associated with this
partial order, and a basis for the neighborhood filters for each element. The
topology can be seen from the table below (we have used that $\Interior{A}
= X\setminus \Closure{(X\setminus A)}$, see~\SetCite{Definition 1.92}):

\begin{center}
  \begin{tabular}{|l|l|l|}\hline
set & closure & interior\\\hline\hline
    $\{1\}$& $\{1\}$& $\emptyset$\\\hline
    $\{2\}$& $\{1, 2\}$&$\emptyset$\\\hline
    $\{3\}$&$\{1, 3\}$&$\emptyset$\\\hline
    $\{6\}$&$\{1, 2, 3, 6\}$&$\{6\}$\\\hline 
    $\{1, 2\}$&$\{1, 2\}$&$\emptyset$\\\hline 
    $\{1, 3\}$&$\{1, 3\}$&$\emptyset$\\\hline
    $\{1, 6\}$&$\{1, 2, 3, 6\}$&$\{6\}$\\\hline
    $\{2, 3\}$&$\{1, 2, 3, 5\}$&$\emptyset$\\\hline
    $\{2, 6\}$&$\{1, 2, 3, 6\}$&$\{2, 6\}$\\\hline
    $\{3, 6\}$&$\{1, 2, 3, 6\}$&$\emptyset$\\\hline
    $\{1, 2, 3\}$&$\{1, 2, 3\}$&$\emptyset$\\\hline
    $\{1, 2, 6\}$&$\{1, 2, 3, 6\}$&$\{2, 6\}$\\\hline
    $\{1, 3, 6\}$&$\{1, 2, 3, 6\}$&$\{3, 6\}$\\\hline
    $\{2, 3, 6\}$&$\{1, 2, 3, 6\}$&$\{2, 3, 6\}$\\\hline
    $\{1, 2, 3, 6\}$&$\{1, 2, 3, 6\}$&$\{1, 2, 3, 6\}$\\\hline
  \end{tabular}
\end{center}
This is the topology:
\begin{equation*}
  \tau = \bigl\{\emptyset, \{6\}, \{2, 6\}, \{3, 6\}, \{2, 3, 6\}, \{1, 2, 3, 6\}\bigr\}.
\end{equation*}
A basis for the respective neighborhood filters is given in this table:
\begin{center}
  \begin{tabular}{|c|l|}\hline
element & basis\\\hline\hline
1 & $\bigl\{\{1, 2, 3, 6\}\bigr\}$\\\hline
2 &  $\bigl\{\{2, 6\}, \{1, 2, 3, 6\}\bigr\}$\\\hline
3 &  $\bigl\{\{3, 6\}, \{2, 3, 6\}, \{1, 2, 3, 6\}\bigr\}$\\\hline
6 &  $\bigl\{\{6\}, \{2, 6\}, \{3, 6\}, \{2, 3, 6\}, \{1, 2, 3, 6\}\bigr\}$\\\hline
  \end{tabular}
\end{center}
\EndExample

The next example deals with topological groups, i.e., topological spaces which have also a group structure rendering multiplication continuous. Here the neighborhood structure is fairly uniform --- if you know the neighborhood filter of the neutral element, you know the neighborhood filter of each element, because you can obtain them by a left shift or a right shift. 

\BeginExample{top-group}
Let $(G, \cdot)$ be a group, and $\tau$ be a topology on $G$ such that the map $\langle x, y\rangle\mapsto xy^{-1}$ is continuous. Then $(G, \cdot, \tau)$ is called a \emph{\index{topology!topological group}topological \index{group!topological}group}.  We will write down a topological group as $G$, the group operations and the topology will  not be mentioned. The neutral element is denoted by $e$, multiplication will usually be omitted. Given a subsets $U$ of $G$, define $gU := \{gh\mid h\in U\}$ and $Ug := \{hg\mid h\in U\}$ for $g\in G$.  

Let us look at the algebraic operations in a group. Put $\zeta(x, y) := xy^{-1}$, then the map $\xi: g \mapsto g^{-1}$ which maps each group element to its inverse is just $\zeta(e,g)$, hence the cut of a continuous map, to it is continous as well. $\xi$ is a bijection with $\xi\circ \xi = id_{G}$, so it is in fact a homeomorphism. We obtain multiplication as $xy = \zeta(x, \xi(y))$, so multiplication is also continuous. Fix $g\in G$, then multiplication $\lambda_{g}: x\mapsto gx$ from the left and $\rho_{g}: x\mapsto xg$ from the right are continuous. Now both $\lambda_{g}$ and $\rho_{g}$ are bijections, and $\lambda_{g}\circ \lambda_{g^{-1}} = \lambda_{g^{-1}}\circ \lambda_{g} = id_{G}$, also $\rho_{g}\circ \rho_{g^{-1}} = \rho_{g^{-1}}\circ \rho_{g} = id_{G}$, thus $\lambda_{g}$ and $\rho_{g}$ are homeomophisms for every $g\in G$. 

Thus we have in a topological group this characterization of the neighborhood filter for every $g\in G$: 
\begin{equation*}
  \upsilon(g)   = \{gU\mid U\in \upsilon(e)\}  = \{Ug\mid U\in\upsilon(e)\}.
\end{equation*}
In fact, let $U$ be a neighborhood of $g$, then $\InvBild{\lambda_{g}}{U} = g^{-1}U$ is a neighborhood of $e$, so is $\InvBild{\rho_{g}}{U} = Ug^{-1}$. Conversely, a neighborhood $V$ of $e$ determines a neighborhood $\InvBild{\lambda_{g^{-1}}}{V} = gV$ resp. $\InvBild{\rho_{g^{-1}}}{V} = Vg$ of $g$. 
\EndExample


\Subsection{Filters and Convergence} 
\label{sec:filters-and-convergence}
  The relationship between topologies and filters turns out to be
  fairly tight, as we saw when discussing the neighborhood filter of a
  point. We saw also that we can actually grow a topology from a
  suitable family of neighborhood filters. This relationship is even
 closer, as we will discuss now when having a look at
  convergence.

  Let $(x_{n})_{n\in \Nat}$ be a sequence in $\Real$ which converges to $x\in\Real$. This means that for any given open neighborhood $U$ of $x$ there exists an index $n\in \Nat$ such that $\{x_{m}\mid m\geq n\}\subseteq U$, so all members of the sequence having an index larger that $n$ are members of $U$. Now consider the filter $\fiF$ generated by the set $\{\{x_{m}\mid m\geq n\}\mid n\in \Nat\}$ of tails. The condition above says exactly that $\upsilon(x)\subseteq \fiF$, if you think a bit about it. This leads to the definition of convergence in terms of filters.

\BeginDefinition{def-convergence-filter}
Let $X$ be a topological space, $\fiF$ a filter on $X$. Then
$\fiF$ \emph{\index{convergence!filter}converges} to a limit $x\in X$ iff
$\upsilon(x)\subseteq \fiF$\MMP{$\fiF\to x$}. This is denoted by \index{${\cal
    F}\to x$}$\fiF\to x$.
\EndDefinition

Plainly, $\upsilon(x)\to x$ for every $x$. Note that the definition above does not force the limit to be uniquely determined. If if two different points $x, y$ share their neighborhood filter, then $\fiF\to x$ iff $\fiF\to y$. Look again at Example~\ref{simple-ordered-set}. There all neighborhood filters are contained in $\upsilon(6)$, so that we have $\upsilon(6)\to t$ for $t\in\{1, 2, 3, 6\}$. It may seem that the definition of convergence through a filter is too involved (after all, being a filter should not be taken on a light shoulder!). In fact, sometimes convergence is defined through a \emph{\index{net}net}\MMP{Net}\label{def-net} as follows. Let $(I, \leq)$ be a directed set, i.e., $\leq$ is a partial order such that, given $i, j\in I$ there exists $k$ with $i\leq k$ and $j\leq k$. An $I$-indexed family $(x_{i})_{i\in I}$ is said to converge\index{convergence!net}\index{net!convergence} to a point $x$ iff, given a neighborhood $U\in\upsilon(x)$ there exists $k\in I$ such that $x_{i}\in U$ for all $i\geq k$. This generalizes the concept of convergence from sequences to index sets of arbitrary size. But look at this. The sets $\bigl\{\{x_{j}\mid j\geq i\}\mid i\in I\bigr\}$ form a filter base, because $(I, \leq)$ is directed. The corresponding filter converges to $x$ iff the net converges to $x$.

But what about the converse? Take a filter $\fiF$ on $X$, then $F_{1}\leq F_{2}$ iff $F_{2}\subseteq F_{1}$ renders $(\fiF, \leq)$ a net. In fact, given $F_{1}, F_{2}\in\fiF$, we have $F_{1}\leq F_{1}\cap F_{2}$ and $F_{2}\leq F_{1}\cap F_{2}$. Now pick $x_{F}\in F$. Then the net $(x_{F})_{F\in\fiF}$ converges to $x$ iff $\fiF\to x$. Assume that $\fiF\to x$; take $U\in\upsilon(x)$, then $U\in\fiF$, thus if $F\in\fiF$ with $F\geq U$, then $F\subseteq U$, hence $x_{F}\in U$ for all such $x_{F}$. Conversely, if each net $(x_{F})_{F\in\fiF}$ derived from $\fiF$ converges to $x$, then for a given $U\in\upsilon(x)$ there exists $F_{0}$ such that $x_{F}\in U$ for $F\subseteq F_{0}$. Since $x_{F}$ has been chosen arbitrarily from $F$, this can only hold if $F\subseteq U$ for $F\subseteq F_{0}$, so that $U\in\fiF$. Because $U\in\upsilon(x)$ was arbitrary, we conclude $\upsilon(x)\subseteq\fiF$.

Hence we find that filters offer a uniform
generalization.

The argument above shows that we may select the elements $x_{F}$ from a
base for $\fiF$. If the filter has a countable base, we construct in this way a
sequence; conversely, the filter constructed from a sequence has a
countable base. Thus the convergence of sequences and the convergence
of filters with a countable base are equivalent concepts.
 
We investigate the characterization of the topological closure in
terms of filters. In order to do this, we need to be able to restrict
a filter to a set\MMP{Trace}, i.e., looking at the footstep the filter leaves on
the set, hence at $$\fiF\cap A := \{F\cap A\mid F\in\fiF\}.$$ This is
what we will do now.

\BeginLemma{localize-filter}
Let $X$ be a set, and $\fiF$ be a filter on $X$. 
Then $\fiF\cap A$ is a filter on $A$ iff $F\cap A\not=\emptyset$ for
all $F\in\fiF$.
\EndLemma

\BeginProof
Since a filter must not contain the empty set, the condition is
necessary. But it is also sufficient, because it makes sure that the
laws of a filter are satisfied.
\EndProof

Looking at $\fiF\cap A$ for an ultrafilter $\fiF$, we know that either
$A\in\fiF$ or $X\setminus A\in \fiF$, so if $F\cap
A\not=\emptyset$ holds for all $F\in \fiF$, then this implies that $A\in\fiF$. Thus we obtain

\BeginCorollary{localize-ultrafilter}
Let $X$ be a set, and $\fiF$ be an ultrafilter on $X$. Then $\fiF\cap
A$ is a filter iff $A\in\fiF$. Moreover, in this case $\fiF\cap A$ is
an ultrafilter on $A$.
\EndCorollary

\BeginProof
It remains to show that $\fiF\cap A$ is an ultrafilter on $A$,
provided, $\fiF\cap A$ is a filter. Let $B\not\in\fiF\cap A$ for some
subset $B\subseteq A$. Since $A\in\fiF$, we conclude $B\not\in\fiF$,
thus $X\setminus B\in\fiF$, since $\fiF$ is an ultrafilter. Thus
$(X\setminus B)\cap A = A\setminus B\in\fiF\cap A$, so $\fiF\cap A$ is
an ultrafilter by~\SetCite{Lemma 1.63}.  
\EndProof

From Lemma~\ref{localize-filter} we obtain a simple and elegant
characterization of the topological closure of a set.

\BeginProposition{char-top-closure}
Let $X$ be a topological space, $A\subseteq X$. Then $x\in\Closure{A}$
iff $\upsilon(x)\cap A$ is a filter on $A$. Thus $x\in \Closure{A}$
iff there exists a filter $\fiF$ on $A$ with $\fiF\to x$.
\EndProposition

\BeginProof
We know from the definition of $\Closure{A}$ that $x\in \Closure{A}$
iff $U\cap A\not=\emptyset$ for all $U\in\upsilon(x)$. This is by
Lemma~\ref{localize-filter} equivalent to $\upsilon(x)\cap A$ being a
filter on $A$.
\EndProof

We know from Calculus that continuous functions preserve convergence,
i.e., if $x_{n}\to x$ and $f$ is continuous, then $f(x_{n})\to
f(x)$. We want to carry this over to the world of filters. For this,
we have to define the image of a filter. Let $\fiF$ be a filter on a
set $X$, and $f: X\to Y$ a map, then $$f(\fiF) := \{B\subseteq Y\mid
\InvBild{f}{B}\in \fiF\}$$ is a filter on $Y$. In fact,\MMP[b]{Image of a filter}
$\emptyset\not\in f(\fiF)$, and, since $f^{-1}$ preserves the Boolean
operations, $f(\fiF)$ is closed under finite intersections. Let $B\in
f(\fiF)$ and $B\subseteq B'$. Since $\InvBild{f}{B}\in\fiF$, and
$\InvBild{f}{B}\subseteq \InvBild{f}{B'}$, we conclude
$\InvBild{f}{B'}\in \fiF$, so that $B'\in f(\fiF)$. Hence $f(\fiF)$ is
also upper closed, so that it is in fact a filter.

This is an easy representation through the direct image.

\BeginLemma{direct-image-filter}
Let $f: X\to Y$ be a map, $\fiF$ a filter on $X$, then $f(\fiF)$
equals the filter generated by $\{\Bild{f}{A}\mid A\in\fiF\}$.
\EndLemma

\BeginProof
Because
$\Bild{f}{A_{1}\cap A_{2}}\subseteq \Bild{f}{A_{1}}\cap
\Bild{f}{A_{2}}$, the set ${\cal G}_{0} := \{\Bild{f}{A}\mid A\in\fiF\}$ is a filter base. Denote by ${\cal G}$ the filter generated by ${\cal G}_{0}$. We claim that $f(\fiF) = {\cal G}$.

``$\subseteq$'': Assume that $B\in f(\fiF)$, hence $\InvBild{f}{B}\in \fiF$ . Since $\Bild{f}{\InvBild{f}{B}} \subseteq B$, we conclude that $B$ is contained in the filter generated by ${\cal G}_{0}$, hence in ${\cal G}$.

``$\supseteq$'': If $B\in{\cal G}_{0}$, we find $A\in\fiF$ with $B = \Bild{f}{A}$, hence $A \subseteq \InvBild{f}{\Bild{f}{A}} = \InvBild{f}{B}\in\fiF$, so that $B\in f(\fiF)$. This implies the desired inclusion, since $f(\fiF)$ is a filter. 
\EndProof

We will see now that not only the filter property is transported
through maps, but also the property of being an ultrafilter.

\BeginLemma{ultrafilter-remains-under-map}
Let $f: X\to Y$ be a map, $\fiF$ an ultrafilter on $X$. Then $f(\fiF)$
is an ultrafilter on $Y$.
\EndLemma

\BeginProof
It it enough to show that if $f(\fiF)$ does not contain a set, it will
contain its complement~\SetCite{Lemma 1.63}. In fact, assume that
$H\not\in f(\fiF)$, so that $\InvBild{f}{H}\not\in\fiF$. Since $\fiF$
is an ultrafilter, we know that $X\setminus\InvBild{f}{H}\in\fiF$; but
$X\setminus\InvBild{f}{H}=\InvBild{f}{Y\setminus H}$, so that
$Y\setminus H\in f(\fiF)$. 
\EndProof

\BeginExample{filter-for-product}
Let $X$ be the product of the topological spaces $(X_{i})_{i\in I}$
with projections $\pi_{i}: X\to X_{i}$. For a filter $\fiF$ on $X$, we
have
$\pi_{j}(\fiF) = \{A_{j}\subseteq X_{j}\mid
A_{j}\times\prod_{i\not=j}X_{i}\in\fiF\}.$
\EndExample

Continuity preserves convergence:

\BeginProposition{cont-pres-conv}
Let $X$ and $Y$ be topological spaces, and $f: X\to Y$ a map.
\begin{enumerate}
\item If $f$ is continuous, and $\fiF$ a filter on $X$, then $\fiF\to x$ implies $f(\fiF)\to
  f(x)$ for all $x\in X$.
\item If $\fiF\to x$
implies $f(\fiF)\to f(x)$ for all $x\in X$ and all filters $\fiF$ on $X$, then $f$ is continuous.
\end{enumerate}
\EndProposition

\BeginProof
Let $V\in \upsilon(f(x))$, then there exists $U\in\upsilon(f(x))$ open
with $U\subseteq V$. Since
$\InvBild{f}{U}\in\upsilon(x)\subseteq\fiF$, we conclude $U\in
f(\fiF)$, hence $V\in f(\fiF)$. Thus $\upsilon(f(x))\subseteq
f(\fiF)$, which means that $f(\fiF)\to f(x)$ indeed. This establishes
the first part.

Now assume that $\fiF\to x$
implies $f(\fiF)\to f(x)$ for all $x\in X$ and an arbitrary filter
$\fiF$ on $X$. Let $V\subseteq Y$ be
open. Given $x\in \InvBild{f}{V}$, we find an open set $U$ with $x\in
U\subseteq\InvBild{f}{V}$ in the following way. Because $x\in
\InvBild{f}{V}$, we know $f(x)\in V$. Since $\upsilon(x)\to x$, we
obtain from the assumption that $f(\upsilon(x))\to f(x)$, thus
$\upsilon(f(x))\subseteq f(\upsilon(x))$. Because
$V\in\upsilon(f(x))$, it follows $\InvBild{f}{V}\in\upsilon(x)$, hence
we find an open set $U$ with $x\in U\subseteq
\InvBild{f}{V}$. Consequently, $\InvBild{f}{V}$ is open in $X$. 
\EndProof

Thus continuity and filters cooperate in a friendly manner.

\BeginProposition{initial-topology}
Assume that $X$ carries the initial topology with respect to a family
$(f_{i}: X\to X_{i})_{i\in I}$ of functions. Then $\fiF\to x$ iff
$f_{i}(\fiF)\to f_{i}(x)$ for all $i\in I$.
\EndProposition

\BeginProof
Proposition~\ref{cont-pres-conv} shows that the condition is
necessary. Assume that $f_{i}(\fiF)\to f_{i}(x)$ for every $i\in I$,
let $\tau_{i}$ be the topology on $X_{i}$. The
sets $$\bigl\{\{\InvBild{f_{i_{1}}}{G_{i_{1}}}\cap\dots\cap\InvBild{f_{i_{k}}}{G_{i_{k}}}\}\mid
i_{1}, \dots, i_{k}\in I, f_{i_{1}}(x)\in G_{i_{1}}\in\tau_{i_{1}},
\dots, f_{i_{k}}(x)\in G_{i_{k}}\in \tau_{i_{k}}, k\in \Nat\bigr\}$$
form a base for the neighborhood filter for $x$ in the initial
topology. Thus, given an open neighborhood $U$ of $x$, we have
$\InvBild{f_{i_{1}}}{G_{i_{1}}}\cap\dots\cap\InvBild{f_{i_{k}}}{G_{i_{k}}}\subseteq
U$ for some suitable finite set of indices. Since $f_{i_{j}}(\fiF)\to
f_{i_{j}}(x)$, we infer $G_{i_{j}}\in f_{i_{j}}(\fiF)$, hence
$\InvBild{f_{i_{j}}}{G_{i_{j}}}\in\fiF$ for $1 \leq j\leq k$, thus
$U\in\fiF$. This means $\upsilon(x)\subseteq\fiF$. Hence $\fiF\to x$, as asserted.
\EndProof

We know that in a product a sequence converges iff its components converge. This is the counterpart for filters:

\BeginCorollary{conv-in-a-product}
Let $X = \prod_{i\in I} X_{i}$ be the product of the topological spaces. Then $\fiF\to (x_{i})_{i\in I}$ in $X$ iff $\fiF_{i}\to x_{i}$ in $X_{i}$ for all $i\in I$, where $\fiF_{i}$ it the $i$-th projection $\pi_{i}(\fiF)$ of $\fiF$. 
\QED
\EndCorollary

The next observation further tightens the connection between topological properties and filters. It requires the existence of ultrafilters, so recall the we assume that the Axiom of Choice holds.

\BeginTheorem{conv-vs-ultrafilter} Let $X$ be a topological space. Then $X$ is compact iff each ultrafilter converges.  \EndTheorem

Thus we tie compactness, i.e., the possibility to extract from each
cover a finite subcover, to the convergence of ultrafilters. Hence an
ultrafilter in a compact space cannot but converge. The proof of
Alexander's Subbase Theorem~\SetCite{Theorem 1.99} indicates already
that there is a fairly close connection between the Axiom of Choice
and topological compactness. This connection is tightened here.

\BeginProof
1.  Assume that $X$ is compact, but that we find an ultrafilter $\fiF$
which fails to converge. Hence we can find for each $x\in X$ an open
neighborhood $U_{x}$ of $x$ which is not contained in $\fiF$. Since
$\fiF$ is an ultrafilter, $X\setminus U_{x}\in \fiF$. Thus
$\{X\setminus U_{x}\mid x\in X\}\subseteq\fiF$ is a collection of
closed sets with $\bigcap_{x\in X}(X\setminus U_{x}) =
\emptyset$. Since $X$ is compact, we find a finite subset $F\subseteq
X$ such that $\bigcap_{x\in F}(X\setminus U_{x}) = \emptyset$. But
$X\setminus U_{x}\in \fiF$, and $\fiF$ is closed under finite
intersections, hence $\emptyset\in\fiF$. This is a contradiction.

2.  Assume that each ultrafilter converges. It is sufficient to show
that each family ${\cal H}$ of closed sets for which every finite
subfamily has a non-empty intersection has a non-empty intersection
itself. Now, the set
$\{\bigcap{\cal H}_{0}\mid {\cal H}_{0}\subseteq{\cal H}\text{
  finite}\}$
of all finite intersections forms the base for a filter $\fiF_{0}$,
which may be extended to an ultrafilter $\fiF$~\SetCite{Theorem
  1.80}. By assumption $\fiF\to x$ for some $x$, hence
$\upsilon(x)\subseteq\fiF$. The point $x$ is a candidate for being a
member in the intersection. Assume the contrary. Then there exists
$H\in {\cal H}$ with $x\not\in H$, so that $x\in X\setminus H$, which
is open. Thus $X\setminus H\in\upsilon(x)\subseteq\fiF$. On the other
hand, $H = \bigcap\{H\}\in\fiF_{0}\subseteq\fiF$, so that
$\emptyset\in\fiF$. Thus we arrive at a contradiction, and
$x\in \bigcap{\cal H}$. Hence $\bigcap{\cal H}\not=\emptyset$.
\EndProof

From Theorem~\ref{conv-vs-ultrafilter} we obtain Tihonov's celebrated
theorem\footnote{``The Tychonoff Product Theorem concerning the
  stability of compactness under formation of topological products may
well be regarded as the single most important theorem of general
topology'' according to H. Herrlich and G. E. Strecker, quoted
from~\cite[p. 85]{Herrlich-Choice}}\index{theorem!Tihonov} as an easy consequence.

\BeginTheorem{tihonov}(Tihonov's Theorem)
The product $\prod_{i\in I}X_{i}$ of topological spaces with
$X_{i}\not= \emptyset$ for all $i\in I$ is compact iff
each space $X_{i}$ is compact. 
\EndTheorem

\BeginProof
If the product $X := \prod_{i\in I}X_{i}$ is compact, then
$\Bild{\pi_{i}}{X} = X_{i}$ is compact by
Proposition~\ref{image-is-compact}. Let, conversely, be $\fiF$ an
ultrafilter on $X$, and assume all $X_{i}$ are compact. Then
$\pi_{i}(\fiF)$ is by Lemma~\ref{ultrafilter-remains-under-map} an
ultrafilter on $X_{i}$ for all $i\in I$, which converges to some
$x_{i}$ by Theorem~\ref{conv-vs-ultrafilter}. Hence $\fiF\to
(x_{i})_{i\in I}$ by Corollary~\ref{conv-in-a-product}. This implies
the compactness of $X$ by another application of Theorem~\ref{conv-vs-ultrafilter}.
\EndProof

According to~\cite[p. 146]{Engelking}, Tihonov established the theorem
for a product of an arbitrary numbers of closed and bounded intervals
of the real line (we know from the Heine-Borel
Theorem~\SetCite{Theorem 1.88} that these intervals are
compact). Kelley~\cite[p. 143]{Kelley} gives a proof of the
non-trivial implication of the theorem which relies on
Alexander's Subbase Theorem~\SetCite{Theorem 1.99}. It goes like
this. It is sufficient to establish that, whenever we have a family of
subbase elements each finite family of which fails to cover $X$, then
the whole family will not cover $X$. The sets
$\bigl\{\InvBild{\pi_{i}}{U}\mid U\subseteq X_{i}\text{ open}, i\in I\bigr\}$
form a subbase for the product topology of $X$. Let ${\cal S}$ be a
family of sets taken from this subbase such that no finite family of
elements of ${\cal S}$ covers $X$. Put ${\cal S}_{i} := \{U\subseteq
X_{i}\mid \InvBild{\pi_{i}}{U}\in {\cal S}\}$, then ${\cal S}_{i}$ is
a family of open sets in $X_{i}$. Suppose ${\cal S}_{i}$ contains sets
$U_{1}, \dots, U_{k}$ which cover $X_{i}$, then
$\InvBild{\pi_{i}}{U_{1}}, \dots, \InvBild{\pi_{i}}{U_{k}}$ are
elements of ${\cal S}$ which cover $X$; this is impossible, hence
${\cal S}_{i}$ fails to contain a finite family which covers
$X_{i}$. Since $X_{i}$ is compact, there exists  a point
$x_{i}\in X_{i}$ with $x_{i}\not\in\bigcup{\cal S}_{i}$. But then $x
:= (x_{i})_{i\in I}$ cannot be a member of $\bigcup{\cal S}$. Hence
${\cal S}$ does not cover $X$. This completes the proof.

Both proof rely heavily on the Axiom of Choice\MMP[t]{Axiom of Choice}, the first one through the existence of an ultrafilter extending a given filter, the second one through Alexander's Subbase Theorem. The relationship of Tihonov's Theorem to the Axiom of Choice is even closer: It can actually be shown that the theorem and the Axiom of Choice are equivalent~\cite[Theorem~4.68]{Herrlich-Choice}; this requires, however, establishing the existence of topological products without any recourse to the infinite Cartesian product as a carrier.

We have defined above the concept of a limit point of a filter. A weaker concept is that of an accumulation point. Taking in terms of sequences, an accumulation point of a sequence has the property that each neighborhood of the point contains infinitely many elements of the sequence. This carries over to filters in the following way.

\BeginDefinition{acc-point-filter}
Given a topological space $X$, the point $x\in X$ is called an
\emph{\index{accumulation point}\index{filter!accumulation
    point}accumulation point} of filter $\fiF$ iff
$U\cap F\not=\emptyset$ for every $U\in\upsilon(x)$ and every
$F\in\fiF$.
\EndDefinition

Since $\fiF\to x$ iff $\upsilon(x)\subseteq\fiF$, it is clear that $x$
is an accumulation point. But a filter may fail to have an
accumulation point at all. Consider the filter $\fiF$ over $\Real$
which is generated by the filter base
$\bigl\{]a, \infty[\mid a\in \Real\bigr\}$; it is immediate that $\fiF$ does not
have an accumulation point. Let us have a look at a sequence
$(x_{n})_{n\in\Nat}$, and the filter $\fiF$ generated by the infinite
tails $\bigl\{\{x_{m}\mid m\geq n\}\mid n\in\Nat\bigr\}$. If $x$ is an
accumulation point of the sequence, $U\cap \{x_{m}\mid m\geq
n\}\not=\emptyset$ for every neighborhood $U$ of $x$, thus $U\cap
F\not=\emptyset$ for all $F\in \fiF$ and all such $U$. Conversely, if
$x$ is an accumulation point for filter $\fiF$, it is clear that the
defining property holds also for the elements of the base for the
filter, thus $x$ is an accumulation point for the sequence. Hence we
have found the ``right'' generalization from sequences to filters.

An easy characterization of the set of all accumulation
points goes like this.

\BeginLemma{all-acc-points}
The set of all accumulation points of filter $\fiF$ is exactly
$\bigcap_{F\in\fiF}\Closure{F}$. 
\EndLemma

\BeginProof
This follows immediately from the observation that $x\in\Closure{A}$
iff $U\cap A\not=\emptyset$ for each neighborhood $U\in\upsilon(x)$. 
\EndProof

The lemma has an interesting consequence for the characterization of
compact spaces through filters.

\BeginCorollary{char-acc-point-ultra}
$X$ is compact iff each filter on $X$ has an accumulation point.
\EndCorollary

\BeginProof
Let $\fiF$ be a filter in a compact space $X$, and assume that $\fiF$
does not have an accumulation point. Lemma~\ref{all-acc-points}
implies that $\bigcap_{F\in\fiF}\Closure{F} = \emptyset$. Since
$X$ is compact, we find $F_{1}, \dots, F_{n}\in\fiF$ with
$\bigcap_{i=1}^{n}\Closure{F_{i}} = \emptyset$. Thus
$\bigcap_{i=1}^{n}F_{i} = \emptyset$. But this set is a member of
$\fiF$, a contradiction. 

Now assume that each filter has an accumulation point. It is by
Theorem~\ref{conv-vs-ultrafilter} enough to
show that every ultrafilter $\fiF$ converges. An accumulation point $x$
for $\fiF$ is a limit: assume that $\fiF\not\to x$, then there exists
$V\in\upsilon(x)$ with $V\not\in\fiF$, hence $X\setminus
V\in\fiF$. But $V\cap F\not=\emptyset$ for all $F\in \fiF$, since $x$
is an accumulation point. This is a contradiction. 
\EndProof

This is a characterization of accumulation points in terms of
converging filters.

\BeginLemma{acc-iff-conv-finer-filter}
In a topological space $X$, the point $x \in X$ is an accumulation
point of filter $\fiF$ iff there exists a filter $\fiF_{0}$ with
$\fiF\subseteq\fiF_{0}$ and $\fiF_{0}\to x$.  
\EndLemma

\BeginProof
Let $x$ be an accumulation point of $\fiF$, then $\{U\cap F\mid
U\in\upsilon(x), F\in\fiF\}$ is a filter base. Let $\fiF_{0}$ be the
filter generated by this base, then $\fiF\subseteq\fiF_{0}$, and
certainly $\upsilon(x)\subseteq\fiF_{0}$, thus $\fiF_{0}\to x$. 

Conversely, let $\fiF\subseteq\fiF_{0}\to x$. Since
$\upsilon(x)\subseteq\fiF_{0}$ holds as well, we conclude
$U\cap F\not=\emptyset$ for all neighborhoods $U$ and all elements
$F\in\fiF$, for otherwise we would have $\emptyset = U\cap F\in \fiF$
for some $U, F\in \fiF$, which contradicts $\emptyset\in\fiF$. Thus
$x$ is indeed an accumulation point of $\fiF$.
\EndProof


\Subsection{Separation Properties}
\label{sec:separation-props} 

We see from Example~\ref{simple-ordered-set} that a filter may
converge to more than one point. This may be undesirable. Think of a
filter which is based on a sequence, and each element of the sequence
indicates an approximation step. Then you want the approximation to
converge, but the result of this approximation process should be
unique. We will have a look at this question, and we will see that
this is actually a special case of separation properties.

\BeginProposition{limits-are-unique}
Given a topological space $X$, the following properties are equivalent
\begin{enumerate}
\item\label{limits-are-unique-1} If $x\not=y$ are different points in $X$, there exists
  $U\in\upsilon(x)$ and $V\in\upsilon(y)$ with $U\cap V=\emptyset$.
\item\label{limits-are-unique-2} The limit of a converging filter is uniquely determined.
\item\label{limits-are-unique-3} $\{x\} = \bigcap\{U\mid U\in\upsilon(x)\text{ is closed}\}$ or
  all points $x$.
\item\label{limits-are-unique-4} The diagonal $\Delta := \{\langle x, x\rangle\mid  x\in X\}$ is
  closed in $X\times X$.   
\end{enumerate}
\EndProposition

\BeginProof

\labelImpl{limits-are-unique-1}{limits-are-unique-2}: If $\fiF\to x$
and $\fiF\to y$ with $x\not= y$, we have $U\cap V\in\fiF$ for all
$U\in\upsilon(x)$ and $V\in\upsilon(y)$, hence
$\emptyset\in\fiF$. This is a contradiction. 

\labelImpl{limits-are-unique-2}{limits-are-unique-3}: 
Let $y\in
\bigcap\{U\mid U\in\upsilon(x)\text{ is closed}\}$, thus $y$ is an
accumulation point of $\upsilon(x)$. Hence there exists a filter
$\fiF$ with $\upsilon(x)\subseteq\fiF\to y$ by
Lemma~\ref{acc-iff-conv-finer-filter}. Thus $x=y$. 

\labelImpl{limits-are-unique-3}{limits-are-unique-4}: 
Let $\langle x,
y\rangle\not\in\Delta$, then there exists a closed neighborhood $W$ of
$x$ with $y\not\in W$. Let $U\in\upsilon(x)$ open with $U\subseteq W$,
and put $V := X\setminus W$, then $\langle x,
y\rangle\in U\times V\cap\Delta=\emptyset$, and $U\times V$ is
open in $X\times X$. 

\labelImpl{limits-are-unique-4}{limits-are-unique-1}: If $\langle x,
y\rangle\in(X\times X)\setminus\Delta$, there exists open sets
$U\in\upsilon(x)$ and $V\in\upsilon(y)$ with $U\times
V\cap\Delta=\emptyset$, hence $U\cap V=\emptyset$. 
\EndProof

Looking at the proposition, we see that having a unique limit for a
filter is tantamount to being able to separate two different points
through disjoint open neighborhoods. Because these spaces are
important, they deserve a special name.

\BeginDefinition{def-hausdorff}
A topological space is called a \emph{Hausdorff
  \index{space!Hausdorff, $T_{2}$}space} iff any two
different points in $X$ can be separated by disjoint open
neighborhoods, i.e., iff condition~(\ref{limits-are-unique-1}) in
Proposition~\ref{limits-are-unique} holds. Hausdorff spaces are also
called \emph{$T_{2}$-spaces}. 
\EndDefinition

\BeginExample{sorgenfrey-line}
Let $X := \Real$, and define a topology through the base $\bigl\{[a,
b[\mid a, b\in\Real, a < b\bigr\}$. Then this is a Hausdorff
space. This space is sometimes called the \emph{\index{Sorgenfrey line}Sorgenfrey line}. 
\EndExample

Being Hausdorff can be discerned from neighborhood filters:

\BeginLemma{neighborhood-filters-hausdorff}
Let $X$ be a topological space. Then $X$ is a Hausdorff space iff each
$x\in X$ has a base $\upsilon_{0}(x)$ for its neighborhood filters
such that for any $x\not= y$ there exists $U\in\upsilon_{0}(x)$ and
$V\in\upsilon_{0}(y)$ with $U\cap V=\emptyset$.
\QED
\EndLemma

It follows a first and easy consequence for maps into a Hausdorff space,
viz., the set of arguments on which they coincide is closed.

\BeginCorollary{eqset-is closed}
Let $X$, $Y$ be topological spaces, and $f, g: X\to Y$  continuous
maps. If $Y$ is a Hausdorff space, then $\{x\in X\mid f(x) =
g(x)\}$ is closed.
\EndCorollary

\BeginProof
The map $t: x\mapsto \langle f(x), g(x)\rangle$ is a continuos map
$X\to Y\times Y$. Since $\Delta\subseteq Y\times Y$ is closed by
Proposition~\ref{limits-are-unique}, the set $\InvBild{t}{\Delta}$ is
closed. But this is just the set in question.
\EndProof

The reason for calling a Hausdorff space a $T_{2}$
space\footnote{\emph{T} stands for German \emph{Trennung}, i.e.,
  separation} will become clear once we have discussed other ways of
separating points and sets; then $T_{2}$ will be a point in a spectrum
denoting separation properties.  For the moment, we introduce two
other separation properties which deal with the possibility of
distinguishing two different points through open sets. Let for this
$X$ be a topological space.
\begin{description}
\item[$T_0$-space:] $X$ is called a \emph{$T_{0}$-space}\MMP[h]{$T_{0}, T_{1}$} iff, given two
  different points $x$ and $y$, there exists an open set $U$ which
  contains exactly one of them.
\item[$T_1$-space:] $X$ is called a
  \emph{$T_{1}$-space\index{space!$T_{0}, T_{1}$}} iff, given two
  different points $x$ and $y$, there exist open neighborhoods $U$ of
  $x$ and $V$ of $y$ with $y\not\in U$ and $x\not\in V$. 
\end{description}

The following examples demonstrate these spaces.

\BeginExample{ex-t0-space}
Let $X := \Real$, and define the topologies on the real numbers
through 
\begin{align*}
  \tau_{<} & := \{\emptyset, \Real\}\cup\bigl\{]-\infty, a[\mid
             a\in\Real\bigr\},\\
\tau_{\leq} & := \{\emptyset, \Real\}\cup\bigl\{]-\infty, a]\mid
             a\in\Real\bigr\}.
\end{align*}
Then $\tau_{<}$ is a $T_{0}$-topology. $\tau_{\leq}$ is a
$T_{1}$-topology which is not $T_{0}$. 
\EndExample

This is an easy characterization of $T_{1}$-spaces.

\BeginProposition{char-t1-spaces}
A topological space $X$ is a $T_{1}$-space iff $\{x\}$ is closed for
all $x\in X$.
\EndProposition

\BeginProof
Let $y\in\Closure{\{x\}}$, then $y$ is in every open neighborhood $U$
of $x$. But this can happen in a $T_{1}$-space only if
$x=y$. Conversely, if $\{x\}$ is closed, and $y\not= x$, then there
exists a neighborhood $U$ of $x$ which does not contain $y$, and $x$
is not in the open set $X\setminus\{x\}$. 
\EndProof

\BeginExample{ex-not-t1-space}
Let $X$ be a set with at least two points, $x_{0}\in X$ be fixed. Put
$\closOp{\emptyset} := \emptyset$ and for
$\closOp{A} := A\cup\{x_{0}\}$ for $A\not=\emptyset$. Then
$\closOp{\cdot}$ is a closure operator, we look at the associated
topology. Since $\{x\}$ is open for $x\not= x_{0}$, $X$ is a $T_{0}$
space, and since $\{x\}$ is not closed for $x\not= x_{0}$, $X$ is not $T_{1}$. 
\EndExample

\BeginExample{finite-ordered-t1}
Let $(D, \leq)$ be a partially ordered set. The topology associated
with the closure operator for this order according to
Example~\ref{finite-ordered-for-closure} is $T_{1}$ iff
$y\leq x \Leftrightarrow x=y$, because this is what
$\closOp{\{x\}} = \{x\}$ says.
\EndExample

\BeginExample{cofinite-not-t2}
Let $X := \Nat$, and $\tau := \{A\subseteq \Nat\mid A\text{ is
  cofinite}\}\cup\{\emptyset\}$. Recall that a cofinite set is defined as having a finite complement. Then $\tau$ is a topology on $X$ such that $X\setminus\{x\}$ is open for each $x\in X$. Hence $X$ is a $T_{1}$-space. But $X$ is not Hausdorff. If $x\not= y$ and $U$ is an open neighborhood of $x$, then $X\setminus U$ is finite. Thus if $V$ is disjoint from $U$, we have $V\subseteq X\setminus U$. But then $V$ cannot be an open set with $y\in V$.  
\EndExample

While the properties discussed so far deal with the relationship of two different points to each other, the next group of axioms looks at closed sets; given a closed set $F$, we call an open set $U$ with $F\subseteq U$ a neighborhood of $F$. Let again $X$ be a topological space.
\begin{description}
\item[$T_3$-space:] $X$ is a \emph{$T_{3}$-space}\MMP{$T_{3},
    T_{3\half}, T_{4}$} iff given a point $x$ and a closed set $F$, which does not contain $x$, there exist disjoint open neighborhoods of $x$ and of $F$. 
\item[$T_{3\half}$-space:] $X$ is a
  \emph{$T_{3\half}$-space} iff given a
  point $x$ and a closed set $F$ with $x\not\in F$ there exists a
  continuous function $f: X\to \Real$ with $f(x) = 1$ and $f(y) = 0$
  for all $y\in F$.
\item[$T_4$-space:] $X$ is a \emph{$T_{4}$-space} \index{space!$T_{3},
    T_{3\half}, T_{4}$} iff two disjoint closed sets have disjoint open neighborhoods. 
\end{description}

$T_{3}$ and $T_{4}$ deal with the possibility of separating a closed
set from a point resp. another closed set. $T_{3\half}$ is squeezed-in
between these axioms. Because $\{x\in X\mid f(x) < 1/2\}$ and
$\{x\in X\mid f(x) > 1/2\}$ are disjoint open sets, it is clear that
each $T_{3\half}$-space is a $T_{3}$-space. It is also clear that the
defining property of $T_{3\half}$ is a special property of $T_{4}$,
provided singletons are closed. The relationship and further
properties will be explored now. 

It might be noted that continuous
functions play now an important r\^ole here in separating
objects. $T_{3\half}$ entails among others that there are ``enough''
continuous functions. Engel\-king~\cite[p. 29 and 2.7.17]{Engelking}
mentions that there are spaces which satisfy $T_{3}$ but have only
constant continuous functions, and comments ``they are, however,
fairly complicated ...'' (p. 29), Kuratowski~\cite[p. 121]{Kuratowski}
makes a similar remark. So we will leave it at that and direct the
reader, who want to know more, to these sources and the papers quoted
there.

We look at some examples.

\BeginExample{t3-nott2-nott1}
Let $X := \{1, 2, 3, 4\}$.  
\begin{enumerate}
\item With the indiscrete topology $\{\emptyset, X\}$, $X$ is a $T_{3}$
  space, but it is neither $T_{2}$ nor $T_{1}$.
\item Take the topology $\bigl\{\{1\}, \{1, 2\}, \{1, 3\}, \{1, 2,
  3\}, X, \emptyset\bigr\}$, then two closed sets are only disjoint
  when one of them is empty, because all of them contain the point $4$
  (with the exception of $\emptyset$, of course). Thus the space is
  $T_{4}$. The point $1$ and
  the closed set $\{4\}$ cannot be separated by a open sets, thus the
  space is not $T_{3}$.
\end{enumerate}
\EndExample

The next example displays a space which is $T_{2}$ but not $T_{3}$.

\BeginExample{t2-but-not-t3}
Let $X := \Real$, and put $Z := \{1/n\mid n\in \Nat\}$. Define in
addition for $x\in \Real$ and $i\in\Nat$ the sets
$
B_{i}(x) := ]x-1/i, x+1/i[.
$
Then $\upsilon_{0}(x) := \{B_{i}(x) \mid i\in \Nat\}$ for $x\not=0$,
and $\upsilon_{0}(0) := \{B_{i}(0)\setminus Z\mid i\in \Nat\}$ define
neighborhood filters for a Hausdorff space by
Lemma~\ref{neighborhood-filters-hausdorff}. But this is is not a
$T_{3}$-space. One notes first that $Z$ is closed: if $x\not\in Z$ and
$x\not\in[0, 1]$, one certainly finds $i\in \Nat$ with $B_{i}(x)\cap
Z=\emptyset$, and if $0< x \leq 1$, there exists $k$ with $1/(k+1)
< x < 1/k$, so taking $1/i$ less than the minimal distance of $x$ to
$1/k$ and $1/(k+1)$, one has $B_{i}(x)\cap Z = \emptyset$. If $x=0$,
each neighborhood contains an open set which is disjoint from $Z$.  Now each
open set $U$ which contains $Z$ contains also $0$, so we cannot
separate $0$ from $Z$. 
\EndExample

Just one positive message: the reals satisfy $T_{3\half}$. 

\BeginExample{reals-are-t3half}
Let
$F\subseteq \Real$ be non-empty, then 
\begin{equation*}
  f(t) := \inf_{y\in F}\frac{|t-y|}{1+|t-y|}
\end{equation*}
defines a continuous function $f: \Real\to [0, 1]$ with $z\in
F\Leftrightarrow f(z) = 0$. Thus, if $x\not\in F$, we have $f(x)>0$,
so that
$
  y \mapsto f(y)/f(x)
$
is a continuous function with the desired properties. Thus the reals with the usual topology are a $T_{3\half}$-space. 
\EndExample

The next proposition is a characterization of $T_{3}$-spaces in terms of open
neighborhoods, motivated by the following observation. Take a point
$x\in \Real$ and an open set $G\subseteq\Real$ with $x\in G$. Then
there exists $r>0$ such that the open interval $]x-r, x+r[$ is
entirely contained in $G$. But we can say more: by making this open
interval a little bit smaller, we can actually fit a closed interval
around $x$ into the given neighborhood as well, so, for example,
$x\in\ ]x-r/2, x+r/2[\ \subseteq\ [x-r/2, x+r/2]\ \subseteq\ ]x-r,
x+r[\ \subseteq G$.
Thus we find for the given neighborhood another neighborhood the
closure of which is entirely contained in it.

\BeginProposition{char-t3-ngbh}
Let $X$ be a topological space. Then the following are equivalent.
\begin{enumerate}
\item\label{char-t3-ngbh-1} $X$ is a $T_{3}$-space.
\item\label{char-t3-ngbh-2} For every point $x$ and every open neighborhood $U$ of $x$ there
  exists an open neighborhood $V$ of $x$ with $\Closure{V}\subseteq U$.
\end{enumerate}
\EndProposition

\BeginProof
\labelImpl{char-t3-ngbh-1}{char-t3-ngbh-2}: Let $U$ be an open
neighborhood of $x$, then $x$ is not contained in the closed set
$X\setminus U$, so by $T_{3}$ we find disjoint open sets $U_{1}, U_{2}$ with $x\in U_{1}$
and $X\setminus U\subseteq U_{2}$, hence $X\setminus U_{2}\subseteq U$. Because $U_{1}\subseteq X\setminus
U_{2}\subseteq U$, and $X\setminus U_{2}$ is closed, we conclude
$\Closure{U}_{1}\subseteq U$.

\labelImpl{char-t3-ngbh-2}{char-t3-ngbh-1}: Assume that we have a
point $x$ and a closed set $F$ with $x\not\in F$. Then $x\in
X\setminus F$, so that $X\setminus F$ is an open neighborhood of
$x$. By assumption, there exists an open neighborhood $V$ of $x$ with
$x\in \Closure{V}\subseteq X\setminus F$, then $V$ and
$X\setminus(\Closure{V})$ are disjoint open neighborhoods of $x$
resp. $F$.
\EndProof

This characterization can be generalized to $T_{4}$-spaces (roughly,
by replacing the point through a closed set) in the following way.

\BeginProposition{char-t4-ngbh}
Let $X$ be a topological space. Then the following are equivalent.
\begin{enumerate}
\item\label{char-t4-ngbh-1} $X$ is a $T_{4}$-space.
\item\label{char-t4-ngbh-2} For every closed set $F$ and every open neighborhood $U$ of $F$ there
  exists an open neighborhood $V$ of $F$ with $F\subseteq V\subseteq\Closure{V}\subseteq U$.
\end{enumerate}
\EndProposition

The proof of this proposition is actually nearly a copy of the
preceding one, \emph{mutatis mutandis}. 

\BeginProof
\labelImpl{char-t4-ngbh-1}{char-t4-ngbh-2}:
Let $U$ be an open neighborhood of the closed set $F$, then the closed
set $F' :=
X\setminus U$ is disjoint to $F$, so that we can
find disjoint open neighborhoods $U_{1}$ of $F$ and $U_{2}$ of $F'$,
thus $U_{1}\subseteq X\setminus U_{2} \subseteq X\setminus F' = U$, so
$V := U_{1}$ is the open neighborhood we are looking for.

\labelImpl{char-t4-ngbh-2}{char-t4-ngbh-1}: 
Let $F$ and $F'$ be
disjoint closed sets, then $X\setminus F'$ is an open neighborhood for
$F$. Let $V$ be an open neighborhood for $F$ with $F\subseteq
V\subseteq \Closure{V}\subseteq X\setminus F'$, then $V$ and
$U := X\setminus(\Closure{V})$ are disjoint open neighborhoods of $F$ and
$F'$. 
\EndProof
 
We mentioned above that the separation axiom $T_{3\half}$ makes sure
that there are enough continuous functions on the space. Actually,
the continuous functions even determine the topology in this case, as
the following characterization shows.

\BeginProposition{char-t_3half}
Let $X$ be a topological space, then the following statements are
equivalent.
\begin{enumerate}
\item\label{char-t_3half-1} $X$ is a $T_{3\half}$-space.
\item\label{char-t_3half-2} $\beta := \bigl\{\InvBild{f}{U}\mid  f: X\to \Real \text{ is
    continuous}, U\subseteq \Real\text{ is open}\bigr\}$ constitutes a
  basis for the topology of $X$.
\end{enumerate}
\EndProposition

\BeginProof
The elements of $\beta$ are open sets, since they are comprised of
inverse images of open sets under continuous functions.

\labelImpl{char-t_3half-1}{char-t_3half-2}:  
Let $G\subseteq X$ be an open set with $x\in G$. We show that we can
find $B\in \beta$ with $x\in B\subseteq B$. In fact, since $X$ is
$T_{3\half}$, there exists a continuous function $f: X\to \Real$ with
$f(x) = 1$ and $f(y) = 0$ for $y\in X\setminus G$. Then $B := \{x\in
X\mid -\infty < x < 1/2\} = \InvBild{f}{]-\infty, 1/2[}$ is a suitable
element of $\beta$.

\labelImpl{char-t_3half-2}{char-t_3half-1}:
Take $x\in X$ and a closed set $F$ with $x\not\in F$. Then $U :=
X\setminus F$ is an open neighborhood $x$. Then we can find
$G\subseteq \Real$ open and $f: X\to \Real$ continuous with $x\in
\InvBild{f}{G}\subseteq U$. Since $G$ is the union of open intervals,
we find an open interval $I := ]a, b[\ \subseteq G$ with
$f(x)\in I$. Let $G: \Real\to \Real$ be a continuous with $g(f(x)) =
1$ and $g(t) = 0$, if $t\not\in I$; such a function exists since
$\Real$ is a $T_{3\half}$-space (Example~\ref{reals-are-t3half}). Then
$g\circ f$ is a continuous function with the desired
properties. Consequently, $X$ is a $T_{3\half}$-space.
\EndProof

The separation axioms give rise to names for classes of spaces. We
will introduce there traditional names now.

\BeginDefinition{traditional names}
Let $X$ be a topological space, then $X$ is called
\begin{itemize}
\item \emph{\index{space!regular}regular} iff $X$ satisfies $T_{1}$ and $T_{3}$,
\item \emph{\index{space!completely regular}completely regular}, iff $X$ satisfies $T_{1}$ and $T_{3\half}$,
\item \emph{\index{space!normal}normal}, iff $X$ satisfies $T_{1}$ and $T_{4}$. 
\end{itemize}
\EndDefinition

The reason $T_{1}$ is always included is that one wants to have 
every singleton as a closed set, which, as the examples show, is not
always the case. Each regular space is a Hausdorff space, each regular
space is completely regular, and each normal space is regular. We will
obtain as a consequence of Urysohn's Lemma that that each normal space
is completely regular as well (Corollary~\ref{normal-arecompletely-reg}).

In a completely regular space we can separate a point $x$ from a
closed set not containing $x$ through a continuous function. It turns
out that normal spaces have an analogous property: Given two disjoint
closed sets, we can separate these sets through a continuous
function. This is what \emph{\index{Urysohn's Lemma}Urysohn's Lemma} says,
a famous result from the beginnings of set-theoretic topology. To be precise:

\BeginTheorem{urysohns-lemma}(Urysohn)
Let $X$ be a  normal space. Given disjoint closed sets $F_{0}$ and
$F_{1}$, there exists a continuous function $f: X\to \Real$ such that
$f(x) = 0$ for $x\in F_{0}$ and $f(x) = 1$ for $x\in F_{1}$. 
\EndTheorem

We need some technical preparations for proving
Theorem~\ref{urysohns-lemma}; this gives also the opportunity to
introduce the concept of a dense set.

\BeginDefinition{set-is-dense}
A subset $D\subseteq X$ of a topological space $X$ is called
\emph{\index{dense set}dense} iff
$\Closure{D} = X$. 
\EndDefinition

Dense sets are fairly practical when it comes to compare continuous
functions for equality: if suffices that the functions coincide on a
dense set, then they will be equal. Just for the record:

\BeginLemma{equal-on-dense}
Let $f, g: X\to Y$ be continuous maps with $Y$ Hausdorff, and assume
that $D\subseteq X$ is dense. Then $f = g$
iff $f(x) = g(x)$ for all $x\in D$. 
\EndLemma

\BeginProof
Clearly, if $f = g$, then $f(x) = g(x)$ for all $x\in D$. So we have
to establish the other direction.

Because $Y$ is a Hausdorff space, $\Delta_{Y}:= \{\langle y,
y\rangle\mid y\in Y\}$ is closed (Proposition~\ref{limits-are-unique}), and because $f\times g: X\times X\to
Y\times Y$ is continuous, $\InvBild{(f\times g)}{\Delta_{Y}}\subseteq
X\times X$ is closed as well. The latter set contains $\Delta_{D}$, hence its
closure $\Delta_{X}$. 
\EndProof

It is immediate that if $D$ is dense, then $U\cap D\not=\emptyset$ for
each open set $U$, so in particular each neighborhood of a point meets
the dense set $D$. To provide an easy example, both $\Rational$ and
$\Real\setminus\Rational$ are dense subsets of $\Real$ in the usual
topology. Note that $\Rational$ is countable, so $\Real$ has even a
countable dense set.

The first lemma has a family of subsets indexed by a dense subset of
$\Real$ exhaust a given set and provides a useful real function.

\BeginLemma{exhaust-family-1}
Let $M$ be set, $D\subseteq\pReal$ be dense, and $(E_{t})_{t\in D}$ be
a family of subsets of $M$ with these properties:
\begin{itemize}
\item if $t < s$, then $E_{t}\subseteq E_{s}$,
\item $M = \bigcup_{t\in D}E_{t}$.
\end{itemize}
Put $f(m) := \inf\{t\in D\mid m\in E_{t}\}$, then we have for all
$s\in \Real$
\begin{enumerate}
\item $\{m\mid f(m) < s\} = \bigcup\{E_{t}\mid t\in D, t < s\}$,
\item $\{m\mid f(m)\leq s\} = \bigcap\{E_{t}\mid  t\in D, t > s\}$. 
\end{enumerate}
\EndLemma

\BeginProof
1.
Let us work on the first equality. If $f(m) < s$, there exists $t < s$
with $m\in E_{t}$. Conversely, if $m\in E_{t}$ for some $t < s$, then
$f(m) = \inf\{r\in D\mid m\in E_{r}\} \leq t < s$.

2.
For the second equality, assume $f(m) \leq s$, then we can find for each $r > s$ some $t < r$
with $m\in E_{t}\subseteq E_{r}$. To establish the other inclusion,
assume that $f(m) \leq t$ for all $t > s$. If $f(m)  = r > s$, we can
find some $t'\in D$ with $r > t' > s$, hence $f(m) \leq t'$. This is a
contradiction, hence $f(m) \leq s$. 
\EndProof
 
This lemma, which does not assume a topology on $M$, but requires only a plain
set, is extended now for the topological scenario in which we will use
it. We assume that each set $E_{t}$ is open, and we assume that
$E_{t}$ contains the closures of its predecessors. Then it will turn
out that the function we just have defined is continuous,
specifically:

\BeginLemma{exhaust-family-2}
Let $X$ be a topological space, $D\subseteq\pReal$ a dense subset, and
assume that $(E_{t})_{t\in D}$ is a family of open sets with these
properties
\begin{itemize}
\item if $t < s$, then $\Closure{E}_{t}\subseteq E_{s}$,
\item $X = \bigcup_{t\in D}E_{t}$.
\end{itemize}
Then $f: x \mapsto inf\{t\in D\mid x\in E_{t}\}$ defines a continuous
function on $X$.
\EndLemma

\BeginProof
0.  
Because a subbase for the topology on $\Real$ is comprised of the
intervals $]-\infty, x[$ resp. $]x, +\infty[$, we see from
Lemma~\ref{cont-for-subbase} that it is sufficient to show that for
any $s\in \Real$ the sets $\{x\in X \mid f(x) < s\}$ and
$\{x\in X\mid f(x) > s\}$ are open, since they are the corresponding
inverse images under $f$. For the latter set we show that its
complement $\{x\in X\mid f(x) \leq s\}$ is closed. Fix $s\in \Real$.

1.
We obtain from Lemma~\ref{exhaust-family-1} that $\{x\in X \mid f(x) <
s\}$ equals $\bigcup \{E_{t}\mid t\in D, t < s\}$; since all sets
$E_{t}$ are open, their union is. Hence $\{x\in X \mid f(x) <
s\}$ is open.

2.
We obtain again from Lemma~\ref{exhaust-family-1} that $\{x\in X\mid
f(x) \leq s\}$ equals $\bigcap\{E_{t}\mid  t\in D, t > s\}$, so if we
can show that $ \bigcap\{E_{t}\mid  t\in D, t > s\} =
\bigcap\{\Closure{E}_{t}\mid  t\in D, t > s\}$, we are done. In
fact, the left hand side is contained in the right hand side, so
assume that $x$ is an element of the right hand side. If $x$ is not
contained in the left hand side, we find $t'>s$ with $t'\in D$ such that
$x\not\in E_{t'}$. Because $D$ is dense, we find some $r$ with $s < r <
t'$ with $\Closure{E}_{r}\subseteq E_{t'}$. But then $x\not\in
\Closure{E}_{r}$, hence $x\not\in\bigcap
\{\Closure{E}_{t}\mid t\in D, t > s\}$, a contradiction. Thus both
sets are equal, so that $\{x\in X \mid f(x)\geq s\}$ is closed. 
\EndProof

We are now in a position to establish Urysohn's Lemma. The idea of the
proof rests on this observation for a $T_{4}$-space $X$: suppose that
we have open sets $A$ and $B$ with $A\subseteq\Closure{A}\subseteq B$. Then we
can find an open set $C$ such that $\Closure{A}\subseteq C \subseteq
\Closure{C}\subseteq B$, see Proposition~\ref{char-t4-ngbh}.
{\def\ssub{\sqsubseteq^*}
Denote just for the proof for open sets $A, B$ the fact that
$\Closure{A}\subseteq B$ by $A\ssub B$. Then we may express the idea
above by saying that $A\ssub B$ implies the existence of an open set
$C$ with $A\ssub C \ssub B$, so $C$ may be squeezed in. But now we
have $A\ssub C$ and $C\ssub B$, so we find open sets $E$ and $F$ with
$A\ssub E \ssub C$ and $C \ssub F \ssub B$, arriving at the chain $A
\ssub E \ssub C \ssub F \ssub B$. But why stop here? 

The proof makes this argument systematic and constructs in this way a
continuous function.

\BeginProof
1.
Let $D := \{p/2^{q}\mid p, q\text{ non-negative integers}\}$. These
are all dyadic numbers, which are dense in $\pReal$. We are about to
construct a family $(E_{t})_{t\in D}$ of open sets $E_{t}$ indexed by $D$ in the following way.

2.
Put $E_{t} := X$ for $t>1$, and let $E_{1} := X\setminus F_{1}$,
moreover let $E_{0}$ be an open set containing $F_{0}$ which is
disjoint from $E_{1}$. We now construct open sets $E_{p/2^{n}}$ by
induction on $n$ in the following way. Assume that we have already
constructed open sets $$E_{0}\ssub E_{\frac{1}{2^{n-1}}}\ssub
E_{\frac{2}{2^{n-1}}}\dots\ssub E_{\frac{2^{n-1}-1}{2^{n-1}}} \ssub
E_{1}.$$ Let $t = \frac{2m+1}{2^{n}}$, then we find an open set
$E_{t}$ with $E_{\frac{2m}{^{n}}}\ssub E_{t}\ssub
E_{\frac{2m+2}{2^{n}}}$; we do this for all $m$ with $0 \leq m \leq 2^{n-1}-1$. 

3.
Look as an illustration at the case $n=3$. We have found already
the open sets $E_{0}\ssub E_{1/4}\ssub E_{1/2}\ssub E_{3/4}\ssub
E_{1}$. Then the construction goes on with finding open sets $E_{1/8},
E_{3/8}, E_{5/8}$ and $E_{7/8}$ such that after the
step is completed, we obtain this chain.
\begin{equation*}
 E_{0}\ssub E_{1/8}\ssub E_{1/4}\ssub E_{3/8}\ssub E_{1/2}\ssub
 E_{5/8}\ssub E_{3/4}\ssub E_{7/8}\ssub E_{1}. 
\end{equation*}

4.
In this way we construct a family $(E_{t})_{t\in D}$ with the
properties requested by Lemma~\ref{exhaust-family-2}. It yields a
continuous function $f: X\to \Real$ with $f(x) = 0$ for all $x\in F_{0}$
and $f(1)(x) = 1$ for all $x\in F_{1}$. 
\EndProof   
}

Urysohn's Lemma is used to prove the Tietze Extension Theorem, which
we will only state, but not prove.

\BeginTheorem{tietze-extension}
Let $X$ be a $T_{4}$-space, and $f: A\to \Real$ be a function which is
continuous on a closed subset $A$ of $X$. Then $f$ can be extended to
a continuous function $f^{*}$ on all of $X$. 
\QED
\EndTheorem

We obtain as an immediate consequence of Urysohn's Lemma

\BeginCorollary{normal-arecompletely-reg}
A normal space is completely regular. 
\EndCorollary

We have obtained a hierarchy of spaces through gradually tightening
the separation properties, and found that continuous functions help
with the separation. The question arises, how compactness fits into
this hierarchy. It turns out that a compact Hausdorff space is normal;
the converse obviously does not hold: the reals with the Euclidean
topology are normal, but by no means compact. 

We call a subset $K$ in a topological space $X$ compact iff it is
compact as a subspace, i.e., a compact topological space in its own
right. This is a first and fairly straightforward observation.

\BeginLemma{closed-is-compact}
A closed subset $F$ of a compact space $X$ is
compact.
\EndLemma

\BeginProof
Let $(G_{i}\cap F)_{i\in I}$ be an open cover of
$F$ with $G_{i}\subseteq X$ open, then $\{F\}\cup\{G_{i}\mid i\in I\}$
is an open cover of $X$, so we can find a finite subset $J\subseteq I$
such that $\{F\}\cup\{G_{j}\mid i\in J\}$ covers $X$, hence
$\{G_{i}\cap F\mid i\in J\}$ covers $F$. 
\EndProof

In a Hausdorff space, the
converse holds as well:

\BeginLemma{compact-is-closed}
Let $X$ be a Hausdorff space, and $K\subseteq X$ compact, then 
\begin{enumerate}
\item Given $x\not\in K$, there exist disjoint open neighborhoods $U$ of
  $x$ and $V$ of $K$.
\item $K$ is closed.
\end{enumerate}
\EndLemma

\BeginProof
Given $x\not\in K$, we want to find $U\in\upsilon(x)$ with $U\cap
K=\emptyset$ and $V\supseteq K$ open with $U\cap V = \emptyset$. 

Let's see, how to do that. There exists for $x$ and any element $y\in K$
disjoint open neighborhoods $U_{y}\in\upsilon(x)$ and
$W_{y}\in\upsilon(y)$, because $X$ is Hausdorff. Then $(W_{y})_{y\in
  Y}$ is an open cover of $K$, hence by compactness there exists a
finite subset $W_{0}\subseteq W$ such that $\{W_{y}\mid y\in W_{0}\}$
covers $K$. But then $\bigcap_{y\in W_{0}} U_{y}$ is an open
neighborhood of $x$ which is disjoint from $V := \bigcup_{y\in W_{0}}
W_{y}$, hence from $K$. $V$ is the open neighborhood of $K$ we are
looking for. This establishes the first part, the second follows as an
immediate consequence.
\EndProof

Look at the reals as an illustrative example.

\BeginCorollary{compact-subset-real}
$A\subseteq\Real$ is compact iff it is closed and bounded. 
\EndCorollary

\BeginProof If $A \subseteq \Real$ is compact, then it is closed by Lemma~\ref{compact-is-closed}, since $\Real$ is a Hausdorff space. Since $A$ is compact, it is also bounded. If, conversely, $A\subseteq\Real$ is closed and bounded, then we can find a closed interval $[a, b]$ such that $A\subseteq [a, b]$. We know from the Heine-Borel Theorem~\SetCite{Theorem 1.88} that this interval is compact, and a closed subset of a compact space is compact by Lemma~\ref{closed-is-compact}.  \EndProof

This has yet another, frequently used consequence, viz., that a continuous real valued function on a compact space assume its minimal and its maximal value. Just for the record:

\BeginCorollary{compact-yields-extrema}
Let $X$ be a compact Hausdorff space, $f: X\to \Real$ a continuous
map. Then there exist $x_{*}, x^{*}\in X$ with $f(x_{*}) =
\min\Bild{f}{X}$ and $f(x^{*})  = \max \Bild{f}{X}$. 
\QED
\EndCorollary

But ---~after travelling an interesting side path~--- let us return to the problem of establishing that a compact
Hausdorff space is normal. We know now that we can separate a point
from a compact subset through disjoint open neighborhoods. This is but
a small step from establishing the solution to the above problem.

\BeginProposition{compact-is-normal}
A compact Hausdorff space is normal.
\EndProposition

\BeginProof
Let $X$ be compact, $A$ and $B$ disjoint closed subsets. Since $X$ is
Hausdorff, $A$ and $B$ are compact as well. Now the rest is an easy
application of Lemma~\ref{compact-is-closed}. Given $x\in B$, there
exist disjoint open neighborhoods $U_{x}\in\upsilon(x)$ of $x$ and
$V_{x}$ of $A$. Let $B_{0}$ be a finite subset of $B$ such that
$U := \bigcup\{U(x)\mid b\in B_{0}\}$ covers $B$ and
$V := \bigcap\{V_{x}\mid x\in B_{0}\}$ is an open neighborhood of
$A$. $U$ and $V$ are disjoint.
\EndProof

From the point of view of separation, to be compact is for a topological space a stronger property than being normal. The example $\Real$ shows that this is a strictly stronger property. We will show now that $\Real$ is just one point apart from being compact by investigating locally compact spaces.


\Subsection{Local Compactness and Compactification} 
\label{sec:local-comp-comp}

We restrict ourselves in this section to Hausdorff spaces. Sometimes a space is not compact but has enough compact subsets, because each point has a compact neighborhood. These spaces are called locally compact, and we investigate properties they share with and properties they distinguish them from compact spaces. We show also that a locally compact space misses being compact by just one point. Adding this point will make it compact, so we have an example here where we embed a space into one with a desired property. While we are compactifying spaces, we also provide another one, named after Stone and $\check{\mathrm{C}}$ech, which requires the basic space to be completely regular. We establish also another classic, the Baire Theorem, which states that in a locally compact $T_{3}$ space the intersection of a countable number of open dense sets is dense again; applications will later on capitalize on this observation. 

\BeginDefinition{local-compact}
Let $X$ be a Hausdorff space. $X$ is called \emph{\index{space!locally compact}\index{compact!locally compact}locally compact} iff for each $x\in X$ and each open neighborhood $U\in \upsilon(x)$ there exists a neighborhood $V\in\upsilon(x)$ such that $\Closure{V}$ is compact and $\Closure{V}\subseteq U$. 
\EndDefinition

Thus the compact neighborhoods form a basis for the neighborhood filter for each point. This implies that we can find for each compact subset an open neighborhood with compact closure. The proof of this property gives an indication of how to argue in locally compact spaces.

\BeginProposition{compact-closure-for-open}
Let $X$ be a locally compact space, $K$ a compact subset. Then there exists an open neighborhood $U$ of $K$ and a compact set $K'$ with $K\subseteq U\subseteq K'$. 
\EndProposition

\BeginProof
Let $x\in K$, then we find an open neighborhood $U_{x}\in\upsilon(x)$ with $\Closure{U}_{x}$ compact. Then $(U_{x})_{x\in K}$ is a cover for $K$, and there exists a finite subset $K_{0}\subseteq K$ such that $(U_{x})_{x\in K_{0}}$ covers $K$. Put $U := \bigcup_{x\in K_{0}}$, and note that this open set has a compact closure. 
\EndProof

So this is not too bad: We have plenty of compact sets in a locally compact space. Such a space is very nearly compact. We add to $X$ just one point, traditionally called $\infty$ and define the neighborhood for $\infty$ in such a way that the resulting space is compact. The obvious way to do that is to make all complements of compact sets a neighborhood of $\infty$, because it will then be fairly easy to construct from a cover of the new space a finite subcover. This is what the compactification which we discuss now will do for you. We carry out the construction in a sequence of lemmas, just in order to render the process a bit more transparent.

\BeginLemma{l1-alexandrov-one-point}
Let $X$ be a Hausdorff space with topology $\tau$, $\infty\not\in X$
be a distinguished new point. Put $X^{*} := X\cup\{\infty\}$, and
define\MMP{One point extension}
\begin{equation*}
  \tau^{*} := \{U\subseteq X^{*}\mid U\cap X\in\tau\}\cup\{U\subseteq X^{*}\mid \infty\in U, X\setminus U\text{ is compact}\}.
\end{equation*}
Then $\tau^{*}$ is a topology on $X^{*}$, and the identity $i_{X}: X\to X^{*}$ is $\tau$-$\tau^{*}$-continuous. 
\EndLemma

\BeginProof
$\emptyset$ and $X^{*}$ are obviously members of $\tau^{*}$; note that $X\setminus U$ being compact entails $U\cap X$ being open. Let $U_{1}, U_{2}\in\tau^{*}$. If $\infty\in U_{1}\cap U_{2}$, then $X\setminus (U_{1}\cap U_{2})$ is the union of two compact sets in $X$, hence is compact, If $\infty\not\in U_{1}\cap U_{2}$, $X\cap (U_{1}\cap U_{2})$ is open in $X$. Thus $\tau^{*}$ is closed under finite intersections. Let $(U_{i})_{i\in I}$ be a family of elements of $\tau^{*}$. The critical case is that $\infty\in U := \bigcup_{i\in I}U_{i}$, say, $\infty\in U_{j}$. But then $X\setminus U \subseteq X\subseteq U_{j}$, which is compact, so that $U\in\tau^{*}$. Continuity of $i_{X}$ is now immediate.  
\EndProof

We find $X$ in this new construction as a subspace.

\BeginCorollary{x-is-a-subspace}
$(X, \tau)$ is a dense subspace of $(X^{*}, \tau^{*})$. 
\EndCorollary

\BeginProof
We have to show that $\tau =\tau^{*}\cap X$. But this is obvious from the definition of $\tau^{*}$. 
\EndProof

Now we can state and prove the result which has been announced above.

\BeginTheorem{alexandrov-one-point}
Given a Hausdorff space $X$, the one point extension $X^{*}$ is a compact space, in which $X$ is dense. If $X$ is locally compact, $X^{*}$ is a Hausdorff space.  
\EndTheorem

\BeginProof
It remains to show that $X^{*}$ is compact, and that it is a Hausdorff space, whenever $X$ is locally compact. 

Let $(U_{i})_{i\in I}$ be an open cover of $X^{*}$, then $\infty\in U_{j}$ for some $j\in I$, thus $X\setminus U_{j}$ is compact and is covered by $(U_{i})_{i\in I, i\not= j}$. Select an finite subset $J\subseteq I$ such that $(U_{i})_{i\in J}$ covers $X\setminus U_{j}$, then ---~voilà~--- we have found a finite cover $(U_{i})_{i\in J\cup\{j\}}$ of $X^{*}$.

Since the given space is Hausdorff, we have to separate the new point $\infty$ from a given point $x\in X$, provided $X$ is locally compact. But take a compact neighborhood $U$ of $x$, then $X^{*}\setminus U$  is an open neighborhood of $\infty$. 
\EndProof 

$X^{*}$ is called the \emph{\index{compactification!Alexandrov one point}Alexandrov one point compactification} of $X$. The new point is sometimes called the \emph{infinite point}. It is not difficult to show that two different one point compactifications are homeomorphic, so we may talk about \emph{the} (rather than \emph{a}) one-point compactification. 

Looking at the map $i_{X}: X\to X^{*}$, which permits looking at elements of $X$ as elements of $X^{*}$, we see that $i_{X}$ is injective and has the property that $\Bild{i_{X}}{G}$ is an open set in the image $\Bild{i_{X}}{X}$ of $X$ in $X^{*}$, whenever $G\subseteq X$ is open. These properties will be used for characterizing compactifications. Let us first define embeddings, which are of interest independently of compactifications.

\BeginDefinition{embedding-map}
The continuous map $f: X\to Y$ between the topological spaces $X$ and $Y$ is called an \emph{\index{embedding}embedding} iff 
\begin{itemize}
\item $f$ is injective,
\item $\Bild{f}{G}$ is open in $\Bild{f}{X}$, whenever $G\subseteq X$ is open.
\end{itemize}
\EndDefinition
 
So if $f: X\to Y$ is an embedding, we may recover a true image of $X$
from its image $\Bild{f}{X}$, so that $f: X\to \Bild{f}{X}$ is a
homeomorphism.

Let us have a look at the map $[0, 1]^{N}\to [0, 1]^{M}$, which is induced by a map $f: M\to N$ for sets $M$ and $N$, and which we delt with in Lemma~\ref{into-unitcube}. We will put this map to good use in a moment, so it is helpful to analyze it a bit more closely.

\BeginExample{onto-cube-is-embedding}
Let $f: M\to N$ be a surjective map. Then
$f^{*}: [0, 1]^{N}\to [0, 1]^{M}$, which sends $g: N\to [0, 1]$ to
$g\circ f: M\to [0, 1]$ is an embedding. We have to show that $f^{*}$
is injective, and that it maps open sets into open sets in the
image. This is done in two steps:

\begin{description}
\item[$f^{*}$ is injective:] In fact, if $g_{1}\not= g_{2}$, we
  find $n\in N$ with $g_{1}(n)\not=g_{2}(n)$, and because $f$ is onto,
  we find $m$ with $n = f(m)$, hence
  $f^{*}(g_{1})(m) = g_{1}(f(m)) \not=g_{2}(f(m)) =
  f^{*}(g_{2})(m)$. Thus $f^{*}(g_{1})\not=f(g_{2})$ (an alternative
  and more general proof is proposed in~\CategCite{Proposition 1.23}).
\item[Open sets are mapped to open sets:] We know already from Lemma~\ref{into-unitcube} that $f^{*}$ is
  continuous, so we have to show that the image $\Bild{f}{G}$ of an
  open set $G\subseteq [0, 1]^{N}$ is open in the subspace
  $\Bild{f}{[0, 1]^{M}}$. Let $h\in \Bild{f}{G}$, hence $h = f^{*}(g)$
  for some $g\in G$. $G$ is open, thus we can find a subbase element
  $H$ of the product topology with $g\in H\subseteq G$, say,
  $H = \bigcap_{i=1}^{k}\InvBild{\pi_{N, n_{i}}}{H_{i}}$ for some
  $n_{1}, \dots, n_{k}\in N$ and some open subsets
  $H_{1}, \dots, H_{k}$ in $[0, 1]$. Since $f$ is onto,
  $n_{1} = f(m_{1}), \dots, n_{k} = f(m_{k})$ for some
  $m_{1}, \dots, m_{k}\in M$. Since
  $h\in \InvBild{\pi_{N, n_{i}}}{H_{i}}$ iff
  $f^{*}(h)\in\InvBild{\pi_{M, m_{i}}}{H_{i}}$, we obtain
  \begin{equation*}
    h = f^{*}(g) \in \Bild{f^{*}}{\bigcap_{i=1}^{k}\InvBild{\pi_{N, n_{i}}}{H_{i}}}
    = \bigl(\bigcap_{i=1}^{k}\InvBild{\pi_{M, m_{i}}}{H_{i}}\bigr)\cap\Bild{f^{*}}{[0, 1]^{N}}
  \end{equation*}
  The latter set is open in the image of $[0, 1]^{N}$ under $f^{*}$,
  so we have shown that the image of an open set is open relative to
  the subset topology of the image.
\end{description}
These proofs will serve as patterns later on.
\EndExample

Given an embedding, we define the compactification of a space.

\BeginDefinition{compfact-space}
A pair $(e, Y)$ is said to be a
\emph{\index{topology!compactification}\index{compactification}compactification}
of a topological space $X$ iff $Y$ is a compact topological space, and
if $e: X\to Y$ is an embedding.
\EndDefinition

The pair $(i_{X}, X^{*})$ constructed as the Alexandrov one-point
compactification is a compactification in the sense of
Definition~\ref{compfact-space}, provided the space $X$ is locally
compact. We are about to construct another important compactification
for a completely regular space $X$. Define for $X$ the space $\beta X$
as follows\footnote{It is a bit unfortunate that there appears to be
  an ambiguity in notation, since we denote the basis of a topological
  space by $\beta$ as well. But tradition demands this
  compactification to be called $\beta X$, and from the context it
  should be clear what we have in mind.}: Let $F(X)$ be all continuous
maps $X\to [0, 1]$, and map $x$ to its evaluations from $F(X)$, so
construct $e_{X}: X\ni x\mapsto (f(x))_{f\in F(X)}\in[0,
1]^{F(X)}$.
Then $\beta X := \Closure{(\Bild{e_{X}}{X})}$, the closure being taken
in the compact space $[0, 1]^{F(X)}$. We claim that $(e_{X}, \beta X)$
is a compactification of $X$. Before delving into the proof, we note
that we want to have a completely regular space, since there we have
enough continuous functions, e.g., to separate points, as will become
clear shortly. We will first show that this is a compactification
indeed, and then investigate an interesting property of it.

\BeginProposition{beta-x-is-a-compactification}
$(e_{X}, \beta X)$ is a compactification of the completely regular space $X$.
\EndProposition

\BeginProof
1.
We take the closure in the Hausdorff space $[0, 1]^{F(X)}$, which is compact by Tihonov's Theorem~\ref{tihonov}. Hence $\beta X$ is a compact Hausdorff space by Lemma~\ref{closed-is-compact}. 

2.  $e_{X}$ is continuous, because we have $\pi_{f}\circ e_{X} = f$
for $f\in F(X)$, and each $f$ is continuous. $e_{X}$ is also
injective, because we can find for $x\not=x'$ a map $f\in F(X)$ such
that $f(x) \not= f(x')$; this translates into
$e_{X}(x) (f) \not= e_{X}(x')(f)$, hence $e_{X}(x) \not= e_{X}(x')$.

3.  The image of an open set in $X$ is open in the image. In fact, let
$G\subseteq X$ be open, and take $x\in G$. Since $X$ is completely
regular, we find $f\in F(X)$ and an open set $U\subseteq[0, 1]$ with
$x\in \InvBild{f}{U}\subseteq G$; this is so because the inverse
images of the open sets in $[0, 1]$ under continuous functions form a
basis for the topology (Proposition~\ref{char-t_3half}). But
$x\in \InvBild{f}{U}\subseteq G$ is equivalent to
$x\in \InvBild{(\pi_{f}\circ e_{X})}{U}\subseteq G$. Because
$e_{X}: X\to \Bild{e_{X}}{X}$ is a bijection, this implies
$x\in\InvBild{\pi_{f}}{U}\subseteq \Bild{e_{X}}{G}\cap
\Bild{e_{X}}{X}\subseteq
\Bild{e_{X}}{G}\cap\Closure{(\Bild{e_{X}}{X})}$.
Hence $\Bild{e_{X}}{G}$ is open in $\beta X$.
\EndProof

If the space we started from is already compact, then we obtain nothing new:

\BeginCorollary{compact-homeom-to-beta}
If $X$ is a compact Hausdorff space, $e_{X}: X\to \beta X$ is a homeomorphism.
\EndCorollary

\BeginProof
A compact Hausdorff space is normal, hence completely regular by
Proposition~\ref{compact-is-normal} and
Corollary~\ref{normal-arecompletely-reg}, so we can construct the
space $\beta X$ for $X$ compact. The assertion then follows from
Exercise~\ref{ex-compact-homeom}.
\EndProof

This kind of compactification is important, so it deserves a distinguishing name.

\BeginDefinition{stone-cech}
The compactification $(e_{X}, \beta X)$ is called the
\emph{Stone-$\check{\mathrm{C}}$ech
  \index{compactification!Stone-$\check{\mathrm{C}}$ech}compactification}
of the regular space $X$.
\EndDefinition

This compactification permits the extension of continuous maps in the
following sense: suppose that $f: X\to Y$ is continuous with $Y$
compact, then there exists a continuous extension $\beta X\to Y$. This
statement is slightly imprecise, because $f$ is not defined on
$\beta X$, so we want really to extend
$f\circ e_{X}^{-1}: \Bild{e_{X}}{X}\to Y$ ~---~since $e_{X}$ is a
homeomorphism from $X$ onto its image, one tends to identify both
spaces.

\BeginTheorem{extension-of-stone-cech}
Let $(e_{X}, \beta X)$ be the Stoch-$\check{\mathrm{C}}$ech
compactification of the completely regular space $X$. Then, given a
continuous map $f: X\to Y$ with $Y$ compact, there exists a continuous
extension $f_{!}: \beta X\to Y$ to $f\circ e_{X}^{-1}$.
\EndTheorem

The idea of the proof is to capitalize on the compactness of the
target space $Y$, because $Y$ and $\beta Y$ are homeomorphic. This
means that $Y$ has a topologically identical copy in $[0, 1]^{F(Y)}$,
which may be used in a suitable fashion. The proof is adapted
from~\cite[p. 153]{Kelley}; Kelley calls it a ``mildly intricate
calculation''.

\BeginProof
1.  Define $\phi_{f}: F(Y)\to F(X)$ through $h\mapsto f\circ h$, then
this map induces a map $\phi^{*}_{f}: [0, 1]^{F(X)}\to [0, 1]^{F(Y)}$
by sending $t: F(X)\to [0, 1]$ to $t\circ \phi_{f}$. Then
$\phi^{*}_{f}$ is continuous according to Lemma~\ref{into-unitcube}.

2.
Consider this diagram
\begin{equation*}
\xymatrix{
\Bild{e_{X}}{X}\ar[rr]^{\subseteq}&&[0, 1]^{F(X)}\ar[rr]^{\phi^{*}_{f}}&&[0, 1]^{F(Y)}&& \beta Y\ar[ll]_{\supseteq}\\
X\ar[u]_{e_{X}} \ar[rrrrrr]_{f}&&&&&&Y\ar[u]^{e_{Y}}
}
\end{equation*}
We claim that $\phi^{*}_{f}\circ e_{X} = e_{Y}\circ f$. In fact, take $x\in X$ and $h\in F(Y)$, then
\begin{align*}
  (\phi_{f}^{*}\circ e_{X})(x)(h) 
& = 
(e_{X}\circ \phi_{f})(h) \\
& = 
e_{X}(x)(h\circ f) \\
& = 
(h\circ f)(x) \\
& = 
e_{Y}(f(x))(h)\\
& =
(e_{Y}\circ f)(x)(h).
\end{align*}

3.
Because $Y$ is compact, $e_{Y}$ is a homeomorphism by Exercise~\ref{ex-compact-homeom}, and since $\phi^{*}_{f}$ is continuous, we have 
\begin{equation*}
\Bild{\phi_{f}^{*}}{\beta X} = 
\Bild{\phi_{f}^{*}}{\Closure{\Bild{e_{X}}{X}}} \subseteq 
\Closure{\bigl(\Bild{\phi_{f}^{*}}{\Bild{e_{X}}{X}}\bigr)} \subseteq
\beta Y.
\end{equation*}
Thus $e_{X}^{-1}\circ \phi_{f}^{*}$ is a continuous extension to $f\circ e_{X}$. 
\EndProof

It is immediate from Theorem~\ref{extension-of-stone-cech} that a
Stone-$\check{\mathrm{C}}$ech compactification is uniquely determined,
up to homeomorphism. This justifies the probably a bit prematurely used
characterization as \emph{the} Stone-$\check{\mathrm{C}}$ech compactification
above. 

Baire's Theorem, which we will establish now, states a property of
locally compact spaces which has a surprising range of
applications~---~it states that the intersection of dense open sets in
a locally compact $T_{3}$-space is dense again. This applies of course
to compact Hausdorf spaces as well. The theorem has a counterpart
for complete pseudometric spaces, as we will see below. For stating
and proving the theorem we lift the
assumption of working in a Hausdorff space, because it is really not necessary here.

\BeginTheorem{baire-locally-compact}
Let $X$ be a locally compact $T_{3}$-space. Then the intersection of
dense open sets is dense.\index{theorem!Baire!locally compact}
\EndTheorem

\BeginProof
Let $\Folge{D}$ be a sequence of dense open sets. Fix a non-empty open
set $G$, then we have to show that
$G\cap\bigcap_{n\in\Nat}D_{n}\not=\emptyset$. Now $D_{1}$ is dense and
open, hence we find an open set $V_{1}$ such that $\Closure{V}_{1}$ is
compact and $\Closure{V}_{1}\subseteq D_{1}\cap G$ by
Proposition~\ref{char-t3-ngbh}, since $X$ is a $T_{3}$-space. We
select inductively in this way a sequence of open sets $\Folge{V}$
with compact closure such that
$\Closure{V}_{n+1}\subseteq D_{n}\cap V_{n}$. This is possible since
$D_{n}$ is open and dense for each $n\in\Nat$.

Hence we have a decreasing sequence $\Closure{V}_{2} \supseteq \dots
\Closure{V}_{n} \supseteq \dots$ of closed sets in the compact set
$\Closure{V}_{1}$, thus $\bigcap_{n\in\Nat}\Closure{V}_{n} =
\bigcap_{n\in\Nat}V_{n}$ is not empty, which entails
$G\cap\bigcap_{n\in\Nat}D_{n}$ not being empty.
\EndProof

Just for the record:

\BeginCorollary{baire-compact}
The intersection of a sequence of dense open sets in a compact
Hausdorff space is dense.
\EndCorollary

\BeginProof
A compact Hausdorff space is normal by
Proposition~\ref{compact-is-normal}, hence regular by Proposition~\ref{char-t4-ngbh}, thus the assertion
follows from Theorem~\ref{baire-locally-compact}.
\EndProof

We give an example from Boolean algebras.

\BeginExample{bool-alg-dense}
Let $B$ be a Boolean algebra with $\wp_{B}$ as the set of all prime
ideals. Let $X_{a} := \{I\in \wp_{B}\mid a\not\in I\}$ be all prime
ideals which do not contain a given element $a\in B$, then
$\{X_{a}\mid a\in B\}$ is the basis for a compact Hausdorff topology
on $\wp_{B}$, and $a\mapsto X_{a}$ is a Boolean algebra isomorphism, see~\SetCite{Example 1.98}. 

Assume that we have a countable family $S$ of elements of $B$ with $a = \sup\ S\in B$, then we say that the prime ideal $I$ \emph{preserves the supremum} of $S$ iff $\Klasse{a}{\sim_{I}} = \sup_{s\in S}\ \Klasse{s}{\sim_{I}}$ holds.  Here $\sim_{I}$ is the equivalence relation induced by $I$, i.e., $\isEquiv{b}{b'}{\sim_{I}} \Leftrightarrow b\ominus b'\in I$ with $\ominus$ as the symmetric difference in $B$ (\SetCite{Sect. 1.5.7}).

We claim that the set $R$ of all prime ideals, which do \emph{not} preserve the supremum of this family, is closed and has an empty interior. Well, $R = X_{a}\setminus \bigcup_{k\in
  K}X_{a_{k}}$. Because the sets $X_{a}$ and $X_{a_{k}}$ are clopen, $R$ is closed. Assume that the interior of $R$ is not empty, then we find $b\in B$ with $X_{b}\subseteq R$, so that $X_{a_{k}}\subseteq X_{a}\setminus X_{b} = X_{a\wedge-b}$ for all $k\in K$. Since $a\mapsto X_{a}$ is an isomorphism, this means $a_{k}\leq a\wedge-b$, hence $\sup_{k\in K}\ a_{k}\leq a\wedge-b$ for all $k\in K$, thus $a = a\wedge -b$, hence $a\leq -b$. But then $X_{b} \subseteq X_{a}\subseteq X_{-b}$, which is certainly a contradiction. Consequently, the set of all prime ideal preserving this particular supremum is open and dense in $\wp_{B}$.

If we are given for each $n\in\Nat$ a family $S_{n}\subseteq B$ and $a_{0}\in B$ such
that 
\begin{itemize}
\item $a_{0}\not=\top$, the maximal element of $B$,
\item $a_{n} := \sup_{s\in S_{n}} s$ is an element of $B$ for each $n\in\Nat$.
\end{itemize}
Then we claim that there exists a prime ideal $I$ which contains $a_{0}$ and which
preserves all the suprema of $S_{n}$ for $n\in\Nat$. 

Let $P$ be the set of all prime
ideals which preserve all the suprema of the families above, then
\begin{equation*}
  P = \bigcap_{n\in\Nat}P_{n},
\end{equation*}
where $P_{n}$ is the set of all prime ideals which preserve the
supremum $a_{n}$, which is dense and open by the discussion above. Hence $P$ is dense by
Baire's Theorem (Corollary~\ref{baire-compact}). Since $X_{-a_{0}} = \wp_{B}\setminus X_{a_{0}}$ is
open and not empty, we infer that $P\cap X_{-a_{0}}$ is not empty,
because $P$ is dense. Thus we can select an arbitrary prime ideal from
this set. 
\EndExample

This example, which is taken from~\cite[Sect. 5]{Rasiowa-Sikorski-I}, will help in establishing Gödel's Completeness Theorem, see Section~\ref{sec:goedel}. The approach is typical for an application of Baire's Theorem~---~it is used to show that a set $P$, which is obtained from an intersection of countably many open and dense sets in a compact space, is dense, and that the object of one's desire is a member of $P$ intersecting an open set, hence this object must exist.

Having been carried away by Baire's Theorem, let us make some general
remarks. We have seen that local compactness is a somewhat weaker
property than compactness. Other notions of compactness have been
studied; an incomplete list for Hausdorff space $X$ includes
\begin{description}
\item[countably compact:] $X$ is called \emph{\index{compact!countably compact}countably compact} iff each
  countable open cover contains a finite subcover.
\item[Lindelöf space:] $X$ is a \emph{\index{compact!Lindelöf space}Lindelöf space} iff each open cover contains
  a countable subcover. 
\item[paracompactness:]  $X$ is said to be \emph{\index{compact!paracompact}paracompact}
  iff each open cover has a locally finite refinement. This explains
  it:
  \begin{itemize}
  \item An open cover ${\cal B}$ is a
    \emph{\index{compact!paracompact!refinement}refinement} of an open
    cover ${\cal A}$ iff each member of ${\cal B}$ is the subset of a
    member of ${\cal A}$.
  \item An open cover ${\cal A}$ is called \emph{\index{compact!paracompact!locally finite}locally finite} iff each
    point has a neighborhood which intersects a finite number of
    elements of ${\cal A}$.
  \end{itemize}
\item[sequentially compact:] $X$ is called \emph{\index{compact!sequentially
  compact}sequentially
  compact} iff each sequence has a convergent subsequence (we will
deal with this when discussing compact pseudometric spaces, see Proposition~\ref{seq-comp-equiv-comp}).
\end{description}
The reader is referred to~\cite[Chapter 3]{Engelking} for a
penetrating study.

\Subsection{Pseudometric and Metric Spaces} 
\label{sec:metric-and-pseudometric}

We turn to a class of spaces now in which we can determine the distance between any two points. This gives rise to a topology, declaring a set as open iff we can construct for each of its points an open ball which is entirely contained in this set. It is clear that this defines a topology, and it is also clear that having such a metric gives the space some special properties, which are not shared by general topological spaces. It also adds a sense of visual clearness, since an open ball is conceptually easier to visualize that an abstract open set. We will study the topological properties of these spaces now, starting with pseudometrics, with which we may measure the distance between two objects, but if the distance is zero, we cannot necessarily conclude that the objects are identical. This is a situation which occurs quite frequently when modelling an application, so it is sometimes more adequate to deal with pseudometric rather than metric spaces.

\BeginDefinition{def-pseudometric}
A map $d: X\times X\to \pReal$ is called a
\emph{\index{pseudometric}pseudometric on $X$} iff these
conditions hold
\begin{description}
\item[identity:] $d(x, x) = 0$ for all $x\in X$.
\item[symmetry:] $d(x, y) = d(y, x)$ for all $x, y\in X$,
\item[triangle inequality:] $d(x, y) \leq d(x, z) + d(z, y)$ for all
  $x, y, z\in X$. 
\end{description}
Then $(X, d)$ is called a \emph{\index{space!pseudometric}pseudometric space}. 
If, in addition, we have 
\begin{equation*}
  d(x, y) = 0 \Leftrightarrow x = y,
\end{equation*}
then $d$ is called a \emph{metric on $X$}; accordingly, $(X, d)$ is called
a \emph{\index{space!metric}metric space}.
\EndDefinition

The non-negative real number $d(x, y)$ is called the distance of the elements
$x$ and $y$ in a pseudometric space $(X, d)$. It is clear that one
wants to have that each point does have distance $0$ to itself, and
that the distance between two points is determined in a symmetric
fashion. The triangle inequality is intuitively clear as well:
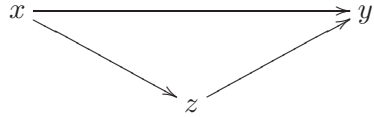
\begin{equation*}
\xymatrix{
x\ar[drr]\ar[rrrr] &&&& y\\
&&z\ar[rru]
}
\end{equation*}
Before proceeding, let us have a look at some examples. Some of them
will be discussed later on in greater detail. 

\BeginExample{for-metric-spaces}
\begin{enumerate}
\item Define for $x, y\in \Real$ the distance as $|x - y|$, hence as
  the absolute value of their difference. Then this defines a
  metric. Define, similarly, 
  \begin{equation*}
    d(x, y) := \frac{|x-y|}{1+|x-y|},
  \end{equation*}
then $d$ defines also a metric on $\Real$ (the triangle inequality
follows from the observation that $a\leq b \Leftrightarrow a/(1+a)\leq b/(1+b)$
holds for non-negative numbers $a$ and $b$). 

\item Given $x, y\in \Real^{n}$ for $n\in \Nat$, then 
  \begin{align*}
    d_{1}(x, y) & := \max_{1\leq i \leq n}|x_{i}-y_{i}|,\\
d_{2}(x, y) & := \sum_{i=1}^{n}|x_{i}-y_{i}|,\\
d_{3}(x, y) & := \sqrt{\sum_{i=1}^{n}(x_{i}-y_{i})^{2}}
  \end{align*}
define all metrics an $\Real^{n}$. Metric $d_{1}$ measures the maximal
distance between the components, $d_{2}$ gives the sum of the
distances, and $d_{3}$ yields the Euclidean, i.e., the
geometric, distance of the given points. The crucial property to be
established is in each case the triangle inequality. It follows for
$d_{1}$ and $d_{2}$ from the triangle inequality for the absolute
value, and for $d_{3}$ by direct computation.  
\item Given a set $X$, define
  \begin{equation*}
d(x, y) := 
    \begin{cases}
      0, & \text{ if } x = y\\
1, & \text{ otherwise}
    \end{cases}
  \end{equation*}
Then $(X, d)$ is a metric space, $d$ is called the
\emph{\index{metric!discrete}discrete metric}. Different points are
assigned the distance $1$, while each point has distance $0$ to
itself. 

\item Let $X$ be a set, ${\cal B}(X)$ be the set of all bounded maps
  $X\to \Real$. Define
  \begin{equation*}
    d(f, g) := \sup_{x\in X}|f(x) - g(x)|.
  \end{equation*}
Then $({\cal B}(X), d)$ is a metric space; the distance between
functions $f$ and $g$ is just their maximal difference.  

\item Similarly, given a set $X$, take a set ${\cal E}\subseteq{\cal
    B}(X)$ of bounded real valued functions as a set of evaluations and determine the distance of two
  points in terms of their evaluations:
  \begin{equation*}
    e(x, y) := \sup_{f\in {\cal F}}|f(x) - f(y)|.
  \end{equation*}
So two points are similar if their evaluations on terms of all
elements of ${\cal F}$ are close. This is a pseudometric on $X$, which
is not a metric if ${\cal F}$ does not separate points. 
\item Denote by $\Cont[[0, 1]]$ the set of all continuous real valued functions
  $[0, 1]\to \Real$, and measure the distance between $f, g\in \Cont[[0, 1]]$
  through
  \begin{equation*}
    d(f, g) := \sup_{0\leq x \leq 1}\ |f(x) - g(x)|.
  \end{equation*}
Because a continuous function on a compact space is bounded, $d(f, g)$
is always finite, and since for each $x\in[0, 1]$ the inequality
$|f(x) -g(x)| \leq |f(x) - h(x)| + |h(x)-g(x)|$ holds, the triangle
inequality is satisfied. Then $(\Cont[[0, 1]], d)$ is a metric space,
because $\Cont[[0, 1]]$ separates points. 
\item Define for the Borel sets $\Borel{[0, 1]}$ on the unit interval
  this distance:
  \begin{equation*}
    d(A, B) := \lambda(A\Delta B)
  \end{equation*}
with $\lambda$ as Lebesgue measure. Then $\lambda(A\Delta B) =
\lambda((A\Delta C)\Delta(C\Delta B) \leq \lambda(A\Delta C) + \lambda(C\Delta
B)$ implies the triangle inequality, so that $(\Borel{[0,
  1]}, d)$ is a pseudometric space. It is no metric space, however,
because $\lambda(\Rational\cap[0, 1]) = 0$, hence $d(\emptyset,
\Rational\cap[0, 1]) = 0$, but the latter set is not empty.
\item Given a non-empty set $X$ and a ranking function $r: X\to \Nat$,
  define the closeness $c(A, B)$ of two subset $A, B$ of $X$ as
  \begin{equation*}
    c(A, B) :=
    \begin{cases}
      +\infty, & \text{ if } A = B,\\
\inf\ \{r(w)\mid w\in A\Delta B\}, & \text{ otherwise}
    \end{cases}
  \end{equation*}
If $w\in A\Delta B$, then $w$ can be interpreted as a witness that $A$
and $B$ are different, and the closeness of $A$ and $B$ is just the minimal
rank of a witness. We observe these properties:
\begin{itemize}
\item $c(A, A) = +\infty$, and $c(A, B) = 0$ iff $A = B$ (because $A =
  B$ iff $A\Delta B =\emptyset$).
\item $c(A, B) = c(B, A)$,
\item $c(A, C) \geq \min\ \{c(A, B), c(B, C)\}$. If $A = C$, this is
  obvious; assume otherwise that $b\in A\Delta C$ is a witness of
  minimal rank. Since $A\Delta C = (A\Delta B)\Delta (B\Delta C)$, $b$
  must be either in $A\Delta B$ or $B\Delta C$, so that $r(b) \geq
  c(A, C)$ or $r(b) \geq c(B, C)$. 
\end{itemize}
Now put $d(A, B) := 2^{-c(A, B)}$ (with $2^{-\infty} := 0$). Then $d$
is a metric on $\PowerSet{X}$. This metric satisfies even
$d(A, B) \leq \max\ \{d(A, C), d(B, C)\}$ for an arbitrary $C$, hence
$d$ is an
\emph{\index{metric!ultrametric}\index{ultrametric}ultrametric}.
\item A similar construction is possible with a decreasing sequence of
  equivalence relation on a set $X$. In fact, let
  $(\rho_{n})_{n\in\Nat}$ be such a sequence, and put $\rho_{0} :=
  X\times X$. Define
  \begin{equation*}
    c(x, y) := 
    \begin{cases}
+\infty, &\text{ if }\langle x, y\rangle\in\bigcap_{n\in\Nat}\rho_{n}\\
      \max\ \{n\in\Nat\mid \langle x, y\rangle\in\rho_{n}\}, & \text{ otherwise}
    \end{cases}
  \end{equation*}
Then it is immediate that $c(x, y) \geq \min\ \{c(x, z), c(z,
y)\}$. Intuitively, $c(x, y)$ gives the degree of similarity of $x$
and $y$ --- the larger this value, the more similar $x$ and $y$ are. 
Then 
\begin{equation*}
  d(x, y) := 
  \begin{cases}
0, & \text{ if }c(x, y) = \infty\\
   2^{-c(x, y)}, & \text{ otherwise} 
  \end{cases}
\end{equation*}
defines a pseudometric. $d$ is a metric iff
$\bigcap_{n\in\Nat}\rho_{n} = \{\langle x, x\rangle\mid x \in X\}$. 
\end{enumerate}
\EndExample

Given a pseudometric space $(X, d)$, define for $x\in X$ and $r>0$ the
\emph{\index{$B(x, r)$}open ball $B(x, r)$}\MMP{$B(x, r)$} with center $x$ and radius $r$ as
\begin{equation*}
  B(x, r) := \{y\in X\mid d(x, y) < r\}.
\end{equation*}
The \emph{\index{$S(x, r)$}closed ball $S(x, r)$} is defined similarly as 
\begin{equation*}
  S(x, r) := \{ y\in X \mid d(x, y)\leq r\}.
\end{equation*}
If necessary, we indicate the pseudometric explicitly with $B$ and
$S$. Note that $B(x, r)$ is open, and $S(x, r)$ is closed, but that the closure
$\Closure{B(x, r)}$ of $B(x, r)$ may be properly contained in the closed ball
$S(x, r)$ (let $d$ be the discrete metric, then $B(x, 1) = \{x\} =
\Closure{B(x, 1)}$, but $S(x, 1) = X$, so both closed set do not
coincide if $X$ has more than one point).

Call $G\subseteq X$ open iff we can find for each $x\in G$ some $r>0$ such
that $B(x, r)\subseteq G$. Then this defines the \emph{pseudometric topology}
on $X$. It has the set $\beta := \{B(x, r) \mid x\in X, r > 0\}$ of
open balls as a basis. Let us have a look at the properties a base is
supposed to have. Assume that $x\in B(x_{1}, r_{1})\cap B(x_{2},
r_{2})$, and select $r$ with $0 < r < \min\{r_{1}-d(x, x_{1}),
r_{2}-d(x, x_{2})\}$. Then $B(x, r)\subseteq B(x_{1}, r_{1})\cap
B(x_{2}, r_{2})$, because we have for $z\in B(x, r)$ 
\begin{equation}
\label{triang-equ}
  d(z, x_{1}) \leq d(z, x) + (x, x_{1}) < r + d(x, x_{1}) \leq (r_{1}-d(x, x_{1})) +
  d(x, x_{1}) = r_{1},
\end{equation}
by the triangle inequality; similarly, $d(x, x_{2}) < r_{2}$. Thus it
follows from  Proposition~\ref{when-is-a-base} that $\beta$ is in fact
a base. 

Call two pseudometrics on $X$
\emph{\index{pseudometrics!equivalent}equivalent} iff they generate
the same topology. An equivalent formulation goes like this. Let
$\tau_{i}$ be the topologies generated from pseudometrics $d_{i}$ for
$i = 1, 2$, then $d_{1}$ and $d_{2}$ are equivalent iff the identity
$(X, \tau_{1})\to (X, \tau_{2})$ is a homeomorphism. These are two
common methods to construct equivalent pseudometrics.

\BeginLemma{gen-equiv-pseudo-metrics}
Let $(X, d)$ be a pseudometric space. Then 
\begin{align*}
  d_{1}(x, y) & := \max\{d(x, y), 1\},\\
d_{2}(x, y) & := \frac{d(x, y)}{1 + d(x, y)}
\end{align*}
both define pseudometrics which are equivalent to $d$. 
\EndLemma

\BeginProof
It is clear that both $d_{1}$ and $d_{2}$ are pseudometrics (for
$d_{2}$, compare Example~\ref{for-metric-spaces}). Let $\tau,
\tau_{1}, \tau_{2}$ be the respective topologies, then it is immediate
that $(X, \tau)$ and $(X, \tau_{1})$ are homeomorphic. Since $d_{2}(x,
y) < r$ iff $d(x, y) < r/(1-r)$, provided $0 < r < 1$, we obtain also
that $(X, \tau)$ and $(X, \tau_{2})$ are homeomorphic. 
\EndProof

These pseudometrics have the advantage that they are bounded, which is
sometimes quite practical for establishing topological
properties. Just as a point in case:

\BeginProposition{countable-product-is pseudo}
Let $(X_{n}, d_{n})$ be a pseudometric space with associated topology
$\tau_{n}$.  Then the topological product $\prod_{n\in\Nat}(X_{n},
\tau_{n})$ is a pseudometric space again.
\EndProposition

\BeginProof
1.
We may assume that each $d_{n}$ is bounded by $1$, otherwise we select
an equivalent pseudometric with this property
(Lemma~\ref{gen-equiv-pseudo-metrics}). Put
\begin{equation*}
  d\bigl((x_{n})_{n\in\Nat}, (y_{n})_{n\in\Nat}\bigr) :=
  \sum_{n\in\Nat}2^{-n}\cdot d_{n}(x_{n}, y_{n}).
\end{equation*}
We claim that the product topology is the topology induced by the
pseudometric $d$ (it is obvious that $d$ is one). 

2.
Let $G_{i}\subseteq X_{i}$ open for $1\leq i\leq k$, and assume that
$x\in G := G_{1}\times \dots\times G_{k}\times\prod_{n>k}X_{n}$. We can
find for $x_{i}\in G_{i}$ some positive $r_{i}$ with $B_{d_{i}}(x_{i},
r_{i})\subseteq G_{i}$. Put $r := \min\{r_{1}, \dots, r_{k}\}$, then
certainly $B_{d}(x, r) \subseteq G$. This implies that each element of
the base for the product topology is open with respect to $d$.

3.
Given the sequence $x$ and $r> 0$, take $y\in B_{d}(x, r)$. Put $t :=
r - d(x, y) > 0$. Select $m \in \Nat$ with $\sum_{n>m}2^{-n}< t/2$, and
let $G_{n} := B_{d_{n}}(y_{n}, t/2)$ for $n\leq m$. If $z\in U :=
G_{1}\times \dots\times G_{n}\times\prod_{k>m}X_{k}$, then 
\begin{align*}
  d(x, z) & \leq 
d(x, y) + d(y, z)\\
& \leq 
r - t + \sum_{n=1}^{m}2^{-n}d_{n}(y_{n}, z_{n}) + \sum_{n>m}2^{-n}\\
& <
r - t + t/2 + t/2\\
& = 
r,
\end{align*}
so that $U\subseteq B_{d}(x, r)$. Thus each open ball is open in the
product topology.
\EndProof

One sees immediately that the pseudometric $d$ constructed above is a
metric, provided each $d_{n}$ is one. Thus

\BeginCorollary{countable-metric-is metric}
The countable product of metric spaces is a metric space in the
product topology.
\QED
\EndCorollary

One expects that each pseudometric space can be made a metric space by
identifying those elements which cannot be separated by the
pseudometric. Let's try:

\BeginProposition{pseudo-to-metric}
Let $(X, d)$ be a pseudometric space, and define $\isEquiv{x}{y}{\sim}$
iff $d(x, y) = 0$ for $x, y\in X$. Then the factor space $\Faktor{X}{\sim}$ is a
metric space with metric $D(\Klasse{x}{\sim}, \Klasse{y}{\sim}) := d(x, y)$. 
\EndProposition

\BeginProof
1.
Because $d(x, x') = 0$ and $d(y, y') = 0$ implies $d(x, y) = d(x',
y')$, $D$ is well-defined, and it is clear that it has all the
properties of a pseudometric. $D$ is also a metric, since
$D(\Klasse{x}{\sim}, \Klasse{y}{\sim}) = 0$ is equivalent to $d(x, y)
= 0$, hence to $\isEquiv{x}{y}{\sim}$, thus to $\Klasse{x}{\sim} =
\Klasse{y}{\sim}$. 

2.
The metric topology is the final topology with respect to the factor
map $\fMap{\sim}$. To establish this, take a map $f:
\Faktor{X}{\sim}\to Z$ with a topological space $Y$. Assume that
$\InvBild{(f\circ \fMap{\sim})}{G}$ is open for $G\subseteq Y$
open. If $\Klasse{x}{\sim}\in \InvBild{f}{G}$, we have
$x\in\InvBild{(f\circ \fMap{\sim})}{G}$, thus there exists
$r>0$ with $B_{d}(x, r)\subseteq
\InvBild{\fMap{\sim}}{\InvBild{f}{G}}$. But this means that $B_{D}(\Klasse{x}{\sim},
r)\subseteq \InvBild{f}{U}$, so that the latter set is open. Thus if
$f\circ \fMap{\sim}$ is continuous, $f$ is. The converse is
established in the same way. This implies that the metric topology is
final with respect to the factor map $\fMap{\sim}$, cp. Proposition~\ref{initial-final}. 
\EndProof

We want to show that a pseudometric space satisfies the $T_{4}$-axiom
(hence that a metric space is normal). So we take two disjoint
closed sets and need to produce two disjoint open sets, each of which
containing one of the closed sets. The following construction is
helpful.

\BeginLemma{distance-is-continuous}
Let $(X, d)$ be a pseudometric space\MMP{$d(x, A)$}. Define the distance of point
$x\in X$ to $\emptyset\not=A\subseteq X$ through
\begin{equation*}
  d(x, A) := \inf_{y\in A}d(x, y).
\end{equation*}
Then $d(\cdot, A)$ is continuous.
\EndLemma

\BeginProof
Let $x, z\in X$, and $y\in A$, then $d(x, y) \leq d(x, z) + d(z,
y)$. Now take lower bounds, then $d(x, A) \leq d(x, z) + d(z,
A)$. This yields $d(x, A) - d(z, A) \leq d(x, z)$. Interchanging the
r\^oles of $x$ and $z$ yields  $d(z, A) - d(x, A) \leq d(z, x)$, thus
$|d(x, A) - d(z, A)|\leq d(x, z)$. This implies continuity of
$d(\cdot, A)$. 
\EndProof

Given a closed set $A\subseteq X$, we find that $A = \{x\in X\mid d(x, A) =
0\}$; we can say a bit more: 

\BeginCorollary{closure-through-distance}
Let $X, A$ be as above, then $\Closure{A} = \{x\in X\mid d(x, A) =
0\}$. 
\EndCorollary

\BeginProof
Since $\{x\in X\mid d(x, A) = 0\}$ is closed, we infer
that $\Closure{A}$ is contained in this set. If, in the other hand, $x\not\in
\Closure{A}$, we find $r>0$ such that $B(x, r)\cap A=\emptyset$, hence $d(x, A)
\geq r$. Thus the other inclusion holds as well.
\EndProof

Armed with this observation, we can establish now

\BeginProposition{pseudo-metric-is t4}
A pseudometric space $(X, d)$ is a $T_{4}$-space. 
\EndProposition

\BeginProof
Let $F_{1}$ and $F_{2}$ be disjoint closed subsets of $X$. Define 
\begin{equation*}
  f(x) := \frac{d(x, F_{1})}{d(x, F_{1}) + d(x, F_{2})},
\end{equation*}
then Lemma~\ref{distance-is-continuous} shows that $f$ is continuous,
and Corollary~\ref{closure-through-distance} indicates that the denominator will not vanish,
since $F_{1}$ and $F_{2}$ are disjoint. It is immediate that
$F_{1}$ is contained in the open set  $\{x\mid f(x) < 1/2\}$, that
$F_{2}\subseteq \{x\mid f(x) > 1/2\}$, and that these open sets are disjoint.
\EndProof

Note that a pseudometric $T_{1}$-space is already a metric space
(Exercise~\ref{ex-pseudo-t1-is-metric}). 

Define for $r>0$ the $r$-neighborhood $A^{r}$ of set $A\subseteq X$ as\MMP{$A^{r}$} 
\begin{equation*}
  A^{r} := \{x\in X\mid d(x, A) < r\}.
\end{equation*}
This makes of course only sense if $d(x, A)$ is finite. Using the
triangle inequality, one calculates $(A^{r})^{s}\subseteq
A^{r+s}$. This observation will be helpful when we look at the next
example.

\BeginExample{vietoris-pseudometric}
Let $(X, d)$ be a pseudometric space, and let 
\begin{equation*}
  \allCompact(X) := \{C\subseteq X\mid C\text{ is compact and not empty}\}
\end{equation*}
be the set of all compact and not empty subsets of $X$. Define 
\begin{equation*}
  \delta_{H}(C, D) := \max\ \{\max_{x\in C}\ d(x, D),
  \max_{x\in D}\ d(x, C)\}
\end{equation*}
for $C, D\in\allCompact(X)$\MMP{$\delta_{H}$}. We claim that $\delta_{H}$ is a pseudometric
on $\allCompact(X)$, which is a metric if $d$ is a metric on $X$. 

One notes first that 
\begin{equation*}
  \delta_{H}(C, D) = \inf\ \{r > 0 \mid C\subseteq D^{r}, D\subseteq C^{r}\}.
\end{equation*}
This follows easily from $C\subseteq D^{r}$ iff $\max_{x\in C}\ d(x,
D) < r$. Hence we obtain that $\delta_{H}(C, D) \leq r$ and
$\delta_{H}(D, E) \leq s$ together imply $\delta_{H}(C, E) \leq r +
s$, which implies the triangle inequality. The other laws for a
pseudometric are obvious. $\delta_{H}$ is called the \emph{\index{metric!Hausdorff}Hausdorff pseudometric}. 

Now assume that $d$ is a metric, and assume $\delta_{H}(C, D) =
0$. Thus $C\subseteq \bigcap_{n\in\Nat}D^{1/n}$ and $D\subseteq
\bigcap_{n\in\Nat}C^{1/n}$. Because $C$ and $D$ are closed, and $d$ is
a metric, we obtain $C = D$, thus $\delta_{H}$ is a metric, which is
accordingly called the \emph{Hausdorff metric}. 
\EndExample

Let us take a magnifying glass and have a look at what happens
locally in a point of a pseudometric space. Given $U\in\upsilon(x)$,
we find an open ball $B(x, r)$ which is contained in $U$, hence we
find even a rational number $q$ with $B(x, q)\subseteq B(x, r)$. But
this means that the open balls with rational radii form a basis for
the neighborhood filter of $x$. This is sometimes also the case in
more general topological spaces, so we define this and two related
properties for topological rather than pseudometric spaces.

\BeginDefinition{space-separable}
A topological space 
\begin{enumerate}
\item satisfies the \emph{first axiom of countability} (and the space is called in this
  case \emph{\index{topology!first countable}first countable}) iff the neighborhood filter of each point has a
  countable base of open sets,
\item satisfies the \emph{second axiom of countability} (the space is called in this
  case \emph{\index{topology!second countable}second countable}) iff
  the topology has a countable base,
\item is \emph{\index{topology!separable}separable} iff it has a countable dense subset.
\end{enumerate}
\EndDefinition

The standard example for a separable topological space is of course
$\Real$, where the rational numbers $\Rational$ form a countable dense
subset.

This is a trivial consequence of the observation just made.

\BeginProposition{pseudo-is-first-count}
A pseudometric space is first countable. 
\QED
\EndProposition

In a pseudometric space separability and satisfying the second axiom
of countability coincide, as the following observation shows.

\BeginProposition{separable-iff-2ndcountable}
A pseudometric space $(X, d)$ is second countable iff it has a
countable dense subset.
\EndProposition

\BeginProof
1.
Let $D$ be a countable dense subset, then 
\begin{equation*}
  \beta := \{B(x, r)\mid x\in D, 0<r\in\Rational\}
\end{equation*}
is a countable base for the topology. For, given $U\subseteq X$ open,
there exists $d\in D$ with $d\in U$, hence we can find a rational
$r>0$ with $B(d, r)\subseteq U$. On the other hand, one shows exactly
as in the argumentation leading to Eq.~(\ref{triang-equ}) on
page~\pageref{triang-equ} that $\beta$ is a base. 

2.
Assume that $\beta$ is a countable base for the topology, pick from
each $B\in\beta$ an element $x_{B}$. Then $\{x_{B}\mid B\in\beta\}$ is
dense: given an open $U$, we find $B\in \beta$ with $B\subseteq U$,
hence $x_{B}\in U$. This argument does not require $X$ being a
pseudometric space (but the Axiom of Choice).
\EndProof 

We know from Exercise~\ref{ex-closure-filter} that a point $x$ in a
topological space is in the closure of a set $A$ iff there exists a
filter $\fiF$ with $i_{A}(\fiF)\to x$ with $i_{A}$ as the injection
$A\to X$. In a first countable space, in particular in a pseudometric
space, we can work with sequences rather than filters, which is
sometimes more convenient.

\BeginProposition{sequences-are-enough}
Let $X$ be a first countable topological space, $A\subseteq X$. Then
$x\in\Closure{A}$ iff there exists a sequence $(x_{n})_{n\in\Nat}$ in
$A$ with $x_{n}\to x$. 
\EndProposition

\BeginProof
If there exists a sequence $(x_{n})_{n\in\Nat}$ which converges to $x$
such that $x_{n}\in A$ for
all $n\in\Nat$, then the corresponding filter converges to $x$, so we
have to establish the converse statement.

Now let $(U_{n})_{n\in\Nat}$ be the basis of the neighborhood filter
of $x\in\Closure{A}$, and $\fiF$ be a filter with $i_{A}(\fiF)\to
x$. Put $V_{n} := U_{1}\cap\dots\cap U_{n}$, then $V_{n}\cap A \in
i_{A}(\fiF)$. The sequence $(V_{n})_{n\in\Nat}$ decreases, and forms
a basis for the neighborhood filter of $x$. Pick from each $V_{n}$ an
element $x_{n}\in A$, and take a neighborhood $U\in\upsilon(x)$. Since
there exists $n$ with $V_{n}\subseteq U$, we infer that $x_{m}\in U$
for all $m\leq n$, hence $x_{n}\to x$. 
\EndProof

A second countable normal space $X$ permits the following remarkable
construction. Let $\beta$ be a countable base for $X$, and define
${\cal A} := \{\langle U, V\rangle\mid  U, V\in \beta,
\Closure{U}\subseteq V\}$. Then ${\cal A}$ is countable as well, and
we can find for each pair $\langle U, V\rangle\in {\cal A}$ a
continuous map $f: X\to [0, 1]$ with $f(x) = 0$ for all $x\in U$ and
$f(x) = 1$ for all $x\in X\setminus V$. This is a consequence of
Urysohn's Lemma (Theorem~\ref{urysohns-lemma}). The collection ${\cal
  F}$ of all
these functions is countable, because ${\cal A}$ is countable. Now
define the embedding map
\begin{equation*}
e:
  \begin{cases}
   X & \to [0, 1]^{{\cal F}}\\
x & \mapsto (f(x))_{f\in {\cal F}} 
  \end{cases}
\end{equation*}
We endow the space $[0, 1]^{{\cal F}}$ with the product topology,
i.e., with the initial topology with respect to all projections
$\pi_{f}: x \mapsto f(x)$. Then we observe these properties
\begin{enumerate}
\item The map $e$ is continuous. This is so because $\pi_{f}\circ e =
  f$, and $f$ is continuous, hence we may infer continuity from
  Proposition~\ref{initial-final}.
\item The map $e$ is injective. This follows from Urysohn's Lemma
  (Theorem~\ref{urysohns-lemma}), since two distinct points constitute
  two disjoint closed sets.
\item If $G\subseteq X$ is open, $\Bild{e}{G}$ is open in
  $\Bild{e}{X}$. In fact, let $e(x)\in \Bild{e}{G}$. We find an open
  neighborhood $H$ of $e(x)$ in $[0, 1]^{{\cal F}}$ such that
  $\Bild{e}{X}\cap H\subseteq \Bild{e}{G}$ in the following way: we
  infer from the construction that we can find a map $f\in {\cal F}$
  such that $f(x) = 0$ and $f(y) = 1$ for all $y\in X\setminus G$,
  hence $f(x)\not\in\Closure{\Bild{f}{X\setminus G}}$; hence the 
  set $H := \{y\in [0, 1]^{{\cal F}}\mid y_{f}\not\in\Bild{f}{X\setminus G}\}$ is
  open in $[0, 1]^{{\cal F}}$, and $H\cap \Bild{e}{X}$ is contained in
  $\Bild{e}{G}$.
\item $[0, 1]^{{\cal F}}$ is a metric space by
  Corollary~\ref{countable-metric-is metric}, because the unit
  interval $[0, 1]$ is a metric space, and because ${\cal F}$ is
  countable.
\end{enumerate}
Summarizing, $X$ is homeomorphic to a subspace of $[0, 1]^{{\cal
    F}}$. This is what \emph{\index{theorem!Urysohn's
    Metrization}Urysohn's Metrization Theorem} says.

\BeginProposition{normal-embed-metric}
A second countable normal topological space is metrizable.
\QED.
\EndProposition

The problem of metrization of topological spaces is non-trivial, as
one can see from Proposition~\ref{normal-embed-metric}. The reader who
wants to learn more about it may wish to consult Kelley's
textbook~\cite[p. 124 f]{Kelley} or Engelking's treatise~\cite[4.5,
5.4]{Engelking}. 

\Subsubsection{Completeness}
\label{sec:completeness}

Fix in this section a pseudometric space $(X, d)$. A \emph{\index{Cauchy sequence}Cauchy \index{sequence!Cauchy}sequence}
$\Folge{x}$ is defined in $X$ just as in $\Real$: Given $\epsilon>0$,
there exists an index $n\in \Nat$ such that $d(x_{m}, x_{m'}) <
\epsilon$ holds for all $m, m'\geq n$. 

Thus we have a Cauchy sequence, when we know that eventually the
members of the sequence will be arbitrarily close; a converging
sequence is evidently a Cauchy sequence. But a sequence which
converges requires the knowledge of its limit; this is sometimes a
disadvantage in applications. It would be helpful if we could conclude
from the fact that we have a Cauchy sequence that we also have a point
to which it converges. Spaces for which this is always guaranteed are
called complete; they will be introduced next, examples show that
there are spaces which are not complete; note, however, that we can
complete each pseudometric space. This will be considered in some detail later on. 

\BeginDefinition{pseudo-is-complete}
The pseudometric space is said to be
\emph{\index{space!pseudometric!complete}complete} iff each Cauchy
sequence has a limit. 
\EndDefinition

Compare in a pseudometric space the statement
$\lim_{n\to \infty} x_{n} = x$ with the statement that $\Folge{x}$ is
a Cauchy sequence. The former requires the knowledge of the limit
point, while the latter is derived from observing the members of the
sequence, but without knowing a limit. Hence we know in a complete
space that a limit will exist, without being obliged to identify
it. This suggests that complete pseudometric spaces are important.

It is well known that the rational numbers are not complete, which is
usually shown by showing that $\sqrt{2}$ is not rational. Another
instructive example proposed by Bourbaki~\cite[II.3.3]{Bourbaki} is the
following.

\BeginExample{rationals-are-not-complete}
The rational numbers $\Rational$ are not complete in the usual
metric. Take 
\begin{equation*}
  x_{n} := \sum_{i=0}^{n}2^{-i\cdot (i+1)/2}.
\end{equation*}
Then $\Folge{x}$ is a Cauchy sequence in $\Rational$: if $m> n$, then 
$
  |x_{m}-x_{n}| \leq 2^{-(n+3)/2}
$
(this is shown easily through the well known identity $\sum_{i=0}^{p}i
= p\cdot (p+1)/2$). Now assume that the sequence converges to
$a/b\in\Rational$, then we can find an integer $h_{n}$ such that
\begin{equation*}
  \bigl|\frac{a}{b} - \frac{h_{n}}{2^{n\cdot (n+1)/2}}\bigr| \leq \frac{1}{2^{n\cdot (n+3)/2}},
\end{equation*}
yielding
\begin{equation*}
  |a\cdot 2^{n\cdot (n+1)/2} - b\cdot h_{n}| \leq \frac{b}{2^{n}}
\end{equation*}
for all $n\in\Nat$. The left hand side of this inequality is a whole
number, the right side is not, once $n>n_{0}$ with $n_{0}$ so large that
$b<2^{n}$. This means that the left hand side must be zero, so that
$a/b = x_{n}$ for $n>n_{0}$. This is a contradiction. 
\EndExample

We know that $\Real$ is complete with the usual metric, the rationals
are not. But there is a catch: if we change the metric, completeness
may be lost.

\BeginExample{loose-completeness}
The half open interval $]0, 1]$ is not complete under the usual metric
$d(x, y) :=  |x - y|$. But take the metric 
\begin{equation*}
 d'(x, y) := \bigl|\frac{1}{x}-\frac{1}{y}\bigr| 
\end{equation*}
Because $a < x < b$ iff $1/b < 1/x < 1/a$ holds for $0< a \leq b \leq
1$, the metrics $d$ and $d'$ are equivalent on $]0, 1]$. Let
$\Folge{x}$ be a $d'$-Cauchy sequence, then $(1/x_{n})_{n\in\Nat}$ is
a Cauchy sequence in $(\Real, |\cdot |)$, hence it converges, so
that $\Folge{x}$ is $d'$-convergent in $]0, 1]$. 

The trick here is to make sure that a Cauchy sequence avoids the
region around the critical value $0$. 
\EndExample

Thus we have to carefully stick to the given metric, and changing the
metric always entails checking completeness properties, if they are
relevant. 

\BeginExample{cont-is-complete}
Endow the set $\Cont[[0, 1]]$ of continuous functions on the unit
interval with the metric
$d(f, g) := \sup_{0\leq x\leq 1}\ |f(x)-g(x)|$, see
Example~\ref{for-metric-spaces}. We claim that this metric space is
complete. In fact, let $\Folge{f}$ be a $d$-Cauchy sequence in
$\Cont[[0, 1]]$. Because we have for each $x\in [0, 1]$ the inequality
$ |f_{n}(x) - f_{m}(x)| \leq d(f_{n}, f_{m}), $ we conclude that
$\bigl(f_{n}(x)\bigr)_{n\in\Nat}$ is a Cauchy sequence for each
$x\in[0, 1]$, which converges to some $f(x)$, since $\Real$ is
complete. We have to show that $f$ is continuous, and that
$d(f, f_{n})\to 0$.

Let $\epsilon>0$ be given, then there exists $n\in\Nat$ such that
$d(f_{m}, f_{m'})< \epsilon/2$ for $m, m'\geq n$; hence we have
$|f_{m}(x) - f_{m'}(x')| \leq |f_{m}(x) - f_{m}(x')| +
|f_{m}(x')-f_{m'}(x')| \leq |f_{m}(x) - f_{m}(x')| + d(f_{m},
f_{m'})$. Choose $\delta>0$ so that $|x-x'|<\delta$ implies
$|f_{m}(x)-f_{m}(x')| < \epsilon/2$, then $|f_{m}(x) - f_{m'}(x')| <
\epsilon$ for $m, m'\geq n$. But this means $|x-x'|<\delta$ implies
$|f(x)-f(x')|\leq \epsilon$. Hence $f$ is continuous. Since $\bigl(\{x\in
[0, 1]\mid |f_{n}(x)-f(x)|\leq \epsilon\}\bigr)_{n\in\Nat}$ constitutes an open cover
of $[0, 1]$, we find a finite cover given by $n_{1}, \dots, n_{k}$;
let $n'$ be the smallest of these numbers, then $d(f, f_{n})\leq
\epsilon$ for all $n\geq n'$, hence $d(f, f_{n})\to 0$. 
\EndExample

The next example is inspired by an observation in~\cite{MacQueen+Plotkin+Sethi}. 
\BeginExample{ranking-complete}
Let $r: X\to \Nat$ be a ranking function, and denote the (ultra-) metric on
$\PowerSet{X}$ constructed from it by $d$, see
Example~\ref{for-metric-spaces}. Then $(\PowerSet{X}, d)$ is
complete. In fact, let $\Folge{A}$ be a Cauchy sequence, thus we find
for each $m\in \Nat$ an index $n\in\Nat$ such that $c(A_{k}, A_{\ell})
\geq m$, whenever $k, \ell\geq n$. We claim that the sequence
converges to
\begin{equation*}
  A := \bigcup_{n\in\Nat}\bigcap_{k\geq n}A_{k},
\end{equation*}
which is the set of all elements in $X$ which are contained in all but a finite
number of sequence elements. Given $m$, fix $n$ as above; we show that
$c(A, A_{k}) > m$, whenever $k>n$. Take an element $x\in A\Delta
A_{k}$ of minimal rank. 
\begin{itemize}
\item If $x\in A$, then there exists $\ell$ such that $x\in A_{t}$ for
  all $t\geq \ell$, so take $t\geq \max\ \{\ell, n\}$, then $x\in A_{t}\Delta
  A_{k}$, hence $c(A, A_{k}) = r(x) \geq c(A_{t}, A_{k}) > m$.
\item If, however, $x\not\in A$, we conclude that $x\not\in A_{t}$ for
  infinitely many $t$, so $b\not\in A_{t}$ for some $t>n$. But since
  $x\in A\Delta A_{k}$, we conclude $x\in A_{k}$, hence $x\in
  A_{k}\Delta A_{t}$, thus $c(A, A_{k}) = r(x) \geq c(A_{k}, A_{t}) >
  m$. 
\end{itemize}
Hence $A_{n}\to A$ in $(\PowerSet{X}, d)$. 
\EndExample

This observation is trivial, but sometimes helpful.

\BeginLemma{closed-is-complete}
A closed subset of a complete pseudometric space is complete.
\QED
\EndLemma

If we encounter a pseudometric space which is not complete, we may
complete it through the following construction. Before discussing it,
we need a simple auxiliar statement, which says that we can check
completeness already on a dense subset.

\BeginLemma{check-on-dense}
Let $D\subseteq X$ be dense. Then the space is complete iff each
Cauchy sequence on $D$ converges.
\EndLemma

\BeginProof
If each Cauchy sequence from $X$ converges, so does each such sequence
from $D$, so we have to establish the converse. Let $\Folge{x}$ be a
Cauchy sequence on $X$. Given $n\in\Nat$, there exists for $x_{n}$ an
element $y_{n}\in D$ such that $d(x_{n}, y_{n})< 1/n$. Because
$\Folge{x}$ is a Cauchy sequence, $\Folge{y}$ is one as well, which
converges by assumption to some $x\in X$; the triangle inequality
shows that $\Folge{x}$ converges to $x$ as well.
\EndProof

This helps in establishing that each pseudometric space can be
embedded into a complete pseudometric space. The approach may be
described as \index{Charly Brown's device}Charly Brown's device ---
``If you can't beat them, join them''. So we take all Cauchy sequences
as our space into which we embed $X$, and ---~intuitively~--- we flesh out
from a Cauchy sequence of these sequences the diagonal sequence,
which then will be a Cauchy sequence as well, and which will be a
limit of the given one. This sounds more complicated than it is, however,
because fortunately Lemma~\ref{check-on-dense} makes life easier, when it comes to
establish completeness. Here we go.

\BeginProposition{embed-into-complete-space}
There exists a complete pseudometric space $(X^{*}, d^{*})$ into which
$(X, d)$ may be embedded isometrically as a dense subset.
\EndProposition

\BeginProof
0.
This  is the line of attack: We define $X^{*}$ and $d^{*}$, show
that we can embed $X$ isometrically into it as a dense subset, and
then we establish completeness with the help of Lemma~\ref{check-on-dense}\MMP{Fairly direct approach}.

1.
Define
\begin{equation*}
  X^{*}  := \{\Folge{x}\mid \Folge{x} \text{ is a $d$-Cauchy sequence
          in }X\},
\end{equation*}
and put 
\begin{equation*}
d^{*}\bigl(\Folge{x}, \Folge{y}\bigr) := \lim_{n\to \infty}d(x_{n}, y_{n})
\end{equation*}
Before proceeding, we should make sure that the limit in question
exists. In fact, given $\epsilon>0$, there exists $n\in\Nat$ such that
$d(x_{m'}, x_{m}) < \epsilon/2$ and $d(y_{m'}, y_{m}) < \epsilon/2$ for
$m, m'\geq n$, thus, if $m, m'\geq n$, we obtain
$$
d(x_{m}, y_{m}) \leq d(x_{m}, x_{m'}) + d(x_{m'}, y_{m'}) + d(y_{m'},
y_{m}) < d(x_{m'}, y_{m'}) + \epsilon,
$$
interchanging the r\^oles of $m$ and $m'$ yields
$$
|d(x_{m}, y_{m}) - d(x_{m'}, y_{m'})| < \epsilon
$$
for $m, m'\geq n$. Hence $(d(x_{n}, y_{n}))_{n\in\Nat}$ is a Cauchy sequence in
$\Real$, which converges by completeness of $\Real$. 

2.  
Given $x\in X$, the sequence $(x)_{n\in\Nat}$ is a Cauchy
sequence, so it offers itself as the image of $x$; let $e:X\to X^{*}$
be the corresponding map, which is injective, and it preserves the
pseudometric. Hence $e$ is continuous. We show that $\Bild{e}{X}$ is
dense in $X^{*}$: take a Cauchy sequence $\Folge{x}$ and
$\epsilon>0$. Let $n\in\Nat$ be selected for $\epsilon$, and assume
$m\geq n$. Then
\begin{equation*}
  D(\Folge{x}, e(x_{m})) = \lim_{n\to \infty}d(x_{n}, x_{m}) < \epsilon.
\end{equation*}

3.
The crucial point is completeness. An appeal to
Lemma~\ref{check-on-dense} shows that it is sufficient to show that a
Cauchy sequence in $\Bild{e}{X}$ converges in $(X^{*}, d^{*})$, because
$\Bild{e}{X}$ is dense. But this is trivial. 
\EndProof

Having  the completion $X^{*}$ of a pseudometric space $X$ at one's disposal, one might be
tempted to extend a continuous map $X\to Y$ to a continuous
map $X^{*}\to Y$ for example in the case that $Y$ is complete. This is
usually not possible, for example, not every continuous function
$\Rational\to \Real$ has a continuous extension. We will deal with
this problem when discussing uniform continuity below, but we will
state and prove here a condition which is sometime helpful when one
wants to extend a function not to the whole completion, but to a
domain which is somewhat larger than the given one. Define the \index{diameter}\emph{diameter} \index{$\mathsf{diam}(A)$}$\mathsf{diam}(A)$
of a set $A$ as\MMP{$\mathsf{diam}(A)$}
\begin{equation*}
  \mathsf{diam}(A) := \sup\ \{d(x, y)\mid x, y\in A\}
\end{equation*}
(note that the diameter may be infinite). It is easy to see that
$\mathsf{diam}(A) = \mathsf{diam}(\Closure{A})$ using
Proposition~\ref{sequences-are-enough}. Now assume that $f: A\to Y$ is
given, then we measure the discontinuity of $f$ at point $x$ through
the \label{ref:oscillation}\emph{oscillation} ${\o}_f(x)$\index{oscillation}\MMP{Oscillation} of $f$ at $x \in
\Closure{A}$, which is defined as the smallest
diameter of the image of an open neighborhood of $x$, formally,
\begin{equation*}
{\o}_f(x) := \inf\{\mathsf{diam}(\Bild{f}{A \cap V}) \mid x \in V, V \text{ open}\}.
\end{equation*}
If $f$ is continuous on $A$, we have ${\o}_f(x) = 0$ for each
element $x$ of $A$. In fact, let $\epsilon>0$ be given, then there
exists $\delta>0$ such that
$\mathsf{diam}(\Bild{f}{A \cap V}) < \epsilon$, whenever $V$ is a
neighborhood of $x$ of diameter less than $\delta$. Thus
${\o}_{f}(x) < \epsilon$; since $\epsilon>0$ was chosen to be
arbitrary, the claim follows.

\BeginLemma{Kuratowski}
Let $Y$ be a complete metric space, $X$ a pseudometric space, then a continuous
map $f: A \rightarrow Y$ can be extended to a continuous map
$
f_*: G \rightarrow Y,
$
where 
$
G := \{x \in \Closure{A} \mid {\o}_f(x) = 0\}
$
has these properties:\MMP{Extension}
\begin{enumerate}
\item $A \subseteq G \subseteq \Closure{A}$,
\item $G$ can be written as the intersection of countably many open sets.
\end{enumerate}
\EndLemma

The basic idea for the proof is rather straightforward. Take an
element in the closure of $A$, then there exists a sequence in $A$
converging to this point\MMP{Idea for the proof}. If the oscillation at that point is zero, the
images of the sequence elements must form a Cauchy sequence, so we
extend the map by forming the limit of this sequence. Now we have to
show that this map is well defined and continuous. 

\BeginProof
1.
We may and do assume that the complete metric $d$ for $Y$ is bounded
by $1$. Define $G$ as above,
then $A \subseteq G \subseteq \Closure{A}$, and $G$ can be written as
the intersection of a sequence of open sets. In fact,
represent $G$ as
\begin{equation*}
G = \bigcap_{n \in \Nat} \{x \in \Closure{A} \mid {\o}_f(x) < \frac{1}{n}\},
\end{equation*}
so we have to show that $\{x \in \Closure{A} \mid {\o}_f(x) < q\}$ is open in $\Closure{A}$ for any $q>0$. But we have 
\begin{equation*}
\{x \in \Closure{A} \mid {\o}_f(x) < q\}
=
\bigcup \{V \cap \Closure{A} \mid \mathsf{diam}(\Bild{f}{V \cap A}) < q\}.
\end{equation*}
This is the union of sets open in $\Closure{A}$, hence is an open set itself. 

2.
Now take an element $x \in G \subseteq \Closure{A}$. Then there exists a sequence
$\Folge{x}$ of elements $x_n \in A$ with $x_n \rightarrow x$. Given $\epsilon > 0$,
we find a neighborhood $V$ of $x$ with $\mathsf{diam}(\Bild{f}{A \cap V}) < \epsilon$, since the oscillation of $f$ at $x$ is $0$. Because
$x_n \rightarrow x$, we know that we can find an index $n_{\epsilon}\in\Nat$ such that $x_m \in V \cap A$ for all $m > n_\epsilon$. This implies that
the sequence $(f(x_n))_{n \in \Nat}$ is a Cauchy sequence in $Y$. It converges
because $Y$ is complete. Put
\begin{equation*}
f_*(x) := \lim_{n \rightarrow \infty} f(x_n).
\end{equation*}

3.
We have to show now that 
\begin{itemize}
\item $f_{*}$ is well-defined.
\item $f_{*}$ extends $f$.
\item $f_{*}$ is continuous.
\end{itemize}
Assume that we can find $x\in G$ such that $\Folge{x}$ and $\Folge{x'}$ are sequences in $A$ with $x_{n}\to x$ and $x'_{n}\to x$, but $\lim_{n\to \infty}f(x_{n}) \not= \lim_{n\to \infty}f(x'_{n})$. Thus we find some $\eta>0$ such that $d(f(x_{n}), f(x'_{n})) \geq \eta$ infinitely often. Then the oscillation of $f$ at $x$ is at least $\eta>0$, a contradiction. This implies that $f_{*}$ is well-defined, and it implies also that $f_{*}$ extends $f$. Now let $x\in G$. If $\epsilon>0$ is given, we find a neighborhood $V$ of $x$  with $\mathsf{diam}(\Bild{f}{A\cap V})<\epsilon$. Thus, if $x'\in G\cap V$, then $d(f_{*}(x), f_{*}(x')) < \epsilon$. Hence $f_{*}$ is continuous. 
\EndProof

A characterization of complete spaces in terms of sequences of
closed sets with decreasing diameters is given below.

\BeginProposition{diam-to-zero-compl}
These statements are equivalent
\begin{enumerate}
\item\label{diam-to-zero-compl-1} $X$ is complete.
\item\label{diam-to-zero-compl-2} For each decreasing sequence $\Folge{A}$ of
  non-empty closed sets the diameter of which tends to zero there
  exists $x\in X$ such that $\bigcap_{n\in\Nat}A_{n} = \Closure{\{x\}}$.
\end{enumerate}
In particular, if $X$ is a metric space, then $X$ is complete iff
each decreasing sequence of non-empty closed sets the diameter of
which tends to zero has exactly one point in common.  
\EndProposition

\BeginProof
The assertion for the metric case follows immediately from the general
case, because $\Closure{\{x\}} = \{x\}$, and because there can be not
more than one element in the intersection.

\labelImpl{diam-to-zero-compl-1}{diam-to-zero-compl-2}: Let
$\Folge{A}$ be a decreasing sequence of non-empty closed sets with
$\mathsf{diam}(A_{n})\to 0$, then we have to show that
$\bigcap_{n\in\Nat}A_{n} = \Closure{\{x\}}$ for some $x\in X$. Pick
from each $A_{n}$ an element $x_{n}$, then $\Folge{x}$ is a Cauchy
sequence which converges to some $x$, since $X$ is complete. Because
the intersection of closed sets is closed again, we conclude
$\bigcap_{n\in\Nat}A_{n} = \Closure{X}$.

\labelImpl{diam-to-zero-compl-2}{diam-to-zero-compl-1}: Take a Cauchy sequence
$\Folge{x}$, then $A_{n} :=\Closure{\{x_{m}\mid
  m\geq n\}}$ is a decreasing sequence of closed sets the diameter of
which tends to zero. In fact, given $\epsilon>0$ there exists $n\in
\Nat$ such that $d(x_{m}, x_{m'}) < \epsilon$ for all $m, m'\geq n$,
hence $\mathsf{diam}(A_{n}) < \epsilon$, and it follows that this
holds also for all $k\geq n$. Then it is
obvious that $x_{n}\to x$ whenever $x\in\bigcap_{n\in\Nat}A_{n}$. 
\EndProof

We mention all too briefly a property of complete spaces which renders them
most attractive, viz., Banach's Fixpoint Theorem. 

\BeginDefinition{def-contraction}
Call $f: X\to X$ a \emph{\index{contraction}contraction} iff there exists $\gamma$ with
$0<\gamma<1$ such that $d(f(x), f(y)) \leq \gamma\cdot
d(x, y)$ holds for all $x, y\in X$.
\EndDefinition

Then one shows 

\BeginTheorem{banach-fixed-point}
Let $f: X\to X$ be a contraction with $X$ complete. Then there exists
$x\in X$ with $f(x) = x$. If $f(y) = y$ holds as well, then $d(x, y) =
0$. In particular, if $X$ is a metric space, then there exists a
unique fixed point for $f$.\MMP{Banach's Fixpoint Theorem}
\EndTheorem

The idea is just to start with an arbitrary element of $X$, and to
iterate $f$ on it. This yields a sequence of elements of $X$. Because
the elements become closer and closer, completeness kicks in and makes
sure that there exists a limit. This limit is independent of the
starting point.

\BeginProof
Define the $n$-th iteration $f^{n}$ of $f$ through $f^{1} := f$ and
$f^{n+1} := f^{n} \circ f$. Now let $x_{0}$ be an arbitrary element of
$X$, and define $x_{n} := f^{n}(x_{0})$. Then $d(x_{n}, x_{n+m})\leq
\gamma^{n}\cdot d(x_{0}, x_{m})$, so that $\Folge{x}$ is a Cauchy
sequence which converges to some $x\in X$, and $f(x) = x$. If $f(y) =
y$, we have $d(x, y) = d(f(x), f(y))\leq \gamma\cdot d(x, y)$, thus
$d(x, y) = 0$. This implies uniqueness of the fixed point as well.  
\EndProof

The Banach Fixed Point Theorem has a wide range of applications, and
it used for  iteratively approximating the solution of equations,
e.g., for implicit functions. The following example permits a glance
at Google's\MMP{Google} page rank algorithm, it
follows~\cite{Rousseau} (the linear algebra behind it is explored in,
e.g.,~\cite{Langville+Meyer, Keener}).

\BeginExample{how-google-works}
Let 
$
S := \{\langle x_{1}, \dots, x_{n}\rangle\mid x_{i}\geq 0,
x_{1}+\dots+x_{n}=1\}
$
be the set of all discrete probability distributions over $n$ objects, and $P:
\Real^{n}\to \Real^{n}$ be a stochastic matrix; this means that $P$
has non-negative entries and the rows all add up to $1$.  The set $\{1, \dots, n\}$ is usually
interpreted as the state space for some random experiment,  entry $p_{i, j}$ is then
interpreted as the probability for the change of state $i$ to state
$j$. We have in particular $P: S\to S$, so a probability distribution
is transformed into another probability distribution. We assume that
$P$ has an eigenvector $v_{1}\in S$ for the eigenvalue $1$, and  that the
other eigenvalues are in absolute value not greater than 1 (this is
what the classic Perron-Frobenius Theorem says, see~\cite{Langville+Meyer, Keener}); moreover we assume
that we can find a base $\{v_{1}, \dots, v_{n}\}$ of eigenvectors, all
of which may be assumed to be in $S$; let $\lambda_{i}$ be the
eigenvector for $v_{i}$, then $\lambda_{1}=1$, and $|\lambda_{i}|\leq
1$ for $i\geq 2$. Such a matrix is called a \emph{regular transition matrix}; these matrices
are investigated in the context of stability of finite Markov transition
chains. 

Define for the distributions $p = \sum_{i=1}^{n}p_{i}\cdot v_{i}$ and
$q = \sum_{i=1}^{n}q_{i}\cdot v_{i}$ their distance  through
\begin{equation*}
  d(p, q) := \half\cdot \sum_{i=1}^{n}|p_{i}-q_{i}|.
\end{equation*}
Because $\{v_{1}, \dots, v_{n}\}$ are linearly independent, $d$ is a
metric. Because this set forms a basis, hence is given through a
bijective linear maps from the base given by the unit vectors, and
because the Euclidean metric is complete, $d$ is complete as well. 

Now define $f(x) := P\cdot x$, then this is a contraction $S\to S$:
\begin{equation*}
  d(P\cdot x, P\cdot y) = \half\cdot \sum_{i=1}^{n}|x_{i}\cdot
  P(v_{i})- y_{i}\cdot P(v_{i})| \leq
  \half\sum_{i=1}^{n}|\lambda_{i}\cdot (x_{i}-y_{i})| \leq \half\cdot
  d(x, y). 
\end{equation*}
Thus $f$ has a fixed point, which must be $v_{1}$ by uniqueness. 

Now assume that we have a (very litte) Web universe with only five
pages. The links are given as in the diagram.
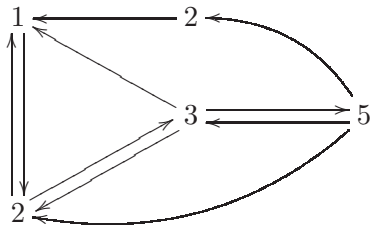
\begin{equation*}
\xymatrix{
1\ar@<.5ex>[dd] && 2\ar[ll]\\
&&3\ar[ull]\ar@<.5ex>[rr]\ar@<.5ex>[dll] && 5\ar@<.5ex>[ll]\ar@/_1pc/[ull]\ar@/^1.5pc/[dllll]\\
2\ar@<.5ex>[uu]\ar@<.5ex>[urr]
}
\end{equation*}
The transitions between pages are at random, the matrix below
describes such a random walk
\def\thrd{\frac{1}{3}}
\begin{equation*}
  P := \left(
    \begin{matrix}
      0 & 1 & 0 & 0& 0\\
\half & 0 & \half & 0 & 0\\
\thrd & \thrd & 0 & 0 & \thrd\\
1 & 0 & 0 & 0 & 0\\
0 & \thrd & \thrd & \thrd & 0
    \end{matrix}
\right)
\end{equation*}
It says that we make a transition from state $2$ to state $1$ with
$p_{2, 1} = \half$, also $p_{2, 3} = \half$, the transition from state
$2$ to state $3$. From state $1$ one  goes with probability one to state
$2$, because $p_{1, 2} = 1$. Iterating $P$ quite a few times will
yield a solution which does not change much after $32$ steps, one
obtains
\begin{equation*}
  P^{32} = \left(
\begin{matrix}
.293 & .390 & .220 & .024 & .073\\
.293 & .390 & .220 & .024 & .073\\
.293 & .390 & .220 & .024 & .073\\
.293 & .390 & .220 & .024 & .073\\
.293 & .390 & .220 & .024 & .073\\
\end{matrix}
\right)
\end{equation*}
The eigenvector $p$ for the eigenvalue $1$ looks like this: $p =
\langle.293, .390, .220, .024, .073\rangle$, so this yields a stationary
distribution. 

In terms of web searches\MMP{Web search}, the importance of the pages is
ordered according this stationary distribution as $2, 1, 3, 5, 4$, so
this is the ranking one would associate with these pages. 

This is the basic idea behind \index{Google}Google's page ranking
algorithm. Of course, there are many practical considerations which
have been eliminated from this toy example. It may be that the matrix
does not follow the assumptions above, so that it has to me modified
accordingly in a preprocessing step; size is a problem, of course, since handling the
extremely large matrices occurring in web searches.  
\EndExample
 
Compact pseudometric spaces are complete. This will be a byproduct
of a more general characterization of compact spaces.  We show
first that compactness and sequential compactness are the same
for these spaces. This is sometimes helpful in those situations in
which a sequence is easier to handle than an open cover, or an
ultrafilter. 

Before discussing this, we introduce\MMP{$\epsilon$-net}
\emph{\index{$\epsilon$-net}$\epsilon$-nets} as a cover of $X$ through
a \emph{finite} family $\{B(x, \epsilon)\mid x\in A\}$ of open balls
of radius $\epsilon$. $X$ may or may not have an $\epsilon$-net for
any given $\epsilon>0$. For example, $\Real$ does not have an
$\epsilon$-net for any $\epsilon > 0$, in contrast to $[0, 1]$ or
$]0, 1]$.

\BeginDefinition{totally-bounded}
The pseudometric space $X$ is \emph{\index{totally bounded}totally bounded} iff there exists for
each $\epsilon>0$ an $\epsilon$-net for $X$. A subset of a
pseudometric space is totally bounded iff it is a totally bounded
subspace. 
\EndDefinition

Thus $A\subseteq X$ is totally bounded iff $\Closure{A}\subseteq X$ is
totally bounded.

We see immediately

\BeginLemma{compact-tot-bounded}
A compact pseudometric space is totally bounded.
\QED
\EndLemma

Now we are in a position to establish this equivalence, which will
help characterize compact pseudometric spaces.

\BeginProposition{seq-comp-equiv-comp}
The following properties are equivalent for the pseudometric space
$X$:
\begin{enumerate}
\item\label{seq-comp-equiv-comp-1} $X$ is compact.
\item\label{seq-comp-equiv-comp-2} $X$ is sequentially compact. 
\end{enumerate}
\EndProposition

\BeginProof
\labelImpl{seq-comp-equiv-comp-1}{seq-comp-equiv-comp-2}: Assume that
the sequence $\Folge{x}$ does not have a convergent subsequence, and
consider the set $F := \{x_{n}\mid n\in\Nat\}$. This set is closed,
since, if $y_{n}\to y$ and $y_{n}\in F$ for all $n\in\Nat$, then $y\in
F$, since the sequence $\Folge{y}$ is eventually constant. $F$ is also
discrete, since, if we could find for some $z\in F$ for each
$n\in\Nat$ an element in $F\cap B(z, 1/n)$ different from $z$, we
would have a convergent subsequence. Hence $F$ is a closed discrete
subspace of $X$ which contains infinitely many elements, which is
impossible. This contradiction shows that each sequence has a
convergent subsequence.

\labelImpl{seq-comp-equiv-comp-2}{seq-comp-equiv-comp-1}\MMP{Plan of attack}: 
Before we enter into the second and harder part of the proof, we
have a look at its plan. Given an open cover for the sequential
compact space  $X$, we have to construct a finite cover from it. If we succeed in constructing
for each $\epsilon>0$ a finite net so that we can fit each ball into
some element of the cover, we are done, because in this case we may
take just these elements of the cover, obtaining a finite cover. That
this fitting in is possible is shown in the first part of the
proof. We construct under the assumption that it is not possible a
sequence, which has a converging subsequence, and the limit of this
subsequence will be used as kind of a flyswatter\index{flyswatter}. 

The second part of the proof is then just a simple application of the
net so constructed.

Let $(U_{i})_{i\in I}$ be a finite cover of $X$. We claim that we can find
for this cover some $\epsilon>0$ such that, whenever $\mathsf{diam}(A)<
\epsilon$, there exists $i\in I$ with $A\subseteq U_{i}$. Assume that
this is wrong, then we find for each $n\in\Nat$ some $A_{n}\subseteq
X$ which is not contained in one single $U_{i}$. Pick from each
$A_{n}$ an element $x_{n}$, then $\Folge{x}$ has a convergent
subsequence, say $\Folge{y}$, with $y_{n}\to y$. There exists a member
$U$  of the cover with $y\in U$, and there exists $r>0$ with $B(y,
r)\subseteq U$. Now we catch the fly. Choose $\ell\in\Nat$ with $1/\ell <
r/2$, then $y_{m}\in B(y, r/2)$ for $m\geq n_{0}$ for some suitable
chosen $n_{0}\in\Nat$, hence, because $\Folge{y}$ is a subsequence of
$\Folge{x}$, there are infinitely many $x_{k}$ contained in $B(y,
r/2)$. But since $\mathsf{diam}(A_{\ell})<1/\ell$, this implies
$A_{\ell}\subseteq B(y, r)\subseteq U$, which is a contradiction. 

Now select for the cover $\epsilon>0$ as above, and let the finite set
$A$ be the set of centers for an $\epsilon/2$-net, say, $A = \{a_{1},
\dots, a_{k}\}$. Then we can find for each $a_{j}\in A$ some member
$U_{i_{j}}$ of this cover with $B(a_{j}, \epsilon/2)\subseteq
U_{i_{j}}$ (note that $\mathsf{diam}(B(x, r) < 2\cdot r$). This
yields a finite cover $\{U_{i_{j}}\mid 1 \leq j \leq k\}$ of $X$. 
\EndProof

This proof was conceptually a little complicated, since we had to make the step
from a sequence (with a converging subsequence) to a cover (with the
goal of finding a finite cover). Both are not immediately related.
The missing link turned out to be measuring the size of a set through
its diameter, and capturing limits through suitable sets. 

Using the last equivalence, we are in a position to characterize
compact pseudometric spaces.

\BeginTheorem{compact-is-tot-bound-compl}
A pseudometric space is compact iff it is totally bounded and complete.
\EndTheorem

\BeginProof
1.
Let $X$ be compact. We know already from Lemma~\ref{compact-tot-bounded} that a compact
pseudometric space is totally bounded. Let $\Folge{x}$ be a Cauchy
sequence, then we know that it has a converging subsequence, which,
being a Cauchy sequence, implies that it converges itself. 

2.
Assume that $X$ is totally bounded and complete. In view of
Proposition~\ref{seq-comp-equiv-comp} it is enough to show that $X$ is
sequentially compact. Let $\Folge{x}$ be a sequence in $X$. Since $X$
is totally bounded, we find a subsequence $(x_{n_{1}})$ which is
entirely contained in an open ball of radius less that $1$. Then we may
extract from this sequence a subsequence $(x_{n_{2}})$ which is
contained in an open ball of radius less than $1/2$. 
Continuing inductively we find a subsequence $(x_{n_{k+1}})$ of $(x_{n_{k}})$ the
members of which are completely contained in an open ball of radius less
than $2^{-(k+1)}$. Now define $y_{n} := x_{n_{n}}$, hence $\Folge{y}$
is the diagonal sequence in this scheme. 

We claim that $\Folge{y}$ is a Cauchy sequence. In fact, let
$\epsilon>0$ be given, then there exists $n\in\Nat$ such that
$\sum_{\ell>n}2^{-\ell}<\epsilon/2$. Then we have for $m>n$
\begin{equation*}
  d(y_{n}, y_{m}) \leq 2\cdot \sum_{\ell=n}^{m}2^{-\ell} < \epsilon.
\end{equation*}
By completeness, $y_{n}\to y$ for some $y\in X$. Hence we have found a
converging subsequence of the given sequence $\Folge{x}$, so that $X$ is sequentially
compact. 
\EndProof

It\MMP[t]{Shift of emphasis} might be noteworthy to observe the shift of emphasis between
finding a finite cover for a given cover, and admitting an
$\epsilon$-net for each $\epsilon>0$. While we have to select a finite
cover from an arbitrarily given cover beyond our control, 
in the case of a totally bounded space we can construct for each $\epsilon>0$
a cover of a certain size, hence we may be in a position to influence the shape of
this special cover. Consequently, the characterization of compact spaces in
Theorem~\ref{compact-is-tot-bound-compl} is very helpful and handy,
but, alas, it works only in the restricted calls of pseudometric
spaces. 

We apply this characterization to $(\allCompact(X), \delta_{H})$, the
space of all non-empty compact subsets of $(X, d)$ with the Hausdorff
metric $\delta_{H}$, see Example~\ref{vietoris-pseudometric}.

\BeginProposition{vietoris-is-complete}
$(\allCompact(X), \delta_{H})$ is complete, if $X$ is a complete
pseudometric space.
\EndProposition

\BeginProof
We fix for the proof a Cauchy sequence $\Folge{C}$ of elements of
$\allCompact(X)$. 

0.
Let us pause a moment and discuss the approach to the proof\MMP{Plan} first. We show in a first step that
$\Closure{(\bigcup_{n\in\Nat}C_{n})}$ is compact by showing that it is
totally bounded and complete. Completeness is trivial, since the space
is complete, and we are dealing with a closed subset, so we focus on
showing that the set is totally bounded. Actually, it is sufficient to show that
$\bigcup_{n\in\Nat}C_{n}$ is totally bounded, because a set is totally
bounded iff its closure is. 

Then compactness of $\Closure{(\bigcup_{n\in\Nat}C_{n})}$ implies that
$C := \bigcap_{n\in\Nat}\Closure{(\bigcup_{k\geq n}C_{k})}$ is compact as
well, moreover, we will argue that $C$ must be non-empty. Then it is shown that $C_{n}\to
C$ in the Hausdorff metric. 

1.
Let $D := \bigcup_{n\in\Nat}C_{n}$, and let $\epsilon>0$ be given. We
will construct an $\epsilon$-net for $D$. Because $\Folge{C}$ is
Cauchy, we find for $\epsilon$ an index $\ell$ so that
$\delta_{H}(C_{n}, C_{m}) < \epsilon/2$ for $n, m\geq \ell$. When
$n\geq \ell$ is fixed, this
means in particular that $C_{m}\subseteq C_{n}^{\epsilon/2}$ for all
$m\geq \ell$, thus $d(x, C_{n}) < \epsilon/2$ for all $x\in C_{m}$ and
all $m\geq \ell$. We will use this observation in a moment.

Let $\{x_{1}, \dots, x_{t}\}$ be an $\epsilon/2$-net for
$\bigcup_{j=1}^{n}C_{j}$, we claim that this is an $\epsilon$-net for
$D$. In fact, let $x\in D$. If $x\in\bigcup_{j=1}^{\ell}C_{j}$, then
there exists some $k$ with $d(x, x_{k})<\epsilon/2$. If $x\in C_{m}$
for some $m>n$, $x\in C_{n}^{\epsilon/2}$, so that we find
$x'\in C_{n}$ with $d(x, x')<\epsilon/2$, and for $x'$ we find $k$
such that $d(x_{k}, x')<\epsilon/2$. Hence we have
$d(x, x'_{k}) < \epsilon$, so that we have shown that
$\{x_{1}, \dots, x_{t}\}$ is an $\epsilon$-net for $D$.  Thus
$\Closure{D}$ is totally bounded, hence compact.

2.
From the first part it follows that $\Closure{(\bigcup_{k\geq
  n}C_{k})}$ is compact for each $n\in\Nat$. Since these sets form
a decreasing sequence of non-empty closed subsets to the compact set
given by $n=1$, their intersection cannot be empty, hence 
$
C := \bigcap_{n\in\Nat}\Closure{(\bigcup_{k\geq n}C_{k})}
$
is compact and non-empty, hence a member of $\allCompact(X)$. 

We claim
that $\delta_{H}(C_{n}, C)\to 0$, as $n\to \infty$. Let $\epsilon>0$
be given, then we find $\ell\in\Nat$ such that $\delta_{H}(C_{m},
C_{n})\leq/2$, whenever $n, m\geq \ell$. We show that
$\delta_{H}(C_{n}, C)< \epsilon$ for all $n\geq \ell$. Let $n\geq \ell$. The proof is
subdivided into showing that $C\subseteq C_{n}^{\epsilon}$ and
$C_{n}\subseteq C^{\epsilon}$. 

Let us work on the first inclusion.
Because $D := \Closure{(\bigcup_{i\geq n}C_{i})}$ is totally
bounded, there exists a $\epsilon/2$-net, say, $\{x_{1}, \dots,
x_{t}\}$, for $D$. If $x\in C\subseteq D$, then there exists $j$ such
that $d(x, x_{j})< \epsilon/2$, so that we can find $y\in C_{n}$
with $d(y, x_{j})<\epsilon/2$. Consequently, we find for $x\in C$ some
$y\in C_{n}$ with $d(x, y)<\epsilon$. Hence $C\subseteq
C_{n}^{\epsilon}$. 

Now for the second inclusion. Take $x\in C_{n}$. Since
$\delta_{H}(C_{m}, C_{n})<\epsilon/2$ for $m\geq \ell$, we have
$C_{n}\subseteq C_{m}^{\epsilon/2}$, hence find
$x_{m}\in C_{m}$ with $d(x, x_{m})<\epsilon/2$. The sequence
$(x_{k})_{k\geq m}$ consists of members of the compact set $D$, so it
has converging subsequence which converges to some $y\in D$. But it
actually follows from the construction that $y\in C$, and $d(x, y)
\leq d(x, x_{m}) + d(x_{m}, y) < \epsilon$ for $m$ taken sufficiently
large from the subsequence. This yields $x\in C^{\epsilon}$. 

Taking these inclusions together, they imply $\delta_{H}(C_{n}, C) <
\epsilon$ for $n>\ell$. This shows that $(\allCompact(X), \delta_{H})$
is a complete pseudometric space, if $(X, d)$ is one.
\EndProof

The topology induced by the Hausdorff metric can be defined in a way
which permits a generalization to arbitrary topological spaces, where
it is called the \emph{\index{topology!Vietoris}Vietoris topology}. It has been studied with
respect to finding continuous selections, e.g., by
Michael~\cite{Michael}, see also~\cite{Jayne+Rogers, Castaing-Valadier}. The reader is also referred
to~\cite[§33]{Kuratowski}, and to~\cite[p. 120]{Engelking} for a study of topologies on subsets. 

We will introduce uniform continuity now and discuss this concept briefly
here. Uniform spaces will turn out to be the proper scenario for the
more extended discussion in Section~\ref{sec:uniform-spaces}. As a motivating example, assume that the
pseudometric $d$ on $X$ is bounded, take a
subset $A\subseteq X$ and look at the function $x\mapsto d(x,
A)$. Since 
\begin{equation*}
  |d(x, A) - d(y, A) | \leq d(x, y),
\end{equation*}
we know that this map is continuous. This means that, given $x\in X$,
there exists $\delta>0$ such that $d(x, x')<\delta$ implies  $|d(x, A)
- d(x', A) |< \epsilon$. We see from the inequality above that the
choice of $\delta$ does only depend on $\epsilon$, but not on
$x$. Compare this with the function $x\mapsto 1/x$ on $]0, 1]$. This
function is continuous as well, but the choice of $\delta$ depends on
the point $x$ you are considering: whenever $0<\delta<\epsilon\cdot
x^{2}/(1+\epsilon\cdot x)$, we may conclude that $|x'-x|\leq\delta$
implies $|1/x'-1/x|\leq\epsilon$. In fact, we may easily infer from
the graph of the function that a uniform choice of $\delta$ for a
given $\epsilon$ is not possible. 

This leads to the definition of uniform continuity in a pseudometric
space: the choice of $\delta$ for a given $\epsilon$ does not depend
on a particular point, but is rather, well, uniform.

\BeginDefinition{unif-continuity}
The map $f: X\to Y$ into the pseudometric space $(Y, d')$ is called
\emph{uniformly \index{continuous!uniformly}continuous} iff given
$\epsilon>0$ there exists $\delta>0$ such that $d'(f(x),
f(x'))<\epsilon$ whenever $d(x, x')<\delta$. 
\EndDefinition

Doing\MMP[t]{Continuity vs. uniform continuity} a game of quantifiers, let us just point out the difference between
uniform continuity and continuity. 
\begin{enumerate}
\item Continuity says
  \begin{equation*}
    \forall \epsilon>0\underline{\forall x\in X\exists\delta>0}\forall
    x'\in X: d(x, x')< \delta \Rightarrow d'(f(x), f(x'))<\epsilon.
  \end{equation*}
\item Uniform continuity says 
  \begin{equation*}
    \forall \epsilon>0\underline{\exists\delta>0\forall x\in X}\forall
    x'\in X: d(x, x')< \delta \Rightarrow d'(f(x), f(x'))<\epsilon.
  \end{equation*}
\end{enumerate}

The formulation suggests that uniform continuity depends on the chosen
metric. In contrast to continuity, which is a property depending on
the topology of the underlying spaces, uniform continuity is a
property of the underlying uniform space, which will be discussed
below. We note that the composition of uniformly continuous maps is
uniformly continuous again. 

A uniformly continuous map is continuous. The converse is not true,
however. 

\BeginExample{cont-not-unif-cont}
Consider the map $f: x\mapsto
x^{2}$, which is certainly continuous on $\Real$. Assume that $f$ is
uniformly continuous, and fix $\epsilon>0$, then there exists
$\delta>0$ such that $|x-y|<\delta$ always implies
$|x^{2}-y^{2}|<\epsilon$. Thus we have for all $x$, and for all $r$
with $0<r\leq\delta$ that $|x^{2}-(x+r)^{2}| = |2\cdot x\cdot
r+r^{2}|<\epsilon$ after Binomi's celebrated theorem. But this would
mean $|2\cdot x+r| < \epsilon/r$ for all $x$, which is not
possible. In general, a very similar argument shows that polynomials
$\sum_{i=1}^{n}a_{i}\cdot x^{i}$ with $n>1$ and $a_{n}\not=0$ are not
uniformly continuous. 
\EndExample

A continuous function on a compact pseudometric space, however, is uniformly
continuous. This is established through an argument constructing a
cover of the space, compactness will then permit us to extract a
finite cover, from which we will infer uniform continuity.

\BeginProposition{compact-cont-unif-cont}
Let $f: X\to Y$ be a continuous map from the compact pseudometric space
$X$ to the pseudometric space $(Y, d')$. Then $f$ is uniformly continuous.
\EndProposition

\BeginProof
Given $\epsilon>0$, there exists for each $x\in X$ a positive
$\delta_{x}$ such that $\Bild{f}{B(x, \delta_{x})} \subseteq
B_{d'}(f(x), \epsilon/3)$. Since $\{B(x, \delta_{x}/3)\mid x\in X\}$ is an
open cover of $X$, and since $X$ is compact, we find $x_{1}, \dots,
x_{n}\in X$ such that $B(x_{1}, \delta_{x_{1}}/3), \dots, B(x_{n},
\delta_{x_{n}}/3)$ cover $X$. Let $\delta$ be the smallest among
$\delta_{x_{1}}, \dots, \delta_{x_{n}}$. If $d(x, x')< \delta/3$, then
there exist $x_{i}, x_{j}$ with $d(x, x_{i})< \delta/3$ and $d(x',
x_{j})< \delta/3$, so that $d(x_{i}, x_{j})\leq d(x_{i}, x) + d(x, x')
+ d(x', x_{j}) < \delta$, hence $d'(f(x_{i}), f(x_{j})) < \epsilon/3$,
thus $d'(f(x), f(x')) \leq d'(f(x), f(x_{i})) + d'(f(x_{i}), f(x_{j}))
+ d'(f(x_{j}), f(x')) < 3\cdot \epsilon/3 = \epsilon.$
\EndProof

One of the most attractive features of uniform continuity is that it permits
extending a function --- given a uniform continuous map $f: D\to Y$ with
$D\subseteq X$ dense and $Y$ complete metric, we can extend $f$ to a
uniformly continuous map $F$ on the whole space\MMP{Idea for a proof}. This extension is
necessarily unique (see Lemma~\ref{equal-on-dense}). The basic idea is
to define $F(x) := \lim_{n\to \infty}f(x_{n})$, whenever $x_{n}\to x$
is a sequence in $D$ which converges to $x$. This requires that the
limit exists, and that it is in this case unique, hence it demands the
range to be a metric space which is complete.

\BeginProposition{extend-unif-cont-maps}
Let $D\subseteq X$ be a dense subset, and assume that $f: D\to Y$ is
uniformly continuous, where $(Y, d')$ is a complete metric space. Then
there exists a unique uniformly continuous map $F: X\to Y$ which
extends $f$. 
\EndProposition

\BeginProof
0.
We have already argued that an extension must be unique, if it
exists. So we have to construct it, and to show that it is uniformly
continuous. We will generalize the argument from above referring
to a limit by considering the oscillation at each point\MMP{Outline~---~use the oscillation}.  A glimpse at
the proof of Lemma~\ref{Kuratowski} shows indeed that we argue with a
limit here, but are able to look at the whole set of points which
makes this possible. 

1.
Let us have a look at the oscillation ${\o}_f(x)$ of $f$ at a point
$x\in X$ (see page~\pageref{ref:oscillation}), and we may assume that
$x\not\in D$. We claim that
${\o}_f(x) = 0$. In fact, given $\epsilon>0$, there exists $\delta>0$
such that $d(x', x'') < \delta$ implies $d'(f(x'), f(x'')) <
\epsilon/3$, whenever $x', x''\in D$. Thus, if $y', y''\in
\Bild{f}{D\cap B(x, \delta/2)}$, we find $x', x''\in D$ with $f(x') =
y', f(x'') = y''$ and $d(x',
x'') \leq d(x, x') + d(x'', x) < \delta$, hence $d'(y', y'') =
d'(f(x'), f(y'')) < \epsilon$. This means that $\mathsf{diam}(\Bild{f}{D\cap
  B(x, \delta/2)}) < \epsilon$. 

2.
Lemma~\ref{Kuratowski} tells us that there exists a continuous
extension $F$ of $f$ to the set $\{x\in X \mid {\o}_f(x) = 0\} =
X$. Hence it remains to show that $F$ is \emph{uniformly}
continuous. Given $\epsilon>0$, we choose the same $\delta$ as above,
which did not depend on the choice of the points we were considering
above. Let $x_{1}, x_{2}\in X$ with $d(x_{1}, x_{2})< \delta/2$, then
there exists $v_{1}, v_{2}\in D$ such that $d(x_{1}, v_{1}) <\delta/4$
with $d'(F(x_{1}), f(v_{1}))\leq \epsilon/3$
and $d(x_{2}, v_{2}) < \delta/4$ with $d'(F(x_{2}), f(v_{2}))\leq
\epsilon/3$. We see as above that $d(v_{1}, v_{2}) < \delta$, thus
$d'(f(v_{1}), f(v_{2})) < \epsilon/3$, consequently, 
\begin{equation*}
  d'(F(x_{1}, x_{2})) \leq d'(F(x_{1}), f(v_{1})) + d'(f(v_{1}),
  f(v_{2})) + d'(f(v_{2}), F(x_{2})) < 3\cdot \epsilon/3 = \epsilon.
\end{equation*}
But this means that $F$ is uniformly continuous. 
\EndProof

Looking at $x\mapsto 1/x$ on $]0, 1]$ shows that uniform continuity is
indeed necessary to obtain a continuous extension. 

\Subsubsection{Baire's Theorem and a Game}
\label{sec:baire+game}

The technique of constructing a shrinking sequence of closed sets with a
diameter tending to zero used for establishing Proposition~\ref{diam-to-zero-compl} is helpful in establishing Baire's Theorem~\ref{baire-locally-compact}
also for complete pseudometric spaces;\index{theorem!Baire!complete pseudometric} completeness then makes sure
that the intersection is not empty. The proof is essentially a blend
of this idea with the proof given above
(page~\pageref{baire-locally-compact}). We will then give an
interpretation of Baire's Theorem in terms of the game \emph{Angel vs. Demon} introduced
in~\SetCite{Section 1.7}. We show that Demon has a winning strategy
iff the space is the
countable union of nowhere dense sets (the space is then called to be
of the \emph{first category}). This is done for a subset of
the real line, but can easily generalized. 
 
This is the version of Baire's Theorem in a complete pseudometric
space.

\BeginTheorem{baire-complete}
Let $X$ be a complete pseudometric space, then the intersection of a
sequence of dense open sets is dense again\MMP{Baire's Theorem}. 
\EndTheorem

\BeginProof
Let $\Folge{D}$ be a sequence of dense open sets. Fix a non-empty open
set $G$, then we have to show that
$G\cap\bigcap_{n\in\Nat}D_{n}\not=\emptyset$. Now $D_{1}$ is dense and
open, hence we find an open set $V_{1}$ and $r > 0$ such that
$\mathsf{diam}(\Closure{V}_{1}) \leq r$ and $\Closure{V}_{1}\subseteq D_{1}\cap G$. We
select inductively in this way a sequence of open sets $\Folge{V}$
with $\mathsf{diam}(\Closure{V}_{n}) < r/n$ such that
$\Closure{V}_{n+1}\subseteq D_{n}\cap V_{n}$. This is possible since
$D_{n}$ is open and dense for each $n\in\Nat$.

Hence we have a decreasing sequence
$\Closure{V}_{1} \supseteq \dots \Closure{V}_{n} \supseteq \dots$ of
closed sets with diameters tending to $0$ in the complete space $X$. Thus
$\bigcap_{n\in\Nat}\Closure{V}_{n} = \bigcap_{n\in\Nat}V_{n}$ is not
empty by Proposition~\ref{diam-to-zero-compl}, which entails
$G\cap\bigcap_{n\in\Nat}D_{n}$ not being empty.
\EndProof

Kelley~\cite[p. 201]{Kelley} remarks that there is a slight
incongruence with this theorem, since the assumption of completeness
is non-topological in nature (hence a property which may get lost when
switching to another  pseudometric, see
Example~\ref{loose-completeness}), but we draw a
topological conclusion. He suggests that the assumption on space $X$
should be reworded to $X$ being a topological space for which there
exists a complete pseudometric. But, alas, the formulation above is
the usual one, because it is pretty suggestive after all.

\BeginDefinition{def-nowhere-dense}
Call a set $A\subseteq X$ \emph{\index{nowhere dense}nowhere dense}
iff $\Closure{\Interior{A}} = \emptyset$, i.e., the closure of the
interior is empty, equivalently, iff the open set
$X\setminus\Closure{A}$ is dense. The space $X$ is said to be of the
\emph{\index{space!first category}first category} iff it can be
written as the countable union of nowhere dense sets.
\EndDefinition

Then Baire's Theorem can be reworded
that the countable union of nowhere dense sets in a complete pseudometric space
is nowhere dense. This is an important example for a nowhere dense set:

\BeginExample{cantor-ternary-nowhere-dense}
Cantor's ternary set $C$ (see~\SetCite{Example 1.104}) can be written as\MMP[t]{Cantor's ternary set}
\begin{equation*}
  C = \bigl\{\sum_{i=1}^{\infty}a_{i}3^{-i}\mid a_{i}\in\{0, 2\}\text{ for
    all  }i\in \Nat\bigr\}.
\end{equation*}
This is seen as follows: Define 
$
[a, b]' := [a + (b-a)/3] \cup [a+2\cdot (b-a)/3]
$
for an interval $[a, b]$, and $(A_{1}\cup\dots\cup A_{\ell})' :=
A_{1}'\cup\dots\cup A_{\ell}'$, then 
$
  C = \bigcap_{n\in\Nat}C_{n}
$
with the inductive definition
$
  C_{1}:= [0, 1]'
$
and
$
C_{n+1} := C_{n}'.
$
It is shown easily by induction that 
\begin{equation*}
  C_{n} = \{\sum_{i=1}^{\infty}a_{i}\cdot 3^{-i}\mid a_{i}\in\{0,
  2\}\text{ for $i \leq n$ and }a_{i}\in\{0, 1, 2\}\text{ for }i>n\}.
\end{equation*}
The representation above implies that the interior of $C$ is empty, so
that $C$ is in fact nowhere dense in the unit interval.
\EndExample

Cantor's ternary set is a helpful device in investigating the
structure of complete metric spaces which have a countable dense
subset, i.e., in Polish spaces.

We will give now a game theoretic interpretation of spaces of the
first category through a game which is attributed to Banach and
Mazur\index{Banach-Mazur game}, tieing the existence of a winning
strategy for Demon to spaces of the first category. For
simplicity, we discuss it for a closed interval of the real line. We
do not assume that the game is determined; determinacy is not
necessary here (and its assumption would bring us into serious
difficulties with the assumption of the validity of the Axiom of
Choice, see~\SetCite{Prop. 1.7.6}).

Let a subset $S$ of a closed interval $L_{0}\subseteq \Real$ be
given; this set is assigned to Angel, its adversary Demon is assigned its complement
$T := L_{0}\setminus S$. The game is played in this way\MMP[t]{Rules of
  the Banach-Mazur game}: 
\begin{itemize}
\item Angel chooses a closed interval $L_{1}\subseteq L_{0}$,
\item Demon reacts with choosing a closed interval $L_{2}\subseteq
  L_{1}$,
\item Angel chooses then ---~knowing the moves $L_{0}$ and $L_{1}$~--- a
  closed interval $L_{2}\subseteq L_{1}$,
\item and so on: Demon chooses the intervals with even numbers, Angel
  selects the intervals with the odd numbers, each interval is closed
  and contained in the previous one, both Angel and Demon have
  complete information about the game's history, when making a move. 
\end{itemize}
Angel wins iff $\bigcap_{n\in\Nat}L_{n}\cap S \not=\emptyset$,
otherwise Demon wins. 

We focus on Demon's behavior. Its strategy for the n-th move is
modelled as a map $f_{n}$ which is defined on $2\cdot n$-tuples $\langle
L_{0}, \dots, L_{2\cdot n-1}\rangle$ of closed intervals with
$L_{0}\supseteq L_{1}\supseteq\dots\supseteq L_{2\cdot n-1}$, taking a
closed interval $L_{2\cdot n}$ as a value with 
\begin{equation*}
  L_{2\cdot n} = f_{n}(L_{0}, \dots, L_{2\cdot n-1})\subseteq L_{2\cdot n-1}.
\end{equation*}
The sequence $\Folge{f}$ will be a \emph{winning strategy} for Demon iff
$\bigcap_{n\in\Nat}L_{n}\subseteq T$, when $\Folge{L}$ is chosen
according to these rules.

The following theorem relates the existence of a winning strategy for
Demon with $S$ being of first category.

\BeginTheorem{baire-game}
There exists a strategy for Demon to win iff $S$ is of the first category.
\EndTheorem

We divide the proof into two parts~---~we show first that we can find
a strategy for Demon, if $S$ is of the first category. The converse is
technically somewhat more complicated, so we delay it and do the
necessary constructions first.

\BeginProof (First part)
Assume that $S$ is of the first category, so that we can write $S =
\bigcup_{n\in\Nat} S_{n}$ with $S_{n}$ nowhere dense for each
$n\in\Nat$. Angel starts with a closed interval $L_{1}$, then demon
has to choose a closed interval $L_{2}$; the choice will be so that
$L_{2}\subseteq L_{1}\setminus S_{1}$. We  have to be sure that such a
choice is possible; our assumption implies that $L_{1}\cap
\Closure{S}_{1}$ is open and dense in $L_{1}$, thus contains an open
interval. In the inductive step, assume that Angel has chosen the
closed interval $L_{2\cdot n-1}$ such that $L_{2\cdot
  n-1}\subseteq\dots\subseteq L_{2}\subseteq L_{1}\subseteq
L_{0}$. Then Demon will select an interval $L_{2\cdot n}\subseteq
L_{2\cdot n-1}\setminus(S_{1}\cup \dots \cup S_{n})$. For the same
reason as above, the latter set contains an open interval. This
constitutes Demon's strategy, and evidently
$\bigcap_{n\in\Nat}L_{n}\cap S = \emptyset$, so Demon wins.
\EndProof

The proof for the second part requires some technical
constructions. We assume that $f_{n}$ assigns to each $2\cdot n$-tuple
of closed intervals
$I_{1}\supseteq I_{2}\supseteq\dots\supseteq I_{2\cdot n}$ a closed
interval $f_{n}(I_{1}, \dots, I_{2\cdot n})\subseteq I_{2\cdot n}$,
but do not make any further assumptions, for the time being, that
is. We are given a closed interval $L_{0}$ and a subset
$S\subseteq L_{0}$.

In a first step we define a sequence $\Folge{J}$ of closed intervals
with these properties:
\begin{itemize}
\item $J_{n}\subseteq L_{0}$ for all $n\in\Nat$,
\item $K_{n} := f_{1}(L_{0}, J_{n})$ defines a sequence $\Folge{K}$ of
  mutually disjoint closed intervals,
\item $\bigcup_{n\in\Nat}\Interior{K}_{n}$ is dense in $L_{0}$. 
\end{itemize}

Let's see how to do this. Define ${\cal F}$ as the sequence of all
closed intervals with rational endpoints that are contained in
$\Interior{L}_{0}$. Take $J_{1}$ as the first element of ${\cal
  F}$. Put $K_{1} := f_{1}(L_{0}, J_{1})$, then $K_{1}$ is a closed
interval with $K_{1}\subseteq J_{1}$ by assumption on $f_{1}$. Let
$J_{2}$ be the first element in ${\cal F}$ which is contained in
$L_{0}\setminus K_{1}$, put $K_{2} := f_{1}(L_{0},
J_{2})$. Inductively, select $J_{i+1}$ as the first element of ${\cal
  F}$ which is contained in $L_{0}\setminus\bigcup_{t=1}^{i}K_{t}$,
and set $K_{i+1} := f_{1}(L_{0}, J_{i+1})$. It is clear from the
construction that $\Folge{K}$ forms a sequence of mutually disjoint
closed intervals with $K_{n}\subseteq J_{n}\subseteq L_{0}$ for each
$n\in\Nat$. Assume that $\bigcup_{n\in\Nat}\Interior{K}_{n}$ is not
dense in $L_{0}$, then we find $x\in L_{0}$ which is not contained in
this union, hence we find an interval $T$ with rational endpoints
which contains $x$ but
$T\cap\bigcup_{n\in\Nat}\Interior{K}_{n}=\emptyset$. So $T$ occurs
somewhere in ${\cal F}$, but it is never the first interval to be
considered in the selection process. Since this is impossible, we
arrive at a contradiction. 

We repeat this process for $\Interior{K}_{i}$ rather than $L_{0}$ for
some $i$, hence we will define a sequence $(J_{i, n})_{n\in\Nat}$ of
closed intervals $J_{i, n}$ with these properties:
\begin{itemize}
\item $J_{i, n}\subseteq \Interior{K}_{i}$ for all $n\in\Nat$,
\item $K_{i, n} := f_{2}(L_{0}, J_{i}, K_{i}, J_{i, n})$ defines a
  sequence $(K_{i, n})_{n\in\Nat}$ of mutually disjoint closed
  intervals,
\item $\bigcup_{n\in\Nat}\Interior{K}_{i, n}$ is dense in $K_{i}$.
\end{itemize}
It is immediate that $\bigcup_{i, j}\Interior{K}_{i, j}$ is dense in
$L_{0}$.

Continuing inductively, we find for each $\ell\in \Nat$ two families
$J_{i_{1}, \dots, i_{\ell}}$ and $K_{i_{1}, \dots, i_{\ell}}$ of
closed intervals with these properties
\begin{itemize}
\item $K_{i_{1}, \dots, i_{\ell}} = f_{\ell}(L_{0}, J_{i_{1}},
  K_{i_{1}}, J_{i_{1}, i_{2}}, K_{i_{1}, i_{2}}, \dots, J_{i_{1},
    \dots, i_{\ell}})$,
\item $J_{i_{1}, \dots, i_{\ell+1}}\subseteq \Interior{K}_{i_{1},
    \dots, i_{\ell}}$, 
\item the intervals
  $(K_{i_{1}, \dots, i_{\ell-1}, i_{\ell}})_{i_{\ell}\in \Nat}$ are
  mutually disjoint for each $i_{1}, \dots, i_{\ell-1}$,
\item $\bigcup\{\Interior{K}_{i_{1}, \dots, i_{\ell-1}, i_{\ell}} \mid \langle
  i_{1}, \dots, i_{\ell-1}, i_{\ell}\rangle\in\Nat^{\ell}\}$ is dense
  in $L_{0}$. 
\end{itemize}
Note\MMP[t]{Relax NOW!} that this sequence depends on the chosen sequence $\Folge{f}$ of
functions that represents the strategy for Demon.

\BeginProof (Second part)
Now assume that Demon has a winning strategy $\Folge{f}$; hence no
matter how Angel plays, Demon will win. For proving the assertion, we
have to construct a sequence of nowhere dense subsets the union of
which is $S$. In the first move, Angel chooses a closed interval
$L_{1} := J_{i_{1}}\subseteq L_{0}$ (we refer here to the enumeration given by
${\cal F}$ above, so the interval chosen by Angel has index
$i_{1}$). Demon's answer is then $$L_{2} := K_{i_{1}} := f_{1}(L_{0},
L_{1}) = f_{1}(L_{0}, J_{i_{1}}),$$ as constructed
above. In the next step, Angel selects $L_{3} := J_{i_{1}, i_{2}}$
among those closed intervals which are eligible, i.e., which are
contained in $\Interior{K}_{i_{1}}$ and have rational endpoints,
Demon's countermove is $$L_{4} := K_{i_{1}, i_{2}} := f_{2}(L_{0},
L_{1}, L_{2}, L_{3}) = f_{2}(L_{0}, J_{i_{1}}, K_{i_{1}}, J_{i_{1},
  i_{2}}).$$

In the n-th step, Angel selects $L_{2\cdot n-1} := J_{1_{1}, \dots,
  i_{n}}$ and Demon selects $L_{2\cdot n} := K_{i_{1}, \dots,
  i_{n}}$. Then we see that the sequence $L_{0}\supseteq
L_{1}\dots\supseteq L_{2\cdot n-1}\supseteq L_{2\cdot n} \dots$
decreases and $L_{2\cdot n} = f_{n}(L_{0}, L_{1}, \dots, L_{2\cdot
  n-1})$ holds, as required. 

Put $T := S\setminus L_{0}$ for convenience, then $\bigcap_{n\in\Nat}
L_{n}\subseteq T$ by assumption (after all, we assume that Demon
wins), put
\begin{equation*}
  G_{n} := \bigcup_{\langle i_{1}, \dots, i_{n}\rangle\in
    \Nat^{n}}\Interior{K}_{i_{1}, \dots, i_{n}}.
\end{equation*}
Then $G_{n}$ is open. Let $E := \bigcap_{n\in\Nat}G_{n}$. Given $x\in
E$, there exists a unique sequence $\Folge{i}$ such that $x\in
K_{i_{1}, \dots, i_{n}}$ for each $n\in\Nat$. Hence $x\in
\bigcap_{n\in\Nat}L_{n}\subseteq T$, so that $E\subseteq T$. But then
we can write $$S = L_{0}\setminus T\subseteq L_{0}\setminus T =
\bigcup_{n\in\Nat}(L_{0}\setminus G_{n}).$$
Because $\bigcup\{\Interior{K}_{i_{1}, \dots, i_{n-1}, i_{n}} \mid \langle
  i_{1}, \dots, i_{n-1}, i_{n}\rangle\in\Nat^{n}\}$ is dense
  in $L_{0}$ for each $n\in\Nat$ by construction, we conclude that
  $L_{0}\setminus G_{n}$ is nowhere dense, so $S$ is of the first category. 
\EndProof

Games are an interesting tool for proofs, as we can see in this example; we have shown already that games may be used for other purposes, e.g., demonstrating the each subset of $[0, 1]$ is Lebesgue measurable under the Axiom of Determinacy~\SetCite{Section 1}. Further examples for using games to derive properties in a metric space can be found, e.g., in Kechris' book~\cite{Kechris}.

\Subsection{A Gallery of Spaces and Techniques} 
\label{sec:top-gallery}
The discussion of the basic properties and techniques suggest that we now have a powerful collection of methods at our disposal. Indeed, we set up a small gallery of show cases, in which we demonstrate some approaches and methods. 

We first look at the use of topologies in logics from two different angles. The more conventional one is a direct application of the important Baire Theorem, which permits the construction of a model in a countable language of first order logic. Here the application of the theorem lies at the heart of the application, which is a proof of Gödel's Completeness Theorem. The other vantage point starts from a calculus of observations and develops the concept of topological systems from it, stressing an order theoretic point of view by perceiving topologies as complete Heyting algebras, when considering them as partially ordered subset of the power set of their carrier. Since partial orders may generate topologies on the set they are based on, this yields an interesting interplay between order and topology, which is reflected here in the Hofmann-Mislove Theorem. 

Then we return to the green pastures of classic applications and give a proof of the Stone-Weierstraß Theorem, one of the true classics. It states that a subring of the space of continuous functions on a compact Hausdorff space, which contains the constants, and which separates points is dense in the topology of uniform convergence. We actually give two proofs for this. One is based on a covering argument in a general space, it has a wide range of applications, of course. The second proof is no less interesting. It is essentially based on Weierstraß' original proof and deals with polynomials over $[0, 1]$ only; here concepts like elementary integration and uniform continuity are applied in a very concise and beautiful way.  

Finally, we deal with uniform spaces; they are a generalization of pseudometric spaces, but more specific than topological spaces. We argue that the central concept is closeness of points, which is, however, formulated in conceptual rather than quantitative terms. It is shown that many concepts which appear specific to the metric approach like uniform continuity or completeness may be carried into this context. Nevertheless, uniform spaces are topological spaces, but the assumption on having a uniformity available has some consequences for the associated topology. 

The reader probably misses Polish spaces in this little gallery. We deal with these spaces in depth, but since most of our applications of them are measure theoretic in nature, we deal with them in the context of a discussion of measures as a kind of natural habitat~\MeasCite. 

\Subsubsection{Gödel's Completeness Theorem}
\label{sec:goedel}
Gödel's Completeness Teorem states that a set of sentences of first order logic is consistent iff it has a model. The crucial part is the construction of a model for a consistent set of sentences. This is usually done through Henkin's approach, see, e.g.,~\cite[4.2]{Shoenfield}, \cite[Chapter 2]{Chang+Keisler} or~\cite[5.1]{Srivastava-Logic}. Rasiowa and Sikorski~\cite{Rasiowa-Sikorski-I} followed a completely different path in their topological proof by making use of Baire's Category Theorem and using the observation that in a compact topological space the intersection of a sequence of open and dense sets is dense again. The compact space is provided by the clopen sets of a Boolean algebra which in turn is constructed from the formulas of the first order language upon factoring. The equivalence relation is induced by the consistent set under consideration. 

We present the fundamental ideas of their proof in this section, since it is an unexpected application of a combination of the topological version of Stone's Representation Theorem for Boolean algebras and Baire's Theorem, hinted at already in Example~\ref{bool-alg-dense}. Since we assume that the reader is familiar with the semantics of first order languages, we do not want to motivate every definition for this area in detail, but we sketch the definitions, indicate the deduction rules, say what a model is, and rather focus on the construction of the model. The references given above may be used to fill in any gaps.

A slightly informal description of the first order language $\theLang$ with identity with which we will be working is given first. For this, we assume that we have a countable set $\{x_{n}\mid n\in \Nat\}$ of variables and countably many constants. Moreover, we assume countably many function symbols and countably many predicate symbols. In particular, we have a binary relation $==$, the identity. Each function and each predicate symbol has a positive arity.

These are the components of our language $\theLang$. 
\begin{description}
\item[Terms.] A variable is a term and a constant symbol is a term. If $f$ is a function symbol of arity $n$, and $t_{1}, \dots, t_{n}$ are terms, then $f(t_{1}, \dots, t_{n})$ is a term. Nothing else is a term.
\item[Atomic Formulas.] If $t_{1}$ and $t_{2}$ are terms, then $t_{1} == t_{2}$ is an atomic formula. If $p$ is a predicate symbol of arity $n$, and $t_{1}, \dots, t_{n}$ are terms, then $p(t_{1}, \dots, t_{n})$ is an atomic formula.
\item[Formulas.] An atomic formula is a formula. If $\phi$ and $\psi$ are formulas, then $\phi\wedge\psi$ and $\neg\phi$ are formulas. If $x$ is a variable and $\phi$ is a formula, then $\forall x.\phi$ is a formula. Nothing else is a formula.
\end{description}
Because there are countably many variables resp. constants, the language has countably many formulas. 

One usually adds parentheses to the logical symbols, but we do without, using them, however, freely, when necessary. We will use also disjunction [$\phi\vee\psi$ abbreviates $\neg(\neg\phi\wedge\neg\psi)$], implication [$\phi\to \psi$ for $\neg\phi\vee\psi$], logical equivalence [$\phi\leftrightarrow\psi$ for $(\phi\to \psi)\wedge(\psi\to \phi)$] and existential quantification [$\exists x.\phi$ for $\neg(\forall x.\neg\phi)$]. Conjunction and disjunction are associative.

We need logical axioms and inference rules as well. We have four groups of axioms
\begin{description}
\item[Propositional Axioms.] Each propositional tautology is an axiom.
\item[Identity Axioms.] $x==x$, when $x$ is a variable.
\item[Equality Axioms.] $y_{1}==z_{1} \to \dots\to y_{n}==z_{n}\to f(y_{1}, \dots, y_{n}) = f(z_{1}, \dots, z_{n})$, whenever $f$ is a function symbol of arity $n$, and $y_{1}==z_{1} \to \dots\to y_{n}==z_{n}\to p(y_{1}, \dots, y_{n}) \to  p(z_{1}, \dots, z_{n})$ for a predicate symbol of arity $n$.
\item[Substitution Axiom.] If $\phi$ is a formula, $\phi_{x}[t]$ is obtained from $\phi$ by freely substituting all free occurrences of variable $x$ by term $t$, then $\phi_{x}[t]\to \exists x.\phi$ is an axiom.
\end{description}

These are the inference rules.
\begin{description}
\item[Modus Ponens.] From $\phi$ and $\phi\to \psi$ infer $\psi$.
\item[Generalization Rule.] From $\phi$ infer $\forall x.\phi$. 
\end{description}

A \emph{\index{sentence}sentence} is a formula without free variables. Let $\Sigma$ be a set of sentences, $\phi$ a formula, then we denote that $\phi$ is deducible from $\Sigma$ by $\Sigma\vdash \phi$\MMP{$\Sigma\vdash \phi$}, i.e., iff there is a proof for $\phi$ in $\Sigma$. $\Sigma$ is called \emph{inconsistent} iff $\Sigma\vdash \bot$, or, equivalently, iff each formula can be deduced from $\Sigma$. If $\Sigma$ is not inconsistent, then $\Sigma$ is called \emph{consistent} or a \emph{theory}.  

Fix a theory $T$, and define 
\begin{equation*}
  \isEquiv{\phi}{\psi}{\sim} \text{ iff } T \vdash \phi\leftrightarrow \psi
\end{equation*}
for formulas $\phi$ and $\psi$, then this defines an equivalence relation on the set of all formulas. Let $\mathbf{B}_{T}$ be the set of all equivalence classes $[\phi]$\MMP{$\textbf{B}_T$}, and define
\begin{align*}
  [\phi]\wedge[\psi] & := [\phi\wedge\psi]\\
  [\phi]\vee[\psi] & := [\phi\vee\psi]\\
-[\phi] & := [\neg\phi].
\end{align*}
This defines a Boolean algebra structure on $\mathbf{B}_{T}$, the \emph{Lindenbaum \index{algebra!Lindenbaum}algebra} of $T$. The maximal element $\top$ of $\mathbf{B}_{T}$\MMP{Lindenbaum algebra} is $\{\phi\mid T\vdash \phi\}$, its minimal element $\bot$ is $\{\phi\mid T\vdash \neg \phi\}$. The proof that $\mathbf{B}_{T}$ is a Boolean algebra follows the lines of~\SetCite{1.5.7} closely, hence it can be safely omitted. It might be noted, however, that the individual steps in the proof require additional properties of $\vdash$, for example, one has to show that $T\vdash \phi$ and $T\vdash \psi$ together imply $T\vdash \phi\wedge\psi$. We trust that the reader is in a position to recognize and accomplish this; \cite[Chapter 4]{Srivastava-Logic} provides a comprehensive catalog of useful derivation rules with their proofs. 

Let $\phi$ be a formula, then denote by $\phi(k/p)$ the formula obtained in this way: 
\begin{itemize}
\item all bound occurrences of $x_{p}$ are replaced by $x_{\ell}$,
  where $x_{\ell}$ is the first variable among $x_{1}, x_{2}, \dots$ which does not occur in $\phi$,
\item all free occurrences of $x_{k}$ are replaced by $x_{p}$.
\end{itemize}
This construction is dependent on the integer $\ell$, so the formula $\phi(k/p)$ is not uniquely determined, but its class is. We have these representations in the Lindenbaum algebra for existentially resp. universally quantified formulas.

\BeginLemma{repr-sup-in-lindenbaum}
Let $\phi$ be a formula in $\theLang$, then we have for every $k\in \Nat$ 
\begin{enumerate}
\item\label{repr-sup-in-1} $\sup_{p\in\Nat}[\phi(k/p)] = [\exists x_{k}.\phi]$,
\item\label{repr-sup-in-2} $\inf_{p\in\Nat}[\phi(k/p)] = [\forall x_{k}.\phi]$.
\end{enumerate}
\EndLemma

\BeginProof
1.
Fix $k\in\Nat$, then we have $T\vdash \phi(k/p)\to \exists x_{k}.\phi$ for each $p\in \Nat$ by the $\exists$ introduction rule. This implies $[\phi(k/p)]\leq [\exists x_{k}.\phi]$ for all $p\in \Nat$, hence $\sup_{p\in\Nat}[\phi(k/p)] \leq [\exists x_{k}.\phi]$, thus $[\exists x_{k}.\phi]$ is an upper bound to $\{[\phi(k/p)]\mid p\in \Nat\}$ in the Lindenbaum algebra. We have to show that it is also the least upper bound, so take a formula $\psi$ such that $[\phi(k/p)]\leq[\psi]$ for all $k\in \Nat$. Let $q$ be an index such that $x_{q}$ does not occur free in $\psi$, then we conclude from $T\vdash \phi(k/p)\to \psi$ for all $p$ that $\exists x_{q}.\phi(k/q)\to \psi$. But $T\vdash \exists x_{k}.\phi \leftrightarrow \exists x_{q}.\phi(k/q)$, hence $T\vdash \exists x_{k}.\phi \to \psi$. This means that $[\exists x_{k}.\phi]$ is the least upper bound to $\{[\phi(k/p)]\mid p\in \Nat\}$, proving the first equality. 

2.
The second equality is established in a very similar way. 
\EndProof

These representations motivate 

\BeginDefinition{preserve-sups-and-infs}
Let $\fiF$ be an ultrafilter on the Lindenbaum algebra $\mathbf{B}_{T}$, $S\subseteq \mathbf{B}_{T}$. 
\begin{enumerate}
\item $\fiF$ \emph{preserves the supremum of $S$} iff $\sup S\in \fiF \Leftrightarrow s\in \fiF$ for some $s\in S$.
\item $\fiF$ \emph{preserves the infimum of $S$} iff $\inf S\in \fiF \Leftrightarrow s\in \fiF$ for all $s\in S$.
\end{enumerate}
\EndDefinition

Preserving the supremum of a set is similar to being inaccessible by joins (see Definition~\ref{scott-open-gen}), but inaccessibility refers to directed sets, while we are not making any assumption on $S$, except, of course, that its supremum exists in the Boolean algebra. Note also that one of the characteristic properties of an ultrafilters is that the join of two elements is in the ultrafilter iff it contains at least one of them. Preserving the supremum of a set strengthens this property \emph{for this particular set only}. 

The de Morgan laws and  $\fiF$ being an ultrafilter make it clear that $\fiF$ preserves $\inf S$ iff it preserves $\sup \{-s\mid  s\in S\}$, resp. that $\fiF$ preserves $\sup S$ iff it preserves $\inf\{-s\mid s\in S\}$. This cuts our work in half.

\BeginProposition{preserving-ultrafilter-exists}
Let $\Folge{S}$ be a sequence of subsets $S_{n}\subseteq \mathbf{B}_{T}$ such that $\sup S_{n}$ exists in $\mathbf{B}_{T}$. Then there exists an ultrafilter $\fiF$ such that $\fiF$ preserves the supremum of $S_{n}$ for all $n\in\Nat$. 
\EndProposition

\BeginProof
This is an application of Baire's Category Theorem~\ref{baire-locally-compact} and is discussed  in Example~\ref{bool-alg-dense}. We find there a prime ideal which does not preserve the supremum for $S_{n}$ for all $n\in\Nat$. Since the complement of a prime ideal in a Boolean algebra is an ultrafilter~\SetCite{Lemma 1.5.38, Lemma 1.5.37}, the assertion follows. 
\EndProof

So much for the syntactic side of our language $\theLang$. We will leave the ultrafilter $\fiF$ alone for a little while and turn to the semantics of the logic.

An \emph{interpretation} of $\theLang$ is given by a carrier set $A$, each constant $c$ is interpreted through an element $c_{A}$ of $A$, each function symbol $f$ with arity $n$ is assigned a map $f_{A}: A^{n}\to A$, and each $n$-ary predicate $p$ is interpreted through an $n$-ary relation $p_{A}\subseteq A^{n}$; finally, the binary predicate $==$ is interpreted through equality on $A$. We also fix a sequence $\{w_{n}\mid n\in\Nat\}$ of elements of $A$ for the interpretation of variables,  set ${\cal A} := (A, \{w_{n}\mid n\in\Nat\})$, and call ${\cal A}$ a \emph{\index{model}model} for the first order language. We then proceed inductively:
\begin{description}
\item[Terms.] Variable $x_{i}$ is interpreted by $w_{i}$. Assume that  the term $f(t_{1}, \dots, t_{n})$ is given. If the terms $t_{1}, \dots, t_{n}$\MMP{${\cal A}\models \phi$} are interpreted, through the respective elements $t_{A, 1}, \dots, t_{A, n}$ of $A$, then $f(t_{1}, \dots, t_{n})$ is interpreted through $f_{A}(t_{A, 1}, \dots, t_{A, n})$.
\item[Atomic Formulas.] The atomic formula $t_{1} == t_{2}$ is interpreted through $t_{A, 1} = t_{A, 2}$. If the $n$-ary predicate $p$ is assigned  $p_{A}\subseteq A^{n}$, then $p(t_{1}, \dots, t_{n})$ is interpreted as $\langle t_{A, 1}, \dots, t_{A, n}\rangle\in p_{A}$.

We denote by ${\cal A}\models \phi$ that the interpretation of the atomic formula $\phi$ yields the value true. We say that $\phi$ holds in ${\cal A}$. 
\item[Formulas.] Let $\phi$ and $\psi$ be formulas, then  ${\cal A}\models \phi\wedge\psi$ iff ${\cal A}\models \phi$ and ${\cal A}\models \psi$, and ${\cal A}\models \neg\phi$ iff ${\cal A}\models \phi$ is false. Let $\phi$  be the formula $\forall x_{i}. \psi$, then ${\cal A}\models \phi$ iff 
${\cal A}\models \psi_{x_{i}\mid a}$ for every $a\in A$, where $\psi_{x\mid a}$ is the formula $\psi$ with each free occurrence of $x$ replaced by $a$.  
\end{description}

Construct the ultrafilter $\fiF$ constructed in Proposition~\ref{preserving-ultrafilter-exists} for all possible suprema arising from existentially quantified formulas according to Lemma~\ref{repr-sup-in-lindenbaum}. There are countably many suprema, because the number of formulas is countable. This ultrafilter and the Lindenbaum algebra $\mathbf{B}_{T}$ will be used now for the construction of a \emph{\index{model}model\MMP{Model}} ${\cal A}$ for $T$ (so that ${\cal A}\models \phi$ holds for all $\phi\in T$). 

We will first need to define the carrier set $A$. Define for the variables $x_{i}$ and $x_{j}$ the equivalence relation $\approx$ through $\isEquiv{x_{i}}{x_{j}}{\approx}$ iff $[x_{i}==x_{j}]\in \fiF$; denote by $\hat{x}_{i}$ the $\approx$-equivalence class of $x_{i}$. The carrier set $A$ is defined as $\{\hat{x}_{n}\mid n\in\Nat\}$. 

Let us take care of the constants now. Given a constant $c$, we know that $\vdash \exists x_{i}.c==x_{i}$ by substitution. Thus $[\exists x_{i}.c==x_{i}] = \top\in \fiF$. But $[\exists x_{i}.c==x_{i}] = \sup_{i\in \Nat} [c==x_{i}]$, and  $\fiF$ preserves suprema, so we conclude that there exists $i$ with $[c==x_{i}]\in \fiF$. We pick this $i$ and define $c_{A} := \hat{x}_{i}$. Note that it does not matter which $i$ to choose. Assume that there is more than one. Since $[c==x_{i}]\in \fiF$ and $[c==x_{j}]\in \fiF$ implies $[c==x_{i}\wedge c==x_{j}]\in \fiF$, we obtain $[x_{i}==x_{j}]\in \fiF$, so the class is well defined. 

Coming to terms, let $t$ be a variable or a constant, so that it has an interpretation already, and assume that $f$ is a unary function. Then $\vdash \exists x_{i}.f(t) == x_{i}$, so that $[\exists x_{i}.f(t) == x_{i}]\in \fiF$, hence there exists $i$ such that $[f(t)==x_{i}]\in \fiF$, then put $f_{A}(t_{A}) := \hat{x}_{i}$. Again, if $[f(t) == x_{i}]\in \fiF$ and $[f(t) == x_{j}]\in\fiF$, then $[x_{i}==x_{j}]\in \fiF$, so that $f_{A}(c_{A})$ is well defined. The argument for the general case is very similar. Assume that terms $t_{1}, \dots, t_{n}$ have their interpretations already, and $f$ is a function with arity $n$, then $\vdash \exists x_{i}.f(t_{1}, \dots, t_{n}) == x_{i}$, hence we find $j$ with $[f(t_{1}, \dots, f(t_{n}) == x_{j}]\in \fiF$, so put $f_{A}(t_{A, 1}, \dots, t_{A, n}) := \hat{x}_{j}$. The same argument as above shows that this is well defined. 

Having defined the interpretation $t_{A}$ for each term $t$, we define for the $n$-ary relation symbol $p$ the relation $p_{A}\subseteq A^{n}$ by
\begin{equation*}
  \langle t_{A, 1}, \dots, t_{A, n}\rangle \in p_{A}\Leftrightarrow [p(t_{1}, \dots, t_{n}]\in\fiF
\end{equation*}
Then $p_{A}$ is well defined by the equality axioms. 

Thus ${\cal A}\models \phi$ is defined for each formula $\phi$, hence we know how to interpret each formula in terms of the Lindenbaum algebra of $T$ (and the ultrafilter $\fiF$). We can show now that a formula is valid in this model iff its class is contained in ultrafilter $\fiF$.
\BeginProposition{valid-iff-in-filter}
${\cal A}\models \phi$ iff $[\phi]\in \fiF$ holds for each formula $\phi$ of $\theLang$.
\EndProposition

\BeginProof
The proof is done by induction on the structure of formula $\phi$ and is straightforward, using the properties of an ultrafilter. For example, 
\begin{align*}
  {\cal A}\models \phi\wedge\psi 
& \Leftrightarrow
  {\cal A}\models \phi \text{ and }  {\cal A}\models \psi&\text{ (definition)} \\
& \Leftrightarrow
[\phi]\in\fiF\text{ and }[\psi]\in\fiF&\text{ (induction hypothesis)}\\
&\Leftrightarrow 
[\phi\wedge\psi]\in\fiF&\text{ ($\fiF$ is an ultrafilter)}
\end{align*}
For establishing the equivalence for universally quantified formulas $\forall x_{i}.\psi$, assume that $x_{i}$ is a  free variable in $\psi$ such that ${\cal A}\models \psi_{x_{i}\mid a} \Leftrightarrow [\psi_{x_{i}\mid a}]\in\fiF$ has been established for all $a\in A$. Then
\begin{align*}
  {\cal A}\models \forall x_{i}.\psi
&\Leftrightarrow 
{\cal A}\models \psi_{x_{i}\mid a}\text{ for all }a\in A&\text{ (definition)}\\
&\Leftrightarrow
[\psi_{x_{i}\mid a}]\in\fiF\text{ for all }a\in A&\text{ (induction hypothesis)}\\
&\Leftrightarrow
\sup_{a\in A}[\psi_{x_{i}\mid a}]\in\fiF&\text{ ($\fiF$ preserves the infimum)}\\
&\Leftrightarrow
[\forall x_{i}.\psi]\in\fiF&\text{ (by Lemma~\ref{repr-sup-in-lindenbaum})}
\end{align*}
This completes the proof. 
\EndProof

As a consequence, we have established this version of Gödel's Completeness Theorem:

\BeginCorollary{a-is-a-model-for}
${\cal A}$ is a model for the consistent set $T$ of formulas.
\QED
\EndCorollary

This approach demonstrates how a topological argument is used at the center of a  construction in logic. It should be noted, however, that the argument is only effective since the universe in which we work is countable. This is so because the Baire Theorem, which enables the construction of the ultrafilter, works for a countable family of open and dense sets. If, however, we work in an uncountable language $\theLang$, this instrument is no longer available (\cite[Exercise 2.1.24]{Chang+Keisler} points to a possible generalization). But even in the countable case one cannot help but note that the construction above depends on the Axiom of Choice, because we require an ultrafilter. The approach in~\cite[Exercise 2.1.22]{Chang+Keisler} resp.~\cite[Theorem 2.21]{Koppelberg} point to the construction of a filter without the help of a topology, but, alas, this filter is extended to an ultrafilter, and here the dreaded axiom is needed again. 

\Subsubsection{Topological Systems or: Topology Via Logic}
\label{sec:topology-via-logic}

This section investigates topological systems. They abstract from topologies being sets of subsets and focus on the order structure imposed by a topology instead. We focus on the interplay between a topology and the base space by considering these objects separately. A topology is considered a complete Heyting algebra, the carrier set is, well, a set of points, both are related through a validity relation $\models $ which mimics the $\in$-relation between a set and its elements. This leads to the definition of a topological system, and the question is whether this separation really bears fruits. It does; for example we may replace the point set by the morphisms from the Heyting algebra to the two element algebras $\Zwei$, giving sober spaces, and we show that, e.g., a Hausdorff space is isomorphic to such a structure. 

The interplay of the order structure of a topology and its topological obligations will be investigated through the Scott topology on a dcpo, a directed complete partial order, leading to the Hofmann-Mislove Theorem which characterizes compact sets that are represented as the intersection of the open sets containing them in terms of Scott open filters. 

Before we enter into a technical discussion, however, we put the following definitions on record.

\BeginDefinition{complete-heyting-algebra}
A partially ordered set $P$ is called a \emph{complete Heyting algebra}\index{algebra!Heyting} iff 
\begin{enumerate}
\item each finite subset $S$ has a join $\bigwedge S$,
\item each subset $S$ has a meet $\bigvee S$,
\item finite meets distribute over arbitrary joins, i.e., 
  \begin{equation*}
    a\wedge\bigvee S = \bigvee\{a\wedge s \mid s\in S\}
  \end{equation*}
holds for $a\in L, S\subseteq L$.
\end{enumerate}
A \emph{\index{algebra!Heyting!morphism}morphism} $f$ between the complete Heyting algebras $P$ and $Q$ is a map $f: P\to Q$ such that 
\begin{enumerate}
\item $f(\bigwedge S) = \bigwedge \Bild{f}{S}$ holds for finite $S\subseteq P$,
\item $f(\bigvee S) = \bigvee \Bild{f}{S}$ holds for arbitrary $S\subseteq P$. 
\end{enumerate}
$\|P, Q\|$\MMP{$\|P, Q\|$} denotes the set of all morphisms $P\to Q$. 
\EndDefinition

The definition of a complete Heyting algebra is a bit redundant, but never mind. Because the join and the meet of the empty set is a member of such an algebra, it contains a smallest element $\bot$ and a largest element $\top$, and $f(\bot) = \bot$ and $f(\top) = \top$ follows. A topology is a complete Heyting algebra with inclusion as the partial order~\SetCite{Exercise 29}. Sometimes, complete Heyting algebras are called \emph{\index{frame}frames}; but since the structure underlying the interpretation of modal logics are also called frames, we stick here to the longer name.

\BeginExample{compl-heyting-alg}
Call a lattice $V$ pseudo-complemented iff given $a, b\in V$, there exists $c\in V$ such that $x\leq c$ iff $x\wedge a \leq b$; $c$ is usually denoted by $a\to b$. A complete Heyting algebra is pseudo-complemented. In fact, let $c :=  \bigvee\{x \in V\mid x\wedge a \leq b\}$, then
\begin{equation*}
  c\wedge a = \bigvee\{x\in V\mid x\wedge a \leq b\}\wedge a = \bigvee\{x\wedge a\in V\mid x\wedge a \leq b\}\leq b
\end{equation*}
by the general distributive law, hence $x \leq c$ implies $x\wedge a \leq b$. Conversely, if $x\wedge a \leq b$, then $x\leq c$ follows. 
\EndExample

\BeginExample{compl-heyting-alg-foll}
Assume that we have a complete lattice $V$ which is pseudo-complemented. Then the lattice satisfies the general distributive law. In fact, given $a\in V$ and $S\subseteq V$,  we have $s\wedge a \leq \bigvee\{a\wedge b\mid b\in S\}$ ,
thus $s \leq a\to \bigvee\{a\wedge b\mid b\in S\}$ for all $s\in S$, from which we obtain $\bigvee S \leq a\wedge\bigvee\{a\wedge b\mid b\in S\}$, which in turn gives $a\wedge \bigvee S \leq a\wedge\bigvee\{a\wedge b\mid b\in S\}$. On the other hand, $\bigvee \{a\wedge b \mid  b\in S\}\leq \bigvee S$, and $\bigvee \{a\wedge b \mid  b\in S\}\leq a$, so that we obtain $\bigvee\{a\wedge b \mid  b\in S\}\leq a\wedge \bigvee S$. 

\EndExample

We note
\BeginCorollary{cha-is-lattice}
A complete Heyting algebra is a complete distributive lattice. 
\QED
\EndCorollary
 
Quite apart from investigating what can be said if open sets are replaced by an element of a complete Heyting algebra, and thus focussing on the order structure, one can argue as follows. Suppose we have observers and events, say, $X$ is the set of observers, and $A$ is the set of events. The observers are not assumed to have any structure, the events have a partial order making them a distributive lattice; an observation may be incomplete, so $a\leq b$ indicates that observing event $b$ contains more information than observing event $a$. If observer $x\in X$ observes event $a\in A$, we denote this as $x\models a$. The lattice structure should be compatible with the observations, that is, we want to have for $S\subseteq A$ that
\begin{align*}
  x\models \bigwedge S & \text{ iff }x\models a\text{ for all }a\in S, S\text{ finite},\\
x \models \bigvee S & \text{ iff }x\models a\text{ for some }a\in S, S\text{ arbitrary}.
\end{align*}
(recall $\bigwedge \emptyset = \top$ and $\bigvee\emptyset = \bot$). Thus our observations should be closed under finite conjunctions and arbitrary disjunctions; replacing disjunctions by intersections and conjunctions by unions, this shows a somewhat topological face. We define accordingly

\BeginDefinition{topol-syst}
A \emph{\index{topological system}topological system} $(\pTX, \oMX, \models)$ has a set $\pTX$ of points, a complete Heyting algebra $\oMX$ of observations, and a satisfaction relation $\models\ \subseteq\ \pTX\times\oMX$ (written as $x\models a$ for $\langle x, a\rangle\in\ \models$) such that we have for all $x\in \oMX$
\begin{itemize}
\item If $S\subseteq \oMX$ is finite, then $x\models \bigvee S$ iff $x\models a$ for all $a\in S$.
\item For $S\subseteq \oMX$ arbitrary, $x\models \bigvee S$ iff $x\models a$ for some $a\in S$.
\end{itemize}
The elements of $\pTX$ are called \emph{\index{topological system!points}points}, the elements of $\oMX$\MMP{$\pTX, \oMX$} are called \emph{\index{topological system!opens}opens}. 
\EndDefinition

We will denote a topological system $X = (\pTX,\oMX)$ usually without writing down the satisfaction relation, which is either explicity defined or understood from the context. 

\BeginExample{top-system-example}
\begin{enumerate}
\item\label{top-system-example-1} The obvious example for a topological system $D$ is a topological space $(X, \tau)$ with $\pTX[D] := X$ and $\oMX[D] := \tau$, ordered through inclusion. The satisfaction relation $\models $ is given by the containment relation $\in$, so that we have $x\models G$ iff $x\in G$ for $x\in \pTX[D]$ and $G\in\oMX[D]$.
\item\label{top-system-example-2} But it works the other way around as well. Given a topological system $X$, define for the open $a\in\oMX$ its \emph{\index{topological system!opens!extension}extension}\MMP{Extension $\eXT[\cdot]$}
  \begin{equation*}
    \eXT := \{x\in \pTX\mid x\models a\}.
  \end{equation*}
Then $\tau := \{\eXT\mid a\in \oMX\}$ is a topology on $\pTX$. In fact, $\emptyset = \eXT[\bot]$, $\pTX =\eXT[\top]$, and if $S\subseteq \tau$ is finite, say, $S = \{\eXT[a_{1}], \dots, \eXT[a_{n}]\}$,  then $\bigcap S = \eXT[\bigwedge_{i=1}^{n}a_{i}]$. Similarly, if $S = \{\eXT[a_{i}]\mid i\in I\} \subseteq \tau$ is an arbitrary subset of $\tau$, then $\bigcup S = \eXT[\bigvee_{i\in I}a_{i}]$. This follows directly from the laws of a topological system.
\item\label{top-system-example-3} Put $\Zwei\index{2\negthickspace2} := \{\bot, \top\}$, then this is a complete Heyting algebra\MMP{$\Zwei$}. Let $\oMX := A$ be another complete Heyting algebra, and put $\pTX := \|\oMX, \Zwei\|$ defining $x\models a$ iff $x(a) = \top$ then yields a topological system. Thus a point in this topological system is a morphism $\oMX\to \Zwei$, and a point satisfies the open $a$ iff it assigns $\top$ to it. 
\end{enumerate}
\EndExample

Next, we want to define morphisms between topological systems. Before we do that, we have another look at topological spaces and continuous maps. Recall that a map $f: X\to Y$ between topological spaces $(X, \tau)$ and $(Y, \theta)$ is $\tau$-$\theta$-continuous iff $\InvBild{f}{H}\in\tau$ for all $H\in\theta$. Thus $f$ spawns a map $f^{-1}: \theta\to \tau$~---~note the opposite direction. We have $x\in \InvBild{f}{H}$ iff $f(x)\in H$, accounting for containment. 

This leads to the definition of a morphism as a pair of maps, one working in the opposite direction of the other one, such that the satisfaction relation is maintained, formally:

\BeginDefinition{def-morph-top-syst}
Let $X$ and $Y$ be topological systems. Then $f: X\to Y$ is a \emph{\index{topological system!c-morphism}c-morphism} iff 
\begin{enumerate}
\item $f$ is a pair of maps $f = (\pTf, \oMf)$ with $\pTf: \pTX\to \pTX[Y]$, and $\oMf\in \|\oMX[Y], \oMX\|$ is a morphism for the underlying algebras.
\item $\pTf(x)\models_{Y} b$ iff $x\models_{X} \oMf(b)$ for all $x\in \pTX$ and all $b\in \oMX[Y]$.\MMP{$\pTf, \oMf$} 
\end{enumerate}
\EndDefinition

We  have indicated above for the reader's convenience in which system the satisfaction relation is considered. It is evident that the notion of continuity is copied from topological spaces, taking the slightly different scenario into account. 

\BeginExample{cont-morph-example}
Let $X$ and $Y$ be topological systems with $f: X\to Y$ a c-morphism. Let $(\pTX, \tau_{\pTX})$ and $(\pTX[Y], \tau_{\pTX[Y]})$ be the topological spaces generated from these systems through the extent of the respective opens, as in Example~\ref{top-system-example}, part~\ref{top-system-example-2}. Then $\pTf: \pTX\to \pTX[Y]$ is $\tau_{\pTX}$-$\tau_{\pTX[Y]}$-continuous. In fact, let $b\in\oMX[Y]$, then 
\begin{align*}
  x\in \InvBild{(\pTf)}{\eXT[b]}
& \Leftrightarrow
\pTf(x)\in\eXT[b]\\
& \Leftrightarrow
\pTf(x)\models b\\
& \Leftrightarrow
x\models \oMf(b),
\end{align*}
thus $\InvBild{(\pTf)}{\eXT[b]} = \eXT[\oMf(b)]\in\tau_{\pTX}$. 
\EndExample

This shows that continuity of topological spaces is a special case of c-morphisms between topological systems, in the same way as topological spaces are special cases of topological systems.

Let $f: X\to Y$ and $g: Y\to Z$ be c-morphisms, then their composition is defined as $g\circ f := (\pTf[g]\circ \pTf, \oMf\circ \oMf[g])$. The identity $id_{X}:X\to X$ is defined through $id_{X} := (id_{\pTX}, id_{\oMX})$.  
If, given the c-morphism $f: X\to Y$, there there is a c-morphisms $g: Y\to X$ with $g\circ f = id_{X}$ and $f\circ g = id_{Y}$, then $f$ is called a \emph{\index{topological system!homeomorphism}homeomorphism}.

\BeginCorollary{top-syst-form-category}
Topological systems for a category $\catTS$, the objects of which are topological systems, with c-morphisms as morphisms. \QED
\EndCorollary
 
Given a topological system $X$, the topological space $(\pTX, \tau_{\pTX})$ with $\tau_{\pTX} := \bigl\{\eXT\mid a\in \oMX\bigr\}$ is called the \emph{\index{topological system!spatialization}spatialization} of $X$ and denoted by $\funSP(X)$. We want to make $\funSP$ a (covariant) functor $\catTS\to \catTop$, the latter one denoting the category of topological spaces with continuous maps as morphisms\MMP{$\catTS, \catTop, \funSP$}. Thus we have to define the image $\funSP(f)$ of a c-morphism $f: X\to Y$. But this is fairly straightforward, since we have shown in Example~\ref{cont-morph-example} that $f$ induces a continuous map $(\pTX, \tau_{\pTX})\to (\pTX[Y], \tau_{\pTX[Y]})$. It is clear now that $\funSP: \catTS\to \catTop$ is a covariant functor. On the other hand, part~\ref{top-system-example-1} of Example~\ref{top-system-example} shows that we have a forgetful functor $\funV: \catTop\to \catTS$ with $\funV(X, \tau) := (\pTX, \oMX)$ with $\pTX := X$ and $\oMX := \tau$, and $\funV(f) := (f, f^{-1})$. These functors are related.

\BeginProposition{sp-is-right-adjoint-to-forgetful}
$\funSP$ is right adjoint to $\funV$. 
\EndProposition

\BeginProof
0.
Given a topological space $X$ and a topological system $A$ we have to find a bijection $\phi_{X, A}: \hom{\catTS}(\funV(X), A)\to \hom{\catTop}(X, \funSP(A))$ rendering these diagrams commutative:
\begin{equation*}
\xymatrix{
\hom{\catTS}(\funV(X), A)\ar[d]_{F_{*}}\ar[rr]^{\phi_{X, A}}
&&\hom{\catTop}(X, \funSP(A))\ar[d]^{(\funSP(F))_{*}}
\\
\hom{\catTS}(\funV(X), B)\ar[rr]_{\phi_{X, B}}
&&\hom{\catTop}(X, \funSP(B))
}
\end{equation*}
and
\begin{equation*}
\xymatrix{
\hom{\catTS}(\funV(X), A)\ar[d]_{(\funV(G))^{*}}\ar[rr]^{\phi_{X, A}}
&&\hom{\catTop}(X, \funSP(A))\ar[d]^{G^{*}}
\\
\hom{\catTS}(\funV(Y), A)\ar[rr]_{\phi_{Y, A}}
&&\hom{\catTop}(Y, \funSP(A))
}
\end{equation*}
where $F_{*} := \hom{\catTS}(\funV(X), F): f\mapsto F\circ f$ for $F: A\to B$ in $\catTS$, and $G^{*} := \hom{\catTop}(G, \funSP(A)): g\mapsto g\circ G$ for $G: Y\to X$ in $\catTop$, see~\CategCite{Section 1.5}. 

We define $\phi_{X, A}(\pTf, \oMf) := \pTf$, hence we focus on the component of a c-morphism which maps points to points. 

1.
Let us work on the first diagram. Take $f = (\pTf, \oMf): \funV(X)\to A$ as a morphism in $\catTS$, and let $F: A\to B$ be a c-morphism, $F = (\pTf[F], \oMf[F])$, then $\phi_{X, B}(F_{*}(f)) = \phi_{X, B}(F\circ f) = \pTf[F]\circ \pTf$, and $(\funSP(F))_{*}(\phi_{X, A}(f)) = \funSP(F) \circ \pTf = \pTf[F]\circ \pTf$. 

2.
Similarly, chasing $f$ through the second diagram for some continuous map $G: Y\to X$ yields
\begin{align*}
\phi_{Y, A}((\funV(G))^{*}(f)) 
& = \phi_{Y, A}((\pTf, \oMf)\circ (G, G^{-1}))\\ 
& = \pTf\circ G \\
& = G^{*}(\pTf)\\
& = G^{*}(\phi_{X, A}(f)). 
\end{align*}
This completes the proof.
\EndProof

Constructing $\funSP$, we went from a topological space to its associated topological system by exploiting the observation that a topology $\tau$ is a complete Heyting algebra. But we can travel in the other direction as well, as we will show now. 

Given a complete Heyting algebra $A$, we take the elements of $A$ as opens, and take all morphisms in $\|A, \Zwei\|$ as points, defining the relation $\models $ which connects the components through 
\begin{equation*}
x\models a\Leftrightarrow x(a) = \top.
\end{equation*} 
This construction was announced already in Example~\ref{top-system-example}, part~\ref{top-system-example-3}. In order to extract a functor from this construction, we have to cater for morphisms. In fact, let $\psi\in\|B, A\|$ a morphism $B\to A$ of the complete Heyting algebras $B$ and $A$, and $p\in\|A, \Zwei\|$ a point of $A$, then $p\circ \psi\in\|B, \Zwei\|$ is a point in $B$. Let $\catCha$ be the category of all complete Heyting algebras with $\hom{\catCha}(A, B) := \|B, A\|$, then we define the functor $\funLoc: \catCha\to \catTS$ through $\funLoc(A) := (\|A, \Zwei\|, A)$, and $\funLoc(\psi) := (\psi_{*}, \psi)$ for $\psi\in\hom{\catCha}(A, B)$ with $\psi^{*}(p) := p\circ \psi$. Thus $\funLoc(\psi): \funLoc(A)\to \funLoc(B)$, if $\psi: A\to B$ in $\catCha$. In fact, let $f := \funLoc(\psi)$, and $p\in\|A, \Zwei\|$ a point in $\funLoc(A)$, then we obtain for $b\in B$\MMP{$\catCha, \funLoc$} 
\begin{align*}
\pTf(p) \models b
& \Leftrightarrow
\pTf(p)(b) = \top\\
& \Leftrightarrow
(p\circ \psi)(b) = \top&\text{ (since }\pTf=p\circ \psi)\\
& \Leftrightarrow
p\models \psi(b)\\
& \Leftrightarrow
p\models \oMf(b)&\text{ (since }\oMf=\psi).
\end{align*}
This shows that $\funLoc(\psi)$ is a morphism in $\catTS$. $\funLoc(A)$ is called the \emph{\index{topological system!localization}localization}\MMP{Localization} of the complete Heyting algebra $A$. The topological system is called \emph{\index{topological system!localic}localic} iff it is homeomorphic to the localization of a complete Heyting algebra.  

We have also here a forgetful functor $\funV: \catTS\to \catCha$, and with a proof very similar to the one for Proposition~\ref{sp-is-right-adjoint-to-forgetful} one shows

\BeginProposition{loc-is-left-adjoint-to-forgetful}
$\funLoc$ is left adjoint to the forgetful functor $\funV$. 
\QED
\EndProposition

In a localic system, the points enjoy as morphisms evidently much more structure than just being flat points without a face,  in an abstract set. 

Before exploiting this wondrous remark, recall these notations, where $(P, \leq)$ is a reflexive and transitive relation:
\begin{align*}
  \uparrow\! p & := \{q\in P\mid q\geq p\},\\
\downarrow\!p & := \{q\in P\mid q \leq p\}.
\end{align*}

The following properties are stated just for the record.

\BeginLemma{upper-is-prime}
Let $a\in A$ with $A$ a complete Heyting algebra. Then $\uparrow\!a$ is a filter, and $\downarrow\!a$ is an ideal in $A$. \QED
\EndLemma

\BeginDefinition{}
Let $A$ be a complete Heyting algebra. 
\begin{enumerate}
\item $a\in A$ is called a \emph{\index{prime!element}prime element} iff $\downarrow\!a$ is a prime ideal.
\item The filter $F\subseteq A$ is called \emph{completely \index{prime!completely}prime} iff $\bigvee S\in F$ implies $s\in F$ for some $s\in S$, where $S\subseteq A$.
\end{enumerate}
\EndDefinition

Thus $a\in A$ is a prime element iff we may conclude from $\bigwedge S\leq a$ that there exists $s\in S$ with $s\leq a$, provided $S\subseteq A$ is finite. Note that a prime filter has the stipulated property for finite $S\subseteq A$, so a completely prime filter is a prime filter by implication.

\BeginExample{compl-prime-open}
Let $(X, \tau)$ be a topological space, $x\in X$, then
\begin{equation*}
  {\cal G}_{x} := \{G\in \tau\mid x\in G\}
\end{equation*}
is a completely prime filter in $\tau$. It is clear that ${\cal G}_{x}$ is a filter in $\tau$, since it is closed under finite intersections, and $G\in {\cal G}_{x}$ and $G\subseteq H$ implies $H\in {\cal G}_{x}$ for $H\in \tau$. Now let $\bigcup_{i\in I}S_{i}\in {\cal G}_{x}$ with $S_{i}\in\tau$ for all $i\in I$, then there exists $j\in I$ such that $x\in S_{j}$, hence $S_{j}\in {\cal G}_{x}$. 
\EndExample

Prime filters in a Heyting complete Heyting algebras have this useful property: if we have an element which is not in the filter, then we can find a prime element not in the filter dominating the given one. The proof of this property requires Zorn's Lemma, hence a variant of the Axiom of Choice.

\BeginProposition{prime-elements-not-in}
Let $F\subseteq A$ be a prime filter in the complete Heyting algebra $A$. Let $a\not\in F$, then there exists a prime element $p\in A$ with $a\leq p$ and $p\not\in F$. 
\EndProposition

\BeginProof
Let $Z := \{b\in A\mid a\leq b \text{ and }b\not\in F\}$, then $Z\not=\emptyset$, since $a\in Z$. We want to show that $Z$ is inductively ordered, hence take a chain $C\subseteq Z$, then $c := \sup\ C\in A$, since $A$ is a complete lattice. Clearly, $a\leq c$; suppose $c\in F$, then, since $F$ is completely prime, we find $c'\in C$ with $c'\in F$, which contradicts the assumption that $C\subseteq Z$. Hence $Z$ contains a maximal element $p$ by Zorn's Lemma. 

Since $p\in Z$, we have $a\leq p$ and $p\not\in F$, so we have to show that $p$ is a prime element. Assume that $x\wedge y\leq p$, then either of $x\vee p$ or $y\vee p$ is not in $F$: if both are in $F$, we have by distributivity $(x\vee p)\wedge (y\vee p) = (x\wedge y)\vee p = p$, so $p\in F$, since $F$ is a filter; this is a contradiction. Assume that $x\vee p\not\in F$, then $a\leq x\vee p$, since $a\leq p$, hence even $x\vee p\in Z$. Since $p$ is maximal, we conclude $x\vee p \leq p$, which entails $x\leq p$. Thus $p$ is a prime element. 
\EndProof
 
The reader might wish to compare this statement to an argument used in the proof of Stone's Representation Theorem, see~\SetCite{Section 1.5.7}. There it is used that we may find in a Boolean algebra for each ideal a prime ideal which contains it. The argumentation is fairly similar, but, alas, one works there in a Boolean algebra, and not in a complete Heyting algebra, as we do presently.

This is a characterization of completely prime filters and prime elements in a complete Heyting algebra in terms of morphisms into $\Zwei$. We will use this characterization later on.

\BeginLemma{prime-things}
Let $A$ be a complete Heyting algebra, then
\begin{enumerate}
\item\label{prime-things-1} $F\subseteq A$ is a completely prime filter iff $F = f^{-1}(\top) := \InvBild{f}{\{\top\}}$ for some $f\in\|A, \Zwei\|$.
\item\label{prime-things-2} $I = f^{-1}(\bot)$ for some $f\in\|A, \Zwei\|$ iff $I=\downarrow\!p$ for some prime element $p\in A$. 
\end{enumerate}
\EndLemma

\BeginProof
1.
Let $F\subseteq A$ be a completely prime filter, and define
\begin{equation*}
  f(a) :=
  \begin{cases}
    \top, & \text{ if } a\in F\\
\bot, & \text{ if }a\not\in F
  \end{cases}
\end{equation*}
Then $f: A\to \Zwei$ is a morphism for the complete Heyting algebras $A$ and $\Zwei$. Since $F$ is a filter, we have $f(\bigwedge S) = \bigwedge_{s\in S} f(s)$ for $S\subseteq A$ finite. Let $S\subseteq A$, then 
\begin{equation*}
  \bigvee_{s\in S}f(s) = \top \Leftrightarrow f(s) = \top \text{ for some }s\in S \Leftrightarrow f(\bigvee S) = \top,
\end{equation*}
since $F$ is completely prime. Thus $f\in \|A, \Zwei\|$ and $F = f^{-1}(\top)$. Conversely, given $f\in\|A, \Zwei\|$, the filter $f^{-1}(\top)$ is certainly completely prime. 

2.
Assume that $I = f^{-1}(\bot)$ for some $f\in \|A, \Zwei\|$, and put 
\begin{equation*}
  p := \bigvee\{a\in A\mid f(a) = \bot\}.
\end{equation*}
Since $A$ is complete, we have $p\in A$, and if $a\leq p$, then $f(a) = \bot$. Conversely, if $f(a) = \bot$, then $a\leq p$, so that $I=\downarrow\!p$; moreover, $I$ is a prime ideal, for $f(a)\wedge f(b) = \bot \Leftrightarrow f(a) = \bot\text{ or } f(b) = \bot$, thus $a\wedge b\in I$ implies $a\in I$ or $b\in I$. Thus $p$ is a prime element. Let, conversely, the prime element $p$ be given, then one shows as in part 1. that 
\begin{equation*}
  f(a) :=
  \begin{cases}
    \bot, & \text{ if } a\leq p\\
\top, & \text{ otherwise}
  \end{cases}
\end{equation*}
defines a member of $\|A, \Zwei\|$ with $\downarrow\!p = f^{-1}(\bot)$. 
\EndProof

Continuing Example~\ref{compl-prime-open}, we see that there exists for a topological space $X := (X, \tau)$ for each $x\in X$ an element $f_{x}\in\|\tau, \Zwei\|$ such that $f_{x}(G) = \top$ iff $x\in G$. Define the map $\Phi_{X}: X\to \|\tau, \Zwei\|$ through $\Phi_{X}(x) := f_{x}$ (so that $\Phi_{X}(x) = f_{x}$ iff ${\cal G}_{x} = f^{-1}_{x}(\top)$). We will investigate $\Phi_{X}$ now in a little greater detail\MMP{$\Phi_{X}$}.

\BeginLemma{phi-is-injective}
$\Phi_{X}$ is injective iff $X$ is a $T_{0}$-space.
\EndLemma

\BeginProof
Let $\Phi_{X}$ be injective, $x\not= y$, then ${\cal G}_{x}\not={\cal G}_{y}$. Hence there exists an open set $G$ which contains one of $x, y$, but not the other. If, conversely, $X$ is a $T_{0}$-space, then we have by the same argumentation ${\cal G}_{x}\not={\cal G}_{y}$ for all $x, y$ with $x\not=y$, so that $\Phi_{X}$ is injective.
\EndProof

Well, that's not too bad, because the representation of elements into $\|\tau, \Zwei\|$ is reflected by a (very basic) separation axiom. Let us turn to surjectivity. For this, we need to transfer reducibility to the level of open or closed sets; since this is formulated most concisely for closed sets, we use this alternative. A closed set is called irreducible iff each of its covers with closed sets implies its being covered already by one of them.

\BeginDefinition{irred-closed-set}
A closed set $F\subseteq X$ is called \emph{\index{irreducible}irreducible} iff  $F\subseteq \bigcup_{i\in I}F_{i}$ implies $F\subseteq F_{i}$ for some $i\in I$ for any family $(F_{i})_{i\in I}$ of closed sets. 
\EndDefinition

Thus a closed set $F$ is irreducible iff the open set $X\setminus F$ is a prime element in $\tau$. Let's see: Assume that $F$ is irreducible, and let $\bigcap_{i\in I} G_{i}\subseteq X\setminus F$ for some open sets $(G_{i})_{i\in I}$. Then $F\subseteq \bigcup_{i\in I}X\setminus G_{i}$ with $X\setminus G_{i}$ closed, thus there exists $j\in I$ with $F\subseteq X\setminus G_{j}$, hence $G_{j}\subseteq X\setminus F$. Thus $\downarrow\!(X\setminus F)$ is a prime ideal in $\tau$. One argues in exactly the same way for showing that if $\downarrow\!(X\setminus F)$ is a prime ideal in $\tau$, then $F$ is irreducible. 

Now we have this characterization of surjectivity of our map $\Phi_{X}$ through irreducible closed sets.

\BeginLemma{phi-is-onto}
$\Phi_{X}$ is onto iff for each irreducible closed set $F$ there exists $x\in X$ such that $F=\Closure{\{x\}}$. 
\EndLemma

\BeginProof
1.
Let $\Phi_{X}$ be onto, $F\subseteq X$ be irreducible. By the argumentation above, $X\setminus F$ is a prime element in $\tau$, thus we find $f\in\|\tau, \Zwei\|$ with $\downarrow\!(X\setminus F) = f^{-1}(\bot)$. Since $\Phi_{X}$ is into, we find $x\in X$ such that $f = \Phi_{X}(x)$, hence we have $x\not\in G \Leftrightarrow f(x)=\bot$ for all open $G\subseteq X$. It is then elementary to show that $F=\Closure{\{x\}}$. 

2.
Let  $f\in \|\tau, \Zwei\|$, then we know that $f^{-1}(\bot) = \downarrow\!G$ for some prime open $G$. Put $F := X\setminus G$, then $F$ is irreducible and closed, hence $F = \Closure{\{x\}}$ for some $x\in X$. Then we infer $f(H) = \top \Leftrightarrow x\in H$ for each open set $H$, so we have indeed $f = \Phi_{X}(x)$. Hence $\Phi_{X}$ is onto.   
\EndProof

Thus, if $\Phi_{X}$ is a bijection, we can recover (the topology on)
$X$ from the morphisms on the complete Heyting algebra $\|\tau, \Zwei\|$.

\BeginDefinition{def-sober}
A topological space $(X, \tau)$  is called
\emph{\index{topology!sober}sober}\footnote{The rumors in the domain
  theory community that a certain \emph{Johann
  Heinrich-Wilhelm Sober} was a skat partner of Hilbert's gardener at
  Göttingen could not be confirmed~---~anyway, what about the third man?} iff $\Phi_{X}: X\to \|\tau, \Zwei\|$ is a
bijection. 
\EndDefinition

Thus we obtain as a consequence this characterization.

\BeginCorollary{cor-def-sober}
Given a topological space $X$, the following conditions are equivalent
\begin{itemize}
\item $X$ is sober.
\item $X$ is a $T_{0}$-space and for each irreducible closed set $F$ there
  exists $x\in X$ with $F = \Closure{\{x\}}$.  
\end{itemize}
\EndCorollary

Exercise~\ref{ex-hausdorff-sober} shows that each Hausdorff space is
sober. This property is, however, seldom made use of the the context
of classic applications of Hausdorff spaces in, say, analysis. 

Before continuing, we generalize the Scott topology, which has been defined in Example~\ref{scott-open} for inductively ordered sets. The crucial property is closedness under joins, and we stated this properts in a linearly ordered set by saying that, if the supremum of a set $S$ is in a Scott open set $G$, then we should find an element $s\in S$ with $s\in G$. This will have to be relaxed somewhat. Let us analyze the argument why the intersection $G_{1}\cap G_{2}$ of two Scott open sets (old version) $G_{1}$ and $G_{2}$ is open by taking a set $S$ such that $\bigvee S\in G_{1}\cap G_{2}$. Because $G_{i}$ is Scott open, we find $s_{i}\in S$ with $s_{i}\in G_{i}$ ($i = 1, 2$), and because we work in a linear ordered set, we know that either $s_{1}\leq s_{2}$ or $s_{2}\leq s_{1}$. Assuming $s_{1}\leq s_{2}$, we conclude that $s_{2}\in G_{1}$, because open sets are upward closed, so that $G_{1}\cap G_{2}$ is indeed open. The crucial ingredient here is that we can find for two elements of $S$ an element which dominates both, and this is the key to the generalization.

We want to be sure that each directed set has an upper bound; this is the case, e.g., when we are working in a complete Heyting algebra. The structure we are defining now, however, is considerably weaker, but makes sure that we can do what we have in mind.

\BeginDefinition{def-dcpo}
A partially ordered set in which every directed subset has an upper bound is called a \emph{directed completed partial ordered set}, abbreviated as \emph{\index{dcpo}dcpo}\MMP{dcpo}. 
\EndDefinition

Evidently, complete Heyting algebras are dcpos, in particular topologies are under inclusion. Sober topological spaces with the specialization order induced by the open sets furnish another example for a dcpo.

\BeginExample{sober-tops-are-dcpos}
Let $X = (X, \tau)$ be a sober topological space. Hence the points in $X$ and the morphisms in $\|\tau, \Zwei\|$ are in a bijective correspondence. Define for $x, x'\in X$ the relation $x\sqsubseteq x'$ iff we have for all open sets $x\models G \Rightarrow x'\models G$ (thus $x\in G$ implies $x'\in G$). If we think that being contained in more open sets means having better information, $x\sqsubseteq x'$ is then interpreted as $x'$ being better informed than $x$; $\sqsubseteq$ is sometimes called the \emph{specialization order}.

Then $(X, \sqsubseteq)$ is a partially ordered set, antisymmetry following from the observation that a sober space is a $T_{0}$-space. But $(X, \sqsubseteq)$ is also a dcpo. Let $S\subseteq X$ be a directed set, then $L := \Bild{\Phi_{X}}{S}$ is directed in $\|\tau, \Zwei\|$. Define 
\begin{equation*}
  p(G) := 
  \begin{cases}
    \top, & \text{ if there exists $\ell\in L$ with }\ell(G) = \top\\
\bot, & \text{ otherwise}
  \end{cases}
\end{equation*}
We claim that $p\in \|\tau, \Zwei\|$. It is clear that $p(\bigvee W) = \bigvee_{w\in W} p(w)$ for $W\subseteq \tau$. Now let $W\subseteq \tau$ be finite, and assume that $\bigwedge \Bild{p}{W} = \top$, hence $p(w) = \top$ for all $w\in W$. Thus we find for each $w\in W$ some $\ell_{w}\in L$ with $\ell_{w}(w)=\top$. Because $L$ is directed, and $W$ is finite, we find an upper bound $\ell \in L$ to $\{\ell_{w}\mid w\in W\}$, hence $\ell(w) = \top$ for all $w\in W$, so that $\ell(\bigwedge W) = \top$, hence $p(\bigwedge W) = \top$. This implies $\bigwedge \Bild{p}{W} = p(\bigwedge W)$. Thus $p\in \|\tau, \Zwei\|$, so that there exists $x\in X$ with $x=\Phi_{X}(p)$. Clearly, $x$ is an upper bound to $S$. 
\EndExample

\BeginDefinition{scott-open-gen}
Let $(P, \leq)$ be a dcpo, then $U\subseteq P$ is called \emph{Scott open\index{open!Scott}} iff
\begin{enumerate}
\item $U$ is upward closed.
\item If $\sup\ S\in U$ for some directed set $S$, then there exists $s\in S$ with $s\in U$. 
\end{enumerate}
\EndDefinition

The second property can be described as \emph{inaccessability through directed joins}: If $U$ contains the directed join of a set, it must contain already one of its elements. The following example is taken from~\cite[p. 136]{Cont-Latt}.

\BeginExample{scott-open-example}
The powerset $\PowerSet{X}$ of a set $X$ is a dcpo under inclusion. The sets $\{{\cal F}\subseteq \PowerSet{X}\mid {\cal F}\text{ is of finite character}\}$ are Scott open (${\cal F}\subseteq\PowerSet{X}$ is of \emph{finite character} iff this condition holds: $F\in {\cal F}$ iff some finite subset of $F$ is in ${\cal F}$). Let ${\cal F}$ be of finite character. Then ${\cal F}$ is certainly upward closed. Now let $S := \bigcup{\cal S}\in{\cal F}$ for some directed set ${\cal S}\subseteq\PowerSet{X}$, thus there exists a finite subset $F\subseteq S$ with $F\in{\cal F}$. Because ${\cal S}$ is directed, we find $S_{0}\in {\cal S}$ with $F\subseteq S_{0}$, so that $S_{0}\in {\cal F}$. 
\EndExample

In a topological space, each compact set gives rise to a Scott open filter as a subset of the topology.

\BeginLemma{compact-yields-scott-open}
Let $(X, \tau)$ be a topological space, and $C\subseteq X$ compact, then
\begin{equation*}
  H(C) := \{U\in \tau\mid C\subseteq U\}
\end{equation*}
is a Scott open filter.
\EndLemma

\BeginProof
Since $H(C)$ is upward closed and a filter, we have to establish that it is not accessible by directed joins. In fact, let ${\cal S}$ be a directed subset of $\tau$ such that $\bigcup{\cal S}\in H(C)$. Then ${\cal S}$ forms a cover of the compact set $C$, hence there exists ${\cal S}_{0}\subseteq{\cal S}$ finite such that $C\subseteq \bigcup{\cal S}_{0}$. But ${\cal S}$ is directed, so ${\cal S}_{0}$ has an upper bound $S\in{\cal S}$, thus $S\in H(C)$. 
\EndProof

Scott opens form in fact a topology, and continuous functions are characterized in a fashion similar to Example~\ref{scott-continuous-map}. We just state and prove these properties for completeness, before entering into a discussion of the Hofmann-Mislove Theorem.

\BeginProposition{scott-open-top-cont}
Let $(P, \leq)$ be a dcpo,
\begin{enumerate}
\item $\{U\subseteq P\mid U\text{ is Scott open}\}$ is a topology on $P$, the Scott topology of $P$\index{topology!Scott}. 
\item $F\subseteq P$ is Scott closed iff $F$ is downward closed ($x\leq y$ and $y\in F$ imply $x\in F$) and closed with respect to suprema of directed subsets.
\item Given a dcpo $(Q, \leq)$, a map $f: P\to Q$ is continuous with respect to the corresponding Scott topologies iff $f$ preserves directed joins (i.e., if $S\subseteq P$ is directed, then $\Bild{f}{S}\subseteq Q$ is directed and $\sup\ \Bild{f}{S} = f(\sup\ S)$). 
\end{enumerate}
\EndProposition

\BeginProof
1.
Let $U_{1}, U_{2}$ be Scott open, and $\sup\ S\in U_{1}\cap U_{2}$ for the directed set $S$. Then there exist $s_{i}\in S$ with $s_{i}\in G_{i}$ for $i = 1, 2$. Since $S$ is directed, we find $s\in S$ with $s\geq s_{1}$ and $s\geq s_{2}$, and since $U_{1}$ and $U_{2}$ both are upward closed, we conclude $s\in U_{1}\cap U_{2}$. Because $U_{1}\cap U_{2}$ is plainly upward closed, we conclude that $U_{1}\cap U_{2}$ is  Scott open, hence the set of Scott opens is closed under finite intersections. The other properties of a topology are evidently satisfied. This establishes the first part.

2.
The characterization of closed sets follows directly from the one for open sets by taking complements.

3.
Let $f: P\to Q$ be Scott-continuous. Then $f$ is monotone: if $x\leq x'$, then $x'$ is contained in the closed set $\InvBild{f}{\downarrow\!f(x')}$, thus $x\in \InvBild{f}{\downarrow\!f(x')}$, hence $f(x)\leq f(x')$. Now let $S\subseteq P$ be directed, then $\Bild{f}{S}\subseteq Q$ is directed by assumption, and $S\subseteq \InvBild{f}{\downarrow\!(\sup_{s\in S}f(s))}$. Since the latter set is closed, we conclude that it contains $\sup\ S$, hence $f(\sup\ S) \leq \sup\ \Bild{f}{S}$. On the other hand, since $f$ is monotone, we know that $f(\sup\ S) \geq \sup\ \Bild{f}{S}$. Thus $f$ preserves directed joins. 

Assume that $f$ preserves directed joins, then, if $x\leq x'$, $f(x') = f(\sup\ \{x, x'\}) = \sup\ \{f(x), f(x')\}$ follows, hence $f$ is monotone. Now let $H\subseteq Q$ be Scott open, then $\InvBild{f}{H}$ is upward closed. Let $S\subseteq P$ be directed, and assume that $\sup\ \Bild{f}{S}\in H$, then there exists $s\in S$ with $f(s) \in H$, hence $s\in \InvBild{f}{H}$, which therefore is Scott open. Hence $f$ is Scott continuous. 
\EndProof

Following~\cite[Chapter II-1]{Cont-Latt}, we show that in a sober space there is an order morphism between Scott open filters and certain compact subsets. Preparing for this, we observe that in a sober space every open subset which contains the intersection of a Scott open filter is already an element of the filter.  This will turn out to be a consequence of the existence of prime elements not contained in a prime filter, as stated in Proposition~\ref{prime-elements-not-in}. 

\BeginLemma{open-is-contained-in-filter}
Let ${\cal F} \subseteq\tau$ be a Scott open filter of open subsets in a sober topological space $(X, \tau)$. If $\bigcap{\cal F}\subseteq U$ for the open set $U$, then $U\in {\cal F}$. 
\EndLemma

\BeginProof
0.
The plan\MMP{Plan} of the proof goes like this: Since ${\cal F}$ is Scott open, it is a prime filter in $\tau$. We assume that there exists an open set which contains the intersection, but which is not in ${\cal F}$. This is exactly the situation in Proposition~\ref{prime-elements-not-in}, so there exists an open set which is maximal with respect to not being a member of ${\cal F}$, and which is prime, hence we may represent this set as $f^{-1}(\bot)$ for some $f\in \|\tau, \Zwei\|$. But now sobriety kicks in, and we represent $f$ through an element $x\in X$. This will then lead us to the desired contradiction. 

1.
Because ${\cal F}$ is Scott open, it is a prime filter in $\tau$. Let $G := \bigcap {\cal F}$, and assume that $U$ is open with $G\subseteq U$ (note that we do not know whether or not $G$ is empty). Assume that $U\not\in{\cal F}$, then we obtain from Proposition~\ref{prime-elements-not-in} a prime open set $V$ which is not in ${\cal F}$, which contains $U$, and which is maximal. Since $V$ is prime, there exists $f\in\|\tau, \Zwei\|$ such that $\{H\in \tau\mid f(H) = \bot\} = \downarrow\!V$ by Lemma~\ref{prime-things}. Since $X$ is sober, we find $x\in X$ such that $\Phi_{X}(x) = f$, hence $X\setminus V = \Closure{\{x\}}$. 

2.
We claim that $\Closure{\{x\}}\subseteq G$. If this is not the case, we have $z\not\in H$ for some $H\in{\cal F}$ and $z\in\Closure{\{x\}}$. Because $H$ is open, this entails $\Closure{\{x\}}\cap H = \emptyset$, thus by maximality of $V$, $H\subseteq V$. Since ${\cal F}$ is a filter, this implies $V\in{\cal F}$, which is not possible. Thus $\Closure{\{x\}}\subseteq G$, hence $G\not=\emptyset$, and $X\setminus V \cap G = \emptyset$. Thus $U\cap G=\emptyset$, contradicting the assumption. 
\EndProof

This is a fairly surprising and strong statement, because we usually cannot conclude from $\bigcap {\cal F}\subseteq U$ that $U\in {\cal F}$ holds, when ${\cal F}$ is an arbitrary filter. But we work here under stronger assumptions: the underlying space is sober, so each point is given by a morphism for the underlying complete Heyting algebra \emph{and vice versa}. In addition we deal with Scott open filters. They have the pleasant property that they are inaccessible by directed suprema. 

But we may even say more, viz., that the intersection of these filters is compact. For, if we have an open cover of the intersection, the union of this cover is open, thus must be an element of the filter by the previous lemma. We may write the union as a union of a directed set of open sets, which then lets us apply the assumption that the filter is inaccessible.

\BeginCorollary{cor-is-contained-in-filter}
Let $X$ be sober, and ${\cal F}$ be a Scott open filter. Then $\bigcap {\cal F}$ is compact and nonempty. 
\EndCorollary

\BeginProof
Let $K := \bigcap {\cal F}$, and ${\cal S}$ be an open cover of $K$. Thus $U := \bigcup {\cal S}$ is open with $K\subseteq S$, hence $U\in {\cal F}$ by Lemma~\ref{open-is-contained-in-filter}. But $\bigcup {\cal S} = \bigcup\bigl\{\bigcup {\cal S}_0\mid {\cal S}_0\subseteq{\cal S}\text{ finite}\bigr\}$, and the latter collection is directed, so there exists ${\cal S}_{0}\subseteq {\cal S}$ finite with $\bigcup{\cal S}_{0}\in{\cal F}$. But this means ${\cal S}_{0}$ is a finite subcover of $K$, which consequently is compact. If $K$ is empty, $\emptyset\in{\cal F}$ by Lemma~\ref{open-is-contained-in-filter}, which is impossible.
\EndProof

This gives a complete characterization of the Scott open filters in a sober space. The characterization involves compact sets which are represented as the intersections of these filters. But we can represent only those compact sets $C$ which are upper sets in the specialization order, i.e., for which holds $x\in C$ and $x\sqsubseteq x'$ implies $x'\in C$. These sets are called \emph{\index{set!saturated}saturated}. Recall that $x\sqsubseteq x'$ means $x\in G \Rightarrow x'\in G$ for all open sets $G$, hence a set is saturated iff it equals the intersection of all open sets containing it. With this in mind, we state the \emph{\index{theorem!Hofmann-Mislove}Hofmann-Mislove Theorem}. 

\BeginTheorem{hofmann-mislove}
Let $X$ be a sober space. Then the Scott open filters are in one-to-one and order preserving correspondence with the non-empty saturated compact subsets of $X$ via ${\cal F}\mapsto \bigcap{\cal F}$.  
\EndTheorem

\BeginProof
We have shown in Corollary~\ref{cor-is-contained-in-filter} that the intersection of a Scott open filter is compact and nonempty; it is saturated by construction. Conversely, Lemma~\ref{compact-yields-scott-open} shows that we may obtain from a compact and saturated subset $C$ of $X$ a Scott open filter, the intersection of which must be $C$. It is clear that the correspondence is order preserving. 
\EndProof

It is quite important for the proof of Lemma~\ref{open-is-contained-in-filter} that the underlying space is sober. Hence it does not come as a surprise that the Theorem of Hofmann-Mislove can be used for a characterization of sober spaces as well~\cite[Theorem II-1.21]{Cont-Latt}.

\BeginProposition{char-sober-hofmann-mislove}
Let $X$ be a $T_{0}$-space. Then the following statements are equivalent.
\begin{enumerate}
\item\label{char-sober-hofmann-mislove-1} $X$ is sober.
\item\label{char-sober-hofmann-mislove-2} Each Scott open filter ${\cal F}$ on $\tau$ consists of all open sets containing $\bigcap{\cal F}$.
\end{enumerate}
\EndProposition

\BeginProof
\labelImpl{char-sober-hofmann-mislove-1}{char-sober-hofmann-mislove-2}: This follows from Lemma~\ref{open-is-contained-in-filter}.

\labelImpl{char-sober-hofmann-mislove-2}{char-sober-hofmann-mislove-1}: Corollary~\ref{cor-def-sober} tells us that it is sufficient to show that each irreducible closed sets is the closure of one point. 

Let $A\subseteq X$ be irreducible and closed. Then ${\cal F} := \{G \text{ open}\mid G\cap A\not=\emptyset\}$ is closed under finite intersections, since $A$ is irreducible. In fact, let $G$ and $H$ be open sets with $G\cap A\not=\emptyset$ and $H\cap A\not=\emptyset$. If $A\subseteq (X\setminus G)\cup (X\setminus H)$, then $A$ is a subset of one of these closed sets, say, $X\setminus G$, but then $A\cap H=\emptyset$, which is a contradiction. This implies that ${\cal F}$ is a filter, and ${\cal F}$ is obviously Scott open.

Assume that $A$ cannot be represented as $\Closure{\{x\}}$ for some $x$. Then $X\setminus\Closure{\{x\}}$ is an open set the intersection of which with $A$ is not empty, hence $X\setminus\Closure{\{x\}}\in{\cal F}$. We obtain from the assumption that $X\setminus A\in {\cal F}$, because with 
$ 
  K := \bigcap{\cal F} \subseteq \bigcap_{x\in X} (X\setminus\Closure{\{x\}}),
$ 
we have $K\subseteq X\setminus A$, and $X\setminus A$ is open. Consequently, $A\cap X\setminus A\not=\emptyset$, which is a contradiction. 

Thus there exists $x\in X$ such that $A = \Closure{\{x\}}$. Hence $X$ is sober. 
\EndProof

These are the first and rather elementary discussions of the interplay between topology and order, considered in a systematic fashion in domain theory. The reader is directed to~\cite{Cont-Latt} or to~\cite{Abramsky+Jung} for further information.  

\Subsubsection{The Stone-Weierstraß Theorem}
\label{sec:stone-weierstrass}

This section will see the classic Stone-Weierstraß Theorem on the approximation of continuous functions on a compact topological space. We need for this a ring of continuous functions, and show that ---~under suitable conditions~--- this ring is dense. This requires some preliminary considerations on the space of continuous functions, because this construction evidently requires a topology.

Denote for a topological space $X$ by \index{$\Cont$}$\Cont$\MMP{$\Cont$} the space of all continuous and bounded functions $f: X\to \Real$. 

The structure of $\Cont $ is algebraically fairly rich; just for the record:

\BeginProposition{cont-is-vector-lattice}
$\Cont$ is a real vector space which is closed under constants, multiplication and under the lattice operations. \QED
\EndProposition

This describes the algebraic properties, but we need a topology on this space, which is provided for by the supremum norm. Define for $f\in \Cont$ 
\begin{equation*}
  \|f\| := \sup_{x\in X}|f(x)|.
\end{equation*}
Then $(\Cont, \|\cdot \|)$ is an example for a normed linear (or vector) space.

\BeginDefinition{vector-space}
Let $V$ be a real vector space. A \emph{\index{norm}norm} $\|\cdot \|: V\to \pReal$ assigns to each vector $v$ a non-negative real number $\|v\|$ with these properties:
\begin{enumerate}
\item $\|v\| \geq 0$, and $\|v\| = 0$ iff $v = 0$.
\item $\|\alpha\cdot v\| = |\alpha|\cdot \|v\|$ for all $\alpha\in\Real$ and all $v\in V$.
\item $\|x+y\|\leq \|x\| + \|y\|$ for all $x, y\in V$.
\end{enumerate}
A vector space with a norm is called a \emph{normed \index{space!normed}space}.
\EndDefinition

It is immediate that a normed space is a metric space, putting $d(v, w) := \|v-w\|$. It is also immediate that $f\mapsto \|f\|$ defines a norm on $\Cont$. But we can say actually a bit more: with this definition of a metric, $\Cont$ is a complete metric space; we have established this for the compact interval $[0, 1]$ in Example~\ref{cont-is-complete} already. Let us have a look at the general case.

\BeginLemma{cx-is-complete}
$\Cont$ is complete with the metric induced by the supremum norm.
\EndLemma

\BeginProof
Let $\Folge{f}$ be a $\|\cdot \|$-Cauchy sequence in $\Cont$, then $(f_{n}(x))_{n\in\Nat}$ is bounded, and $f(x) := \lim_{n\to \infty}f_{n}(x)$ exists for each $x\in X$. Let $\epsilon>0$ be given, then we find $n_{0}\in \Nat$ such that $\|f_{n}-f_{m}\|< \epsilon$ for all $n, m\geq n_{0}$, thus $|f(x)-f_{n}(x)|\geq \epsilon$ for $n\geq n_{0}$. This inequality holds for each $x\in X$, so that we obtain $\|f-f_{n}\|\leq \epsilon$ for $n\geq n_{0}$. It implies also that $f\in\Cont$.  
\EndProof

Normed spaces for which the associated metric space are special, so they deserve their own name.

\BeginDefinition{banach-space}
A normed space $(V, \|\cdot \|)$ which is complete in the metric associated with $\|\cdot \|$ is called a \emph{Banach \index{space!Banach}space}. 
\EndDefinition

The topology induced by the supremum norm is called the \emph{\index{topology!uniform convergence}topology of uniform convergence}, so that we may restate Lemma~\ref{cx-is-complete} by saying the $\Cont$ is closed under uniform convergence. A helpful example is Dini's \index{theorem!Dini}Theorem for uniform convergence on $\Cont$ for compact $X$. It gives a criterion of uniform convergence, provided we know already that the limit is continuous.

\BeginProposition{dini-uniform-convergence}
Let $X$ be a compact topological space, and assume that $\Folge{f}$ is a sequence of continuous functions which increases monotonically to a continuous function $f$. Then $\Folge{f}$ converges uniformly to $f$. 
\EndProposition

\BeginProof
We know that $f_{n}(x) \leq f_{n+1}(x)$ holds for all $n\in\Nat$ and all $x\in X$, and that $f(x) := \sup_{n\in\Nat}f_{n}(x)$ is continuous. Let $\epsilon>0$ be given, then $F_{n} := \{x\in X\mid f(x) \geq f_{n}(x) - \epsilon\}$ defines a closed set with $\bigcap_{n\in\Nat} F_{n} = \emptyset$, moreover, the sequence $\Folge{F}$ decreases. Thus we find $n_{0}\in \Nat$ with $F_{n} = \emptyset$ for $n\geq n_{0}$, hence $\|f-f_{n}\| < \epsilon$ for $n\geq n_{0}$. 
\EndProof

The goal of this section is to show that, given the compact topological space $X$, we can approximate each continuous real function uniformly through elements of a subspace of $\Cont$. It is plain that this subspace has to satisfy some requirements: it should
\begin{itemize}
\item be a vector space itself,
\item contain the constant functions,
\item separate points,
\item be closed under multiplication.
\end{itemize}
Hence it is in particular a subring of the ring $\Cont$. Let $A$ be such a subset, then we want to show that the closure $\Closure{A}$ with respect to uniform convergence equals $\Cont$. We will show first that $\Closure{A}$ is closed under the lattice operations, because we will represent an approximating function as the finite supremum of a finite infimum of simpler approximations. So the first goal will be to establish closure under $\inf$ and $\sup$. Recall that 
\begin{align*}
  f\wedge g & = \frac{1}{2}\cdot (f + g - |f-g|),\\
f\vee g & = \frac{1}{2}\cdot (f + g + |f-g|).
\end{align*}
Now it is easy to see that $\Closure{A}$ is closed under the vector
space operations, if $A$ is. Our first step boils down to showing that
$|f|\in \Closure{A}$ if $f\in A$. Thus, given $f\in A$, we have to
find a sequence $\Folge{f}$ such that $|f|$ is the uniform limit of
this sequence. \MMP{But wait!} It is actually enough to show that
$t\mapsto \sqrt{t}$ can be approximated uniformly on the unit interval
$[0, 1]$, because we know that $|f| = \sqrt{f^{2}}$ holds. It is enough to do this on the unit interval, as we will see below.

\BeginLemma{approx-square-root}
There exists a sequence $\Folge{f}$ in $\Cont[[0, 1]]$ which converges uniformly to the function $t\mapsto \sqrt{t}$. 
\EndLemma

\BeginProof
Define inductively for $t\in [0, 1]$
\begin{align*}
  f_{0}(t) & := 0,\\
f_{n+1}(t) & := f_{n}(t) + \frac{1}{2}\cdot (t - f_{n}^{2}(t)).
\end{align*}
We show by induction that $f_{n}(t) \leq \sqrt{t}$ holds. This is clear for $n = 0$. If we know already that the assumption holds for $n$, then we write
\begin{equation*}
  \sqrt{t} - f_{n+1}(t) = \sqrt{t} - f_{n}(t) - \frac{1}{2}\cdot (t - f_{n}^{2}(t)) = (\sqrt{t}-f_{n}(t))\cdot \bigl((1-\frac{1}{2}\cdot (\sqrt{t}+f_{n}(t))\bigr). 
\end{equation*}
Because $t\in[0, 1]$ and from the induction hypothesis, we have $\sqrt{t}+f_{n}(t)\leq 2\cdot \sqrt{t}\leq 2$, so that $\sqrt{t} - f_{n+1}(t) \geq 0$. 

Thus we infer that $f_{n}(t) \leq \sqrt{t}$ for all $t\in[0, 1]$, and $\lim_{n\to \infty}f_{n}(t) = \sqrt{t}$. From Dini's Proposition~\ref{dini-uniform-convergence} we now infer that the convergence is uniform. 
\EndProof

This is the desired consequence from this construction.

\BeginCorollary{ring-is-lattice}
Let $X$ be a compact topological space, and let $A\subseteq \Cont$ be a ring of continuous functions which contains the constants, and which is closed under uniform convergence. Then $A$ is a lattice.
\EndCorollary

\BeginProof
It is enough to show that $A$ is closed under taking absolute values. Let $f\in A$, then we may and do assume that $0\leq f \leq 1$ holds (otherwise consider $(f-\|f\|)/\|f\|$, which is an element of $A$ as well). Because $|f| = \sqrt{f^{2}}$, and the latter is a uniform limit of elements of $A$ by Lemma~\ref{approx-square-root}, we conclude $|f|\in A$, which entails $A$ being closed under the lattice operations. 
\EndProof

We are now in a position to establish the classic Stone-Weierstraß Theorem, which permits to conclude that a ring of bounded continuous functions on a compact topological space $X$ is dense with respect to uniform convergence in $\Cont$, provided it contains the constants and separates points. The latter condition is obviously necessary, but has not been used in the argumentation so far. It is clear, however, that we cannot do without this condition, because $\Cont$ separates points, and it is difficult to see how a function which separates points could be approximated from a collection which does not.

The polynomials on a compact interval in the reals are an example for a ring which satisfies all these assumptions. This collection shows also that we cannot extend the result to a non-compact base space like the reals. Take $x\mapsto \sin x$ for example; this function cannot be approximated uniformly over $\Real$ by polynomials. For, assume that given $\epsilon>0$ there exists a polynomial $p$ such that $\sup_{x\in \Real}|p(x) - \sin x| < \epsilon$, then we would have 
$
-\epsilon - 1 < p(x) < 1 + \epsilon
$
for all $x\in \Real$, which is impossible, because a polynomial is unbounded.  

Here, then, is the Stone-Weierstraß \index{theorem!Stone-Weierstraß}Theorem for compact topological spaces.

\BeginTheorem{stone-weierstrass}
Let $X$ be a compact topological space, and $A\subseteq \Cont$ be a ring of functions which separates points, and which contains all constant functions. 
\EndTheorem

\BeginProof
0.
Our goal\MMP{Approach} is to find for some given $f\in \Cont$ and an arbitrary $\epsilon>0$ a function $F\in \Closure{A}$ such that $\|f-F\|< \epsilon$. Since $X$ is compact, we will find $F$ through a refined covering argument in the following way. If $a, b\in X$ are given, we find a continuous function $f_{a, b}\in A$ with $f_{a, b} = f(a)$ and $f_{a, b} = f(b)$. From this we construct a cover, using sets like $\{x\mid f_{a, b}(x) < f(x) + \epsilon\}$ and $\{x\mid f_{a, b}(x) > f(x) - \epsilon\}$, extract finite subcovers and construct from the corresponding functions the desired function through suitable lattice operations. 

1.
Fix $f\in \Cont$ and $\epsilon>0$. Given distinct point $a\not= b$, we find a function $h\in A$ with $h(a) \not= h(b)$, thus 
\begin{equation*}
  g(x) := \frac{h(x) - h(a)}{h(b) - h(a)}
\end{equation*}
defines a function  $g\in A$ with $g(a) = 0$ and $g(b) = 1$. Then 
\begin{equation*}
  f_{a, b}(x) := (f(b) - f(a))\cdot g(x) + f(a)
\end{equation*}
is also an element of $A$ with $f_{a, b}(a) = f(a)$ and $f_{a, b}(b) = f(b)$. Now define
\begin{align*}
  U_{a, b} & := \{x\in X \mid  f_{a, b}(x) < f(x) + \epsilon\},\\
V_{a, b} & := \{x\in X \mid  f_{a, b}(x) > f(x) - \epsilon\},
\end{align*}
then $U_{a, b}$ and $V_{a, b}$ are open sets containing $a$ and $b$. 

2.
Fix $b$, then $\{U_{a, b}\mid a\in X\}$ is an open cover of $X$, so we can find points $a_{1}, \dots, a_{k}$ such that $\{U_{a_{1}, b}, \dots, U_{a_{k}, b}\}$ is an open cover of $X$ by compactness. Thus
\begin{equation*}
  f_{b} := \bigwedge_{i=1}^{k}f_{a_{i}, b}
\end{equation*}
defines an element of $\Closure{A}$ by Corollary~\ref{ring-is-lattice}. We have $f_{b}(x) < f(x) + \epsilon$ for all $x\in X$, and we know that $f_{b}(x) > f(x) - \epsilon$ for all $x\in V_{b} := \bigcap_{i=1}^{k}V_{a_{i}, b}$. The set $V_{b}$ is an open neighborhood of $b$, so from the open cover $\{V_{b}\mid  b\in X\}$ we find $b_{1}, \dots, b_{\ell}$ such that $X$ is covered through $\{V_{b_{1}}, \dots , V_{b_{\ell}}\}$. Put 
\begin{equation*}
  F := \bigvee_{i=1}^{\ell}f_{b_{i}},
\end{equation*}
then $f_{\epsilon}\in\Closure{A}$ and $\|f-F\|<\epsilon$. 
\EndProof

This is the example already discussed above.

\BeginExample{stone-weierstrass-polynomials}
Let $X := [0, 1]$ be the closed unit interval, and let $A$ consist of all polynomials $\sum_{i=0}^{n}a_{i}\cdot x^{i}$ for $n\in \Nat$ and $a_{0}, \dots, a_{n}\in\Real$. Polynomials are continuous, they form a vector space and are closed under multiplication. Moreover, the constants are polynomials. Thus we obtain from the Stone-Weierstraß Theorem~\ref{stone-weierstrass} that every continuous function on $[0, 1]$ can be uniformly approximated through a sequence of polynomials. 
\EndExample

It is said that Oscar Wilde\MMP[t]{Oscar Wilde} could resist everything but a temptation. The author concurs. Here is a classic proof of the Weierstraß Approximation Theorem, the original form of Theorem~\ref{stone-weierstrass}, which deals with polynomials on $[0, 1]$ only, and establishes the statement given in Example~\ref{stone-weierstrass-polynomials}. We will give this proof now, based on the discussion in the classic~\cite[§ II.4.1]{Hilbert-Courant}. This proof is elegant and based on the manipulation of specific functions (we are all too often focussed on our pretty little garden of beautiful abstract structures, all too often in danger of loosing the contact to concrete mathematics, and our roots).

As a preliminary consideration, we will show that 
\begin{equation*}
  \lim_{n\to \infty}\frac{\int_{\delta}^{1}(1-v^{2})^{n}\ dv}{\int_{0}^{1}(1-v^{2})^{n}\ dv} = 0
\end{equation*}
for every $\delta\in]0, 1[$. Define for this
\begin{align*}
  J_{n} & := \int_{0}^{1}(1-v^{2})^{n}\ dv,\\
J_{n}^{*} & := \int_{\delta}^{1}(1-v^{2})^{n}\ dv.
\end{align*}
(we will keep these notations for later use). We have
\begin{equation*}
  J_{n}  > \int_{0}^{1}(1-v)^{n}\ dv
 = \frac{1}{n+1}
\end{equation*}
and
\begin{equation*}
J_{n}^{*}  = \int_{\delta}^{1}(1-v^{2})^{n}\ dv
 < (1-\delta^{2})^{n}\cdot (1-\delta)
< (1-\delta^{2})^{n}.
\end{equation*}
Thus 
\begin{equation*}
\frac{J_{n}^{*}}{J_{n}}  < (n+1)\cdot (1-\delta^{2})^{n}\to 0.
\end{equation*}
This establishes the claim.

Let $f: [0, 1]\to \Real$ be continuous. Given $\epsilon>0$, there exists $\delta>0$ such that $|x-y|<\delta$ implies $|f(x) - f(y)|< \epsilon$ for all $x\in[0, 1]$, since $f$ is uniformly continuous by Proposition~\ref{compact-cont-unif-cont}. Thus $0\leq v < \delta$ implies $|f(x+v)-f(x)|< \epsilon$ for all $x\in[0, 1]$.  

Put
\begin{align*}
  Q_{n}(x) & := \int_{0}^{1}f(u)\cdot \bigl(1-(u-x)^{2}\bigr)^{n}\ du,\\
  P_{n}(x) & := \frac{Q_{n}(x)}{2\cdot J_{n}}.
\end{align*}
We will show that $P_{n}$ converges to $f$ in the topology of uniform convergence. 

We note first that $Q_{n}$ is a polynomial of degree $2 n$. In fact, put
\begin{equation*}
  A_{j} := \int_{0}^{1}f(u)\cdot u^{j}\ du
\end{equation*}
for $j \geq 0$, expanding yields the formidable representation
\begin{equation*}
  Q_{n}(x) = \sum_{k=0}^{n}\sum_{j=0}^{2 k}\binom{n}{k}\binom{2 k}{j}(-1)^{n-k+j}A_{j} \cdot x^{2 k-j}.
\end{equation*}

Let us work on the approximation. We fix $x\in[0, 1]$, and note that the inequalities derived below do not depend on the specific choice of $x$. Hence they provide a uniform approximation.

Substitute $u$ by $v+x$ in $Q_{n}$; this yields
\begin{align*}
  \int_{0}^{1}f(u) \bigl(1-(u-x)^{2}\bigr)^{n}\ du 
& =
\int_{\-x}^{1-x}f(v+x)(1-v^{2})^{n}\ dv\\
& = I_{1} + I_{2} + I_{3}
\end{align*}
with
\begin{align*}
I_{1} & := \int_{-x}^{-\delta}f(v+x)(1-v^{2})^{n}\ dv,\\
I_{2} & := \int_{-\delta}^{+\delta}f(v+x)(1-v^{2})^{n}\ dv,\\
I_{3} & := \int_{+\delta}^{1-x}f(v+x)(1-v^{2})^{n}\ dv.
\end{align*}
We work on these integrals separately. Let $M := \max_{0\leq x \leq 1}|f(x)|$, then
\begin{equation*}
  I_{1}  \leq M \int_{-1}^{-\delta}(1-v^{2})^{n}\ dv
 = M\cdot  J_{n}^{*},
\end{equation*}
and
\begin{equation*}
I_{3}  \leq M \int_{\delta}^{1}(1-v^{2})^{n}\ dv
 = M\cdot  J_{n}^{*}.
\end{equation*}
We can rewrite $I_{2}$ as follows:
\begin{align*}
  I_{2} & = f(x) \int_{-\delta}^{+\delta}(1-v^{2})^{n}\ dv + 
\int_{-\delta}^{+\delta}\bigl(f(x+v)-f(x)\bigr) (1-v^{2})^{n}\ dv\\
& = 2 f(x) (J_{n}-J_{n}^{*}) + 
\int_{-\delta}^{+\delta}\bigl(f(x+v)-f(x)\bigr) (1-v^{2})^{n}\ dv.
\end{align*}
From the choice of $\delta$ for $\epsilon$ we obtain
\begin{equation*}
  \bigl|\int_{-\delta}^{+\delta}\bigl(f(x+v)-f(x)\bigr) (1-v^{2})^{n}\ dv\bigr|
\leq \epsilon \int_{-\delta}^{+\delta}(1-v^{2})^{n}\ dv
< \epsilon \int_{-1}^{+1}(1-v^{2})^{n}\ dv = 2 \epsilon\cdot  J_{n}
\end{equation*}
Combining these inequalities, we obtain
\begin{equation*}
  |P_{n}(x) - f(x)| < 2 M\cdot  \frac{J_{n}^{*}}{J_{n}} + \epsilon.
\end{equation*}
Hence the difference can be made arbitrarily small, which means that $f$ can be approximated uniformly through polynomials.

The two approaches presented are structurally very different, it would be difficult to recognize the latter as a precursor of the former. While both make substantial use of uniform continuity, the first one is an existential proof, constructing two covers from which to choose a finite subcover each, and from this deriving the existence of an approximating function. It is non-constructive because it would be difficult to construct an approximating function from it, even if the ring of approximating functions is given by a base for the underlying vector space. The second one, however, starts also from uniform continuity and uses this property to find a suitable bound for the difference of the approximating polynomial and the function proper through integration. The representation of $Q_{n}$ above shows what the constructing polynomial looks like, and the coefficients of the polynomials may be computed (in principle, at least). And, finally, the abstract situation gives us a greater degree of freedom, since we deal with a ring of continuous functions observing certain properties, while the original proof works for the class of polynomials only. 

\Subsubsection{Uniform Spaces}
\label{sec:uniform-spaces}

This section will give a brief introduction to uniform spaces. The objective is to demonstrate in what ways the notion of a metric space can be generalized without arriving at the full generality of topological spaces, but retaining useful properties like completeness or uniform continuity. While pseudometric spaces formulate the concept of two points to be close to each other through a numeric value, an general topological spaces use the concept of an open neighborhood, uniform spaces formulate neighborhoods on the Cartesian product. This concept is truly in the middle: each pseudometric generates neighborhoods, and from a neighborhood we may obtain the neighborhood filter for a point. 

For motivation and illustration, we consider a pseudometric space $(X, d)$ and say that two point are neighbors iff their distance is smaller that $r$ for some fixed $r>0$; the degree of neighborhood is evidently depending on $r$. The set
\begin{equation*}
  V_{d, r} := V_{r} := \{\langle x, y\rangle\mid d(x, y)< r\}
\end{equation*}
is then the collection of all neighbors\MMP[h]{$V_{d, r}; V_{r}$}. We may obtain from $V_{r}$ the neighborhood $B(x, r)$ for some point $x$ upon extracting all $y$ such that $\langle x, y\rangle\in V_{r}$, thus 
\begin{equation*}
  B(x, r) = V_{r}[x] := \{y\in X\mid \langle x, y\rangle\in V_{r}\}. 
\end{equation*}
The collection of all these neighborhoods observes these properties.
\begin{enumerate}
\item The diagonal $\Delta := \Delta_{X}$ is contained in $V_{r}$ for all $r>0$, because $d(x, x) = 0$.
\item $V_{r}$ is ---~as a relation on $X$---~symmetric: $\langle x, y\rangle\in V_{r}$ iff $\langle y, x\rangle\in V_{r}$, thus $V_{r}^{-1} = V_{r}$. This property reflects the symmetry of $d$.
\item $V_{r}\circ V_{s}\subseteq V_{r+s}$ for $r, s>0$; this property is inherited from the triangle inequality for $d$. 
\item $V_{r_{1}}\cap V_{r_{2}} = V_{\min\{r_{1}, r_{2}\}}$, hence this collection is closed under finite intersections.
\end{enumerate} 

It is convenient to consider not only these immediate neighborhoods but rather the filter generated by them on $X\times X$ (which is possible because the empty set is not contained in this collection, and the properties above shows that they form the base for a filter indeed). This leads to this definition of a uniformity. It focusses on the properties of the neighborhoods rather than on that of a pseudometric, so we formulate it for a set in general.

\BeginDefinition{def-uniformity}
Let $X$ be a set. A filter $\filterFont{u}$  on $\PowerSet{X\times X}$ is called a \emph{\index{uniformity}uniformity on $X$} iff these properties are satisfied
\begin{enumerate}
\item $\Delta\subseteq U$ for all $U\in\filterFont{u}$.
\item If $U\in\filterFont{u}$, then $U^{-1}\in\filterFont{u}$.
\item If $U\in\filterFont{u}$, there exists $V\in\filterFont{u}$ such that $V\circ V\subseteq U$.
\item $\filterFont{u}$ is closed under finite intersections.
\item If $U\in\filterFont{u}$ and $U\subseteq W$, then $W\in\filterFont{u}$. 
\end{enumerate}
The pair $(X, \filterFont{u})$ is called a \emph{\index{space!uniform}uniform space}. The elements of $\filterFont{u}$ are called $\filterFont{u}$-\emph{neighborhoods}. 
\EndDefinition
The first three properties are gleaned from those of the pseudometric neighborhoods above, the last two are properties of a filter, which have been listed here just for completeness. 

We will omit $\filterFont{u}$ when talking about a uniform space, if this does not yield ambiguities. The term ``neighborhood'' is used for elements of a uniformity and for the neighborhoods\MMP[t]{Neigh\-borhood, entourage, Nachbarschaft} of a point. There should be no ambiguity, because the point is always attached, when talking about neighborhood in the latter, topological sense. Bourbaki uses the term \emph{\index{entourage}entourage} for a neighborhood in the uniform sense, the German word for this is \emph{\index{Nachbarschaft}Nachbarschaft} (while the term for a neighborhood of a point is \emph{\index{Umgebung}Umgebung}). 

We will need some relational identities; they are listed in Figure~\ref{tab:some-identities} for the reader's convenience. 
 
\begin{figure}[t]
  \centering
  \begin{align*}
    U\circ V & := \{\langle x, z\rangle\mid \exists y: \langle x, y\rangle\in U, \langle y, z\rangle\in V\}\\
U^{-1} & := \{\langle y, x\rangle\mid  \langle x, y\rangle\in U\}\\
U[M] & := \{y\mid \exists x\in M: \langle x, y\rangle\in U\}\\
\index{$U[x]$}U[x] & := U[\{x\}]\\
U\text{ is symmetric } &:\Leftrightarrow U^{-1} = U\\
(U\circ V)\circ W & = U\circ (V\circ W)\\
(U\circ V)^{-1} & = V^{-1}\circ U^{-1}\\
(U\circ V)[M] & = U[V[M]]\\
V\circ U\circ V & = \bigcup_{\langle x, y\rangle \in U}V[x]\times V[y]\text{ ($V$ symmetric)}
  \end{align*}
Here $U, V, W\subseteq X\times X$ and $M\subseteq X$. 
  \caption{\label{tab:some-identities}Some Relational Identities}
\end{figure}

As in the case of topologies, where we do not always specify the entire topology, but occasionally make use of the possibility to define it through a base, we will proceed similarly here, where we deal with filter bases. We have this characterization for the base of a uniformity.

\BeginProposition{char-base-unif}
A family $\emptyset\not=\filterFont{b}\subseteq \PowerSet{X\times X}$ is the base\MMP{Base} for a uniformity iff it has the following properties:
\begin{enumerate}
\item\label{char-base-unif-1} Each member of $\filterFont{b}$ contains the diagonal of $X$.
\item\label{char-base-unif-2} For $U\in\filterFont{b}$ there exists $V\in\filterFont{b}$ with $V\subseteq U^{-1}$.
\item\label{char-base-unif-3} For $U\in\filterFont{b}$ there exists $V\in\filterFont{b}$ with $V\circ V\subseteq U$.
\item\label{char-base-unif-4} For $U, V\in\filterFont{b}$ there exists $W\in\filterFont{b}$ with $W\subseteq U\cap V$. 
\end{enumerate}
\EndProposition

\BeginProof 
Recall that the filter generated by a filter base $\filterFont{b}$ is defined through $\{F\mid U\subseteq F\text{ for some }U\in\filterFont{b}\}$.  With this in mind, the proof is straightforward.  
\EndProof

This permits a description of a uniformity in terms of a base, which is usually easier than giving a uniformity as a whole. Let us look at some examples. 

\BeginExample{example-uniformities}
\begin{enumerate}
\item\label{example-uniformities-1} The uniformity $\{\Delta, X\times X\}$ is called the \emph{\index{uniformity!indiscrete}indiscrete uniformity}, the uniformity $\{A\subseteq X\times X\mid \Delta\subseteq A\}$ is called the \emph{\index{uniformity!discrete}discrete uniformity} on $X$. 
\item\label{example-uniformities-2} Let $V_{r} := \{\langle x, y\rangle\mid x, y\in \Real, |x-y|<r\}$, then $\{V_{r}\mid r>0\}$ is a base for a uniformity on $\Real$. Since it makes use of the structure of $(\Real, +)$ as an additive group, it is called the \emph{additive \index{uniformity!additive}uniformity} on $\Real$.
\item\label{example-uniformities-3} Put $V_{E} := \{\langle x, y\rangle\in\Real^{2}\mid x/y\in E\}$ for some neighborhood $E$ of $1\in\Real\setminus\{0\}$. Then the filter generated by $\{V_{E}\mid E \text{ is a neighborhood of 1}\}$ is a uniformity. This is so because the logarithm function is continuous on $\pReal\setminus\{0\}$. This uniformity nourishes itself from the multiplicative group $(\Real\setminus\{0\}, \cdot)$, so it is called the \emph{multiplicative \index{uniformity!multiplicative}uniformity} on $\Real\setminus\{0\}$. This is discussed in greater generality in part~\ref{example-uniformities-9}.
\item\label{example-uniformities-4} A \emph{\index{partition}partition} $\pi$ on a set $X$ is a collection of non-empty and mutually disjoint subsets of $X$ which covers $X$. It generates an equivalence relation on $X$ by rendering two elements of $X$ equivalent iff they are in the same partition element.  Define 
$
V_{\pi} := \bigcup_{i=1}^{n}(P_{i}\times P_{i})
$
for a finite partition $\pi = \{P_{1}, \dots, P_{k}\}$. Then $$\filterFont{b} := \{V_{\pi}\mid \pi\text{ is a finite partition on }X\}$$ is the base for a uniformity. Let $\pi$ be a finite partition, and denote the equivalence relation generated by $\pi$  by $|\pi|$, hence $\isEquiv{x}{y}{|\pi|}$ iff $x$ and $y$ are in the same element of $\pi$. 
\begin{itemize}
\item $\Delta\subseteq  V_{\pi}$ is obvious, since $|\pi|$ is reflexive.
\item $U^{-1} = U$ for all $U\in V_{\pi}$, since $|\pi|$ is symmetric.
\item Because $|\pi|$ is transitive, we have $V_{\pi}\circ V_{\pi}\subseteq V_{\pi}$.
\item Let $\pi'$ be another finite partition, then 
$
\{A\cap B\mid A\in \pi, B\in\pi', A\cap B\not=\emptyset\}
$
defines a partition $\pi''$ such that $V_{\pi''}\subseteq V_{\pi}\cap V_{\pi'}$. 
\end{itemize}
Thus $\filterFont{b}$ is the base for a uniformity, which is, you guessed it, called the \emph{\index{uniformity!finite partitions}uniformity of finite partitions}. 
\item\label{example-uniformities-5} Let $\emptyset\not={\cal I}\subseteq\PowerSet{X}$ be an ideal (\SetCite{Definition 1.5.32}), and define
\begin{align*}
{\cal A}_{E} & := \{\langle A, B\rangle\mid A\Delta B\in E\}\text{ for }E\in{\cal I},\\
\filterFont{b} & := \{{\cal A}_{E}\mid E\in{\cal I}\}.
\end{align*}
Then $\filterFont{b}$ is a base for a uniformity on $\PowerSet{X}$. In fact, it is clear that $\Delta_{\PowerSet{X}} \subseteq {\cal A}_{E}$ always holds, and that each member of $\filterFont{b}$ is symmetric. Let $A\Delta B\subseteq E$ and $B\Delta C\subseteq F$, then 
$ A\Delta C = (A\Delta B)\Delta (B\Delta C) \subseteq (A\Delta B)\cup (A\Delta C)\subseteq E\cup F$, thus ${\cal A}_{E}\circ {\cal A}_{F} \subseteq {\cal A}_{E\cup F}$, and finally ${\cal A}_{E}\cap{\cal A}_{F} \subseteq {\cal A}_{E\cap F}$. Because ${\cal I}$ is an ideal, it is closed under finite intersections and finite unions, the assertion follows.
\item\label{example-uniformities-7} Let $p$ be a prime, and put $W_{k} := \{\langle x, y\rangle\mid x, y\in\Ganz, p^{k}\text{ divides } x - y\}$. Then $W_{k}\circ W_{\ell} \subseteq W_{\min\{k, \ell\}} = W_{k}\cap W_{\ell}$, thus $\filterFont{b} := \{W_{k}\mid k\in\Nat\}$ is the base for a uniformity $\filterFont{u}_{p}$ on $\Ganz$, the \emph{$p$-adic \index{uniformity!$p$-adic}uniformity}.   
\item\label{example-uniformities-6} Let $A$ be a set, $(X, \filterFont{u})$ a uniform space, and let $F(A, X)$ be the set of all maps $A\to X$. We will define a uniformity on $F(A, X)$; the approach is similar to Example~\ref{ex-topol-bases-weak}. Define for $U\in\filterFont{u}$ the set
  \begin{equation*}
    U_{F} := \{\langle f, g\rangle\in F(A, X) \mid  \langle f(x), g(x)\rangle\in U\text{ for all }x\in X\}.
  \end{equation*}
Thus two maps are close with respect to $U_{F}$ iff all their images are close with respect to $U$. It is immediate that $\{U_{F}\mid U\in\filterFont{u}\}$ forms a uniformity, and that  $\{U_{F}\mid U\in\filterFont{b}\}$ is a base for a uniformity, provided $\filterFont{b}$ is a base for uniformity $\filterFont{u}$. 

If $X=\Real$ is endowed with the additive uniformity, a typical set of the base is given for $\epsilon>0$ through
\begin{equation*}
  \{\langle f, g\rangle\in F(A, \Real)\mid \sup_{a\in A}|f(a)-g(a)|<\epsilon\},
\end{equation*}
hence the images of $f$ and of $g$ have to be uniformly close to each other. 
\item\label{example-uniformities-8} Call a map $f: \Real\to \Real$ \emph{\index{map!affine}affine} iff it can be written as $f(x) = a\cdot x+b$ with $a\not=0$; let $f_{a, b}$ be the affine map characterized by the parameters $a$ and $b$, and define $X := \{f_{a, b}\mid a, b\in \Real, a\not=0\}$ the set of all affine maps. Note that an affine map is bijective, and that its inverse is an affine map again  with $f_{a, b}^{-1} = f_{1/a, -b/a}$; the composition of an affine map is an affine map as well, since $f_{a, b}\circ f_{c, d} = f_{ac, ad+b}$. Define for $\epsilon>0, \delta>0$ the $\epsilon, \delta$-neighborhood $U_{\epsilon, \delta}$ by 
  \begin{equation*}
    U_{\epsilon, \delta} := \{f_{a, b}\in X\mid |a-1|<\epsilon, |b|<\delta\}
  \end{equation*}
Put 
\begin{align*}
  U_{\epsilon, \delta}^{L} & := \{\langle f_{x, y}, f_{a, b}\rangle\in X\times X\mid f_{x, y}\circ f_{a, b}^{-1}\in U_{\epsilon, \delta}\},\\
\filterFont{b}_{L} & := \{U_{\epsilon, \delta}^{L}\mid  \epsilon>0, \delta>0\},\\
U_{\epsilon, \delta}^{R} & := \{\langle f_{x, y}, f_{a, b}\rangle\in X\times X \mid  f_{x, y}^{-1}\circ f_{a, b}\in U_{\epsilon, \delta}\},\\
\filterFont{b}_{R} & := \{U_{\epsilon, \delta}^{R}\mid  \epsilon>0, \delta>0\}.\\
\end{align*}
Then $\filterFont{b}_{L}$ resp. $\filterFont{b}_{R}$ is the base for a uniformity $\filterFont{u}_{L}$ resp. $\filterFont{u}_{R}$ on $X$. Let us check this for $\filterFont{b}_{R}$. Given positive $\epsilon, \delta$, we want to find positive $r, s$ with $\langle f_{m, n}, f_{p, q}\rangle\in V_{r, s}^{R}$ implies $\langle f_{p, q}, f_{m, n}\rangle\in U_{\epsilon, \delta}^{R}$. Now we can find for $\epsilon > 0$ and  $\delta>0$ some $r > 0$ and $s>0$ so that
\begin{align*}
  |\frac{p}{m} - 1| < r & \Rightarrow |\frac{m}{p} - 1| < \epsilon\\
|\frac{q}{m} - \frac{n}{m} | < s & \Rightarrow |\frac{n}{p} - \frac{q}{p}| < \delta
\end{align*}
holds, which is just what we want, since it translates into $V_{r, s}^{R}\subseteq \bigl(U_{\epsilon, \delta}^{R}\bigr)^{-1}$. The other properties of a base are easily seen to be satisfied. One argues similarly for $\filterFont{b}_{L}$. 

Note that $(X, \circ)$ is a topological group with the sets $\{U_{\epsilon, \delta}\mid \epsilon>0, \delta>0\}$ as a base for the neighborhood filter of the neutral element $f_{1, 0}$ (topological groups are introduced in Example~\ref{top-group} on page~\pageref{top-group}). 
\item\label{example-uniformities-9} Let, in general, $G$ be a topological group with neutral element $e$ . Define for $U\in\upsilon(e)$ the sets
  \begin{align*}
    U_{L} & := \{\langle x, y\rangle\mid  xy^{-1}\in U\},\\
U_{R} & := \{\langle x, y\rangle\mid  x^{-1}y\in U\},\\
U_{B} & := U_{L}\cap U_{R}.
   \end{align*}
Then $\{U_{L}\mid  U\in\upsilon(e)\}$, $\{U_{R}\mid  U\in\upsilon(e)\}$ and $\{U_{B}\mid  U\in\upsilon(e)\}$ define bases for uniformities on $G$; it can be shown that they do not necessarily coincide (of, course, they do, if $G$ is Abelian).  
\end{enumerate}
\EndExample

Before we show that a uniformity generates a topology, we derive a sufficient criterion for a family of subsets of $X\times X$ is a subbase for a uniformity.

\BeginLemma{subbase-for-uniformity}
Let $\filterFont{s}\subseteq\PowerSet{X\times X}$, then $\filterFont{s}$ is the subbase\MMP{Subbase} for a uniformity on $X$, provided the following conditions hold.
\begin{enumerate}
\item $\Delta\subseteq S$ for each $S\in\filterFont{s}$.
\item Given $U\in \filterFont{s}$, there exists $V\in\filterFont{s}$ such that $V\subseteq U^{-1}$.
\item For each $U\in\filterFont{s}$ there exists $V\in\filterFont{s}$ such that $V\circ V\subseteq U$.  
\end{enumerate}
\EndLemma

\BeginProof
We have to show that 
\begin{equation*}
  \filterFont{b} := \{U_{1}\cap\dots\cap U_{n}\mid U_{1}, \dots, U_{n}\in\filterFont{s}\text{ for some }n\in\Nat\}
\end{equation*}
constitutes a base for a uniformity. It is clear that every element of $\filterFont{b}$ contains the diagonal. Let $U = \bigcap_{i=1}^{n}U_{i}\in\filterFont{b}$ with $U_{i}\in\filterFont{s}$ for $i = 1, \dots, n$, choose $V_{i}\in\filterFont{s}$ with $V_{i}\subseteq U_{i}^{-1}$ for all $i$, then $V := \bigcap_{i=1}^{n}V_{i}\in\filterFont{b}$ and $V\subseteq U^{-1}$. If we select $W_{i}\in\filterFont{s}$ with $W_{i}\circ W_{i}\subseteq U_{i}$, then $W := \bigcap_{i=1}^{n}W_{i}\in\filterFont{b}$ and $W\circ W\subseteq U$. The last condition of Proposition~\ref{char-base-unif} is trivially satisfied for $\filterFont{b}$, since $\filterFont{b}$ is closed under finite intersections. Thus we conclude that $\filterFont{b}$ is a base for a uniformity on $X$ by Proposition~\ref{char-base-unif}, which in turn entails that $\filterFont{s}$ is a subbase. 
\EndProof

\Subsubsubsection{The Topology Generated by a Uniformity}
A pseudometric space $(X, d)$ generates a topology by declaring a set $G$ open iff there exists for $x\in G$ some $r>0$ with $B(x, r)\subseteq G$; from this we obtained the neighborhood filter $\upsilon(x)$ for a point $x$. Note that in the uniformity associated with the pseudometric the identity
\begin{equation*}
B(x, r) = V_{r}[x]
\end{equation*}
holds. Encouraged by this, we approach the topology for a uniform space in the same way. Given a uniform space $(X, \filterFont{u})$, a subset $G\subseteq X$ is called open iff we can find for each $x\in G$ some neighborhood $U\in\filterFont{u}$ such that $U[x]\subseteq G$. The following proposition investigates this construction.

\BeginProposition{topol-from-unif}
Given a uniform space $(X, \filterFont{u})$, \MMP[t]{From $\filterFont{u}$ to $\tau_{\filterFont{u}}$}for each $x\in X$ the family $\filterFont{u}[x] :=\{U[x]\mid U\in \filterFont{u}\}$ is the base for the neighborhood filter of $x$ for a topology $\tau_{\filterFont{u}}$, which is called the \emph{\index{topology!uniform}\index{space!uniform!topology}uniform topology}. The neighborhoods for $x$ in $\tau_{\filterFont{u}}$ are just $\filterFont{u}[x]$. 
\EndProposition

\BeginProof
It follows from Proposition~\ref{def-through-nbh-filters} that $\filterFont{u}[x]$ defines a topology $\tau_{\filterFont{u}}$, it remains to show that the neighborhoods of this topology are just $\filterFont{u}[x]$. We have to show that $U\in\filterFont{u}$ there exists $V\in\filterFont{u}$ with $V[x]\subseteq U[x]$
and $V[x]\in\filterFont{u}[y]$ for all $y\in V[x]$, then the assertion will follow from Corollary~\ref{cor-def-through-nbh-filters}. For $U\in\filterFont{u}$ there exists $V\in\filterFont{u}$ with $V\circ V\subseteq U$, thus $\langle x, y\rangle\in V$ and $\langle y, z\rangle\in V$ implies $\langle x, z\rangle\in U$. Now let $y\in V[x]$ and $z\in V[y]$, thus $z\in U[x]$, but this means $U[x]\in\filterFont{u}[y]$ for all $x\in V[y]$. Hence the assertion follows. 
\EndProof

These are some illustrative example. They indicate also that different uniformities can generate the same topology.  

\BeginExample{top-for-unif}
\begin{enumerate}
\item The topology obtained from the additive uniformity on $\Real$ is the usual topology. The same holds for the multiplicative uniformity on $\Real\setminus\{0\}$. Both can be shown to be distinct~\cite[Ch. 3, §6]{Bourbaki}.
\item The topology induced by the discrete uniformity is the discrete topology, in which each singleton $\{x\}$ is open. Since $\bigl\{\{x\}, X\setminus\{x\}\bigr\}$ forms a finite partition of $X$, the discrete topology is induced also by the uniformity defined by the finite partitions.
\item Let $F(A, \Real)$ be endowed with the uniformity defined by the sets $\{\langle f, g\rangle\in F(A, \Real)\mid \sup_{a\in A}|f(a)-g(a)|<\epsilon\}$, see Example~\ref{example-uniformities}. The corresponding topology yields for each $f\in F(A, \Real)$ the neighborhood $\{g\in F(A, \Real)\mid \sup_{a\in A}|f(a)-g(a)|< \epsilon\}$. This is the topology of uniform convergence. 
\item Let $\filterFont{u}_{p}$ for a prime $p$ be the $p$-adic uniformity on $\Ganz$, see Example~\ref{example-uniformities}, part~\ref{example-uniformities-7}. The corresponding topology $\tau_{p}$ is called the $p$-adic topology. A basis for the neighborhoods of $0$ is given by the sets $V_{k} := \{x\in \Ganz \mid p^{k}\text{ divides } x\}$. Because $p^{m}\in V_{k}$ for $m\geq k$, we see that $\lim_{n\to \infty}\ p^{n} = 0$ in $\tau_{p}$, but not in $\tau_{q}$ for $q\not= p$, $q$ prime. Thus the topologies $\tau_{p}$ and $\tau_{q}$ differ, hence also the uniformities $\filterFont{u}_{p}$ and $\filterFont{u}_{q}$.  
\end{enumerate}
\EndExample

Now that we know that each uniformity yields a topology on the same space, some questions are immediate:
\begin{itemize}
\item Do the open resp. the closed sets play a particular r\^ole in describing the uniformity?
\item Does the topology have particular properties, e.g., in terms of separation axioms?
\item What about metric spaces --- can we determine from the uniformity that the topology is metrizable?
\item Can we find a pseudometric for a given uniformity?
\item Is the product topology on $X\times X$ somehow related to $\filterFont{u}$, which is defined on $X\times X$, after all? 
\end{itemize}

We will give answers to some of these questions, some will be treated only lightly, with an in depth treatment to be found in the vast literature on uniform spaces, see the Bibliographic Notes in Section~\ref{sec:top-bib-notes}. 

Fix a uniform space $X$ with uniformity $\filterFont{u}$ and associated topology $\tau$. References to neighborhoods and open sets are always to $\filterFont{u}$ resp. $\tau$, unless otherwise stated. 

This is a first characterization of the interior of an arbitrary set. Recall that in a pseudometric space $x$ is an interior point of $A$ iff $B(x, r)\subseteq A$ for some $r>0$; the same description applies here as well, \emph{mutatis mutandis} (of course, this ``mutatis mutandis'' part is the interesting one).

\BeginLemma{unif-top-descr-interior}
Given $A\subseteq X$, $x\in\Interior{A}$ iff there exists a neighborhood $U$ with $U[x]\subseteq A$. 
\EndLemma

\BeginProof 
Assume that
$x\in\Interior{A} = \bigcup\{G\mid G\text{ open and }G\subseteq A\}$,
then it follows from the definition of an open set that we must be
able to find an neighborhood $U$ with $U[x]\subseteq A$.

Conversely, we show that the set
$ B := \{x\in X\mid U[x]\subseteq A\text{ for some neighborhood }U\} $
is open, then this must be the largest open set which is contained in
$A$, hence $B=\Interior{A}$. Let $x\in B$, thus $U[x]\subseteq A$, and
we should find now a neighborhood $V$ such that $V[y]\subseteq B$ for
$y\in V[x]$. But we find a neighborhood $V$ with
$V\circ V\subseteq U$. Let's see whether $V$ is suitable: if
$y\in V[x]$, then $V[y]\subseteq (V\circ V)[x]$ (this is so because
$\langle x, y\rangle\in V$, and if $z\in V[y]$, then
$\langle y, z\rangle\in V$; this implies
$\langle x, z\rangle\in V\circ V$, hence $z\in (V\circ V)[x]$). But
this implies $V[y]\subseteq U[x]\subseteq B$, hence $y\in B$. But this
means $V[x]\subseteq B$, so that $B$ is open. 
\EndProof

This gives us a handy way of describing the base for a neighborhood filter for a point in $X$. It states that we may restrict our attention to the members of a base or of a subbase, when we want to work with the neighborhood filter for a particular element.

\BeginCorollary{descr-nbhd-base}
If $\filterFont{u}$ has base or subbase $\filterFont{b}$, then $\{U[x]\mid U\in\filterFont{b}\}$ is a base resp. subbase for the neighborhood filter for $x$.  
\EndCorollary

\BeginProof
This follows immediately from Lemma~\ref{unif-top-descr-interior} together with Proposition~\ref{char-base-unif} resp. Lemma~\ref{subbase-for-uniformity}
\EndProof

Let us have a look at the topology on $X\times X$ induced by $\tau$. Since the open rectangles generate this topology, and since we can describe the open rectangles in terms of the sets $U[x]\times V[y]$, we can expect that these open sets can also related to the uniformity proper. In fact:

\BeginProposition{product-open-base}
If $U\in\filterFont{u}$, then both $\Interior{U}\in \filterFont{u}$ and $\Closure{U}\in\filterFont{u}$.
\EndProposition

\BeginProof
1.
Let $G\subseteq X\times X$ be open, then $\langle x, y\rangle\in G$ iff there exist neighborhoods $U, V\in\filterFont{u}$ with $U[x]\times V[y]\subseteq G$, and because $U\cap V\in\filterFont{u}$, we may even find some $W\in\filterFont{u}$ such that  $W[x]\times W[y]\subseteq G$. Thus
\begin{equation*}
  G = \bigcup\{W[x]\times W[y]\mid \langle x, y\rangle\in G, W\in\filterFont{u}\}.
\end{equation*}

2.
Let $W\in\filterFont{u}$, then there exists a symmetric $V\in\filterFont{u}$ with $V\circ V\circ V\subseteq W$, and by the identities in Figure~\ref{tab:some-identities} we may write 
\begin{equation*}
  V\circ V\circ V = \bigcup_{\langle x, y\rangle\in V}V[x]\times V[y].
\end{equation*}
Hence $\langle x, y\rangle\in \Interior{W}$ for every $\langle x, y\rangle\in V$, so $V\subseteq \Interior{W}$, and since $V\in\filterFont{u}$, we conclude $\Interior{W}\in\filterFont{u}$ from $\filterFont{u}$ being upward closed.

3.
Because $\filterFont{u}$ is a filter, and $U\subseteq \Closure{U}$, we infer $\Closure{U}\in\filterFont{u}$.
\EndProof

The closure of a subset of $X$ and the closure of a subset of $X\times X$ may be described as well directly through uniformity $\filterFont{u}$. These are truly remarkable representations.

\BeginProposition{closure-and-square}
$\Closure{A} = \bigcap\{U[A]\mid U\in\filterFont{u}\}$ for $A\subseteq X$, and $\Closure{M} = \bigcap\{U\circ M\circ U\mid U\in\filterFont{u}\}$  for $M\subseteq X\times X$.
\EndProposition

\BeginProof
1.
We use the characterization of a point $x$ in the closure through its neighborhood filter from \SetCite{Lemma 1.94}: $x\in\Closure{A}$ iff $U[x]\cap A\not=\emptyset$ for all symmetric $U\in\filterFont{u}$, because the symmetric neighborhoods form a base for $\filterFont{u}$. Now $z\in U[x]\cap A$ iff $z\in A$ and $\langle x, z\rangle\in U$ iff $z\in U[z]$ and $z\in A$, hence $U[x]\cap A \not=\emptyset$ iff $x\in U[A]$, because $U$ is symmetric. But this means $\Closure{A}= \bigcap\{U[A]\mid U\in\filterFont{u}\}$. 

2.
Let $\langle x, y\rangle\in\Closure{M}$, then $U[x]\times U[y]\cap M\not=\emptyset$ for all symmetric neighborhoods $U\in\filterFont{u}$, so that $\langle x, y\rangle\in U\circ M\circ U$ for all symmetric neighborhoods. This accounts for the inclusion from left to right. If $\langle x, y\rangle\in U\circ M\circ U$ for all neighborhoods $U$, then for every $U\in\filterFont{u}$ there exists $\langle a, b\rangle\in M$ with $\langle a, b\rangle\in U[x]\times (U^{-1})[y]$, thus $\langle x, y\rangle\in \Closure{M}$. 
\EndProof

Hence 

\BeginCorollary{clos-symm-are-base}
The closed symmetric neighborhoods form a base for the uniformity.
\EndCorollary

\BeginProof
Let $U\in\filterFont{u}$, then there exists a symmetric $V\in\filterFont{u}$ with $V\circ V\circ V\subseteq U$ with $V\subseteq\Closure{V}\subseteq V\circ V\circ V$ by Proposition~\ref{limits-are-unique}. Hence $W := \Closure{V}\cap(\Closure{V})^{-1})$ is a member of $\filterFont{u}$ which is contained in $U$.  
\EndProof

Proposition~\ref{closure-and-square} has also an interesting consequence when looking at the characterization of Hausdorff spaces in Proposition~\ref{limits-are-unique}. Putting $M = \Delta$, we obtain $\Closure{\Delta} = \bigcap\{U\circ U\mid U\in\filterFont{u}\}$, so that the associated topological space is Hausdorff iff the intersection of all neighborhoods is the diagonal  $\Delta$. Uniform spaces with $\bigcap\filterFont{u}=\Delta$ are called \emph{\index{space!uniform!separated}separated}\MMP{Separated}. 

\Subsubsubsection{Pseudometrization}
We will see shortly that the topology for a separated uniform space is completely regular. First, however, we will show that we can generate pseudometrics from the uniformity by the following idea: suppose that we have a neighborhood $V$, then there exists a neighborhood $V_{2}$ with $V_{2}\circ V_{2}\circ V_{2}\subseteq V_{1} := V$; continuing in this fashion, we find for the neighborhood $V_{n}$ a neighborhood $V_{n+1}$ with $V_{n+1}\circ V_{n+1}\circ V_{n+1}\subseteq V_{n}$, and finally put $V_{0} := X\times X$. Given a pair $\langle x, y\rangle\in X\times X$, this sequence $\Folge{V}$ is now used as a witness to determine how far apart these points are: put $f_{V}(x, y) := 2^{-n}$, iff $\langle x, y\rangle\in V_{n}\setminus V_{n-1}$, and $d_{V}(x, y) := 0$ iff $\langle x, y\rangle\in \bigcap_{n\in\Nat}V_{n}$. Then $f_{V}$ will give rise to a pseudometric $d_{V}$\MMP{$d_{V}$}, the pseudometric associated with $V$, as we will show below. 

This means that many pseudometric spaces are hidden deep inside a uniform space! Moreover, if we need a pseudometric, we construct one from a neighborhood. These observations will turn out to be fairly practical later on. But before we are in a position to make use of them, we have to do some work.

\BeginProposition{construct-pseudometric}
Assume that $\Folge{V}$ is a sequence of symmetric subsets of $X\times X$ with these properties for all $n\in\Nat$:
\begin{itemize}
\item $\Delta\subseteq V_{n}$,
\item $V_{n+1}\circ V_{n+1}\circ V_{n+1}\subseteq V_{n}$.
\end{itemize}
Put $V_{0} := X\times X$. Then there exists a pseudometric $d$ with
\begin{equation*}
V_{n}\subseteq \{\langle x, y\rangle\mid d(x, y)<2^{-n}\}\subseteq V_{n-1}
\end{equation*}
for all $n\in\Nat$.
\EndProposition

\BeginProof
0.
The proof uses the idea outlined above. The main effort will be showing that we can squeeze $\{\langle x, y\rangle\mid d(x, y)<2^{-n}\}$ between $V_{n}$ and $V_{n-1}$.

1.
Put $f(x, y) := 2^{-n}$ iff $\langle x, y\rangle\in V_{n}\setminus V_{n-1}$, and let $f(x, y) := 0$ iff $\langle x, y\rangle\in \bigcap_{n\in\Nat}V_{n}$. Then $f(x, x) = 0$, and $f(x, y) = f(y, x)$, because each $V_{n}$ is symmetric. Define 
\begin{equation*}
  d(x, y) := \inf\bigl\{\sum_{i=0}^{k}f(x_{i}, x_{i+1}) \mid  x_{0}, \dots, x_{k+1}\in X\text{ with }x_{0} = x, x_{k+1} = y, k\in\Nat\bigr\} 
\end{equation*}
So we look at all paths leading from $x$ to $y$, sum the weight of all their edges, and look at their smallest value. Since we may concatenate a path from $x$ to $y$ with a path from $y$ to $z$ to obtain one from $x$ to $z$, the triangle inequality holds for $d$, and since $d(x, y) \leq f(x, y)$, we know that $V_{n}\subseteq \{\langle x, y\rangle\mid d(x, y)<2^{-n}\}$. The latter set is contained in $V_{n-1}$; to show this is a bit tricky and requires an intermediary step.

2.
We show by induction on $n$ that 
\begin{equation*}
  f(x_{0}, x_{n+1}) \leq 2\cdot \sum_{i=0}^{n}f(x_{i}, x_{i+1}),
\end{equation*}
so if we have a path of length $n$, then the weight of the edge connecting their endpoints cannot be greater than twice the  weight on an arbitrary path. If $n = 1$, there is nothing to show. So assume the assertion is proved for all path with less that $n$ edges. We take a path from $x_{0}$ to $x_{n+1}$ with $n$ edges $\langle x_{i}, x_{i+1}\rangle$. Let $w$ be the weight of the path from $x_{0}$ to $x_{n+1}$, and let $k$ be the largest integer such that the path from $x_{0}$ to $x_{k}$ is at most $w/2$. Then the path from $x_{k+1}$ to $x_{n+1}$ has a weight at most $w/2$ as well. Now $f(x_{0}, x_{k}) \leq w$ and $f(x_{k+1}, x_{n+1})\leq w$ by induction hypothesis, and $f(x_{k}, x_{k+1})\leq w$. Let $m\in\Nat$ the smallest integer with $2^{-m}\leq w$, then we have $\langle x_{0}, x_{k}\rangle, \langle x_{k}, x_{k+1}\rangle, \langle x_{k+1}, x_{n+1}\rangle\in V_{m}$, thus $\langle x_{0}, x_{n+1}\rangle\in V_{m-1}$. This implies $f(x_{0}, x_{n+1})\leq 2^{-(m-1)}\leq 2\cdot w = 2\cdot \sum_{i=0}^{n}f(x_{i}, x_{i+1})$. 

3.
Now let $d(x, y) < 2^{-n}$, then $f(x, y) \leq 2^{-(n-1)}$ by part 2., and hence $\langle x, y\rangle\in V_{n-1}$. 
\EndProof

This has a ---~somewhat unexpected~--- consequence because it permits characterizing those uniformities, which are generated by a pseudometric.

\BeginProposition{unif-by-pseudometric}
The uniformity $\filterFont{u}$ of $X$ is generated by a pseudometric iff $\filterFont{u}$ has a countable base.
\EndProposition

\BeginProof
Let $\filterFont{u}$ be generated by a pseudometric $d$, then the sets $\{V_{d, r}\mid 0<r\in\Rational\}$ are a countable basis. Let, conversely, $\filterFont{b} := \{U_{n}\mid n\in\Nat\}$ be a countable base for $\filterFont{u}$. Put $V_{0} := X\times X$ and $V_{1} := U_{1}$, and construct inductively the sequence $\Folge{V}\subseteq\filterFont{b}$ of symmetric base elements with $V_{n}\circ V_{n}\circ V_{n}\subseteq V_{n-1}$ and $V_{n}\subseteq U_{n}$ for $n\in\Nat$. Then $\{V_{n}\mid n\in\Nat\}$ is a base for $\filterFont{u}$. In fact, given $U\in\filterFont{u}$, there exists $U_{n}\in\filterFont{b}$ with $U_{n}\subseteq U$, hence $V_{n}\subseteq U$ as well. Construct $d$ for this sequence as above, then we have $V_{n}\subseteq \{\langle x, y\rangle\mid d(x, y) < 2^{-n}\}\subseteq V_{n-1}$. Thus the sets $V_{d, r}$ are a base for the uniformity 
\EndProof

Note that this does not translate into an observation of the metrizability of the underlying topological space. This space may carry a metric, but the uniform space from which it is derived does not.

\BeginExample{partitions-not-metriz}
Let $X$ be an uncountable set, and let $\filterFont{u}$ be the uniformity given by the finite partitions, see Example~\ref{example-uniformities}. Then we have seen in Example~\ref{top-for-unif} that the topology induced by $\filterFont{u}$ on $X$ is the discrete topology, which is metrizable. 

Assume that $\filterFont{u}$ is generated by a pseudometric, then Proposition~\ref{unif-by-pseudometric} implies that $\filterFont{u}$ has a countable base, thus given a finite partition $\pi$, there exists a finite partition $\pi^{*}$ such that $V_{\pi^{*}}\subseteq V_{\pi}$, and $V_{\pi^{*}}$ is an element of this base. Here $V_{\{P_{1}, \dots, P_{n}\}} := \bigcup_{i=1}^{n}(P_{i}\times P_{i})$ is the basic neighborhood  for $\filterFont{u}$ associated with partition $\{P_{1}, \dots, P_{n}\}$. But for any given partition $\pi^{*}$ we can only form a finite number of other partitions $\pi$ with $V_{\pi^{*}}\subseteq V_{\pi}$, so that we have only a countable number of partitions on $X$. 
\EndExample
 
This is another consequence of Proposition~\ref{construct-pseudometric}: each uniform space satisfies the separation axiom $T_{3\half}$. For establishing this claim, we take a closed set $F\subseteq X$ and a point $x_{0}\not\in F$, then we have to produce a continuous function $f: X\to [0, 1]$ with $f(x_{0}) = 0$ and $f(y) = 1$ for $y\in A$. This is how to do it. Since $X\setminus F$ is open, we find a neighborhood $U\in\filterFont{u}$ with $U[x_{0}]\subseteq X\setminus F$. Let $d_{U}$ be the pseudometric associated with $U$, then $\{\langle x, y\rangle\mid d_{U}(x, y)< 1/2\}\subseteq U$. Clearly, $x\mapsto d_{U}(x, x_{0})$ is a continuous function on $X$, hence 
\begin{equation*}
  f(x) := \max\{0, 1-2\cdot d_{U}(x, x_{0})\}
\end{equation*}
is continuous with $f(x_{0}) = 1$ and $f(y) = 0$ for $y\in F$, thus $f$ has the required properties. Thus we have shown

\BeginProposition{tdreieinhalb-unif}
A uniform space is a  $T_{3\half}$-space; a separated uniform space is completely regular. 
\QED
\EndProposition

\Subsubsubsection{Cauchy Filters}

We generalize the notion of a Cauchy sequence to uniform spaces now. We do this in order to obtain a notion of convergence which includes convergence in topological spaces, and which carries the salient features of a Cauchy sequence with it. 

First, we note that filters are a generalization for sequences. So let us have a look at what can be said, when we construct the filter $\fiF$ for a Cauchy sequence $\Folge{x}$ in a pseudometric space $(X, d)$. $\fiF$ has the sets $\filterFont{c} := \{B_{n}\mid n\in\Nat\}$ with $B_n := \{x_{m}\mid m\geq n\}$ as a base. Being a Cauchy filter says that for each $\epsilon>0$ there exists $n\in\Nat$ such that $B_{n}\times B_{n}\subseteq V_{d, \epsilon}$; this inclusion holds then for all $B_{m}$ with $m\geq n$ as well. Because $\filterFont{c}$ is the base for $\fiF$, and the sets $V_{d, r}$ are a base for the uniformity, we may reformulate that $\fiF$ is a Cauchy filter iff for each neighborhood $U$ there exists $B\in \fiF$ such that $B\times B\subseteq U$. Now this looks like a property which may be formulated for general uniform spaces.

Fix the uniform space $(X, \filterFont{u})$. Given $U\in\filterFont{u}$, the set $M\subseteq X$ is called \emph{\index{set!small}\MMP{Small sets}$U$-small}  iff $M\times M\subseteq U$. A collection ${\cal F}$ of sets is said to contain \emph{small sets} iff given $U\in\filterFont{u}$ there exists $A\in{\cal F}$ which is $U$-small, or, equivalently, given $U\in\filterFont{u}$ there exists $x\in X$ with $A\subseteq U[x]$. 

This helps in formulating the notion of a Cauchy filter.

\BeginDefinition{def-cauchy-filter}
A filter $\fiF$ is called a \emph{\index{Cauchy filter}Cauchy \index{filter!Cauchy}filter} iff it contains small sets.
\EndDefinition

In this sense, a Cauchy sequence induces a Cauchy filter. Convergent filters are Cauchy filters as well:

\BeginLemma{convergent-is-cauchy}
If $\fiF\to x$ for some $x\in X$, then $\fiF$ is a Cauchy filter.
\EndLemma

\BeginProof
Let $U\in\filterFont{u}$, then there exists a symmetric $V\in\filterFont{u}$ with $V\circ V\subseteq U$. Because $\upsilon(x)\subseteq \fiF$, we conclude $V[x]\in\fiF$, and $V[x]\times V[x]\subseteq U$, thus $V[x]$ is a $U$-small member of $\fiF$. 
\EndProof

But the converse does not hold, as the following example shows. 

\BeginExample{ultra-is-cauchy}
Let $\filterFont{u}$ be the uniformity induced by the finite partitions with $X$ infinite. We claim that each ultrafilter $\fiF$ is a Cauchy filter. In fact, let $\pi = \{A_{1}, \dots, A_{n}\}$ be a finite partition, then $V_{\pi} = \bigcup_{i=1}^{n} A_{i}\times A_{i}$ is the corresponding neighborhood, then there exists $i^{*}$ with $A_{i^{*}}\in\fiF$. This is so since if an ultrafilter contains the finite union of  sets, it must contain one of them, see~\SetCite{Lemma 1.5.36}. $A_{i^{*}}$ is $V$-small.

The topology induced by this uniformity is the discrete topology, see Example~\ref{top-for-unif}. This topology is not compact, since $X$ is infinite. By Theorem~\ref{conv-vs-ultrafilter} there are ultrafilters which do not converge. 
\EndExample

If $x$ is an accumulation point of a Cauchy sequence in a pseudometric space, then we know that $x_{n}\to x$; this is fairly easy to show. A similar observation can be made for Cauchy filters, so that we have a partial converse to Lemma~\ref{convergent-is-cauchy}.

\BeginLemma{acc-point-cf-converges}
Let $x$ be an accumulation point of the Cauchy filter $\fiF$, then $\fiF\to x$. 
\EndLemma
 
\BeginProof
Let $V\in\filterFont{u}$ be a closed neighborhood; in view of Corollary~\ref{clos-symm-are-base} is is sufficient to show that $V[x]\in \fiF$, then it will follow that $\upsilon(x)\subseteq\fiF$. Because $\fiF$ is a Cauchy filter, we find $F\in\fiF$ with $F\times F\subseteq V$, because $V$ is closed, we may assume that $F$ is closed as well (otherwise we replace it by its closure). Because $F$ is closed and $x$ is an accumulation point of $\fiF$, we know from Lemma~\ref{all-acc-points} that $x\in F$, hence $F\subseteq V[x]$. This implies $\upsilon(x)\subseteq \fiF$. 
\EndProof

\BeginDefinition{complete-unif-space}
The uniform space $(X, \filterFont{u})$ is called \emph{\index{space!uniform!complete}complete} iff each Cauchy filter converges. 
\EndDefinition

Each Cauchy sequence converges in a complete uniform space, because the associated filter is a Cauchy filter. 

A slight reformulation is given in the following proposition, which is the uniform counterpart to the characterization of complete pseudometric spaces in Proposition~\ref{diam-to-zero-compl}. Recall that a collection of sets is said to have the \emph{finite intersection property} iff each finite subfamily has a non-empty intersection.  

\BeginProposition{equiv-uni-complete}
The uniform space $(X, \filterFont{u})$ is complete iff each family of closed sets which has the finite intersection property and which contains small sets has a non-void intersection. 
\EndProposition

\BeginProof
This is essentially a reformulation of the definition, but let's see.

1.
Assume that $(X, \filterFont{u})$ is complete, and let ${\cal A}$ be a family of closed sets with the finite intersection property, which contains small sets. Hence $\fiF_{0} := \{F_{1}\cap \dots\cap F_{n}\mid n\in\Nat, F_{1}, \dots, F_{n}\in {\cal A}\}$ is a filter base. Let $\fiF$ be the corresponding filter, then $\fiF$ is a Cauchy filter, for ${\cal A}$, hence $\fiF_{0}$ contains small sets. Thus $\fiF\to x$, so that $\upsilon(x)\subseteq \fiF$, thus $x\in\bigcap_{F\in\fiF}\Closure{F}\subseteq \bigcap_{A\in{\cal A}}A$. 

2.
Conversely, let $\fiF$ be a Cauchy filter. Since $\{\Closure{F}\mid F\in\fiF\}$ is a family of closed sets with the finite intersection property which contains small sets, the assumption says that $\bigcap_{F\in\fiF}\Closure{F}$ is not empty and contains some $x$. But then $x$ is an accumulation point of $\fiF$ by Lemma~\ref{all-acc-points}, so $\fiF\to x$ by Lemma~\ref{acc-point-cf-converges}. 
\EndProof

As in the case of pseudometric spaces, compact spaces  are derived from a complete uniformity. 

\BeginLemma{compactyields-complete}
Let $(X, \filterFont{u})$ be a uniform space so that the topology associated with the uniformity is compact. Then the uniform space $(X, \filterFont{u})$ is complete. 
\EndLemma

\BeginProof
In fact, let $\fiF$ be a Cauchy filter on $X$. Since the topology for $X$ is compact, the filter has an accumulation point $x$ by Corollary~\ref{char-acc-point-ultra}. But Lemma~\ref{acc-point-cf-converges} tells us then that $\fiF\to x$. Hence each Cauchy filter converges. 
\EndProof  

The uniform space which is derived from an ideal on the powerset of a set, which has been defined in Example~\ref{example-uniformities} (part~\ref{example-uniformities-5}) is complete. We establish this first for Cauchy nets as the natural generalization of Cauchy sequences, and then translate the proof to Cauchy filters. This will permit an instructive comparison of handling these two concepts.

\BeginExample{example-uniformities-5-net}
Recall the definition of a net on page~\pageref{def-net}. A net $(x_{i})_{i\in N}$ in the uniform space $X$ is called a \emph{Cauchy \index{net!Cauchy}net}\MMP{Cauchy net} iff, given a neighborhood $U\in\filterFont{u}$, there exists $i\in N$ such that $\langle x_{j}, x_{\gamma}\rangle\in U$ for all $j, k\in N$ with $j, k\geq i$. The net converges to $x$ iff given a neighborhood $U$ there exists $i\in N$ such that $\langle x_{j}, x\rangle\in U$ for $j\geq i$. 

Now assume that ${\cal I}\subseteq \PowerSet{X}$ is an ideal; part~\ref{example-uniformities-5} of Example~\ref{example-uniformities} defines a uniformity $\filterFont{u}_{{\cal I}}$ on $\PowerSet{X}$ which has the sets 
$
V_{I} := \{\langle A, B\rangle\mid A, B\in\PowerSet{X}, A\Delta B\subseteq I\}
$
as a base, as $I$ runs through ${\cal I}$. We claim that each Cauchy net $(F_{i})_{i\in N}$ converges to $F := \bigcup_{i\in N}\bigcap_{j\geq i}F_{j}$. 

In fact, let a neighborhood $U$ be given; we may assume that $U = V_{I}$ for some ideal $I\in{\cal I}$. Thus there exists $i\in N$ such that $\langle F_{j}, F_{k}\rangle\in V_{I}$ for all $j, k\geq i$, hence $F_{j}\Delta F_{k}\subseteq I$ for all these $j, k$. Let $x\in F\Delta F_{j}$ for $j\geq i$.
\begin{itemize}
\item If $x\in F$, we find $i_{0}\in N$ such that $x\in F_{k}$ for all $k\geq i_{0}$. Fix $k\in N$ so that $k\geq i$ and $k\geq i_{0}$, which is possible since $N$ is directed. Then $x\in F_{k}\Delta F_{j}\subseteq I.$
\item If $x\not\in F$, we find for each $i_{0}\in N$ some $k\geq i_{0}$ with $x\not\in F_{k}$. Pick $k\geq i_{0}$, then $x\not\in F_{k}$, hence $x\in F_{\gamma}\Delta F_{j}\subseteq I$
\end{itemize}
Thus $\langle F, F_{j}\rangle\in V_{I}$ for $j\geq i$, hence the net converges to $F$. 
\EndExample

Now let's investigate convergence of a Cauchy filter. One obvious obstacle in a direct translation seems to be the definition of the limit set, because this appears to be bound to the net's indices. But look at this. If $(x_{i})_{i\in N}$ is a net, then the sets ${\cal B}_{i} := \{x_{j}\mid j\geq i\}$ form a filter base $\filterFont{B}$, as $i$ runs through the directed set $N$ (see the discussion on page~\pageref{def-net}). Thus we have defined $F$ in terms of this base, viz., $F = \bigcup_{{\cal B}\in\filterFont{B}}\bigcap{\cal B}$. This gives an idea for the filter based case.

\BeginExample{example-uniformities-5-filter}
Let $\filterFont{u}_{{\cal I}}$ be the uniformity on $\PowerSet{X}$ discussed in Example~\ref{example-uniformities-5-net}. Then each Cauchy filter $\fiF$ converges. In fact, let $\filterFont{B}$ be a base for $\fiF$, then $\fiF\to F$ with $F := \bigcup_{{\cal B}\in\filterFont{B}}\bigcap{\cal B}$. 

Let $U$ be a neighborhood in $\filterFont{u}_{{\cal I}}$, and we may assume that $U = V_{I}$ for some $I\in{\cal I}$. Since $\fiF$ is a Cauchy filter, we find ${\cal F}\in\fiF$ which is $V_{I}$-small, hence $F\Delta F'\subseteq I$ for all $F, F'\in{\cal F}$. Let $F_{0}\in{\cal F}$, and consider $x\in F\Delta F_{0}$; we show that $x\in I$ by distinguishing these cases:
\begin{itemize}
\item If $x\in F$, then there exists ${\cal B}\in\filterFont{B}$ such that $x\in\bigcap{\cal B}$. Because ${\cal B}$ is an element of base $\filterFont{B}$, and because $\fiF$ is a filter, ${\cal B}\cap {\cal F}\not=\emptyset$, so we find $G\in{\cal B}$ with $G\in{\cal F}$, in particular $x\in G$. Consequently $x\in G\Delta F_{0}\subseteq I$, since ${\cal F}$ is $V_{I}$-small.
\item If $x\not\in F$, we find for each ${\cal B}\in\filterFont{B}$ some $G\in{\cal B}$ with $x\not\in G$. Since $\filterFont{B}$ is a base for $\fiF$, there exists ${\cal B}\in\filterFont{B}$ with ${\cal B}\subseteq{\cal F}$, so there exists $G\in{\cal F}$ with $x\not\in G$. Hence $x\in G\Delta F_{0}\subseteq I$. 
\end{itemize}
Thus $F\Delta F_{0}\subseteq I$, hence $\langle F, F_{0}\rangle\in V_{I}$. This means ${\cal F}\subseteq V_{I}[F]$, which in turn implies $\upsilon(F)\subseteq \fiF$, or, equivalently, $\fiF\to F$. 
\EndExample

For further investigations of uniform spaces, we define uniform continuity as the brand of continuity which is adapted to uniform spaces.

\Subsubsubsection{Uniform Continuity}

Let $f: X\to X'$ be a uniformly continuous map between the pseudometric spaces $(X, d)$ and $(X', d')$. This means that given $\epsilon>0$ there exists $\delta>0$ such that, whenever $d(x, y)<\delta$, $d'(f(x), f(y))<\epsilon$ follows. In terms of neighborhoods, this means $V_{d, \delta}\subseteq\InvBild{(\zZ{f})}{V_{d', \epsilon}}$, or, equivalently, that $\InvBild{(\zZ{f})}{V}$ is a neighborhood in $X$, whenever $V$ is a neighborhood in $X'$. We use this formulation, which is based only on neighborhoods, and not on pseudometrics, for a formulation of uniform continuity.

\BeginDefinition{unif-cont-unif-sp}
Let $(X, \filterFont{u})$ and $(Y, \filterFont{v})$ be uniform spaces. Then $f: X\to Y$ is called \emph{uniformly \index{continuous!uniformly}continuous} iff $\InvBild{(\zZ{f})}{V}\in\filterFont{u}$ for all $V\in\filterFont{v}$. 
\EndDefinition

\BeginProposition{unif-spac-cat}
Uniform spaces for a category with uniform continuous maps as morphisms.
\EndProposition

\BeginProof
The identity is uniformly continuous, and, since $(\zZ{g})\circ (\zZ{f}) = \zZ{(g\circ f)}$, the composition of uniformly continuous maps is uniformly continuous again.
\EndProof

Introducing something new, one checks whether this has some categorical significance, of course. We also want to see what happens in the underlying topological space. But here nothing unexpected will happen: a uniformly continuous map is continuous with respect to the underlying topologies, formally:

\BeginProposition{unif-cont-is-cont}
If $f: (X, \filterFont{u})\to (Y, \filterFont{v})$ is uniformly continuous, then $f: (X, \tau_{\filterFont{u}})\to (Y, \tau_{\filterFont{v}})$ is continuous. 
\EndProposition

\BeginProof
Let $H\subseteq Y$ be open with $f(x)\in H$. If $x\in \InvBild{f}{H}$, there exists a neighborhood $V\in\filterFont{v}$ such that $V[f(x)]\subseteq H$. Since  $U := \InvBild{(\zZ{f})}{V}$ is a neighborhood in $X$, and $U[x]\subseteq \InvBild{f}{H}$,  it follows that $\InvBild{f}{H}$ is open in $X$. 
\EndProof

The converse is not true, however, as Example~\ref{cont-not-unif-cont} shows. 

Before proceeding, we briefly discuss two uniformities on the same topological group which display quite different behavior, so that the identity is not uniformly continuous.
 
\BeginExample{ex-affine-maps}
Let $X := \{f_{a, b}\mid a, b\in \Real, a\not=0\}$ be the set of all affine maps $f_{a, b}: \Real\to \Real$ with the separated uniformities $\filterFont{u}_{R}$ and $\filterFont{u}_{L}$, as discussed in Example~\ref{example-uniformities},  part~\ref{example-uniformities-8}. 

Let $a_{n} := d_{n} := 1/n$, $b_{n} := -1/n$ and $c_{n} := n$. Put $g_{n} := f_{a_{n}, b_{n}}$ and $h_{n} :=  f_{c_{n}, d_{n}}$, $j_{n} := h_{n}^{-1}$. Now $g_{n}\circ h_{n} = f_{1, 1/n^{2}-1/n}\to f_{1, 0}$, $h_{n}\circ g_{n} = f_{1, -1+1/n}\to f_{1, -1}$. Now assume that $\filterFont{u}_{R} = \filterFont{u}_{L}$. Given $U\in \upsilon(e)$, there exists $V\in\upsilon(e)$ symmetric such that $V^{R}\subseteq U^{L}$. Since $g_{n}\circ h_{n}\to f_{1, 0}$, there exists for $V$ some $n_{0}$ such that $g_{n}\circ h_{n}\in V$ for $n\geq n_{0}$, hence $\langle g_{n}, j_{n}\rangle\in V^{R}$, thus $\langle j_{n}, g_{n}\rangle\in V^{R}\subseteq U^{L}$, which means that $h_{n}\circ g_{n}\in U$ for $n\geq n_{0}$. Since $U\in\upsilon(e)$ is arbitrary, this means that $h_{n}\circ g_{n}\to e$, which is a contradiction.

Thus we find that the left and the right uniformity on a topological group are different, although they are derived from the same topology.  In particular, the identity $(X, \filterFont{u}_{R})\to (X, \filterFont{u}_{L})$ is not uniformly continuous. 
\EndExample

We will construct the initial uniformity for a family of maps now. The approach is similar to the one observed for the initial topology (see Definition~\ref{initial-and-final-tops}), but since a uniformity is in particular a filter with certain properties, we have to make sure that the construction can be carried out as intended. Let ${\cal F}$ be a family of functions $f: X\to Y_{f}$, where $(Y_{f}, \filterFont{v}_{f})$ is a uniform space. We want to construct a uniformity $\filterFont{u}$  on $X$ rendering all $f$ uniformly continuous, so $\filterFont{u}$ should contain 
\begin{equation*}
  \filterFont{s} := \bigcup_{f\in {\cal F}}\{\InvBild{(\zZ{f})}{V}\mid V\in\filterFont{v}_{f}\},
\end{equation*}
and it should be the smallest uniformity on $X$ with this property. For this to work, it is necessary for $\filterFont{s}$ to be a subbase. We check this along the properties from Lemma~\ref{subbase-for-uniformity}: 
\begin{enumerate}
\item Let $f\in {\cal F}$  and $V\in \filterFont{v}_{f}$, then $\Delta_{Y_{f}}\subseteq V$. Since $\Delta_{X} = \InvBild{f}{\Delta_{Y_{f}}}$, we conclude $\Delta_{X} \subseteq \{\InvBild{(\zZ{f})}{V}$. Thus each element of $\filterFont{s}$ contains the diagonal of $X$.
\item Because $\bigl(\InvBild{(\zZ{f})}{V}\bigl)^{-1} = \InvBild{(\zZ{f})}{V^{-1}}$, we find that, given $U\in\filterFont{s}$, there exists $V\in\filterFont{s}$ with $V\subseteq U^{-1}$.
\item Let $U\in\filterFont{b}$, so that $U = \InvBild{(\zZ{f})}{V}$ for some $f\in{\cal F}$ and $V\in\filterFont{v}_{f}$. We find $W\in\filterFont{v}_{f}$ with $W\circ W\subseteq V$; put $W_{0} := \InvBild{(\zZ{f})}{W}$, then $W_{0}\circ W_{0}\subseteq \InvBild{(\zZ{f})}{W\circ W}\subseteq\InvBild{(\zZ{f})}{V} = U$, so that we find for $U\in\filterFont{s}$ an element $W_{0}\in\filterFont{s}$ with $W_{0}\circ W_{0}\subseteq U$.  
\end{enumerate}
Thus $\filterFont{s}$ is the subbase for a uniformity, and we have established

\BeginProposition{initial-unif}
Let ${\cal F}$ be a family of maps $X\to Y_{f}$ with $(Y_{f}, \filterFont{v}_{f})$ a uniform space, then there exists a smallest uniformity $\filterFont{u}_{{\cal F}}$ on $X$ rendering all $f\in {\cal F}$ uniformly continuous. $\filterFont{u}_{{\cal F}}$ is called the \emph{initial \index{uniformity!initial}uniformity} on $X$ with respect to ${\cal F}$. 
\EndProposition

\BeginProof
We know that $\filterFont{s} := \bigcup_{f\in {\cal F}}\{\InvBild{(\zZ{f})}{V}\mid V\in\filterFont{v}_{f}\}$ is a subbase for a uniformity $\filterFont{u}$, which is evidently the smallest uniformity so that each $f\in {\cal F}$ is uniformly continuous. So $\filterFont{u}_{f} := \filterFont{u}$ is the uniformity we are looking for.
\EndProof

Having this tool at our disposal, we can now ---~in the same way as we did with topologies~--- define
\begin{description}
\item[Product] The product \index{uniformity!product}uniformity\MMP{Product} for the uniform spaces $(X_{i}, \filterFont{u}_{i})_{i\in I}$ is the initial uniformity on $X := \prod_{i\in I}X_{i}$ with respect to the projections $\pi_{i}: X\to X_{i}$.
\item[Subspace] The subspace uniformity\index{uniformity!subspace}\MMP{Subspace} $\filterFont{u}_{A}$ is the initial uniformity on $A\subseteq X$ with respect to the embedding $i_{A}: x \mapsto x$.
\end{description}

We can construct dually a final uniformity on $Y$ with respect to a family ${\cal F}$ of maps $f: X_{f}\to Y$ with uniform spaces $(X_{f}, \filterFont{u}_{f})$, for example when investigating quotients. The reader is referred to~\cite[II.2]{Bourbaki} or to~\cite[8.2]{Engelking}. 

This is a little finger exercise for the use of a product uniformity. It takes a pseudometric and shows what you would expect: the pseudometric is uniformly continuous iff it generates neighborhoods. The converse holds as well. We do not assume here that $d$ generates $\filterFont{u}$, rather, it is just an arbitrary pseudometric, of which there may be many.

\BeginProposition{d-unif-cont-iff}
Let $(X, \filterFont{u})$ be a uniform space, $d: X\times X\to\pReal$ a pseudometric. Then $d$ is uniformly continuous with respect to the product uniformity on $X\times X$ iff $V_{d, r}\in\filterFont{u}$ for all $r>0$. 
\EndProposition

\BeginProof
1.
Assume first that $d$ is uniformly continuous, thus we find for each $r>0$ some neighborhood $W$ on $X\times X$ such that $\bigl\langle\langle x, u\rangle, \langle y, v\rangle\bigr\rangle\in W$ implies $|d(x, y) - d(u, v)|<r$. We find a symmetric neighborhood $U$ on $X$ such that $U_{1}\cap U_{2} \subseteq W$, where $U_{i} := \InvBild{(\zZ{\pi_{i}})}{U}$ for $i = 1, 2$, and 
\begin{align*}
  U_{1} & = \{\bigl\langle\langle x, u\rangle, \langle y, v\rangle\bigr\rangle\mid \langle x, y\rangle\in U\},\\
 U_{2} & = \{\bigl\langle\langle x, u\rangle, \langle y, v\rangle\bigr\rangle\mid \langle u, v\rangle\in U\}. 
\end{align*}
Thus if $\langle x, y\rangle\in U$, we have $\bigl\langle\langle x, y\rangle, \langle y, y\rangle\bigr\rangle \in W$, hence $d(x, y)< r$, so that $U\subseteq V_{d, r}$, thus $V_{d, r}\in\filterFont{u}$.

2.
Assume that $V_{d, r}\in \filterFont{u}$ for all $r>0$, and we want to show that $d$ is uniformly continuous in the product. If $\langle x, u\rangle, \langle y, v\rangle\in V_{d, r}$, then 
\begin{align*}
  d(x, y) & \leq d(x, u) + d(u, v) + d(v, y)\\
d(u, v) & \leq d(x, u) + d(x, y) + d(y, v),
\end{align*}
hence $|d(x, y) - d(u, v)| < 2\cdot r$. Thus $\InvBild{(\zZ{\pi_{1}})}{V_{d, r}}\cap\InvBild{(\zZ{\pi_{2}})}{V_{d, r}}$ is a neighborhood on $X\times X$ such that  $\bigl\langle\langle x, u\rangle, \langle y, v\rangle\bigr\rangle\in W$ implies $|d(x, y) - d(u, v)|<2\cdot r$.
\EndProof

Combining Proposition~\ref{construct-pseudometric} with the observation from Proposition~\ref{d-unif-cont-iff}, we have established this characterization of a uniformity through pseudometrics.

\BeginProposition{char-thru-pseudometrics}
The uniformity $\filterFont{u}$ is the smallest uniformity which is generated by all pseudometrics which are uniformly continuous on $X\times X$, i.e., $\filterFont{u}$ is the smallest uniformity containing $V_{d, r}$ for all such $d$ and all $r>0$. 
\EndProposition
  
We fix for the rest of this section the uniform spaces $(X, \filterFont{u})$ and $(Y,  \filterFont{v})$.
Note that for checking uniform continuity it is enough to look at a subbase. The proof is straightforward and hence omitted. 

\BeginLemma{unif-subbase-is-enough}
Let $f: X\to Y$ be a map. Then $f$ is uniformly continuous iff $\InvBild{(\zZ{f})}{V}\in\filterFont{u}$ for all elements of a subbase for $\filterFont{v}$. \QED
\EndLemma

Cauchy filters are preserved through uniformly continous maps (the image of a filter is defined on page~\pageref{direct-image-filter}). 

\BeginProposition{preserve-cauchy-filters}
Let $f: X\to Y$ be uniformly continuous and $\fiF$ a Cauchy filter on $X$. Then $f(\fiF)$ is a Cauchy filter.  
\EndProposition

\BeginProof
Let $V\in\filterFont{v}$ be a neighborhood in $Y$, then $U := \InvBild{(\zZ{f})}{V}$ is a neighborhood in $X$, so that there exists $F\in\fiF$ which is $U$-small, hence $F\times F\subseteq U$, hence $\Bild{(\zZ{f})}{F\times F} = \Bild{f}{F}\times\Bild{f}{F} \subseteq V$. Since $\Bild{f}{F}\in f(\fiF)$ by Lemma~\ref{direct-image-filter}, the image filter contains a $V$-small member. 
\EndProof

A first consequence of Proposition~\ref{preserve-cauchy-filters} is shows that the subspaces induced by closed sets in a complete uniform space are complete again.

\BeginProposition{closed-in-compl-unif}
If $X$ is separated, then a complete subspace is closed. Let $A\subseteq X$ be closed and $X$ be complete, then the subspace $A$ is complete.  
\EndProposition

Note that the first part does not assume that $X$ is complete, and that the second part does not assume that $X$ is separated. 

\BeginProof
1.
Assume that $X$ is a Hausdorff space and $A$ a complete subspace of $X$. We show $\partial A \subseteq A$, from which it will follow that $A$ is closed. Let $b\in \partial A$, then $U\cap A\not=\emptyset$ for all open neighborhoods $U$ of $b$. The trace $\upsilon(b)\cap A$ of the neighborhood filter $\upsilon(b)$ on $A$ is a Cauchy filter. In fact, if $W\in\filterFont{u}$ is a neighborhood for $X$, which we may choose as symmetric, then $\bigl((W[b]\cap A)\times (W[b]\cap A)\bigr) \subseteq W\cap (A\times A)$, which means that $W[b]\cap A)$ is $W\cap (A\times A)$- small. Thus $\upsilon(b)\cap A$ is a Cauchy filter on $A$, hence it converges to, say, $c\in A$. Thus $\upsilon(c)\cap A\subseteq \upsilon(b)\in A$, which means that $b = c$, since $X$, and hence $A$, is Hausdorff as a topological space. Thus $b\in A$, and $A$ is closed by Proposition~\ref{char-top-closure}.

2.
Now assume that $A\subseteq X$ is closed, and that $X$ is complete. Let $\fiF$ be a Cauchy filter on $A$, then $i_{A}(\fiF)$ is a Cauchy filter on $X$ by Proposition~\ref{preserve-cauchy-filters}. Thus $i_{A}(\fiF)\to x$ for some $x\in X$, and since $A$ is closed, $x\in A$ follows.  
\EndProof

We show that a uniformly continuous map on a dense subset  into a complete and separated uniform space can be extended uniquely to a uniformly continuous map on the whole space. This was established in Proposition~\ref{extend-unif-cont-maps} for pseudometric spaces; having a look at the proof displays the heavy use of pseudometric machinery such as the oscillation, and the pseudometric itself. This is not available in the present situation, so we have to restrict ourselves to the tools at our disposal, viz., neighborhoods and filters, in particular Cauchy filters for a complete space. We follow Kelley's elegant proof~\cite[p. 195]{Kelley}.

\BeginTheorem{ext-dense-unif-compl}
Let $A\subseteq X$ be a dense subsets of the uniform space $(X, \filterFont{u})$, and $(Y; \filterFont{v})$ be a complete and separated uniform space. Then a uniformly continuous map $f: A\to Y$ can be extended uniquely to a uniformly continous $F: X\to Y$. 
\EndTheorem

\BeginProof
0.
The proof starts from the graph $\{\langle a, f(a)\rangle\mid  a\in A\}$ of $f$ and investigates the properties of its closure in $X\times Y$\MMP{Plan of the proof}. It is shown that the closure is a relation which has $\Closure{A} = X$ as its domain, and which is the graph of a map, since the topology of $Y$ is Hausdorff. This map is an extension $F$ to $f$, and it is shown that $F$ is uniformly continuous. We also use the observation that the image of a converging filter under a uniform continuous map is a Cauchy filter, so that completeness of $Y$ kicks in when needed. We do not have to separately establish uniqueness, because this follows directly from Lemma~\ref{equal-on-dense}. 

1.  
Let $G_{f} := \{\langle a, f(a)\rangle\mid a\in A\}$ be the graph of $f$. We claim that the closure of the domain of $f$ is the domain of the closure of $G_{f}$. Let $x$ be in the domain of the closure of $G_{f}$, then there exists $y\in Y$ with $\langle x, y\rangle\in \Closure{G}_{f}$, thus we find a filter $\fiF$ on $G_{f}$ with $\fiF\to \langle x, y\rangle$. Thus $\pi_{1}(\fiF)\to x$, so that $x$ is in the closure of the domain of $f$. Conversely, if $x$ is in the closure of the domain of $G_{f}$, we find a filter $\fiF$ on the domain of $G_{f}$ with $\fiF\to x$. Since $f$ is uniformly continuous, we know that $f(\fiF)$ generates a Cauchy filter $\filterFont{G}$ on $Y$, which converges to some $y$. The product filter $\fiF\times\filterFont{G}$ converges to $\langle x, y\rangle$, thus $x$ is in the domain of the closure of $G_{f}$. 

2.
Now let $W\in\filterFont{v}$; we show that there exists a neighborhood $U\in\filterFont{u}$ with this property: if $\langle x, y\rangle, \langle u, v\rangle\in\Closure{G}_{f}$, then $x\in U[u]$ implies $y\in W[v]$. After having established this, we know
\begin{itemize}
\item $\Closure{G}_{f}$ is the graph of a function $F$. This is so because $Y$ is separated, hence its topology is Hausdorff. For, assume there exists $x\in X$ some $y_{1}, y_{2}\in Y$ with $y_{1}\not= y_{2}$ and $\langle x, y_{1}\rangle, \langle x, y_{2}\rangle\in\Closure{G}_{f}$. Choose $W\in\filterFont{v}$ with $y_{2}\not\in W[y_{1}]$, and consider $U$ as above. Then $x\in U[x]$, hence $y_{2}\in W[y_{1}]$, contradicting the choice of $W$.
\item $F$ is uniformly continuous. The property above translates to finding for $W\in\filterFont{v}$ a neighborhood $U\in\filterFont{u}$ with $U\subseteq \InvBild{(\zZ{F})}{W}$. 
\end{itemize}
So we are done after having established the statement above.

3.
Assume that $W\in\filterFont{v}$ is given, and choose $V\in\filterFont{v}$ closed and symmetric with $V\circ V\subseteq W$. This is possible by Corollary~\ref{clos-symm-are-base}. There exists $U\in\filterFont{u}$ open and symmetric with $\Bild{f}{U[x]}\subseteq V[f(x)]$ for every $x\in A$, since $f$ is uniformly continuous. If $\langle x, y\rangle, \langle u, v\rangle\in \Closure{G}_{f}$ and $x\in U[u]$, then $U[x]\cap U[u]$ is open (since $U$ is open), and there exists $a\in A$ with $x, u\in U[a]$, since $A$ is dense. We claim $y\in \Closure{\bigl(\Bild{f}{U[a]}\bigr)}$. Let $H$ be an open neighborhood of $y$, then, since $U[a]$ is a neighborhood of $x$, $U[a]\times H$ is a neighborhood of $\langle x, y\rangle$, thus $G_{f}\cap U[a]\times H \not= \emptyset$. Hence we find $y'\in H$ with $\langle x, y'\rangle\in G_{f}$, which entails $H\cap \Bild{f}{U[a]}\not=\emptyset$. Similarly, $z\in \Closure{\bigl(\Bild{f}{U[a]}\bigr)}$; note $\Closure{\bigl(\Bild{f}{U[a]}\bigr)}\subseteq V[f(a)]$. But now $\langle y, v\rangle\in V\circ V\subseteq W$, hence $y\in W[v]$. This establishes the claim above, and finishes the proof.
\EndProof

Let us just have a look at the idea lest it gets lost. If $x\in X$, we find a filter $\fiF$ on $A$ with $i_{A}(\fiF)\to x$. Then $f(i_{a}(\fiF))$ is a Cauchy filter, hence it converges to some $y\in Y$, which we define as $F(x)$. Then it has to be shown that $F$ is well defined, it clearly extends $f$. It finally has to be shown that $F$ is uniformly continuous. So there is a lot technical ground which is covered. 

We note on closing that also the completion of pseudometric spaces can be translated into the realm of uniform spaces. Here, naturally, the Cauchy filters defined on the space play an important r\^ole, and things get very technical. The interplay between compactness and uniformities yields interesting results as well, here the reader is referred to~\cite[Chapter~II]{Bourbaki} or to~\cite{James-Unif}. 

\Subsection{Bibliographic Notes} 
\label{sec:top-bib-notes}
The towering references in this area are~\cite{Bourbaki, Engelking,
  Kuratowski, Kelley}; the author had the pleasure of taking a course on topology from one of the authors of~\cite{Querenburg}, so this text has been an important source, too. The delightful Lecture Note~\cite{Herrlich-Choice} by Herrlich has a chapter ``Disasters without Choice'' which discusses among others the relationship of the Axiom of Choice and various topological constructions. The discussion of the game related to Baire's Theorem in Section~\ref{sec:baire+game} in taken from Oxtoby's textbook~\cite[Sec. 6]{Oxtoby} on the duality between measure and category (\emph{category} in the topological sense introduced on page~\pageref{baire-complete} above); he attributes the game to Banach and Mazur. Other instances of proofs by games for metric spaces are given, e.g., in~\cite[8.H, 21]{Kechris}. The uniformity discused in Example~\ref{example-uniformities-5-filter} has been considered in~\cite{EED-Polymorphic} in greater detail. The section on topological systems follows fairly closely the textbook~\cite{Vickers} by Vickers, but see also~\cite{Cont-Latt, Abramsky+Jung}, and for the discussion of dualities and the connection to intuitionistic logics,~\cite{Johnstone, Goldblatt-Topoi}. The discussion of Gödel's Completeness Theorem in Section~\ref{sec:goedel} is based on the original paper by Rasiowa and Sikorski~\cite{Rasiowa-Sikorski-I} together with occasional glimpses at~\cite[Chapter 2.1]{Chang+Keisler},\cite[Chapter 4]{Srivastava-Logic} and~\cite[Chapter 1.2]{Koppelberg}. Uniform spaces are discussed in\cite{Bourbaki, Engelking, Kelley, Querenburg}, special treatises include\cite{James-Unif} and \cite{Isbell}, the latter one emphasizing a categorical point of view. 

\Subsection{Exercises} 
\label{sec:tops-exercises}
\Exercise{ex-char-final-top}{ Formulate and prove an analogue of
  Proposition~\ref{initial-final} for the final topology for a family
  of maps.  } 
\Exercise{ex-product-euclidean}{ The Euclidean topology
  on $\Real^{n}$ is the same as the product topology on
  $\prod_{i=1}^{n}\Real$.  } \Exercise{ex-prod-discrete}{ Recall that
  the topological space $(X, \tau)$ is called \emph{discrete} iff
  $\tau=\PowerSet{X}$. Show that the product
  $\prod_{i\in I} (\{0, 1\}, \PowerSet{\{0, 1\}})$ is discrete iff the
  index set $I$ is finite.  } \Exercise{ex-compare-tops}{ Let
  $ L := \{(x_{n})_{n\in \Nat}\} \subseteq \Real^{\Nat}\mid \sum_{n\in
    \Nat}|x_{n}| < \infty\} $
  be all sequences of real numbers which are absolutely
  summable. $\tau_{1}$ is defined as the trace of the product topology
  on $\prod_{n\in\Nat}\Real$ on $L$, $\tau_{2}$ is defined in the
  following way: A set $G$ is $\tau_{2}$-open iff given $x\in G$,
  there exists $r>0$ such that
  $ \{y\in L\mid \sum_{n\in \Nat}|x_{n}-y_{n}| < r\} \subseteq G.  $
  Investigate whether the identity maps
  $(L, \tau_{1})\to (L, \tau_{2})$ and
  $(L, \tau_{2})\to (L, \tau_{1})$ are continuous.  }
\Exercise{ex-equiv-classes}{ Define for $x, y\in \Real$ the
  equivalence relation
  $ \isEquiv{x}{y}{\sim} \text{ iff } x - y \in\Ganz.  $ Show that
  $\Faktor{\Real}{\sim}$ is homeomorphic to the unit circle.
  \textbf{Hint:} Example~\ref{example-quotient-top}.  }
\Exercise{ex-scott-equiv}{ Let $A$ be a countable set. Show that a map
  $q: (\partMap{A}{B})\to (\partMap{C}{D})$ is continuous in the
  topology taken from Example~\ref{partial-maps-top} iff it is 
  continuous, when $\partMap{A}{B}$ as well as $\partMap{C}{D}$ are equipped with the Scott topology.  } 
\Exercise{ex-order-top}{ Let $D_{24}$ be the set of
  all divisors of $24$, including $1$, and define an order
  $\sqsubseteq$ on $D_{24}$ through $x\sqsubseteq y$ iff
  $x \text{ divides } y$. The topology on $D_{24}$ is given through
  the closure operator as in
  Example~\ref{finite-ordered-for-closure}. Write a \texttt{Haskell}
  program listing all closed subsets of $D_{24}$, and determining all
  filters $\fiF$ with $\fiF\to 1$.  \textbf{Hint:} It is helpful to
  define a type \texttt{Set} with appropriate operations first,
  see~\cite[4.2.2]{EED-Haskell-Buch}.  } 
\Exercise{ex-closure-filter}{
  Let $X$ be a topological space, $A\subseteq X$, and $i_{A}: A\to X$
  the injection. Show that $x\in\Closure{A}$ iff there exists a filter
  $\fiF$ on $A$ such that $i_{A}(\fiF)\to x$.  }
\Exercise{ex-reals-are-normal}{ Show \emph{by expanding
    Example~\ref{reals-are-t3half}} that $\Real$ with its usual
  topology is a $T_{4}$-space.  } \Exercise{ex-compact-homeom}{ Given
  a continuous bijection $f: X\to Y$ with the Hausdorff spaces $X$ and
  $Y$, show that $f$ is a homeomorphism, if $X$ is compact.  }
\Exercise{ex-inherit-separation}{ Let $A$ be a subspace of a
  topological space $X$.
\begin{enumerate}
\item If $X$ is a $T_{1}, T_{2}, T_{3}, T_{3\half}$ space, so is $A$.
\item If $A$ is closed, and $X$ is a $T_{4}$-space, then so is $A$.
\end{enumerate}
}
\Exercise{ex-semicontinuous}{
A function $f: X\to \Real$ is called \emph{\index{semicontinuous!lower}lower semicontinuous} iff for each
$c\in \Real$ the set $\{x\in X \mid f(x) < c\}$ is open. If $\{x\in
X\mid f(x) > c\}$ is open, then $f$ is called
\emph{\index{semicontinuous!upper}upper semicontinuous}. If $X$ is
compact, then a lower semicontinuous map assumes on $X$ its maximum,
and an upper semicontinuous map assumes its minimum.  
}
\Exercise{ex-prod-locally-compact}{
Let $X := \prod_{i\in I}X_{i}$ be the product of the Hausdorff space
$(X_{i})_{i\in I}$. Show that $X$ is locally compact in the product
topology iff $X_{i}$ is locally compact for all $i\in I$, and all but
a finite number of $X_{i}$ are compact. 
}
\Exercise{ex-river-metric}{
Given $x, y\in \Real^{2}$, define
  \begin{equation*}
D(x, y) := 
    \begin{cases}
     |x_{2}-y_{2}|, & \text{ if }x_{1} = y_{1}\\
|x_{2}| + |y_{2}| + |x_{1}-y_{1}|,& \text{ otherwise}. 
    \end{cases}
  \end{equation*}
Show that this defines a metric on the plane $\Real^{2}$. Draw the open
ball $\{y\mid D(y, 0) < 1\}$ of radius $1$ with the origin as center.}
\Exercise{ex-pseudo-t1-is-metric}{
Let $(X, d)$ be a pseudometric space such that the induced topology is
$T_{1}$. Then $d$ is a metric.
}
\Exercise{ex-cont-seq-testing}{
Let $X$ and $Y$ be two first countable topological spaces. Show that a
map $f: X\to Y$ is continuous iff $x_{n}\to x$ implies always
$f(x_{n})\to f(x)$ for each sequence $(x_{n})_{n\in\Nat}$ in $X$. 
}
\Exercise{ex-contfncts-int}{
Consider the set $\Cont[[0, 1]]$ of all continuous functions on the unit
interval, and define
\begin{equation*}
  e(f, g) := \int_{0}^{1}|f(x) - g(x)|\ dx.
\end{equation*}
Show that
\begin{enumerate}
\item $e$ is a metric on $C([0, 1])$.
\item $\Cont[[0, 1]]$ is not complete with this metric.
\item The metrics $d$ on $\Cont[[0, 1]]$ from Example~\ref{for-metric-spaces} and $e$ are not equivalent. 
\end{enumerate}
}
\Exercise{ex-ultrametric-space}{
Let $(X, d)$ be an ultrametric space, hence $d(x, z) \leq \max\ \{d(x,
y), d(y, z)\}$ (see Example~\ref{for-metric-spaces}). Show that
\begin{itemize}
\item If $d(x, y) \not= d(y, z)$, then $d(x, z) = \max\ \{d(x, y),
  d(y, z)\}$.
\item Any open ball $B(x, r)$ is both open and closed, and $B(x, r) =
  B(y, r)$, whenever $y\in B(x, r)$.
\item Any closed ball $S(x, r)$ is both open and closed, and $S(x, r) =
  S(y, r)$, whenever $y\in S(x, r)$.
\item Assume that $B(x, r)\cap B(x', r')\not=\emptyset$, then $B(x,
  r)\subseteq B(x', r')$ or $B(x', r')\subseteq B(x, r)$. 
\end{itemize}

}
\Exercise{ex-nowhere-dense}{
Show that the set of all nowhere dense sets in a topological space $X$
forms an ideal. Define a set $A\subseteq X$ as \emph{open modulo nowhere
dense sets} iff there exists an open set $G$ such that the symmetric
difference $A\Delta G$ is nowhere dense (hence both $A\setminus G$ and
$G\setminus A$ are nowhere dense). Show that the open sets modulo
nowhere dense sets form an $\sigma$-algebra.}
\Exercise{ex-oxtoby-game}{
Consider the game formulated in Section~\ref{sec:baire+game}; we use
the notation from there. Show that there exists a strategy that Angel
can win iff $L_{1}\cap B$ is of first category for some interval
$L_{1}\subseteq L_{0}$. 
}
\Exercise{ex-open-not-neighbor}{
Let $\filterFont{u}$ be the additive uniformity on $\Real$ from
Example~\ref{example-uniformities}. Show that $\{\langle x,
y\rangle\mid |x-y|<1/(1+|y|)\}$ is not a member of $\filterFont{u}$. 
}
\Exercise{ex-relat-symm-identity}{
Show that $V\circ U\circ V = \bigcup_{\langle x, y\rangle \in
U}V[x]\times V[y]$ for symmetric $V\subseteq X\times X$ and arbitrary
$U\subseteq X\times X$.
}
\Exercise{ex-unif-base}{
Given a base $\filterFont{b}$ for a uniformity $\filterFont{u}$, show that 
\begin{align*}
  \filterFont{b'} & := \{B\cap B^{-1}\mid B\in\filterFont{b}\},\\
\filterFont{b''} & := \{B^{n}\mid B\in \filterFont{b}\}
\end{align*}
are also bases for $\filterFont{u}$, when $n\in\Nat$ (recall $B^{1} := B$ and $B^{n+1} :=
B\circ B^{n}$). 
}
\Exercise{ex-complete-lattice-unif}{
Show that the uniformities on a set $X$ form a complete lattice with
respect to inclusion. Characterize the initial and the final
uniformity on $X$ for a family of functions in terms of this lattice. 
}
\Exercise{ex-inters-small}{
If two subsets $A$ and $B$ in a uniform space $(X, \filterFont{u})$
are $V$-small then $A\cup B$ is $V\circ V$-small, if $A\cap
B\not=\emptyset$. 
}
\Exercise{ex-ultra-discr}{
Show that a discrete uniform space is complete. \textbf{Hint}: A Cauchy filter
is an ultrafilter based on a point.
}
\Exercise{ex-product-top-vs-unif}{
Let ${\cal F}$ be a family of maps $X\to Y_{f}$ with uniform spaces
$(Y_{f}, \filterFont{v}_{f})$. Show that the initial topology on $X$
with respect to ${\cal F}$ is the topology induced by the product
uniformity. 
}
\Exercise{ex-proj-unif-cont}{
Equip the product $X := \prod_{i\in I}X_{i}$ with the product uniformity for
the uniform spaces $\bigl((X_{i}, \filterFont{u}_{i})\bigr)_{i\in
  I}$, and let $(Y, \filterFont{v})$ be a uniform space. A map $f: Y\to X$ is uniformly continuous iff $\pi_{i}\circ f:
Y\to X_{i}$ is uniformly continuous for each $i\in I$.}
\Exercise{ex-spatial-equiv}{
Let $X$ be a topological system. Show that the following statements
are equivalent
\begin{enumerate}
\item $X$ is homeomorphic to $\funSP(Y)$ for some topological system
  $Y$.
\item For all $a, b\in\oMX$ holds $a = b$, provided we have $x\models
  a \Leftrightarrow x\models b$ for all $x\in \pTX$.
\item For all $a, b\in\oMX$ holds $a \leq b$, provided we have $x\models
  a \Rightarrow x\models b$ for all $x\in \pTX$.
\end{enumerate}
}
\Exercise{ex-hausdorff-sober}{
Show that a Hausdorff space is sober.}
\Exercise{ex-stone-weierstrass}{
Let $X$ and $Y$ be compact topological spaces with their Banach spaces
$\Cont$ resp. $\Cont[Y]$ of real continuous maps. Let $f: X\to Y$ be a continuous map, then 
\begin{equation*}
  f^{*}:
  \begin{cases}
    \Cont[Y]&\to \Cont\\
g&\mapsto g\circ f
  \end{cases}
\end{equation*}
defines a continuous map (with respect to the respective norm
topologies). $f^{*}$ is onto iff $f$ is an injection. $f$ is onto iff
$f^{*}$ is an isomorphism of $\Cont[Y]$ onto a ring $A\subseteq \Cont$
which contains constants.  } 
\Exercise{ex-completeness-propositional}{
Let $\theLang$ be a language for propositional logic with constants
$C$ and $V$ as the set of propositinal variables. Prove that a
consistent theory $T$ has a model, hence a map $h: V\to \Zwei$ such
that each formula in $T$ is assigned the vlaue $\top$. \textbf{Hint:}
Fix an ultrafilter on the Lindenbaum algebra of $T$ and consider the
corresponding morphism into $\Zwei$.}
\Exercise{ex-top-group}{
Let $G$ be a topological group, see Example~\ref{top-group}. 
Given $F\subseteq G$ closed, show that
\begin{enumerate}
\item $gF$ and $Fg$ are closed,
\item $F^{-1}$ is closed,
\item $MF$ and $FM$ are closed, provided $M$ is finite. 
\item If $A\subseteq G$, then 
$ 
\Closure{A} = \bigcap_{U\in\upsilon(e)}AU = \bigcap_{U\in\upsilon(e)}UA = \bigcap_{U\in\tau}AU = \bigcap_{U\in\tau}UA.
$ 
\end{enumerate}}

}

\bibliographystyle{alpha}
\newpage
\newcommand{\etalchar}[1]{$^{#1}$}

\newpage
\addcontentsline{toc}{section}{Index}
\begin{theindex}

  \item $B(x, r)$, 38
  \item $S(x, r)$, 38
  \item $U[x]$, 84
  \item $\Cont$, 76
  \item $\epsilon $-net, 50
  \item $\mathsf{diam}(A)$, 46
  \item $\upsilon(x)$, 13
  \item ${\cal     F}\to x$, 15
  \item 2\negthickspace2, 66

  \indexspace

  \item accumulation point, 19
  \item algebra
    \subitem Heyting, 64
      \subsubitem morphism, 65
    \subitem Lindenbaum, 61

  \indexspace

  \item Banach-Mazur game, 56

  \indexspace

  \item Cauchy filter, 92
  \item Cauchy sequence, 43
  \item Charly Brown's device, 45
  \item closure operator, 11
  \item compact
    \subitem countably compact, 34
    \subitem Lindel\IeC {\"o}f space, 35
    \subitem locally compact, 29
    \subitem paracompact, 35
      \subsubitem locally finite, 35
      \subsubitem refinement, 35
    \subitem sequentially compact, 35
  \item compactification, 31
    \subitem Alexandrov one point, 30
    \subitem Stone-$\mathaccentV {check}014{\mathrm  {C}}$ech, 32
  \item continuous, 6
    \subitem uniformly, 53, 94
  \item contraction, 48
  \item convergence
    \subitem filter, 15
    \subitem net, 15

  \indexspace

  \item dcpo, 72
  \item dense set, 25
  \item diameter, 46

  \indexspace

  \item embedding, 30
  \item entourage, 83

  \indexspace

  \item filter
    \subitem accumulation point, 19
    \subitem Cauchy, 92
    \subitem neighborhood, 13
  \item flyswatter, 50
  \item frame, 65

  \indexspace

  \item Google, 50
  \item group
    \subitem topological, 15

  \indexspace

  \item homeomorphism, 10

  \indexspace

  \item irreducible, 71

  \indexspace

  \item map
    \subitem affine, 85
  \item metric
    \subitem discrete, 36
    \subitem Hausdorff, 41
    \subitem ultrametric, 37
  \item model, 62, 63

  \indexspace

  \item Nachbarschaft, 83
  \item net, 15
    \subitem Cauchy, 93
    \subitem convergence, 15
  \item norm, 76
  \item nowhere dense, 55

  \indexspace

  \item open
    \subitem Scott, 6, 73
  \item oscillation, 46

  \indexspace

  \item partition, 84
  \item prime
    \subitem completely, 69
    \subitem element, 69
  \item pseudometric, 35
  \item pseudometrics
    \subitem equivalent, 38

  \indexspace

  \item semicontinuous
    \subitem lower, 99
    \subitem upper, 99
  \item sentence, 61
  \item sequence
    \subitem Cauchy, 43
  \item set
    \subitem saturated, 75
    \subitem small, 92
  \item Sorgenfrey line, 21
  \item space
    \subitem $T_{0}, T_{1}$, 22
    \subitem $T_{3},     T_{3\half}, T_{4}$, 22
    \subitem Banach, 77
    \subitem completely regular, 25
    \subitem first category, 55
    \subitem Hausdorff, $T_{2}$, 21
    \subitem locally compact, 29
    \subitem metric, 35
    \subitem normal, 25
    \subitem normed, 77
    \subitem pseudometric, 35
      \subsubitem complete, 43
    \subitem regular, 25
    \subitem uniform, 83
      \subsubitem complete, 92
      \subsubitem separated, 89
      \subsubitem topology, 87

  \indexspace

  \item theorem
    \subitem Baire
      \subsubitem complete pseudometric, 55
      \subsubitem locally compact, 33
    \subitem Dini, 77
    \subitem Hofmann-Mislove, 75
    \subitem Stone-Weierstra\IeC {\ss }, 79
    \subitem Tihonov, 19
    \subitem Urysohn's Metrization, 42
  \item topological system, 66
    \subitem c-morphism, 67
    \subitem homeomorphism, 67
    \subitem localic, 69
    \subitem localization, 69
    \subitem opens, 66
      \subsubitem extension, 66
    \subitem points, 66
    \subitem spatialization, 67
  \item topology
    \subitem base, 3
    \subitem compactification, 31
    \subitem final, 8
    \subitem first countable, 41
    \subitem initial, 8
    \subitem product, 9
    \subitem quotient, 9
    \subitem Scott, 6, 73
    \subitem second countable, 41
    \subitem separable, 41
    \subitem sober, 72
    \subitem sum, 9
    \subitem topological group, 15
    \subitem trace, 9
    \subitem uniform, 87
    \subitem uniform convergence, 77
    \subitem Vietoris, 53
    \subitem weak, 5
  \item totally bounded, 50

  \indexspace

  \item ultrametric, 37
  \item Umgebung, 83
  \item uniformity, 83
    \subitem $p$-adic, 85
    \subitem additive, 84
    \subitem discrete, 84
    \subitem finite partitions, 85
    \subitem indiscrete, 84
    \subitem initial, 95
    \subitem multiplicative, 84
    \subitem product, 96
    \subitem subspace, 96
  \item Urysohn's Lemma, 25

\end{theindex}

\end{document}